\documentclass[11pt]{article}

\usepackage{graphicx}
\usepackage{color}
\usepackage{colordvi}
\usepackage{latexsym,amsmath,amsfonts,amscd,subfigure,amssymb,amsthm}
\usepackage{threeparttable}
\newcommand{\V}[1]{\mbox{\boldmath $ #1 $}}

\newcommand{\Rmnum}[1]{\expandafter\@slowromancap\romannumeral #1@}

\def \M{\mathbb{M}}

\usepackage[pdftex]{hyperref}

\newcommand{\bey}{\begin{eqnarray}}
\newcommand{\eey}{\end{eqnarray}}

\newcommand{\beq}{\begin{equation}}
\newcommand{\eeq}{\end{equation}}
\theoremstyle{plain}

\theoremstyle{definition}

\theoremstyle{remark}
\newtheorem{exam}{\hspace{6mm}Example}[section]

\title{Moving Mesh Finite Difference Solution\\ of Non-Equilibrium Radiation Diffusion Equations}

\author{
Xiaobo~Yang\footnote{Department of Mathematics, College of Science, China University of Mining and Technology,
Xuzhou, Jiangsu 221116, China. Email: xwindyb@126.com.},
\and Weizhang~Huang\footnote{Department of Mathematics, the University of Kansas, Lawrence, KS 66045, U.S.A.
Email: whuang@ku.edu.},
\and Jianxian~Qiu\footnote{School of Mathematical Sciences, Xiamen University, Xiamen,
Fujian 361005, China. Email: jxqiu@xmu.edu.cn.}
}

\date{}

\begin{document}
\vskip 1cm
\maketitle

\renewcommand{\thefootnote}{\fnsymbol{footnote}}

\renewcommand{\thefootnote}{\arabic{footnote}}

\begin{abstract}

A moving mesh finite difference method based on the moving mesh partial differential equation
is proposed for the numerical solution of the 2T model for multi-material, non-equilibrium
radiation diffusion equations. The model involves nonlinear diffusion coefficients
and its solutions stay positive for all time when they are positive initially.
Nonlinear diffusion and preservation of solution positivity pose challenges in the numerical
solution of the model. A coefficient-freezing predictor-corrector
method  is used for nonlinear diffusion while a cutoff strategy with a positive threshold
is used to keep the solutions positive. Furthermore, a two-level moving mesh
strategy and a sparse matrix solver are used to
improve the efficiency of the computation.
Numerical results for a selection of examples of multi-material non-equilibrium radiation diffusion show that the method is capable of capturing the profiles and local structures of Marshak waves with adequate mesh concentration. The obtained numerical solutions are in good agreement
with those in the existing literature. Comparison studies are also made between uniform and
adaptive moving meshes and between one-level and two-level moving meshes.
\end{abstract}

{\bf Key words:} moving mesh method; non-equilibrium radiation diffusion;
predictor-corrector; positivity; two-level mesh movement.

{\bf AMS subject classification:} 65M06, 65M50

\section{Introduction}
\label{SEC:intro}
\setcounter{equation}{0} \setcounter{figure}{0}
\setcounter{table}{0}

Radiation transport in astrophysical phenomena and inertial confinement fusion
can often be modeled using a set of coupled diffusion equations when photon mean free paths
are much shorter than characteristic length scales. These equations are highly
nonlinear and exhibit multiple time and space scales \cite{MM84}.
Particularly, steep hot wave fronts, called Marshak waves, typically form during
radiation transport processes. Energy density and material temperature near
the steep fronts can vary dramatically in a short distance. Such complex local structures
make mesh adaptation an indispensable tool for use to improve the efficiency in
the numerical solution of radiation diffusion equations because
the number of mesh points can be prohibitively large when a uniform mesh is used.
Research of radiation diffusion has attracted considerable attentions from engineers
and scientists \cite{Castor2004,Kang2002, LC06,MM84,MKR2000, MK2003,Cai2004,
Olson2007,PP06,PP13,Knoll1999,yuan20092,yuan2009,yuan2011,yhq2015,liyonghai2013}.

In this work  we are interested in the non-equilibrium situation where
the radiation field is not in thermodynamics equilibrium with the material temperature.
Marshak \cite{Marshak1958} develops a time-dependent radiative transfer model,
laying the groundwork for the research area. Pomraning \cite{Pomraning1979} obtains
an analytic solution to a particular Marshak wave problem, which is analyzed
more extensively by Su and Olson \cite{Su1996}.
Numerically, Mousseau et al. \cite{MKR2000, MK2003} present
a physics-based preconditioning Newton-Krylov method involving Jacobian-free Newton-Krylov (JFNK),
operator splitting, and multigrid linear solvers and show that the method can capture the Marshak wave
of the thermal transport front properly.
Kang \cite{Kang2002} proposes a $P1$ nonconforming finite element method for non-equilibrium
radiation transport problems.
Olson \cite{Olson2007} considers a hydrogen-like Saha ionization model for a simplified but physically
plausible heat capacity and uses several types of finite difference (FD) schemes to approximate flux-limiting.
Sheng et al. \cite{yuan20092} construct a monotone finite volume scheme for multi-material, non-equilibrium
radiation diffusion equations and show numerically that their method is better than the standard nine-point
finite difference scheme and preserves the nonnegativity of energy density.

On the other hand, there exist only a few published studies that have employed mesh adaptation
for the numerical solution of radiation diffusion equations.
For example, Lapenta and Chac\'{o}n \cite{LC06} use a fully implicit moving mesh method to solve
a one-dimensional equilibrium radiation diffusion equation. They discretize both the mesh and physical
equations using finite volumes and solve the resulting equation with a preconditioned
inexact-Newton method. Their results show great improvements in cost-effectiveness with mesh adaptation.
Yang et al. \cite{yhq2015} study a moving mesh FD method based on the moving-mesh-partial-differential-equation
(MMPDE) strategy \cite{HRR94b,HR11} for equilibrium radiation diffusion equations
and show that the method capture Marshak waves accurately and efficiently.
Pernice et al. \cite{PP13} use adaptive mesh refinement
to solve three-dimensional non-equilibrium radiation diffusion equations. They use implicit time integration
for stiff multi-physics systems as well as
the JFNK  \cite{KRO1999, KRO2001,MKR2000,Knoll1999}
to solve the resulting nonlinear algebraic equations. They also use an optimal multilevel preconditioner
to provide level-independent solver convergence. Non-equilibrium radiation diffusion equations
are challenging to solve, but through the numerical results, they demonstrate their method can
efficiently capture the local structures of Marshak waves and can give convincing results with good
accuracy.

The objective of this work is to study a moving mesh FD solution of two-dimensional
non-equilibrium radiation diffusion systems. The method is based on the MMPDE
moving mesh approach \cite{HRR94b,HR11}. The MMPDE is used to adaptively move
the mesh around evolving features of the physical solution and is defined as the gradient flow
equation of a meshing functional based on mesh equidistribution and alignment.
The shape, size, and orientation of mesh elements are controlled through a monitor function
\cite{Hua01b} defined through the Hessian of the energy density.
A similar moving mesh FD method has been developed in \cite{yhq2015} for equilibrium
radiation diffusion equations, and the current work can be considered as a generalization of \cite{yhq2015}.
However, this generalization is non-trival. Unlike \cite{yhq2015}, we now need to deal with
a system of two coupled equations for the energy density and material temperature. The diffusion
coefficients depend on both the energy density and material temperature and it is more sensitive
to treat diffusion numerically.
Moreover, the system is stiffer, making it more difficult to integrate in time (with smaller time steps)
and more expensive to solve overall. Furthermore, it is more delicate to preserve the solution
positivity. Like \cite{yhq2015}, we use here the cutoff strategy to maintain the positivity in the computed
solutions. It has been shown in \cite{LuHuVV2012} that the strategy retains the accuracy and
convergence order of FD approximation for parabolic PDEs.
It has been found in \cite{yhq2015} that the strategy with a threshold zero (meaning that the computed solutions are
kept to be nonnegative) works for equilibrium radiation diffusion equations.
For the current situation, on the other hand, we have found that a positive threshold is needed
and an empirical choice depending on the mesh size seems to work well for problems we have tested.
Numerical results for a selection of examples are presented.
They show that the method is capable of capturing the profiles and local structures
of Marshak waves with adequate mesh concentration. The obtained numerical solutions
are in good agreement with those of \cite{Kang2002,yuan20092}.
Comparison studies are also made between uniform and
adaptive moving meshes and between one-level and two-level moving meshes.

The  outline of the paper is as follows. The physical model and governing equations
are described in \S\ref{SEC:problem}. The moving mesh FD method and the treatments of nonlinearity
as well as the cutoff strategy are discussed in \S\ref{SEC:MFD}. In \S\ref{SEC:numerics} numerical results
obtained for a selection of examples of multi-material, multiple spot concentration scenarios. Finally,
conclusions are drawn in \S\ref{SEC:conclusions}.

\section{The 2T model for non-equilibrium radiation diffusion}
\label{SEC:problem}

Under the assumption of an optically thick medium (short mean free path of photons) a
first-principle statement of radiation transport reduces to the radiation diffusion limit.
A particular idealized dimensionless form of the governing system, known as the 2T model,
consists of two equations, the radiation diffusion (gray approximation) equation and material energy
balance equation, that is,
\begin{equation}
\begin{cases}
\frac{\partial E}{\partial t} - \nabla\cdot(D_r\nabla E)=\sigma_a(T^4-E), \\
\frac{\partial T}{\partial t} - \nabla\cdot(D_t\nabla T)=-\sigma_a(T^4-E),
\end{cases}
\label{2T}
\end{equation}
where
\begin{equation}
\sigma_a = \frac{z^3}{T^3},\quad
D_r=\frac{1}{3\sigma_a+\frac{1}{E}|\nabla E|}, \quad D_t=\kappa\, T^{\frac52}.
\label{para1}
\end{equation}
Here, $E$ represents the photon energy, $T$ is the material temperature,
$\sigma_a$ is the opacity, $\kappa$ is the material conductivity, and $z$ is the atomic mass number.
Notice that a limiting term $|\nabla E|/E$ is added to the diffusion coefficient $D_r$ to avoid
a possible unphysical behavior that a flux of energy moves faster than the speed of light in regions
of strong gradient where a simple diffusion theory can fail. Moreover, we use
the form of the material (plasma) conduction diffusion coefficient $D_t$ from Spitzer and Harm
\cite{Harm} and take $\kappa = 0.01$ in our computation.
Furthermore, compared to the equilibrium case, the nonlinear source terms on the right-hand sides
of the equations do not vanish in general, reflecting the transfer of energy between the radiation field
and material temperature. Additional nonlinearities come from the particular form of diffusion
coefficients, which are functions of $E$ and $T$.

\begin{figure}[htb]
\begin{center}
\includegraphics[scale = 0.4]{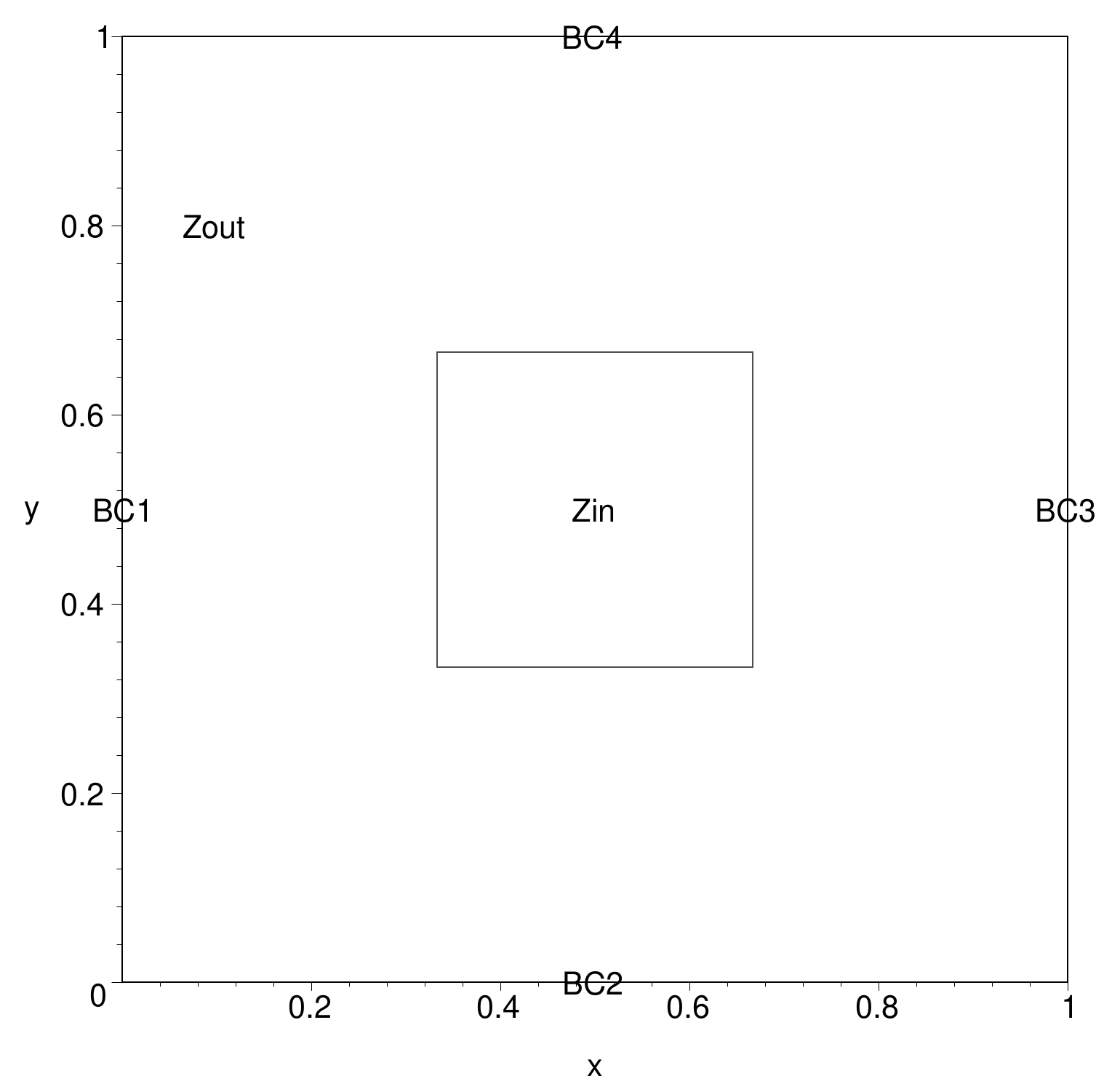}
\end{center}
\caption{The physical domain with a middle inset $(1/3, 2/3)\times (1/3, 2/3)$.
The value of $z$ is $Z_{\textrm{in}}$ and $Z_{\textrm{out}}$ inside and outside the inset, respectively.
BC1 and BC2 are the inflow and outflow boundaries, respectively, while BC3 and BC4 are
perfectly insulated.}
\label{BCD1}
\end{figure}

We consider (\ref{2T}) in two dimensions on the unit square domain, see Fig.~\ref{BCD1}.
Homogenous Neumann boundary conditions are used for boundary segments
$y = 0$ and $y=1$ and inflow and outflow boundary conditions are employed on $x=0$ and
$x=1$, respectively. More specifically, we have
\begin{equation}
\begin{cases}
\frac{\partial E}{\partial y} = 0, \quad \frac{\partial T}{\partial y} = 0, & \quad \text{ on } y = 0 \text{ or } y = 1 \\
\frac{\partial T}{\partial x} = 0, & \quad \text{ on } x = 0 \text{ or } x = 1 \\
\frac14E - \frac{1}{6\sigma_a}\frac{\partial E}{\partial x}=1, & \quad \text{ on } x = 0\\
\frac14E + \frac{1}{6\sigma_a}\frac{\partial E}{\partial x}=0,& \quad \text{ on } x = 1 .
\end{cases}
\label{BC-1}
\end{equation}
The initial conditions are
\begin{equation}
\begin{cases}
E(x,y,0) = (1-\tanh(10 x))(1-10^{-5}) + 10^{-5}, & \quad (x,y) \in \Omega\\
T(x,y,0) = E(x,y,0)^{\frac14}, & \quad (x,y) \in \Omega.
\end{cases}
\label{IC-2}
\end{equation}
Essentially $E(x,y,0)$ is equal to $10^{-5}$ everywhere except on the boundary $x=0$ where it is 1.
A narrow transition between $10^{-5}$ and 1 is used
to avoid a potentially difficult initial start in numerical computation.
(Slightly different initial and boundary conditions are used in the third example in \S\ref{SEC:numerics}.)

\section{The moving mesh FD method}
\label{SEC:MFD}
\setcounter{equation}{0}
\setcounter{figure}{0}
\setcounter{table}{0}

In this section we describe the moving mesh FD method for solving
the initial-boundary value problem (IBVP) (\ref{2T}), (\ref{BC-1}), and (\ref{IC-2}).
We discretize this problem in space using central finite differences and in time using
a Singly Diagonally Implicit Runge-Kutta scheme (SDIRK) \cite{Cas79}. We also discuss
linearization of the equations, preservation of solution positivity, and
adaptive mesh movement.

\subsection{FD discretization on moving meshes}

We denote a curvilinear moving mesh for $\Omega$ by
\begin{equation}
(x_{m,n}(t), \; y_{m,n}(t)), \quad m=1,...,M,\; n=1,...,N
\label{mmesh-1}
\end{equation}
where $M$ and $N$ are positive integers. The generation of such an adaptive moving mesh
will be described in \S\ref{SEC:mmesh}. For the moment, we consider (\ref{mmesh-1}) as the image
of a fixed rectangular mesh under a known coordinate transformation
$x = x(\xi, \eta, t)$, $y = y(\xi, \eta, t)$, i.e.,
\begin{equation}
\label{mmesh-2}
x_{m,n}(t) = x(\xi_m,\eta_n,t),\quad y_{m,n}(t) = y(\xi_m,\eta_n,t), \quad m=1,...,M,\quad n=1,...,N
\end{equation}
where the reference mesh is taken as
\begin{equation}
\label{rmesh-1}
(\xi_m, \eta_n) = ((m-1)\Delta \xi,(n-1)\Delta \eta), \quad m=1,...,M,\quad n=1,...,N
\end{equation}
and $\Delta \xi = 1/(M-1)$ and $\Delta \eta = 1/(N-1)$. The boundary correspondence
between the reference and physical domains is given by
\begin{equation}
x(0, \eta) = 0, \quad x(1, \eta) = 1, \quad y(\xi, 0) = 0, \quad y(\xi, 1) = 1.
\label{mesh-bc-1}
\end{equation}
We let
\[
\hat{E}(\xi,\eta,t) = E(x(\xi,\eta,t),y(\xi,\eta,t),t),\qquad \hat{T}(\xi,\eta,t) = T(x(\xi,\eta,t),y(\xi,\eta,t),t) .
\]

The discretization of the 2T model on the moving mesh (\ref{mmesh-1}) consists of two steps,
its transformation from $(x,y)$ to $(\xi, \eta)$ and discretization on the rectangular reference mesh.
First, using the coordinate transformation, we can transform (\ref{2T})  (e.g., see \cite[\S3.1.4]{HR11})
into the reference domain as
\begin{equation}
\begin{cases}
\hat{E}_t - \V{b}(t) \cdot \hat{\nabla} \hat{E}
= \frac{1}{J(t)} \hat{\nabla} \cdot \left (D_R(\hat{E},\hat{T}) A(t) \hat{\nabla} \hat{E} \right )
+ J(t)\sigma_a \left(\hat{T}^4 - \hat{E} \right ), \\
\hat{T}_t - \V{b}(t) \cdot \hat{\nabla} \hat{T}
= \frac{1}{J(t)} \hat{\nabla} \cdot \left (D_T(\hat{T}) A(t) \hat{\nabla} \hat{T} \right )
-J(t)\sigma_a \left (\hat{T}^4 - \hat{E} \right),
\end{cases}
\label{NRDE-4}
\end{equation}
where
\[
J(t) = x_\xi y_\eta - x_\eta y_\xi,\quad
\hat{\nabla} = \left [ \begin{array}{c} \frac{\partial}{\partial \xi} \\  \frac{\partial}{\partial \eta} \end{array} \right ],
\]
\[
\V{b}(t) = \frac{1}{J(t)} \left [\begin{array}{c} y_\eta x_t - x_\eta y_t \\ -y_\xi x_t + x_\xi y_t \end{array} \right ],\quad
A(t) = \frac{1}{J(t)} \left [\begin{array}{cc} x_\eta^2 + y_\eta^2 & - (x_\xi x_\eta + y_\xi y_\eta) \\
- (x_\xi x_\eta + y_\xi y_\eta) & x_\xi^2 + y_\xi^2 \end{array} \right ] .
\]
Similarly, the boundary condition (\ref{BC-1}) can be transformed into
\begin{equation}
\begin{cases}
\frac{-x_\eta}{x_\xi y_\eta} \frac{\partial \hat{E}}{\partial \xi}
+ \frac{1}{y_\eta} \frac{\partial \hat{E}}{\partial \eta}= 0,
\quad \frac{-x_\eta}{x_\xi y_\eta} \frac{\partial \hat{T}}{\partial \xi}
+ \frac{1}{y_\eta} \frac{\partial \hat{T}}{\partial \eta}= 0,
& \quad \text{ on } \eta = 0 \text{ or } \eta = 1 \\
\frac{-y_\xi}{x_\xi y_\eta} \frac{\partial \hat{T}}{\partial \eta}
+ \frac{1}{x_\xi} \frac{\partial \hat{T}}{\partial \xi}= 0,   & \quad \text{ on } \xi = 0 \text{ or } \xi = 1  \\
\frac14\hat{E} - \frac{1}{6\sigma_a}\left ( \frac{1}{x_\xi} \frac{\partial \hat{E}}{\partial \xi}
+ \frac{-y_\xi}{x_\xi y_\eta} \frac{\partial \hat{E}}{\partial \eta} \right )=1, & \quad \text{ on } \xi = 0\\
\frac14\hat{E} + \frac{1}{6\sigma_a}\left ( \frac{1}{x_\xi} \frac{\partial \hat{E}}{\partial \xi}
+ \frac{-y_\xi}{x_\xi y_\eta} \frac{\partial \hat{E}}{\partial \eta} \right )=0,& \quad \text{ on } \xi = 1
\end{cases}
\label{BC-3}
\end{equation}
where we have used the fact that $x_\eta = 0$ on $\xi = 0$ and $\xi = 1$ and $y_\xi = 0$ on $\eta = 0$ and $\eta = 1$.

The discretization of (\ref{NRDE-4}) and (\ref{BC-3}) on the rectangular reference mesh (\ref{rmesh-1})
using central finite differences is straightforward. To save space, we omit the detail of the derivation
and formulation of the FD approximation here and refer the reader to \cite[\S3.2]{HR11}.
The FD approximation of (\ref{NRDE-4}) can be expressed as
\begin{equation}
\begin{cases}
\frac{d E_h}{d t} - \V{b}_h(t) \cdot \hat{\nabla}_h E_h = \frac{1}{J_h(t)} \hat{\nabla}_h \cdot
\left (D_R(E_h, T_h) A_h(t) \hat{\nabla}_h E_h \right ) + J_h(t)\sigma_a\left (T_h^4 - E_h \right),  \\
\frac{d T_h}{d t} - \V{b}_h(t) \cdot \hat{\nabla}_h T_h = \frac{1}{J_h(t)} \hat{\nabla}_h \cdot
\left (D_T(T_h) A_h(t) \hat{\nabla}_h T_h \right ) - J_h(t)\sigma_a\left (T_h^4 - E_h \right),
\end{cases}
\label{FD-2}
\end{equation}
where $E_h$ and $T_h$ denote the FD approximations of $\hat{E}(\xi,\eta,t)$ and $\hat{T}(\xi,\eta,t)$
on the mesh (\ref{rmesh-1}), respectively.

\subsection{Linearization and predictor-corrector approximation}

Recall that the 2T model (\ref{FD-2}) has nonlinear diffusion coefficients.
Integration of nonlinear radiation diffusion equations have been studied extensively in the past;
e.g., see \cite{Knoll2003,Lowrie2007,Lowrie2004,Cai2004,PP06,Knoll1999}. Generally speaking, there are three
types of method for treating nonlinear diffusion terms \cite{Lowrie2004}, the Beaming-Warming method,
lagged diffusion, and predictor-corrector method.
For the Beaming-Warming method, the diffusion coefficient is expanded up to
linear term of $E$ and $T$ at the pervious time step and it is a second-order approximation to the diffusion equations.
For lagged diffusion, the diffusion coefficient is simply calculated with the energy and material temperature
at the pervious time step and it is only a first-order approximation.
The predictor-corrector method uses the lagged diffusion as the predictor while adding a corrector step
so it gives a second-order approximation.

In this paper, we use the predictor-corrector method for solving non-equilibrium systems since
it is comparable to the Beam-Warming method in terms
of accuracy and stability and to lagged diffusion in terms of simplicity and efficiency.
With the method, the linearized equation of (\ref{FD-2}) reads as
\begin{equation}
\begin{cases}
\frac{d E_h}{d t} - \V{b}_h(t) \cdot \hat{\nabla}_h E_h = \frac{1}{J_h(t)} \hat{\nabla}_h \cdot
\left (D_L(E_h^{*}, T_h^{*}) A_h(t) \hat{\nabla}_h E_h \right ) + J_h(t)\sigma_a\left (T_h^4 - E_h \right),
\quad t_n < t \le t_{n+1}
\\
\frac{d T_h}{d t} - \V{b}_h(t) \cdot \hat{\nabla}_h T_h = \frac{1}{J_h(t)} \hat{\nabla}_h \cdot
\left (D_T(T_h^{*}) A_h(t) \hat{\nabla}_h T_h \right ) - J_h(t)\sigma_a\left (T_h^4 - E_h \right)  ,\quad t_n < t \le t_{n+1}
\\
E_h(t_n) = E_h^n \\
T_h(t_n) = T_h^n
\end{cases}
\label{FD-3}
\end{equation}
where $E_h^n$ and $T_h^n$ are the approximations of the energy density and temperature at $t=t_n$.
During the prediction stage, $E_h^{*}$ and $T_h^{*}$ are taken as the energy density and
material temperature at $t=t_n$, i.e., $E_h^{*} = E_h^n$ and $T_h^{*} = T_h^n$.
This stage is the same as the lagged diffusion method. The solution obtained
in this stage at $t=t_{n+1}$  is used as $E_h^{*}$ and $T_h^{*}$ during the correction stage.
In both stages, the linear equation (\ref{FD-3}) is integrated with a two-stage SDIRK scheme \cite{Cas79}.
The resulting linear systems are solved by the unsymmetric multifrontal sparse LU factorization
package UMFPACK \cite{UMFPACK2004}.

\subsection{Preservation of solution positivity and cutoff}

It is known that the solutions of IBVP (\ref{2T}), (\ref{BC-1}), and (\ref{IC-2}) stay positive for all time.
Unfortunately, the scheme described in the previous subsections does not preserve the solution
positivity and the computed solutions may become zero or even negative at places. Although these
values can be very small in magnitude, they can cause nonphysical oscillations and other problems
such as not-a-number (NaN), divergence of nonlinear iterations, too small time steps, and
even early blowup of computation \cite{yuan2009}. We employ here a cutoff strategy, i.e.,
replace solution values that are below a positive threshold by the threshold. 
Unfortunately, no theory exists so far on how to choose such a threshold.
An empirical formula is $30/((M-1)(N-1))$ (see Table~\ref{mytab}) which has been found to
work well for the examples we consider.  Noticeably, Lu et al. \cite{LuHuVV2012} show
that the cutoff procedure can retain accuracy, convergence order, and stability of finite difference schemes
for linear or nonlinear parabolic PDEs.

\begin{table}[htbp]
 \centering\small
  \caption{The values of the cutoff threshold defined as $\frac{30}{(M-1)(N-1)}$ for various mesh sizes.}
 \medskip
 \begin{tabular}{|c|c|c|c|c|}
 \hline
Mesh $(M\times N)$ & 41$\times$41 \quad& 61$\times$61 \quad& 81$\times$81 \quad& 121$\times$121 \\
\hline
Cutoff threshold  & 1.87e-2 & 8.30e-3 & 4.70e-3  & 2.10e-3 \\
\hline
\end{tabular}
\label{mytab}
\end{table}


\subsection{The MMPDE approach of mesh movement}
\label{SEC:mmesh}

The MMPDE approach \cite{HRR94b,HR97b,HR11} is used here to generate the adaptive moving mesh.
The main idea of the approach is to generate the moving mesh as the image of a fixed, reference mesh
under a time coordinate transformation. Such a coordinate transformation is determined as the solution
of an MMPDE which in turn is defined as the gradient flow equation of a meshing functional.
We use a meshing functional formulated in terms of the inverse coordinate transformation $\xi = \xi(x, y, t)$
and $\eta = \eta(x, y, t)$ and based on mesh equidistribution and alignment.
A monitor function that is symmetric and uniformly positive definite at each point of the domain
is used in the functional to provide the information for the size, shape, and orientation of the mesh elements.
Denote the Hessian of the energy density by
\[
H = \left [ \begin{array}{cc} E_{xx} & E_{xy} \\ E_{xy} & E_{yy} \end{array}\right ] .
\]
Given its eigen-decomposition $H = Q\mbox{diag}(\lambda_1, \lambda_2) Q^T$, we define
$|H| = Q\mbox{diag}(|\lambda_1|, |\lambda_2|) Q^T$.
Then the monitor function is chosen as
\begin{equation}
\M =\textrm{det}(\alpha I + |H|)^{-\frac14}\left [ \frac{}{} \alpha I + |H|\right ],
\label{M-1}
\end{equation}
which is known to be optimal for the $H^1$ norm of the error of linear interpolation \cite{HR11,HS03}.
Here, $\alpha > 0$ is the regularization parameter defined through the equation
\[
\int_{\Omega} \det(\M(\alpha))^{\frac 1 2} d x d y = 2 \int_{\Omega} \det(\M(0))^{\frac 1 2} d x d y,
\]
where $\M(0)$ denotes the monitor function (\ref{M-1}) with $\alpha = 0$.
In practical computation, the Hessian of $E$ is unknown. It is replaced by an approximation
based on $E_h$ (see \S\ref{SEC:soln-procedure} for a more detailed description).
The meshing functional is given by
\begin{align}
I[\xi, \eta] & = 0.1 \int_\Omega \det(\M)^{\frac 1 2} \left (\nabla \xi^T \M^{-1} \nabla \xi
+ \nabla \eta^T \M^{-1} \nabla \eta \right )^2 d x d y
 +  3.2 \int_\Omega \frac{\det(\M)^{\frac 1 2}}{\left (J \det(\M)^{\frac 1 2}\right )^2} d x d y,
\label{functional-1}
\end{align}
where $J = x_\xi y_\eta - x_\eta y_\xi = 1/(\xi_x \eta_y - \xi_y \eta_x)$ is the Jacobian of
the coordinate transformation. This functional is proposed in \cite{Hua01b}
to control mesh equidistribution and alignment.

The MMPDE is defined as the gradient flow equation of the meshing functional, i.e.,
\begin{equation}
\frac{\partial \xi}{\partial t} = - \frac{1}{\tau} \frac{\delta I}{\delta \xi}, \quad
\frac{\partial \eta}{\partial t} = - \frac{1}{\tau} \frac{\delta I}{\delta \eta},
\label{MPDE-1}
\end{equation}
where $\tau$ is a parameter used to control the response of the mesh movement to
the change in the monitor function
and ${\delta I}/{\delta \xi}$ and ${\delta I}/{\delta \eta}$ are the functional derivatives
of $I[\xi, \eta]$. It is not difficult to find that
\begin{align}
 & \frac{\delta I}{\delta \xi} = - 4 \theta \nabla \cdot \left (\det(\M)^{\frac 1 2} \beta \M^{-1} \nabla \xi \right )
 - 8 (1-2\theta) \nabla \cdot \left (\frac{1}{J \det(\M)^{\frac 1 2}} \left [ \begin{array}{r} \eta_y \\ - \eta_x \end{array}\right ]
 \right ),
 \label{funder-1}
 \\
& \frac{\delta I}{\delta \eta} = - 4 \theta \nabla \cdot \left (\det(\M)^{\frac 1 2} \beta \M^{-1} \nabla \eta \right )
 - 8 (1-2\theta) \nabla \cdot \left (\frac{1}{J \det(\M)^{\frac 1 2}} \left [ \begin{array}{r} - \xi_y \\ \xi_x \end{array}\right ]
 \right ) ,
 \label{funder-2}
\end{align}
where
\begin{equation}
\beta = \nabla \xi^T \M^{-1} \nabla \xi + \nabla \eta^T \M^{-1} \nabla \eta.
\label{beta-1}
\end{equation}
By interchanging the roles of independent and dependent variables and after some straightforward but lengthy 
derivations (e.g., see \cite[Chapter 6]{HR11}), we can rewrite the above equation into
\begin{equation}
\frac{\partial}{\partial t} \left [ \begin{array}{c} x \\ y \end{array} \right ]
= \frac{1}{\tau} \left ( A_{11} \frac{\partial^2}{\partial \xi^2} + (A_{12} + A_{21}) \frac{\partial^2}{\partial \xi\partial \eta}
+ A_{22} \frac{\partial^2}{\partial \eta^2} + I_d \; b_1 \frac{\partial}{\partial \xi} +  I_d \; b_2 \frac{\partial}{\partial \eta}
\right ) \left [ \begin{array}{c} x \\ y \end{array} \right ] ,
\label{MPDE-3}
\end{equation}
where $I_d$ is the $2$-by-$2$ identity matrix and the coefficient $A_{ij}$, $b_1$, and $b_2$
can be found in \cite[Chapter~6]{HR11}.

The moving mesh equation (\ref{MPDE-3}) is supplemented with the one-dimensional version of the MMPDE
for the adaptation of boundary points (cf. \cite{Hua01b}). They are discretized in space using central
finite differences and in time by the backward Euler method with coefficients $A_{ij}$ and $b_i$ calculated
at the previous time step. The resulting algebraic systems are solved using the sparse matrix solver
UMFPACK \cite{UMFPACK2004}.

\subsection{The solution procedure}
\label{SEC:soln-procedure}

We now describe the overall solution procedure of the moving mesh FD method.
Assume that the physical solutions $E^n$ and $T^n$,
the mesh $(x^n, y^n)$, and the time step size $\Delta t_n$ are given at $t=t_n$.

\begin{itemize}
\item[] {\bf Step 1. The moving mesh step.}
The monitor function (\ref{M-1}) is computed using $E^n$ and $(x^n, y^n)$ and
smoothed using several sweeps of a low-pass filter. The Hessian of the energy density
used in (\ref{M-1}) is replaced by an approximation obtained using least squares fitting.
More specifically, at any mesh point a local quadratic polynomial is constructed by least squares fitting
of the nodal values of $E^n$ at neighboring mesh points. The approximate Hessian at the given mesh point
is then obtained by differentiating the quadratic polynomial twice.
After the monitor function has been obtained, the mesh equation (\ref{MPDE-3}) is integrated
from $t_n$ to $t_{n+1} = t_n + \Delta t_n$
for the new mesh $(x^{n+1}, y^{n+1})$.

\item[] {\bf Step 2. The predictor of the physical PDE solving step.}
The physical PDE (\ref{FD-2}) is integrated at $t_n$ using the predictor-corrector scheme (\ref{FD-3})
with $E_h^* = E_h^n$ and $T_h^* = T_h^n$. During the integration, the mesh is considered to
move linearly in time, viz.,
\begin{align}
& x(t) = \frac{t-t_n}{\Delta t_n}x^{n+1} + \frac{t_n+\Delta t_n -t}{\Delta t_n}x^n,\qquad
y(t) = \frac{t-t_n}{\Delta t_n}y^{n+1} + \frac{t_n+\Delta t_n -t}{\Delta t_n}y^n .
\label{mesh-2}
\end{align}

\item[]
{\bf Step 3. The corrector of the physical PDE solving step.} The physical PDE (\ref{FD-2})
is integrated from $t_n$ to $t_{n+1}$
using the predictor-corrector scheme (\ref{FD-3}) with $E_h^*$ and  $T_h^*$ being taken as
the solutions obtained in Step 2.
\end{itemize}

\section{Numerical tests}
\label{SEC:numerics}
\setcounter{equation}{0} \setcounter{figure}{0}
\setcounter{table}{0}
\setcounter{secnumdepth}{4}

\begin{figure}
\begin{center}
\hbox{
\hspace{1in}
\begin{minipage}[t]{2.5in}
\centerline{\scriptsize (a):  Example~\ref{Example4.1}}
\includegraphics[scale=0.35]{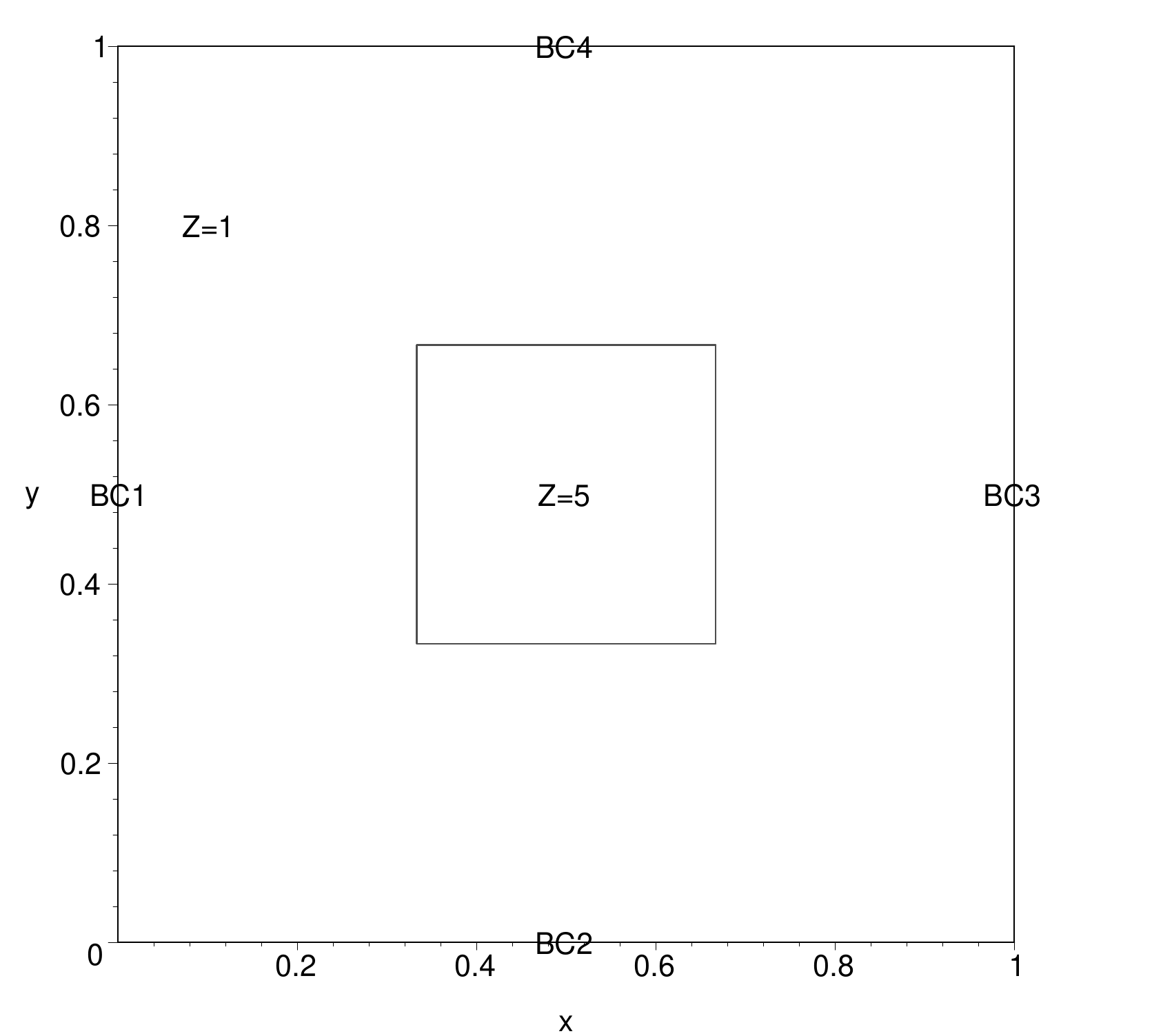}
\end{minipage}
\begin{minipage}[t]{2.5in}
\centerline{\scriptsize (b):  Example~\ref{Example4.2}}
\includegraphics[scale=0.35]{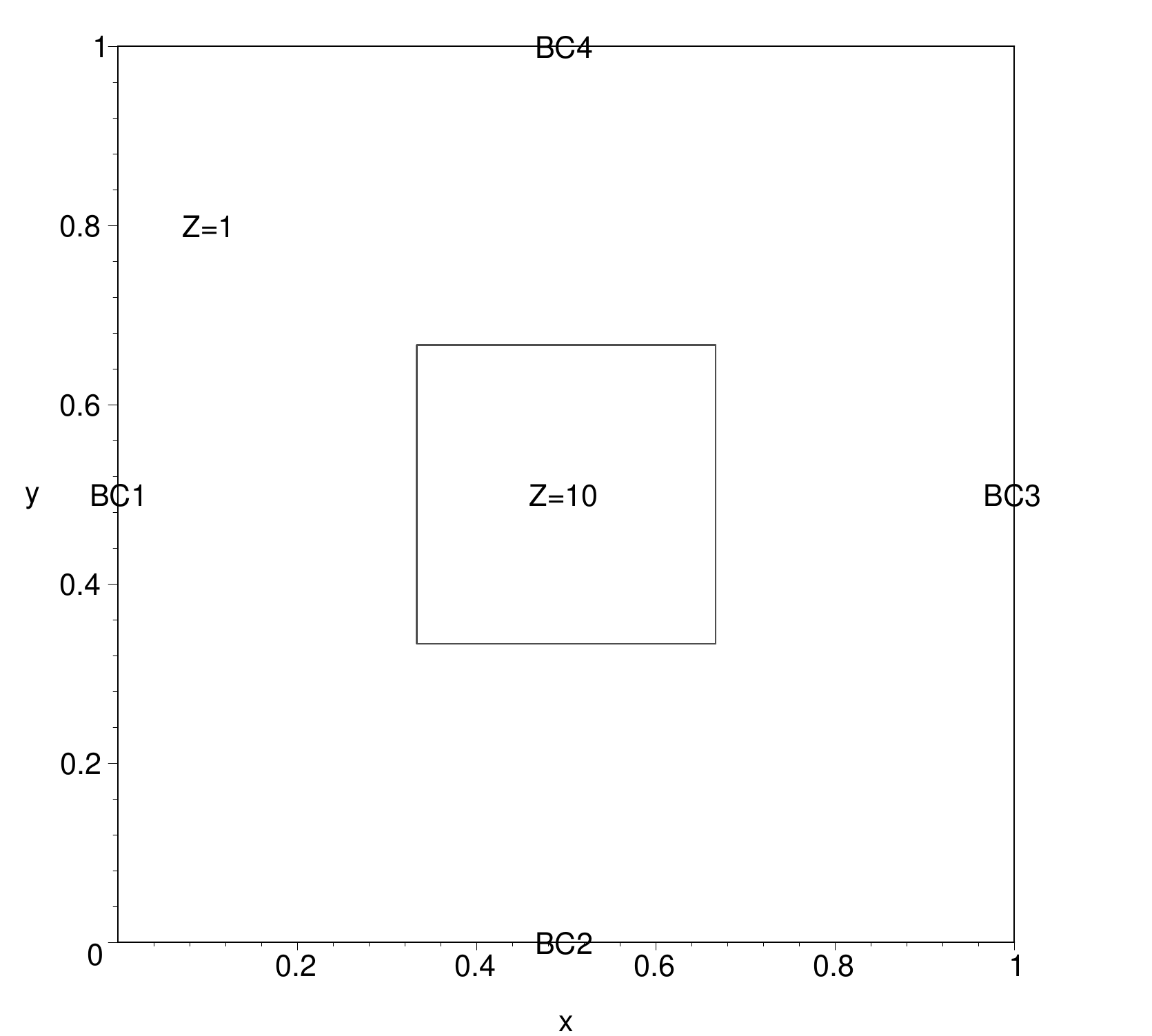}
\end{minipage}
}
\end{center}
\caption{Material configuration for Examples~\ref{Example4.1} and \ref{Example4.2}.}
\label{M1}
\end{figure}

In this section we present numerical results obtained by the moving mesh FD method
described in the previous section for three examples of multi-material radiation diffusion.
The material configuration is given in Fig.~\ref{M1} for the first two examples
and in Fig.~\ref{BCD3} for the third one (which also has a slightly different boundary
condition than (\ref{BC-1})). In the results, MM, MM1, and MM2 stand for moving mesh,
one-level moving mesh, and two-level moving mesh, respectively.

\begin{figure}
\begin{center}
\hbox{
\begin{minipage}[t]{2.3in}
\centerline{\scriptsize (a): t=1.0}
\includegraphics[width=2.3in]{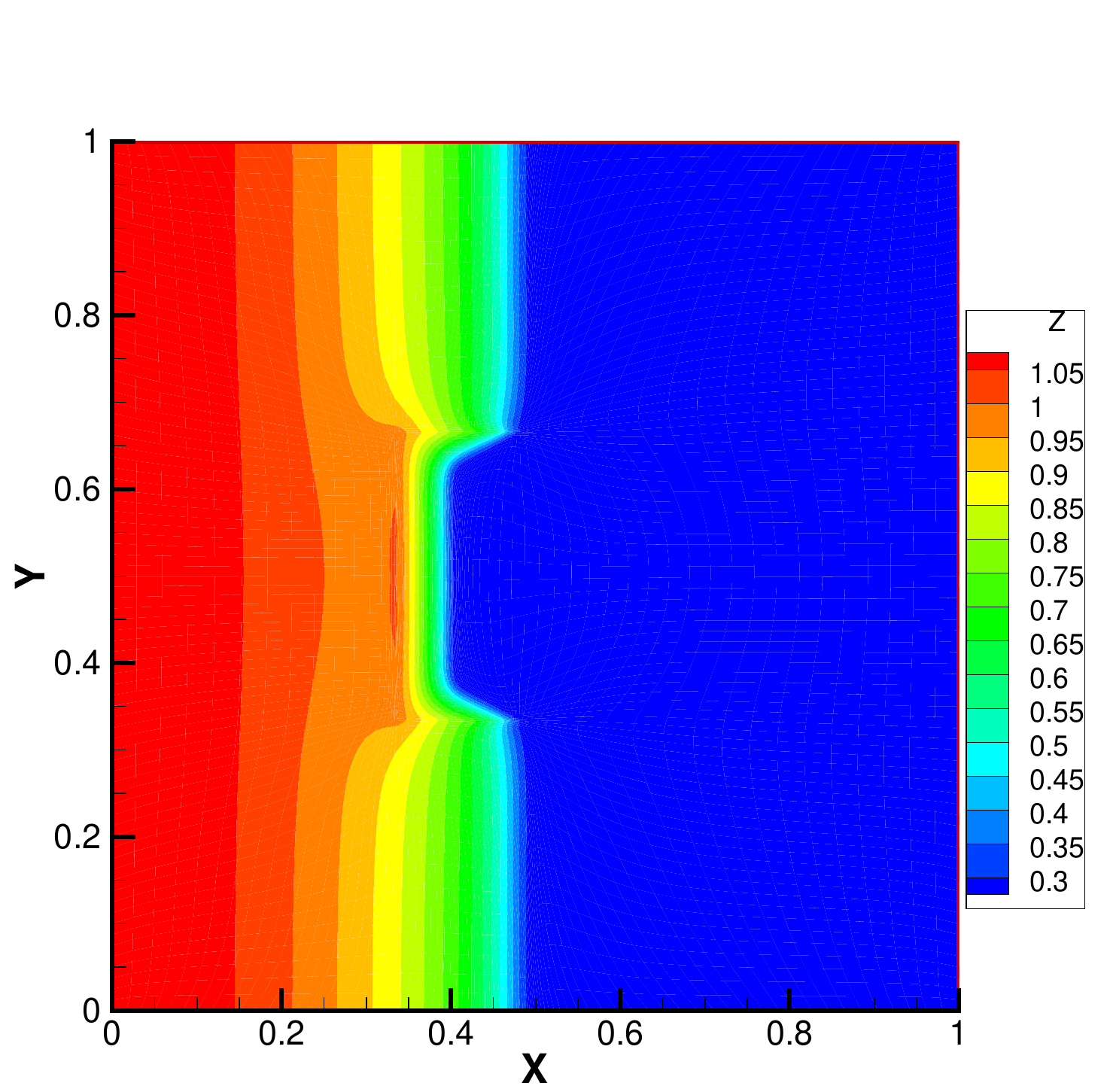}
\end{minipage}
\begin{minipage}[t]{2.3in}
\centerline{\scriptsize (b): t=1.5}
\includegraphics[width=2.3in]{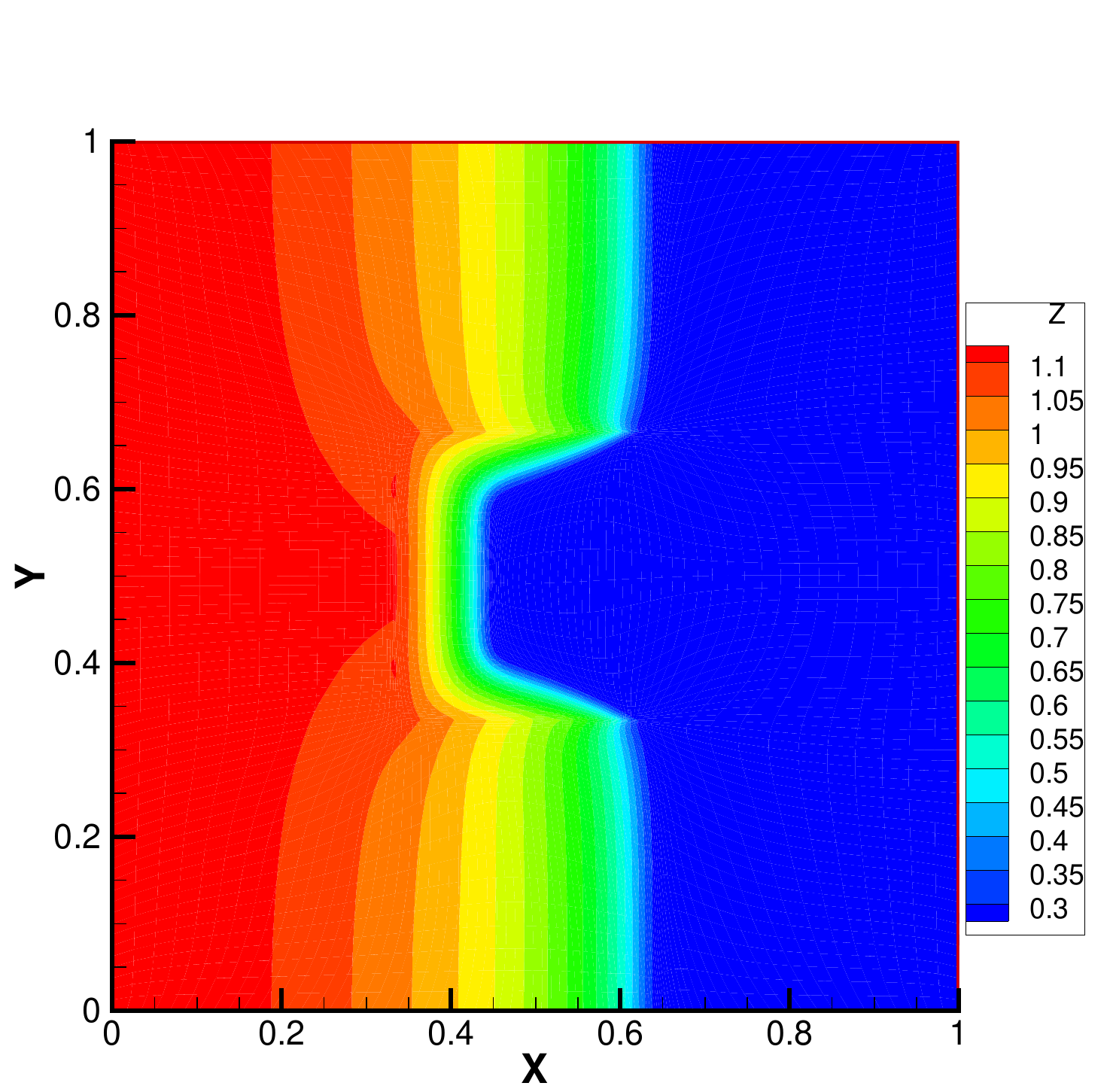}
\end{minipage}
\begin{minipage}[t]{2.3in}
\centerline{\scriptsize (c): t=2.0}
\includegraphics[width=2.3in]{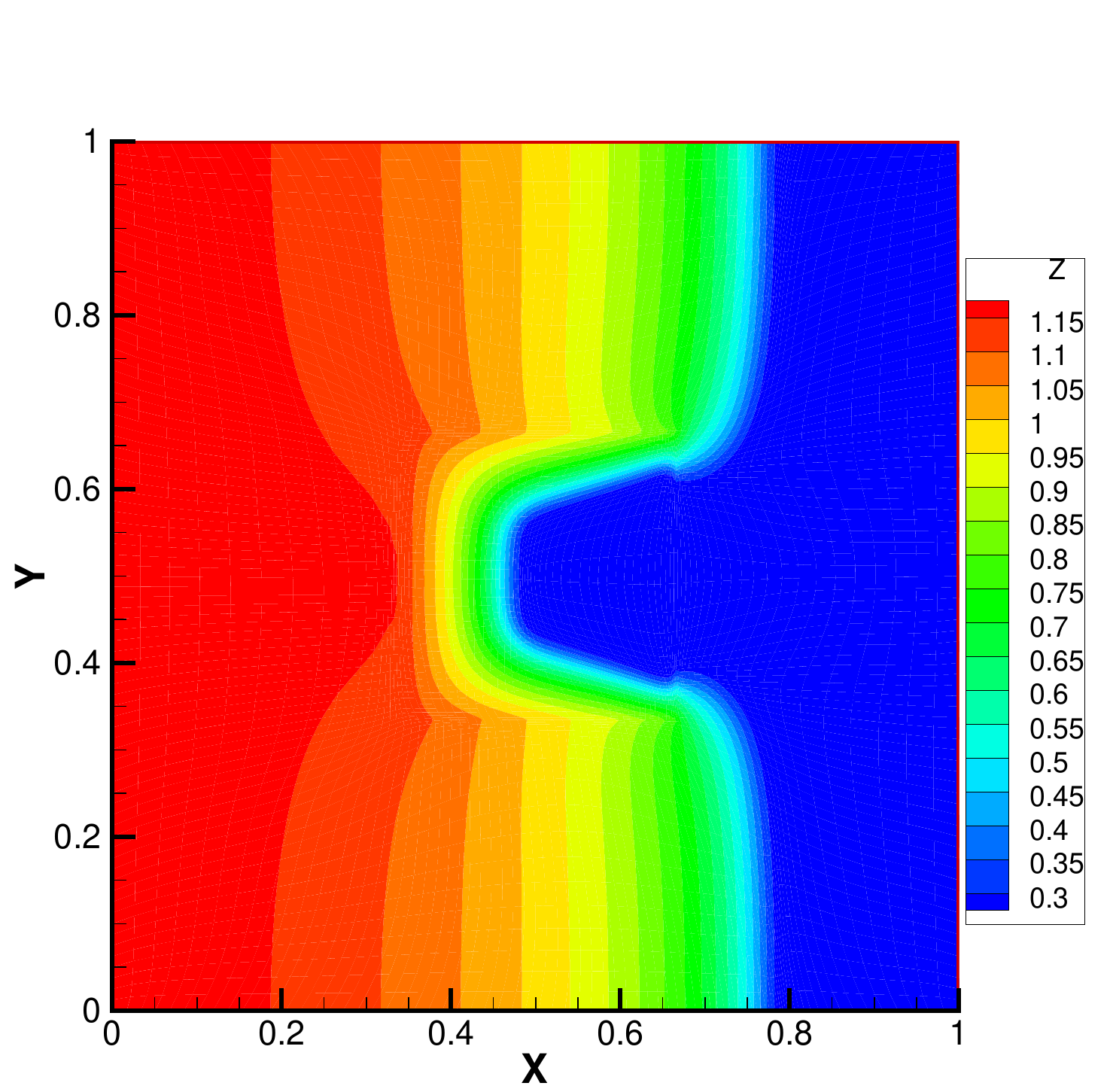}
\end{minipage}
}
\vspace{5mm}
\hbox{
\begin{minipage}[t]{2.3in}
\centerline{\scriptsize (d): t=2.4}
\includegraphics[width=2.3in]{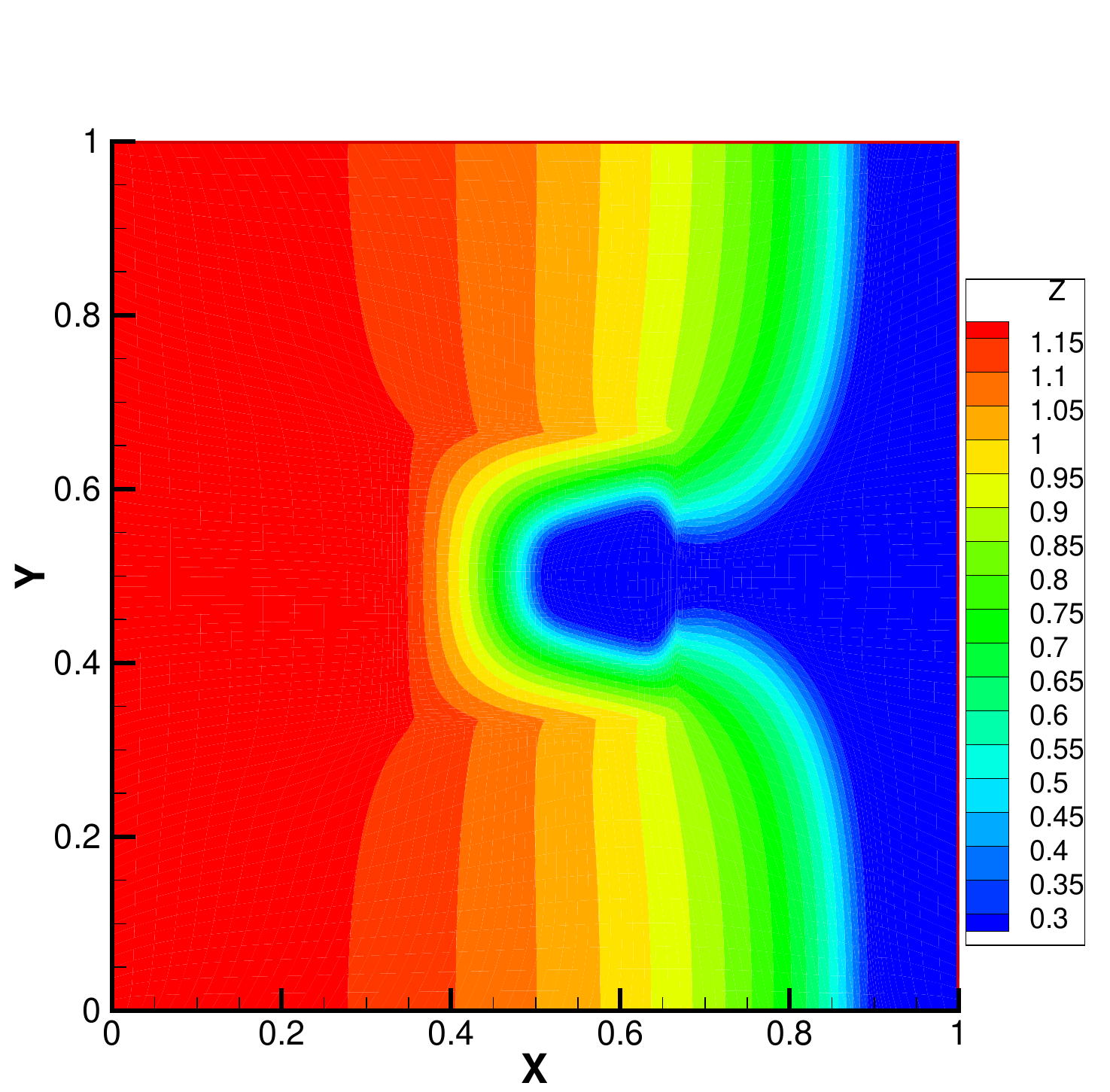}
\end{minipage}
\begin{minipage}[t]{2.3in}
\centerline{\scriptsize (e): t=2.8}
\includegraphics[width=2.3in]{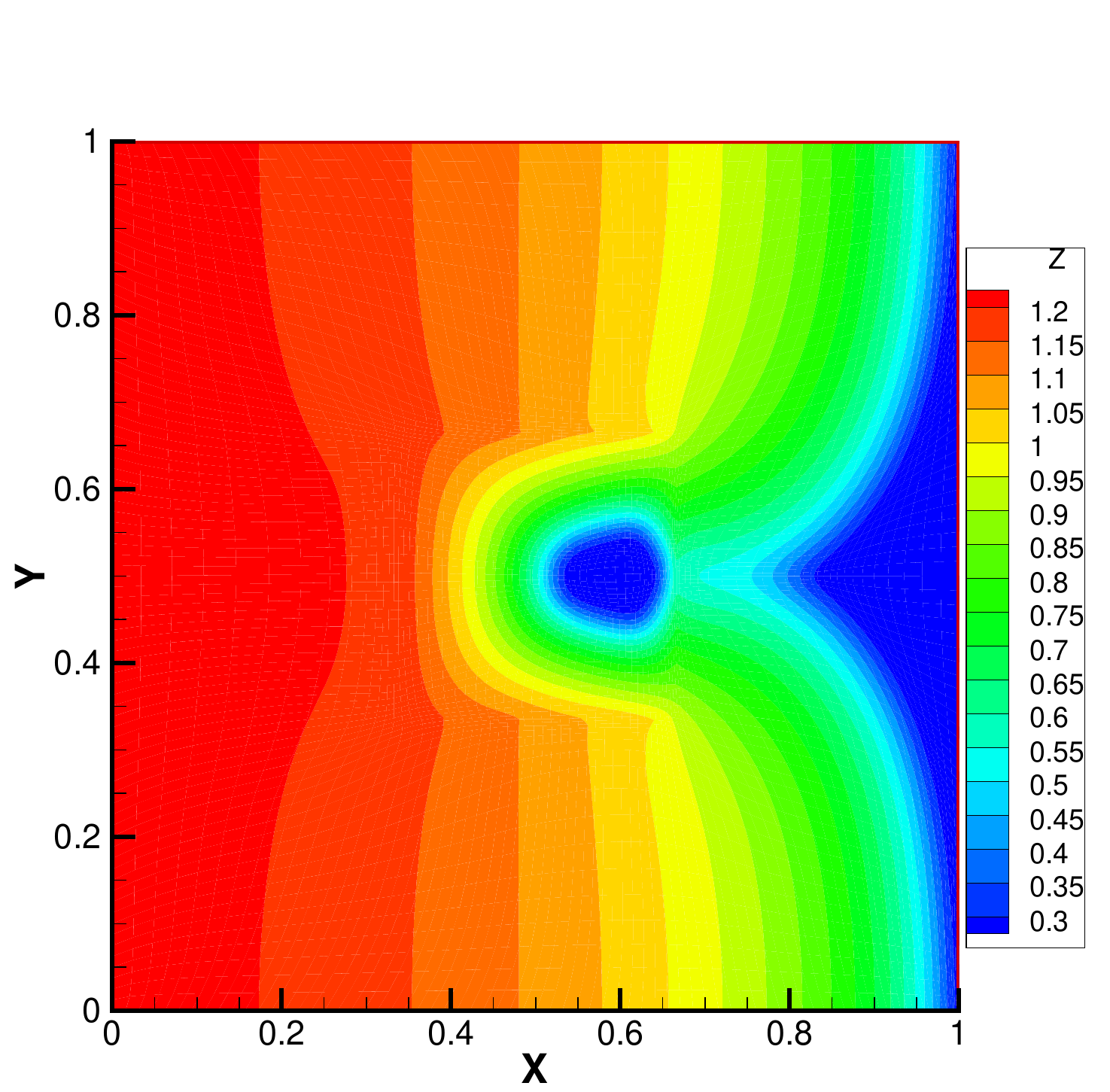}
\end{minipage}
\begin{minipage}[t]{2.3in}
\centerline{\scriptsize (f): t=3.0}
\includegraphics[width=2.3in]{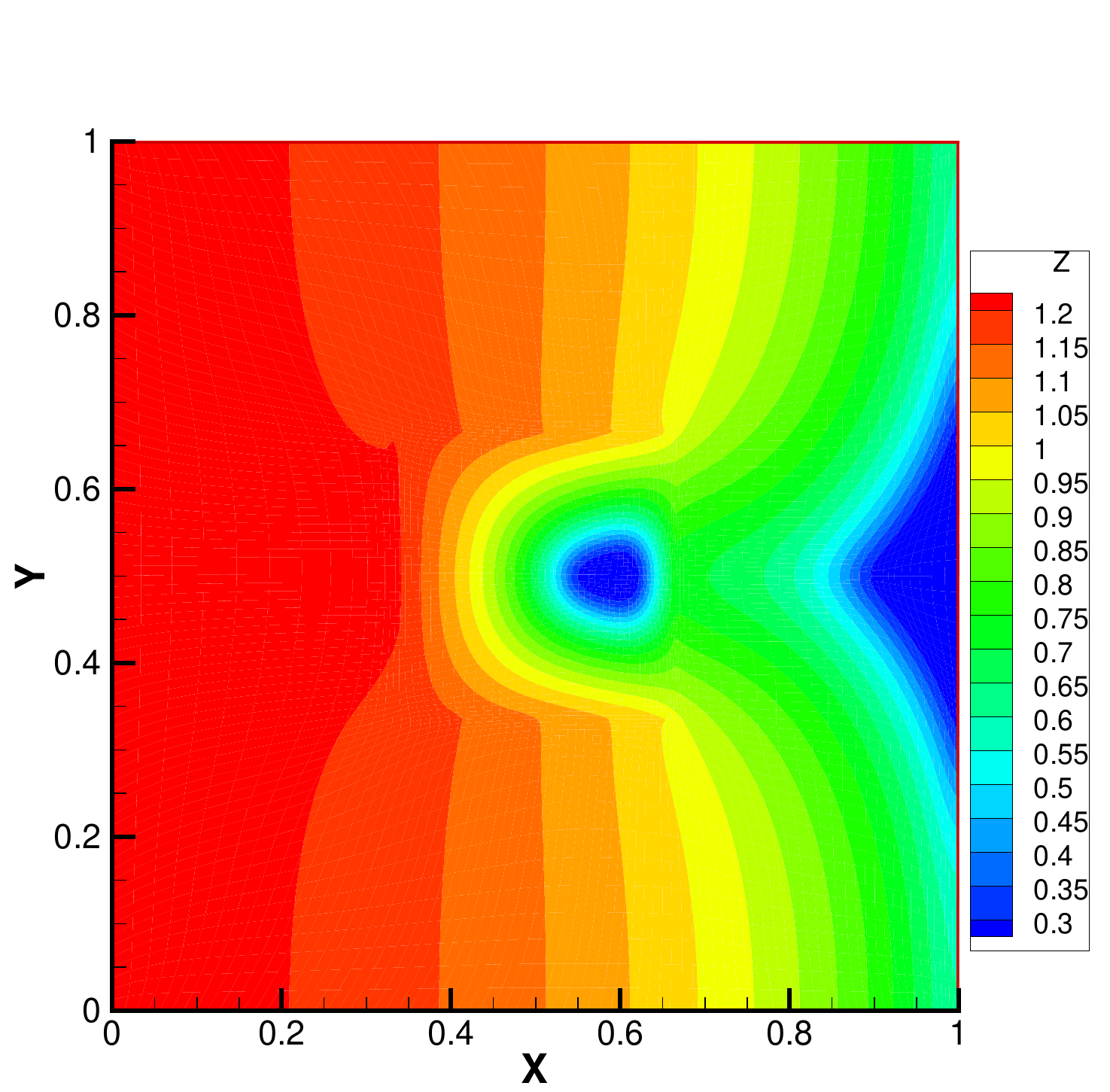}
\end{minipage}
}
\end{center}
\caption{Example~\ref{Example4.1}. The computed solution at $t=1.0, 1.5, 2.0, 2.4, 2.8, 3.0$ is obtained
with a moving mesh of $81\times 81$.}
\label{T5}
\end{figure}

\begin{figure}
\begin{center}
\hbox{
\begin{minipage}[t]{2.3in}
\centerline{\scriptsize (a): t=1.0}
\includegraphics[width=2.3in]{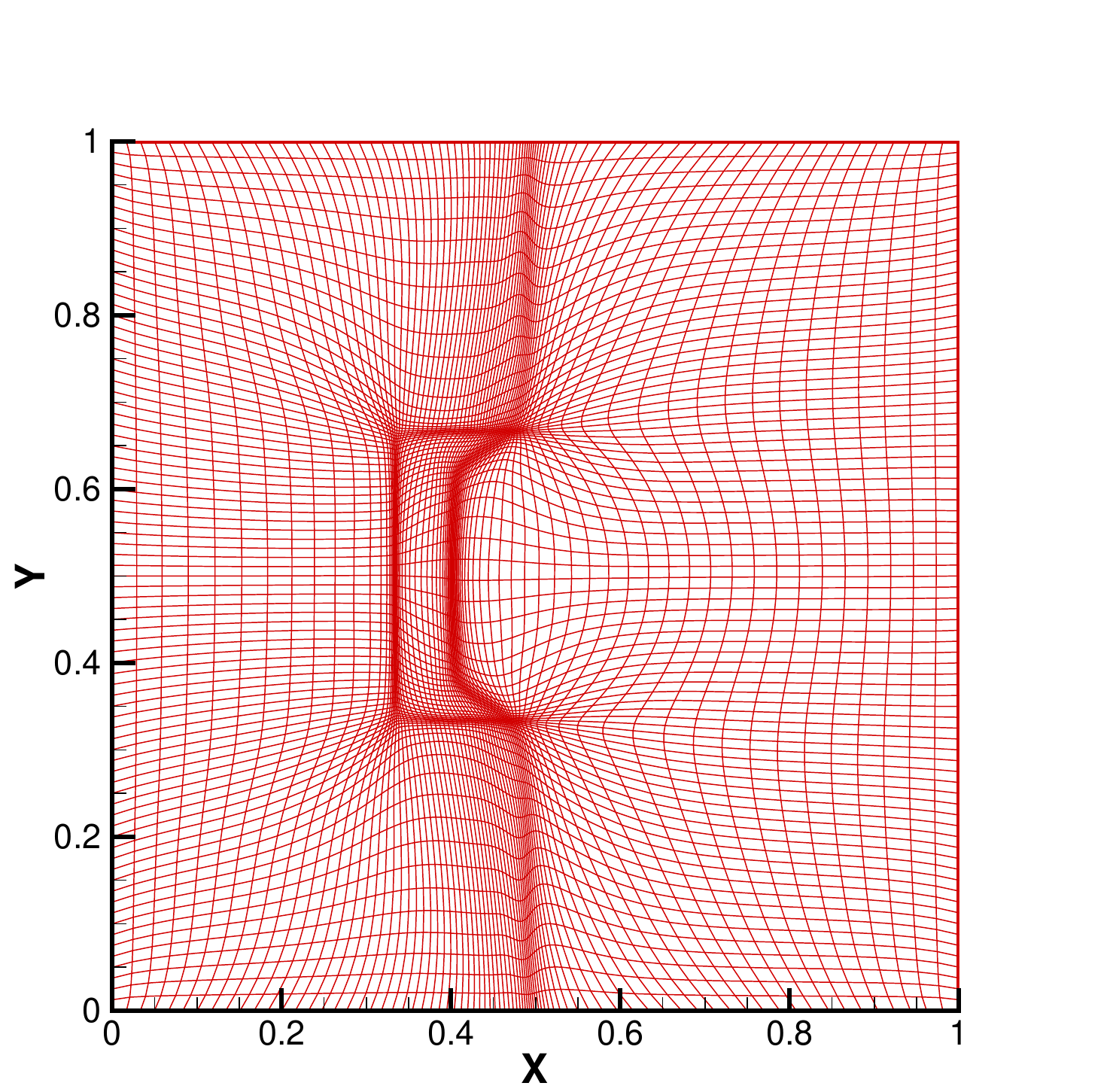}
\end{minipage}
\begin{minipage}[t]{2.3in}
\centerline{\scriptsize (b): t=1.5}
\includegraphics[width=2.3in]{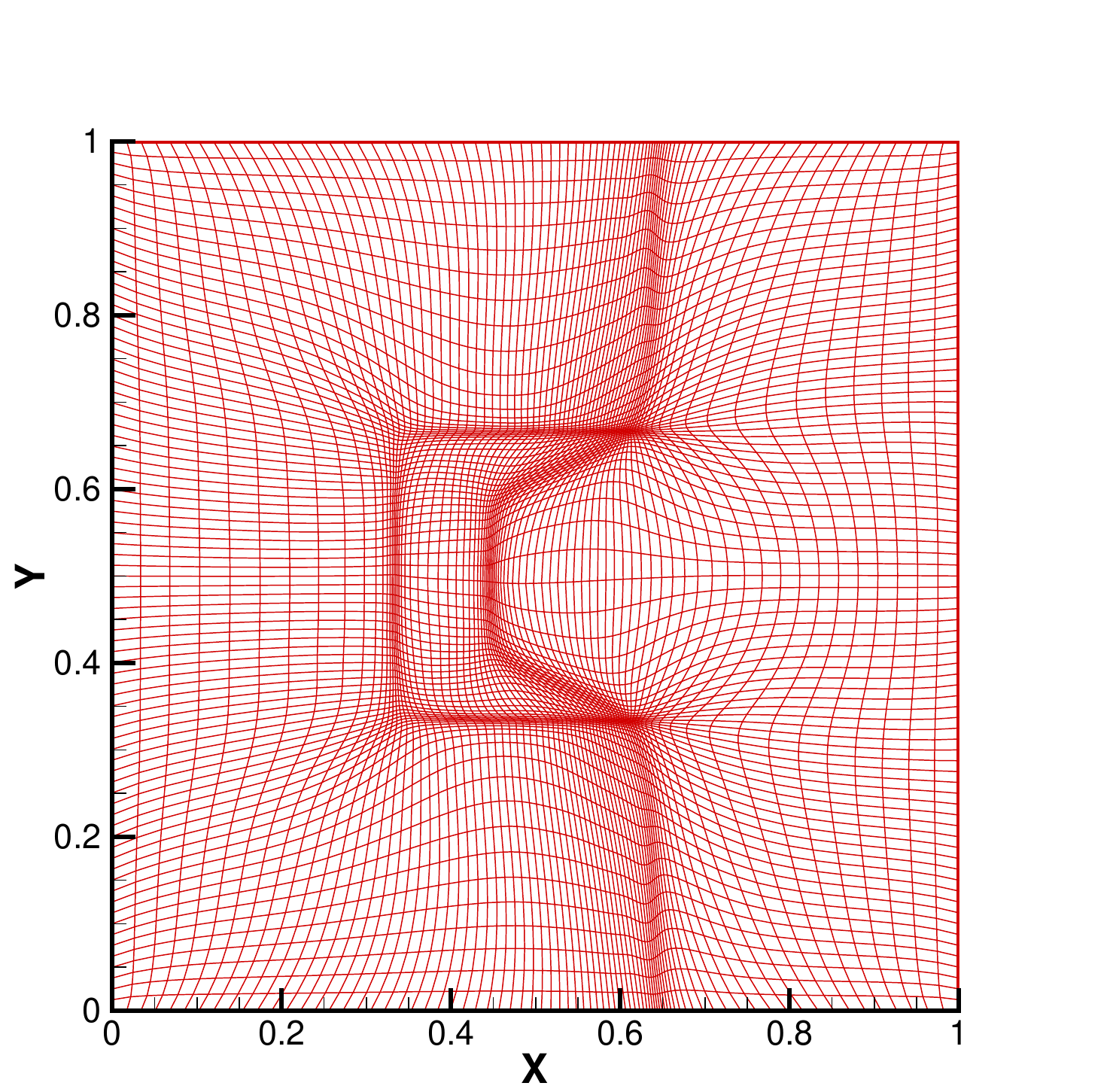}
\end{minipage}
\begin{minipage}[t]{2.3in}
\centerline{\scriptsize (c): t=2.0}
\includegraphics[width=2.3in]{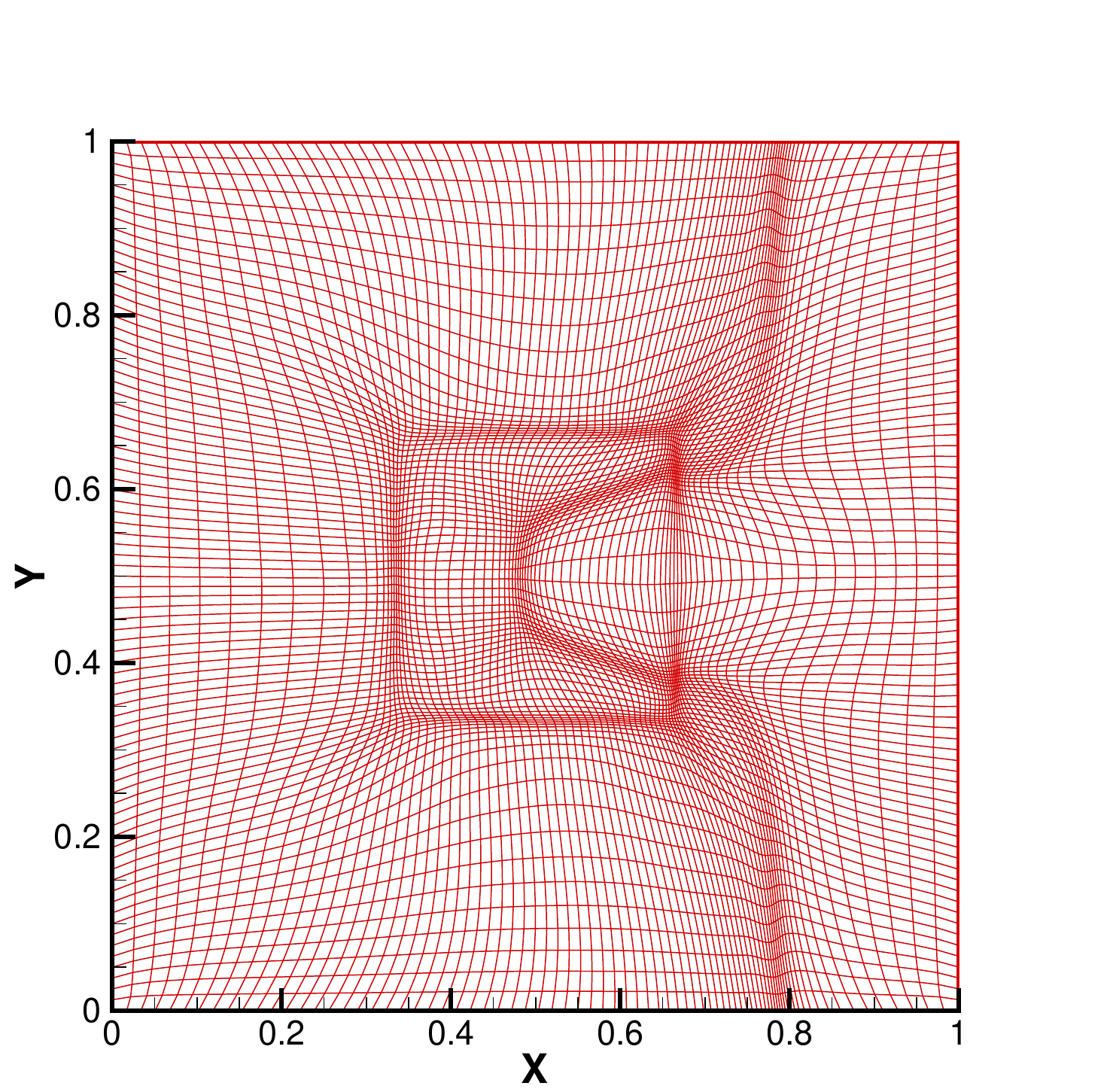}
\end{minipage}
}
\vspace{5mm}
\hbox{
\begin{minipage}[t]{2.3in}
\centerline{\scriptsize (d): t=2.4}
\includegraphics[width=2.3in]{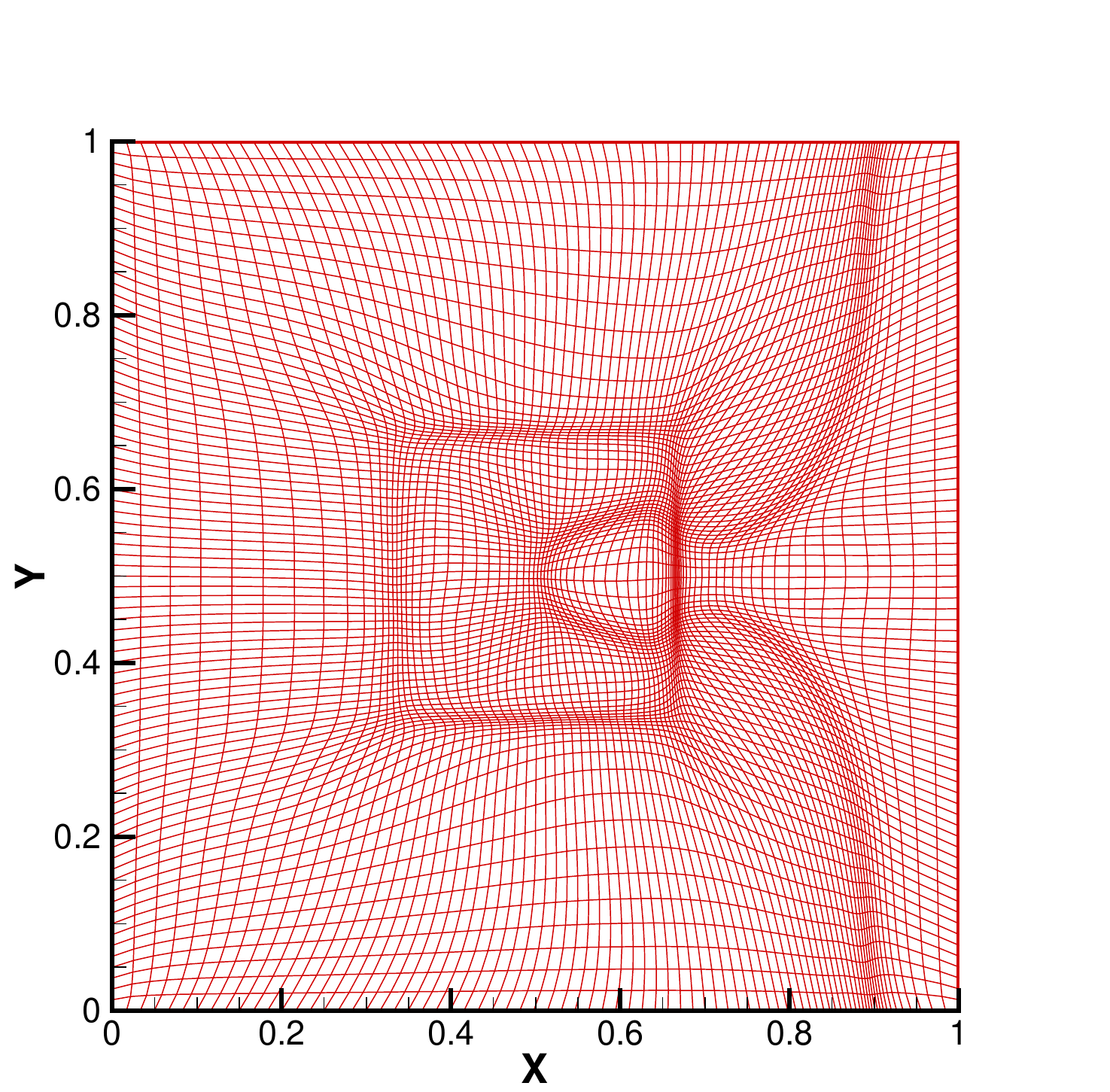}
\end{minipage}
\begin{minipage}[t]{2.3in}
\centerline{\scriptsize (e): t=2.8}
\includegraphics[width=2.3in]{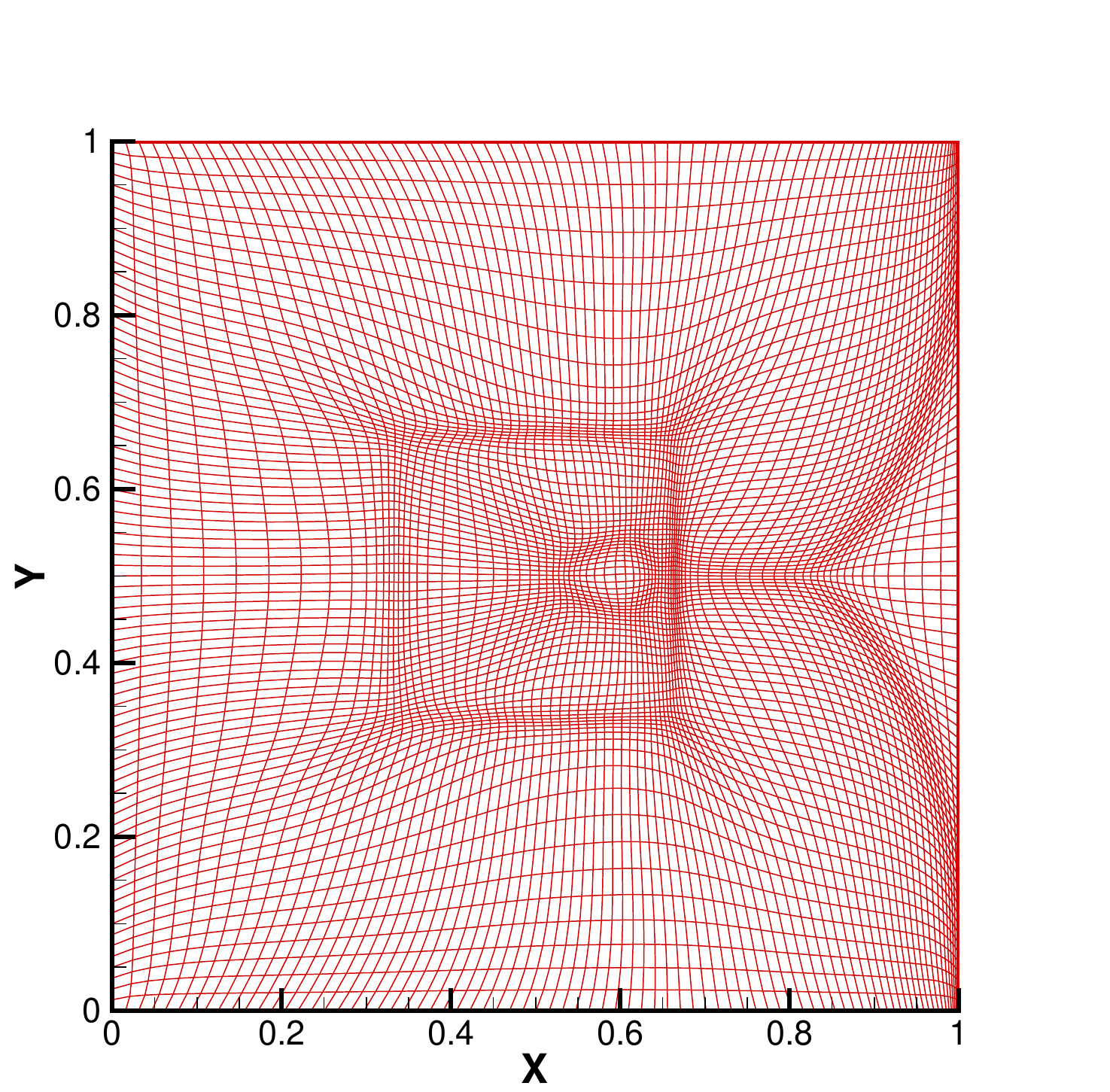}
\end{minipage}
\begin{minipage}[t]{2.3in}
\centerline{\scriptsize (f): t=3.0}
\includegraphics[width=2.3in]{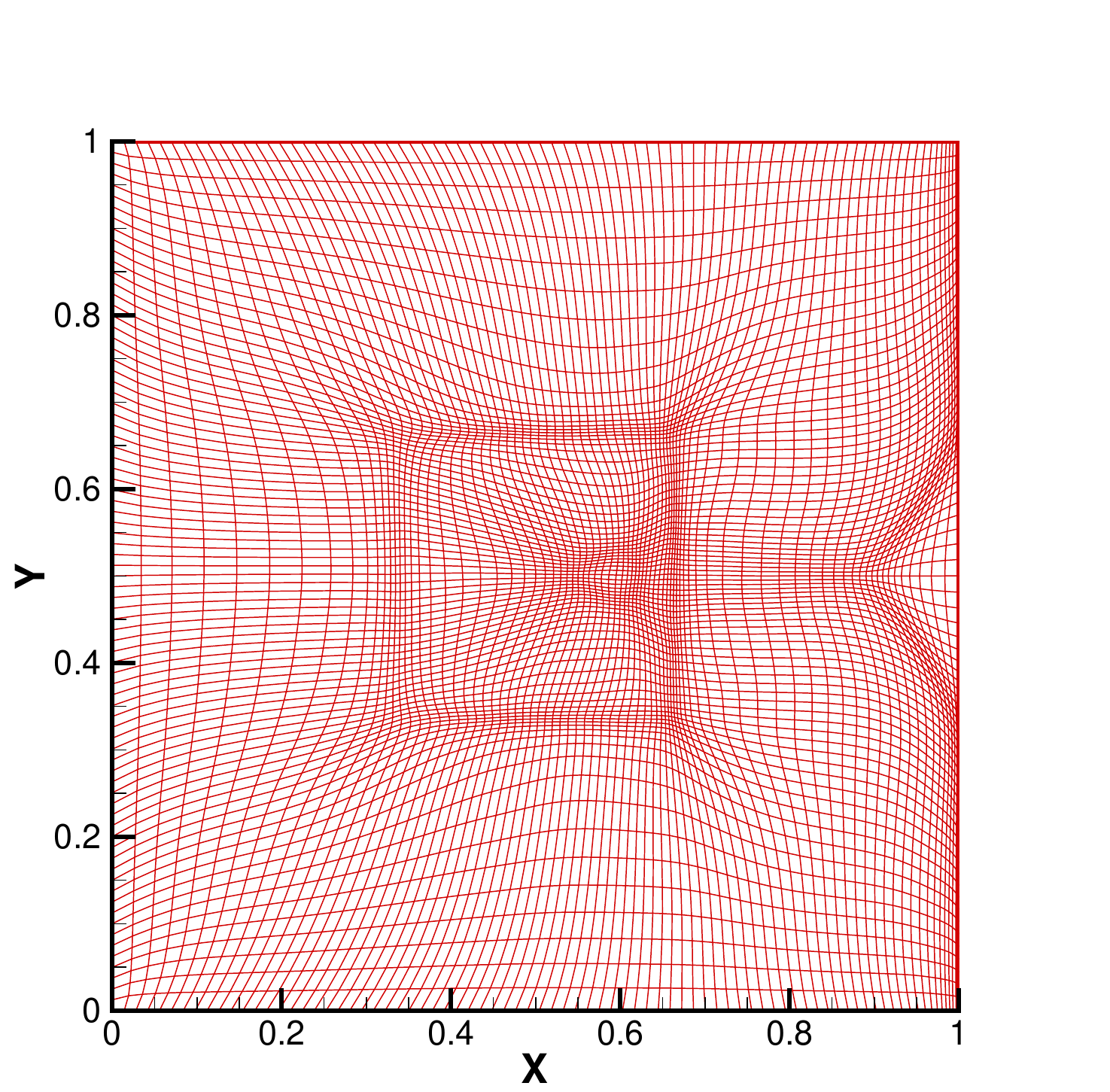}
\end{minipage}
}
\end{center}
\caption{Example~\ref{Example4.1} The moving mesh ($81\times 81$) is shown at $t=1.0, 1.5, 2.0, 2.4, 2.8, 3.0$.}
\label{T6}
\end{figure}

\begin{figure}
\begin{center}
\hbox{
\hspace{1in}
\begin{minipage}[t]{2.0in}
\centerline{\scriptsize (a):  with MM at $t=1.0$}
\includegraphics[width=2.0in]{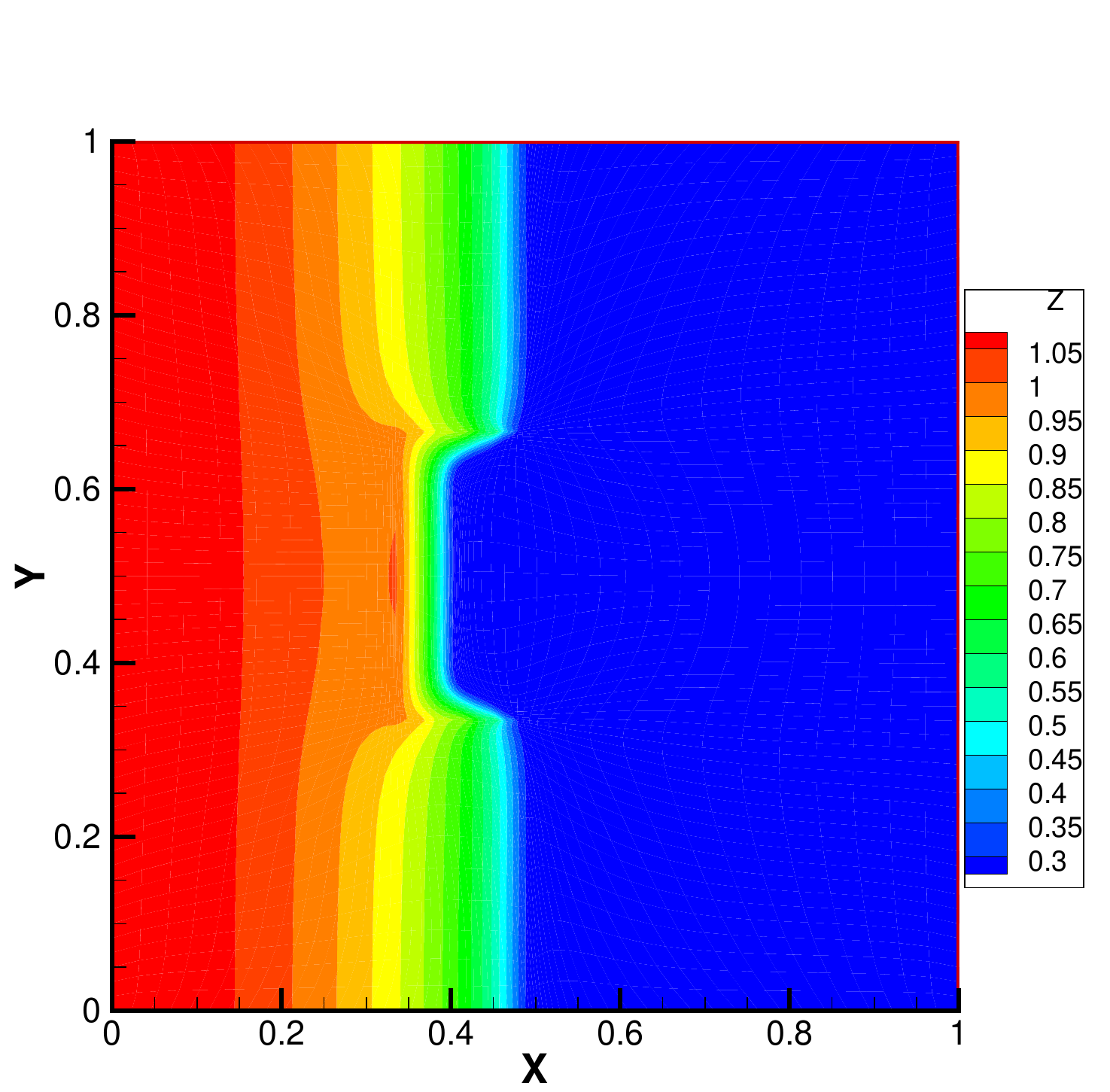}
\end{minipage}
\hspace{0.5in}
\begin{minipage}[t]{2.0in}
\centerline{\scriptsize (b):  with UM at $t=1.0$}
\includegraphics[width=2.0in]{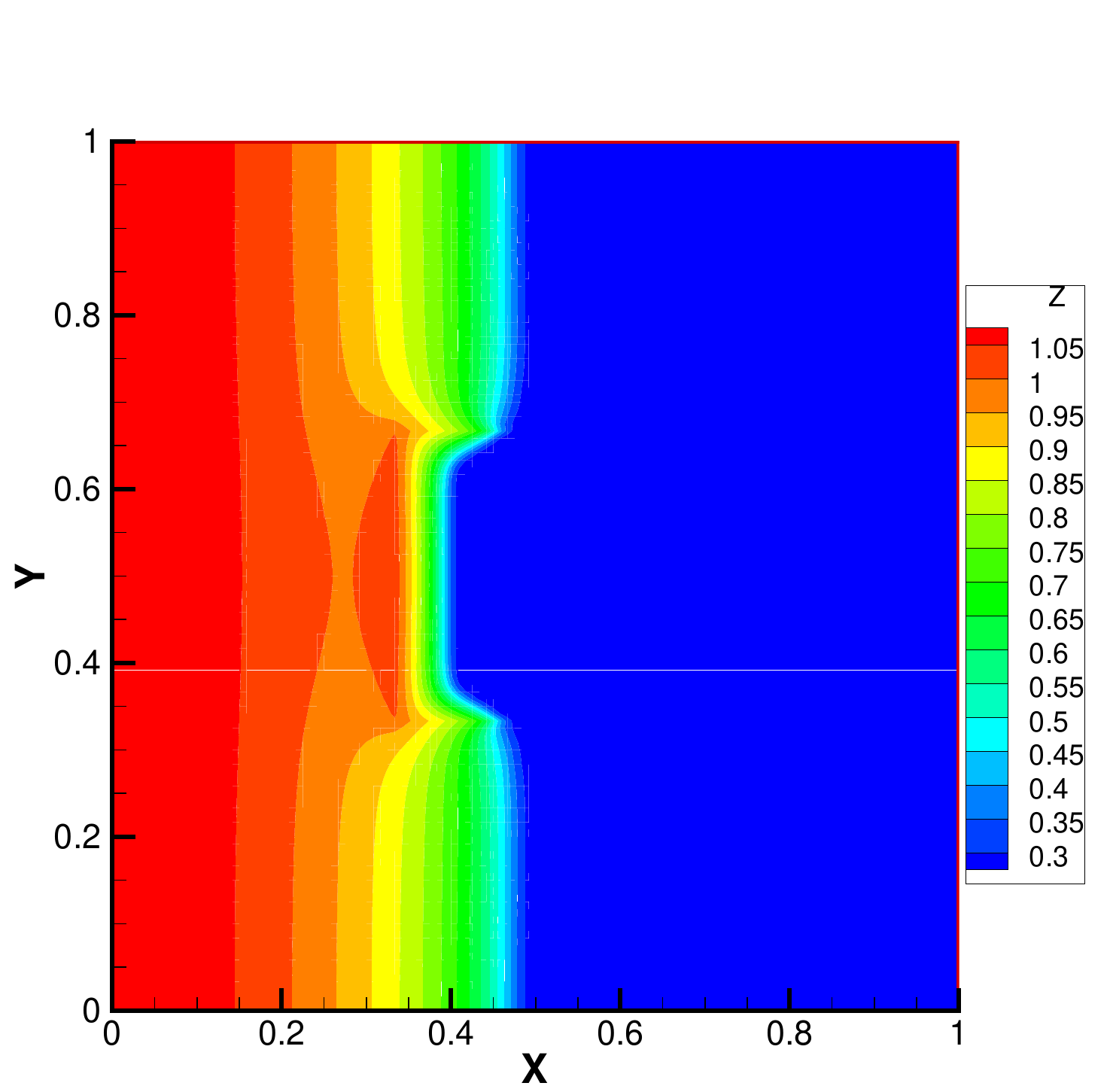}
\end{minipage}
}
\hbox{
\hspace{1in}
\begin{minipage}[t]{2.0in}
\centerline{\scriptsize (c):  with MM at $t=2.0$}
\includegraphics[width=2.0in]{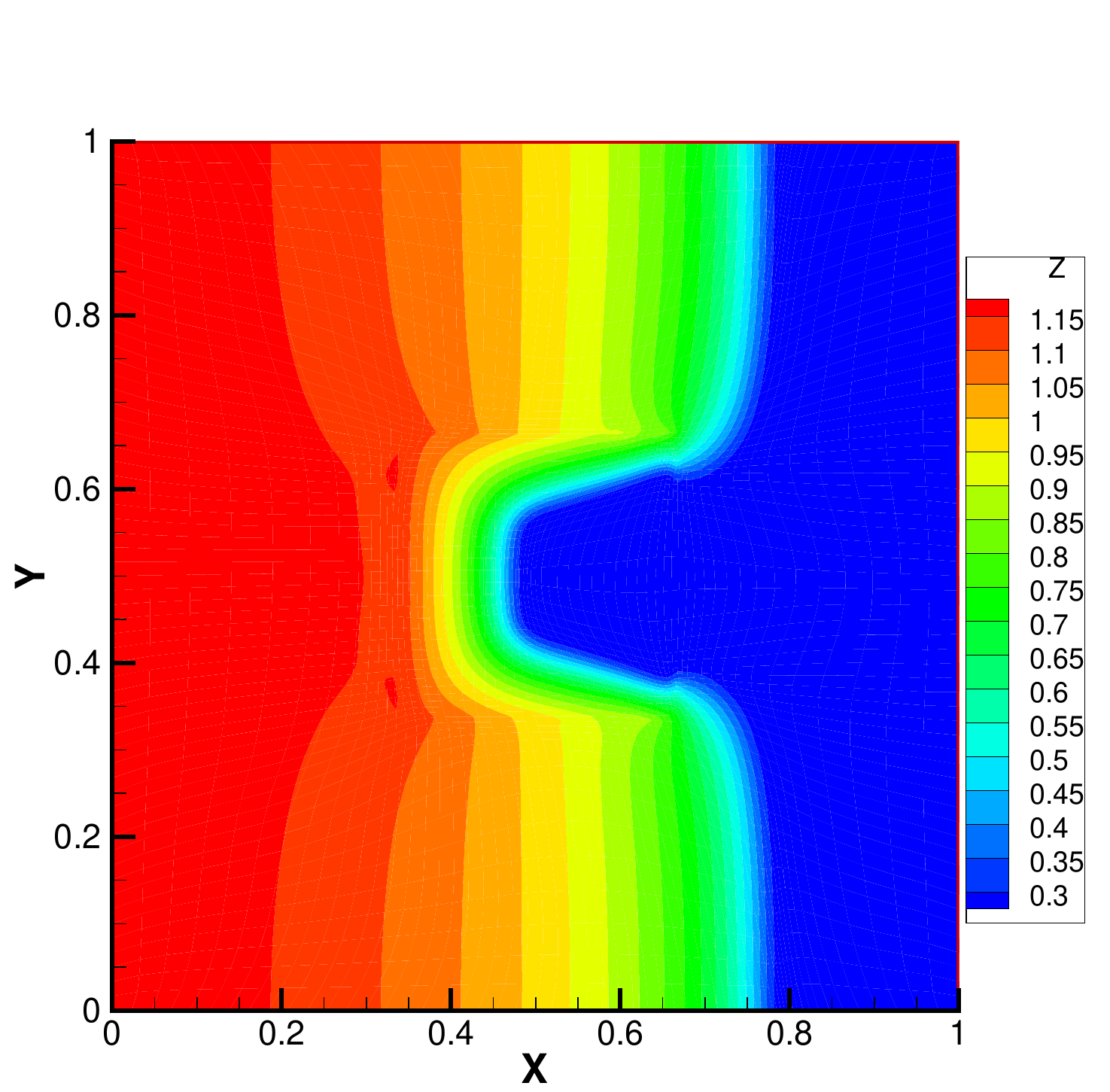}
\end{minipage}
\hspace{0.5in}
\begin{minipage}[t]{2.0in}
\centerline{\scriptsize (d):  with UM at $t=2.0$}
\includegraphics[width=2.0in]{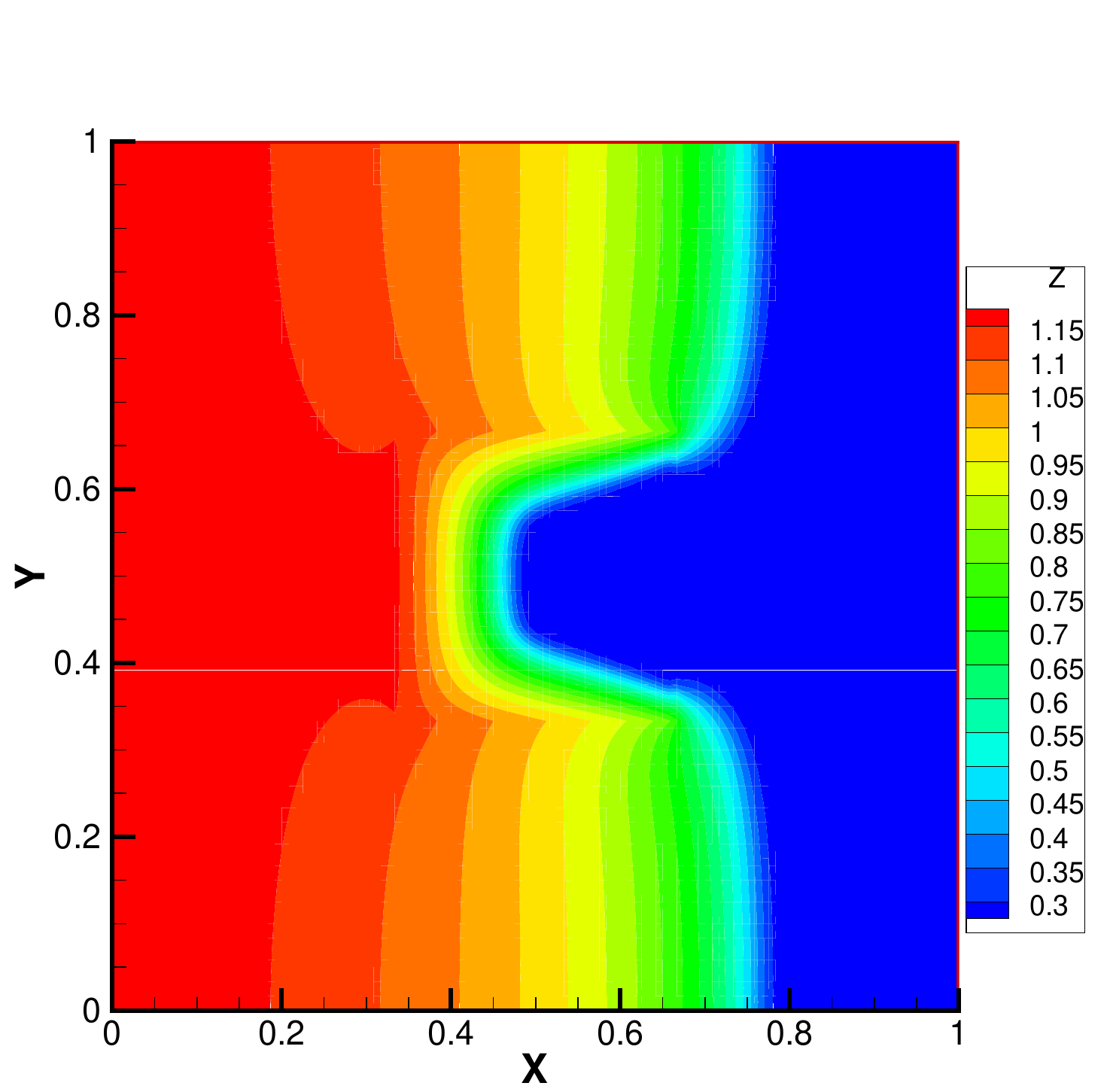}
\end{minipage}
}
\hbox{
\hspace{1in}
\begin{minipage}[t]{2.0in}
\centerline{\scriptsize (e):  with MM at $t=2.5$}
\includegraphics[width=2.0in]{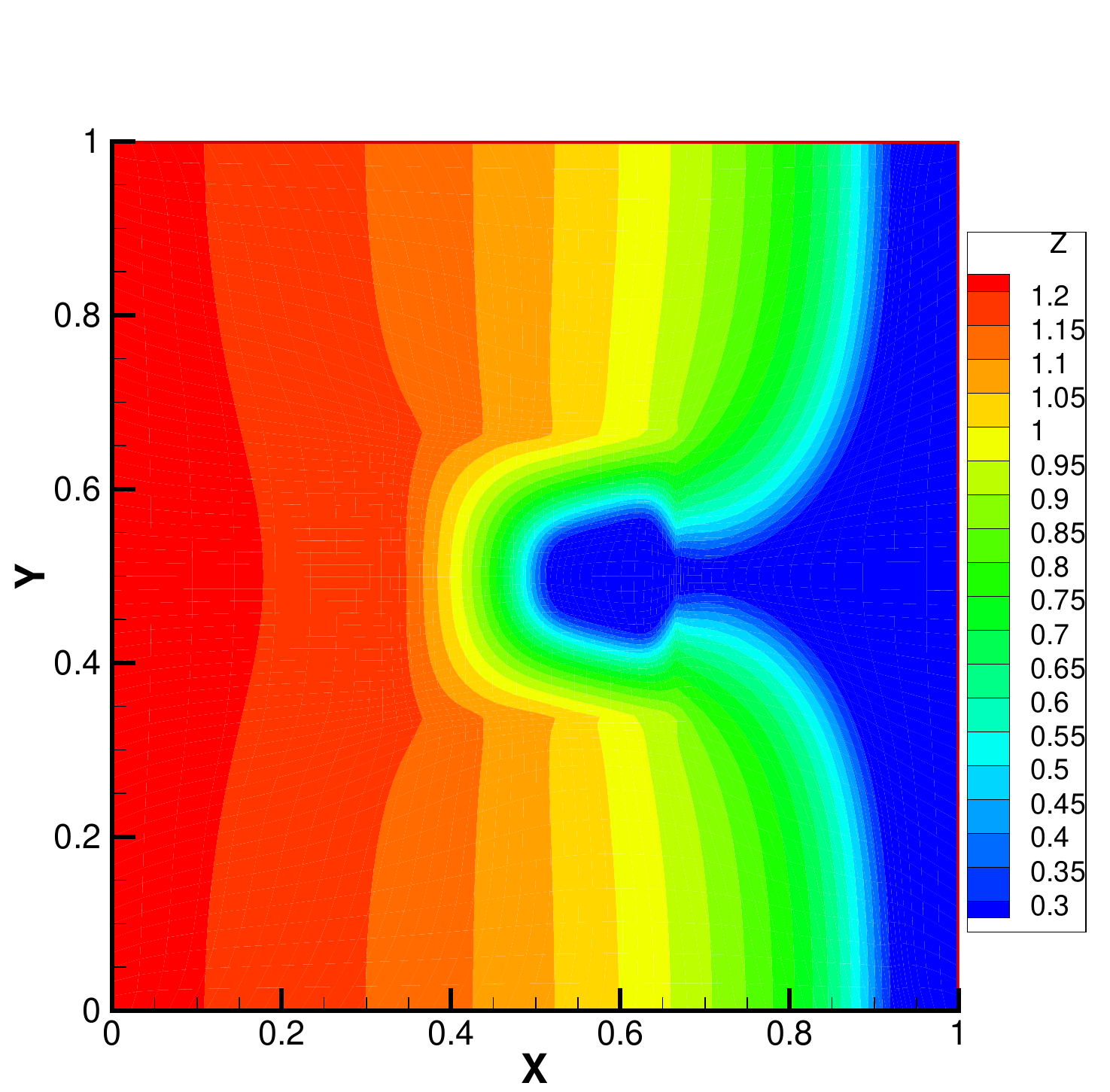}
\end{minipage}
\hspace{0.5in}
\begin{minipage}[t]{2.0in}
\centerline{\scriptsize (f):  with UM at $t=2.5$}
\includegraphics[width=2.0in]{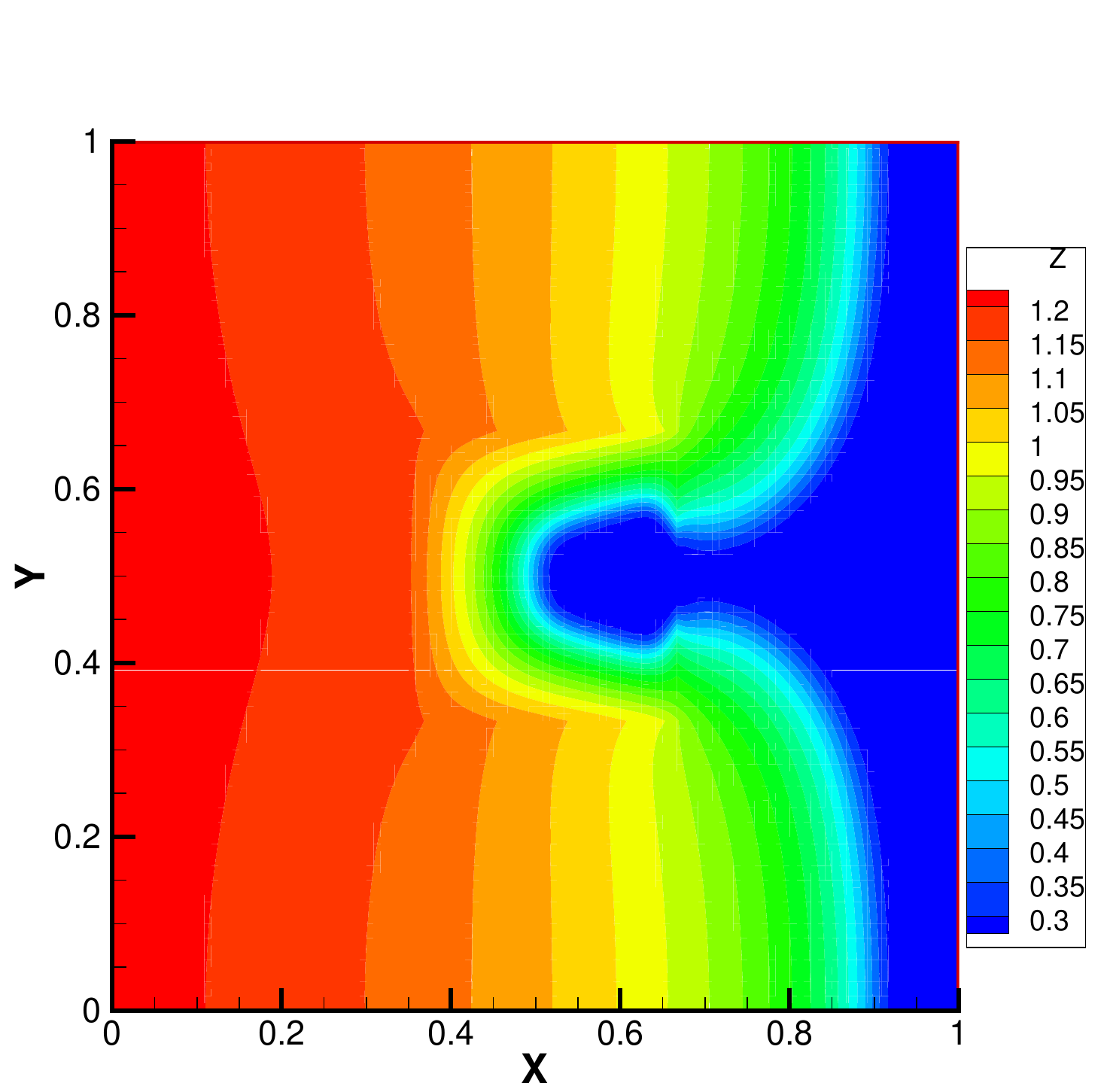}
\end{minipage}
}
\hbox{
\hspace{1in}
\begin{minipage}[t]{2.0in}
\centerline{\scriptsize (g):  with MM at $t=3.0$}
\includegraphics[width=2.0in]{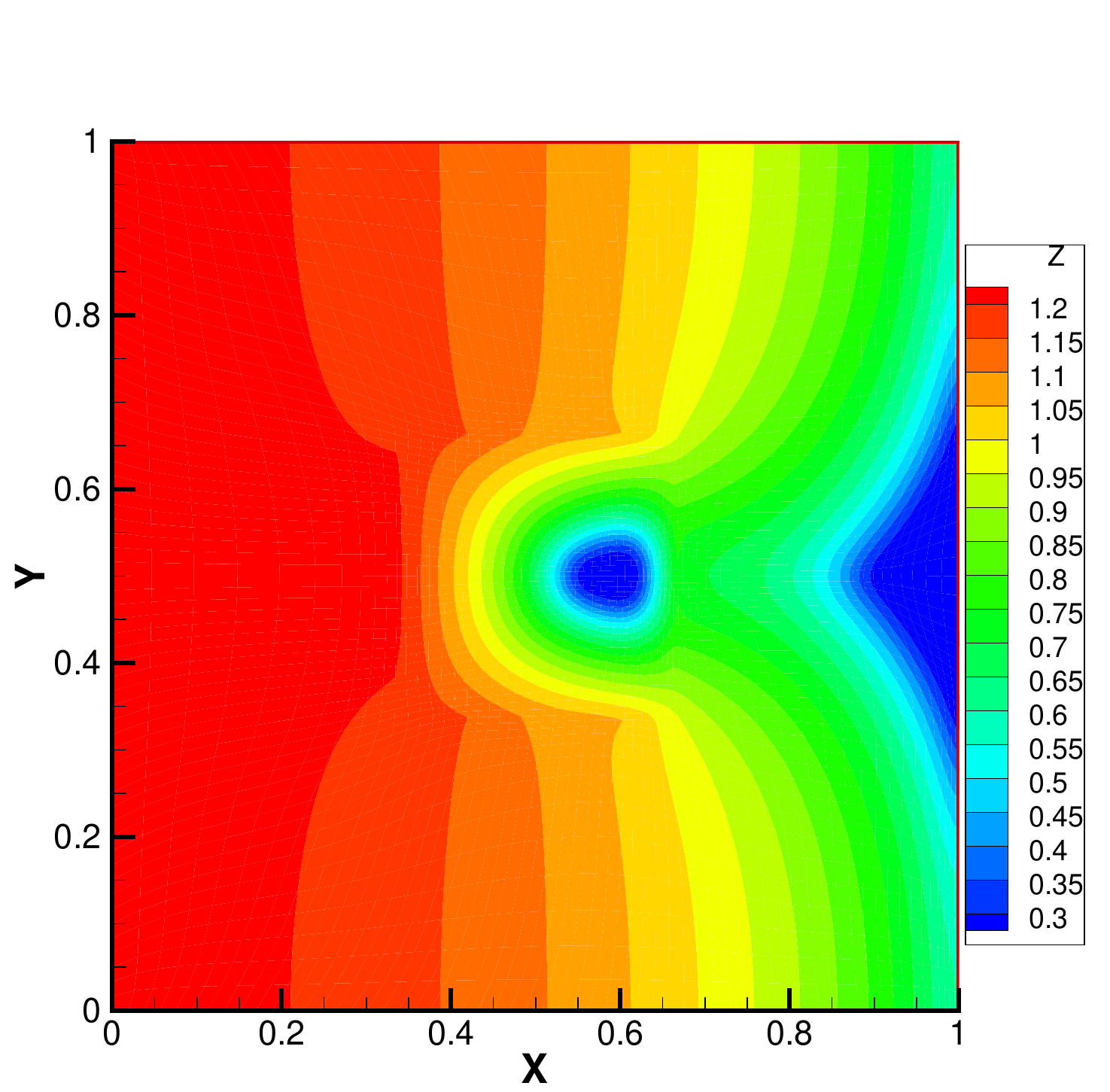}
\end{minipage}
\hspace{0.5in}
\begin{minipage}[t]{2.0in}
\centerline{\scriptsize (h):  with UM at $t=3.0$}
\includegraphics[width=2.0in]{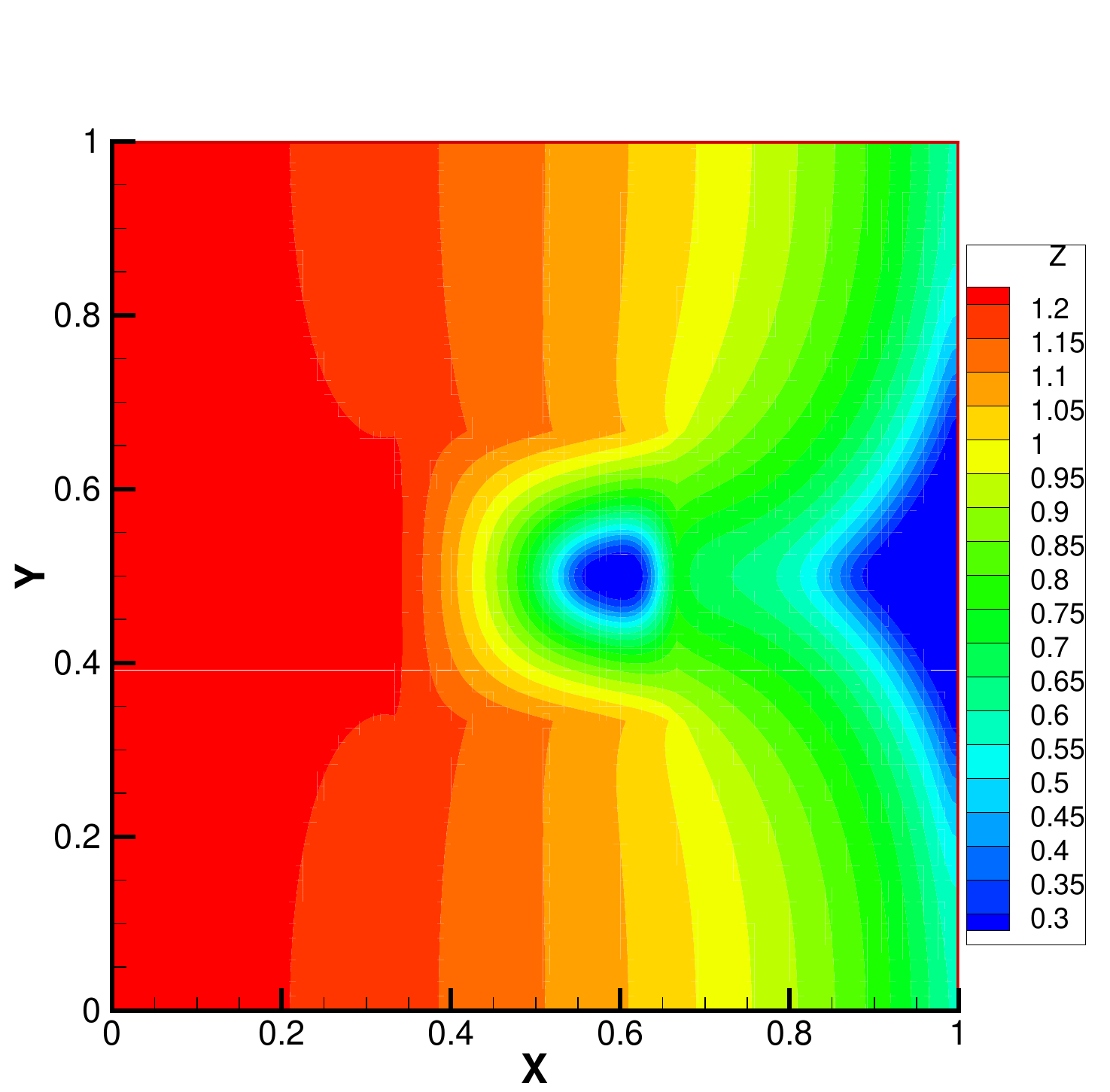}
\end{minipage}
}
\end{center}
\caption{Example~\ref{Example4.1}. The contours of the temperature obtained
with a moving mesh (MM) of size  $61\times 61$ are compared with those
obtained with a uniform mesh (UM) of size $121 \times121$.}
\label{T7}
\end{figure}

\begin{figure}
\begin{center}
\hbox{
\hspace{1in}
\begin{minipage}[t]{2.0in}
\centerline{\scriptsize (a):  with MM1 at $t=1.0$}
\includegraphics[width=2.0in]{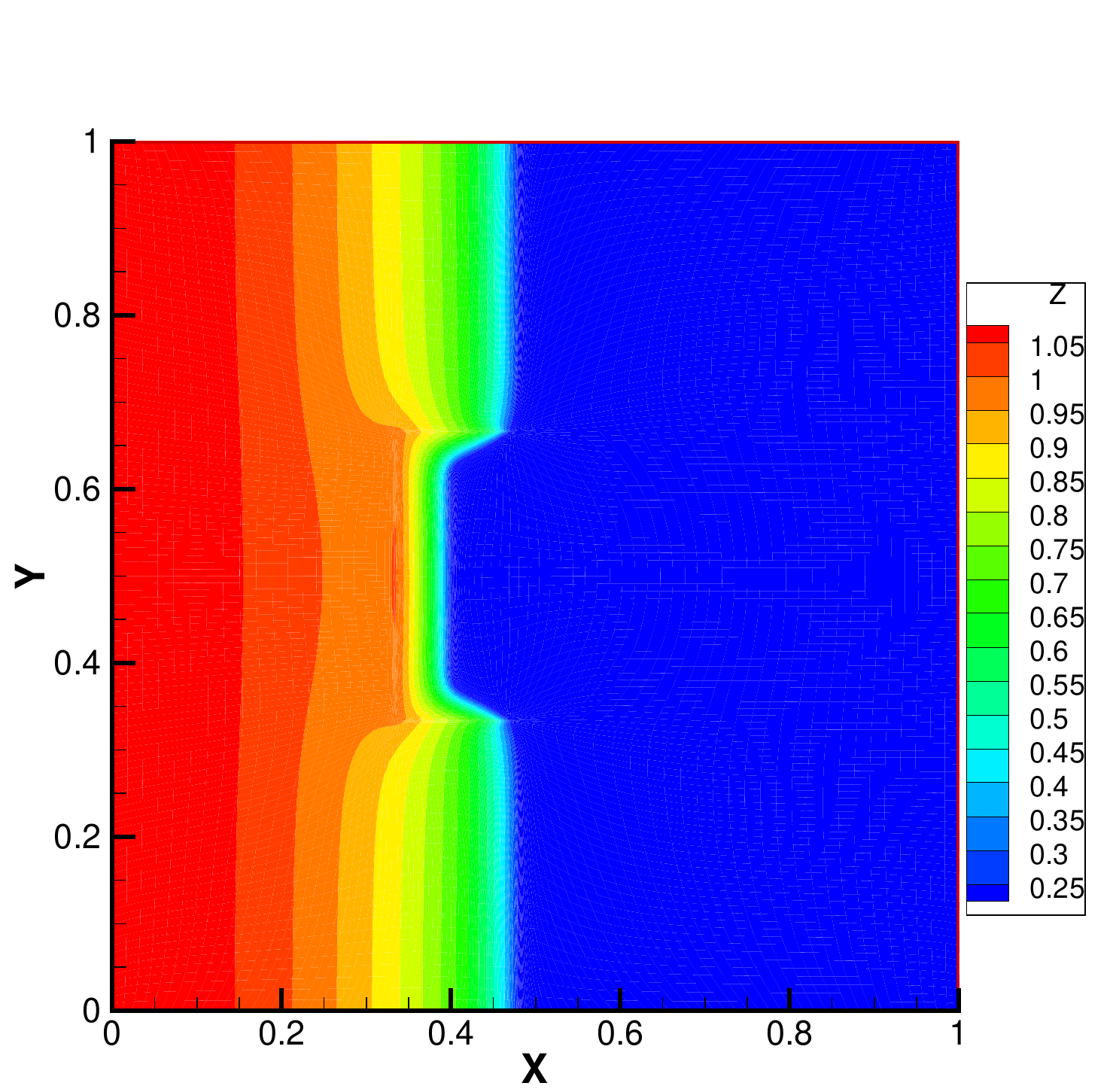}
\end{minipage}
\hspace{0.5in}
\begin{minipage}[t]{2.0in}
\centerline{\scriptsize (b):  with MM2 at $t=1.0$}
\includegraphics[width=2.0in]{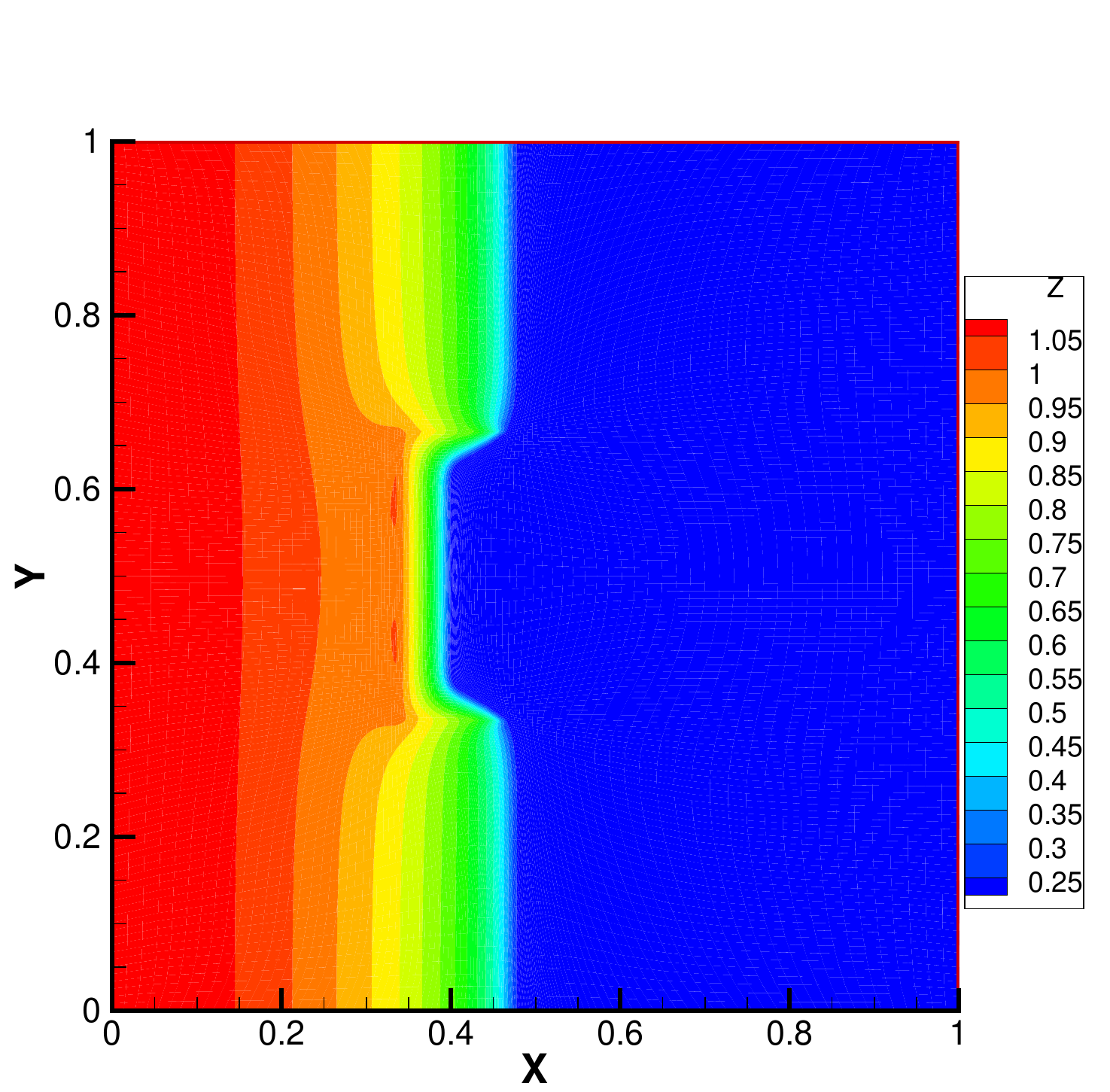}
\end{minipage}
}
\hbox{
\hspace{1in}
\begin{minipage}[t]{2.0in}
\centerline{\scriptsize (c):  with MM1 at $t=2.0$}
\includegraphics[width=2.0in]{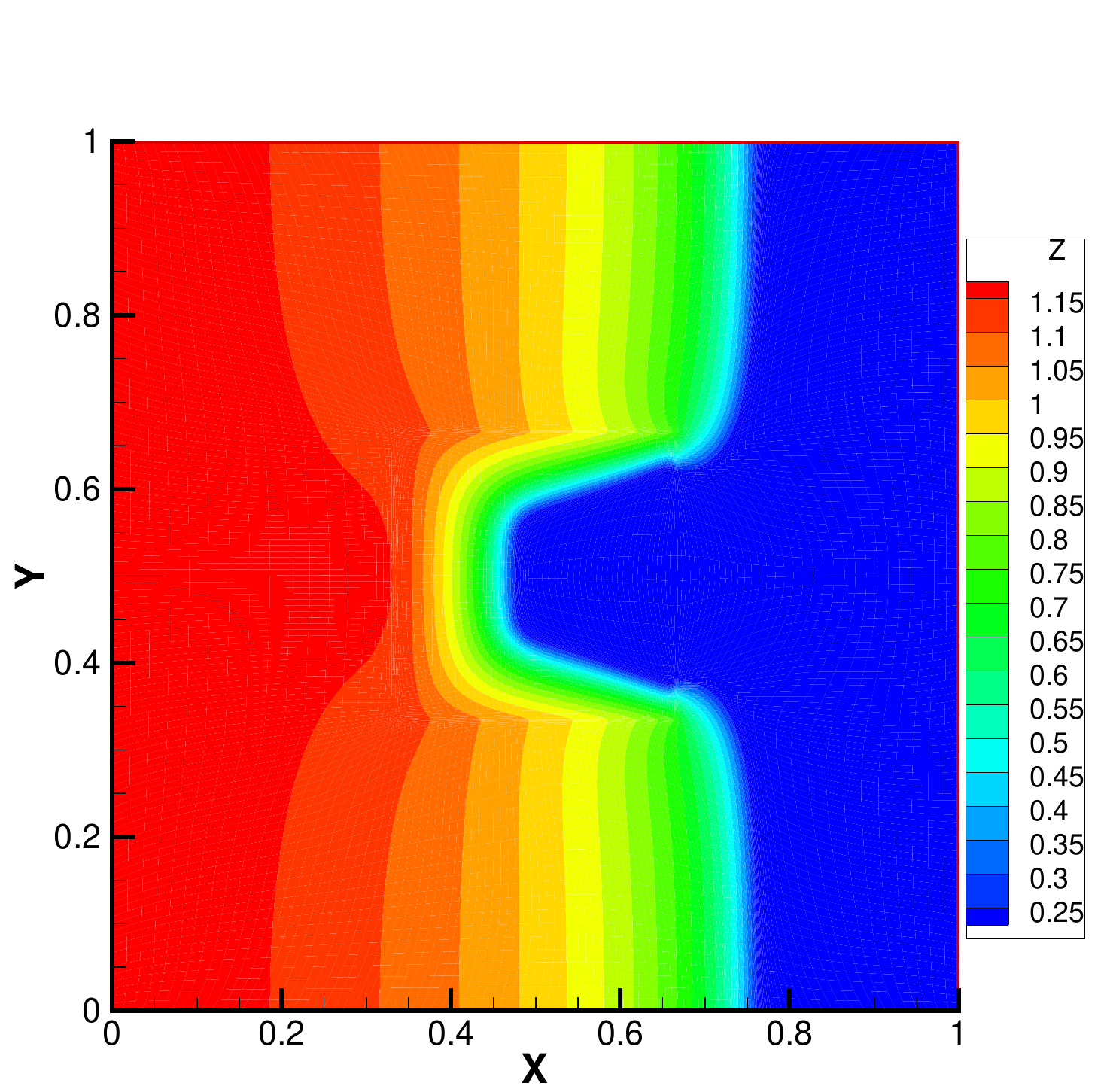}
\end{minipage}
\hspace{0.5in}
\begin{minipage}[t]{2.0in}
\centerline{\scriptsize (d):  with MM2 at $t=2.0$}
\includegraphics[width=2.0in]{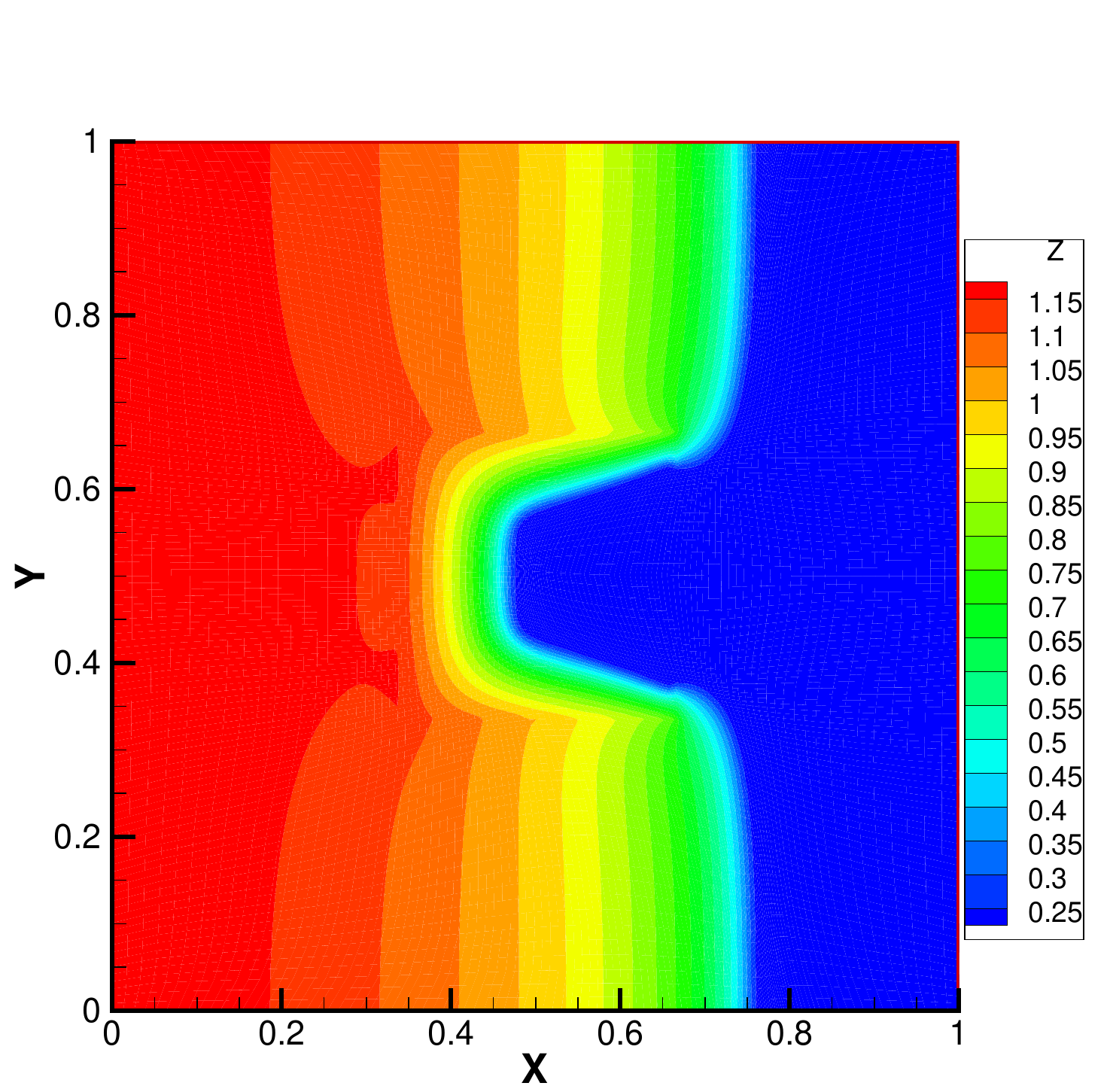}
\end{minipage}
}
\hbox{
\hspace{1in}
\begin{minipage}[t]{2.0in}
\centerline{\scriptsize (e):  with MM1 at $t=2.5$}
\includegraphics[width=2.0in]{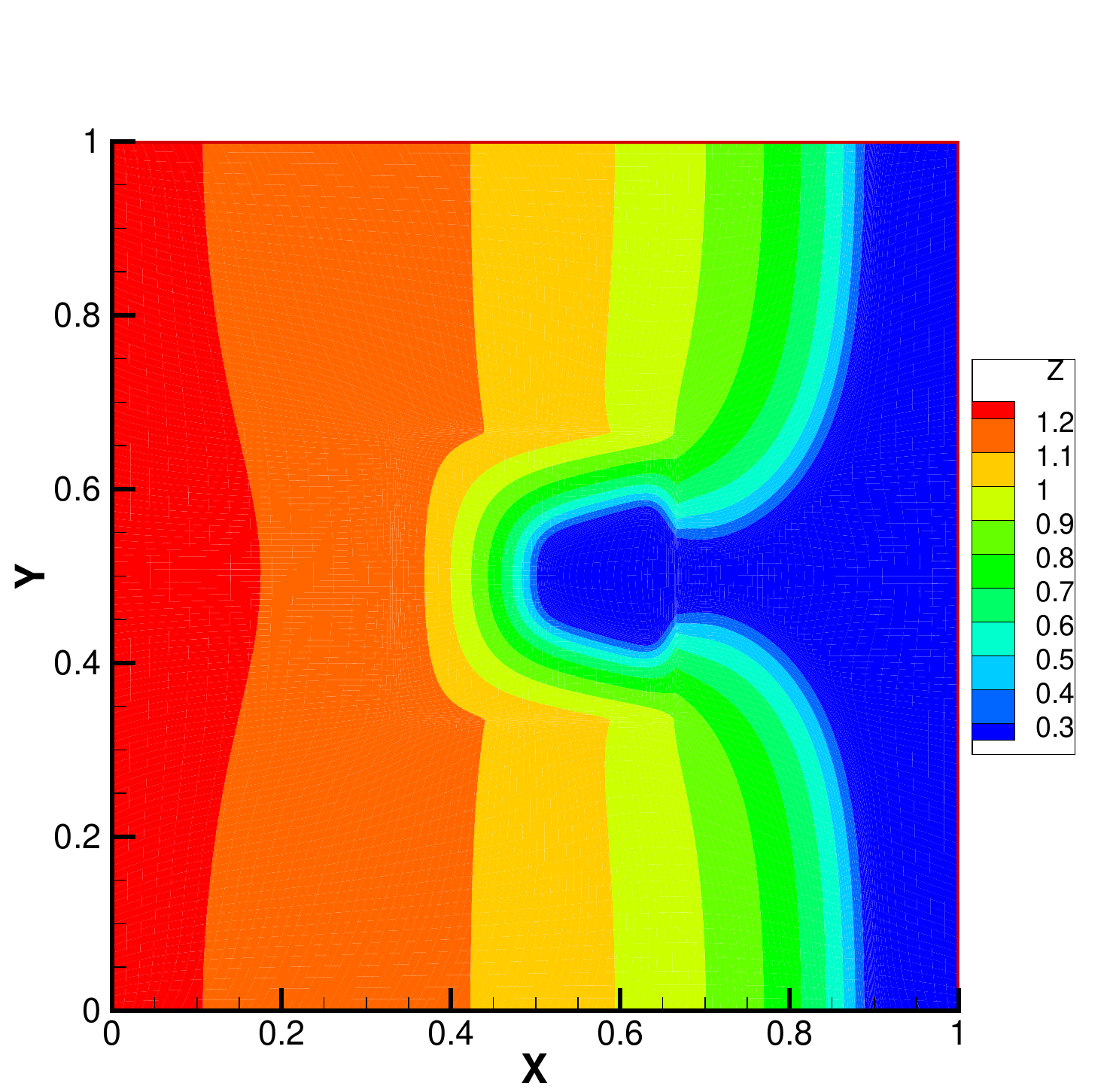}
\end{minipage}
\hspace{0.5in}
\begin{minipage}[t]{2.0in}
\centerline{\scriptsize (f):  with MM2 at $t=2.5$}
\includegraphics[width=2.0in]{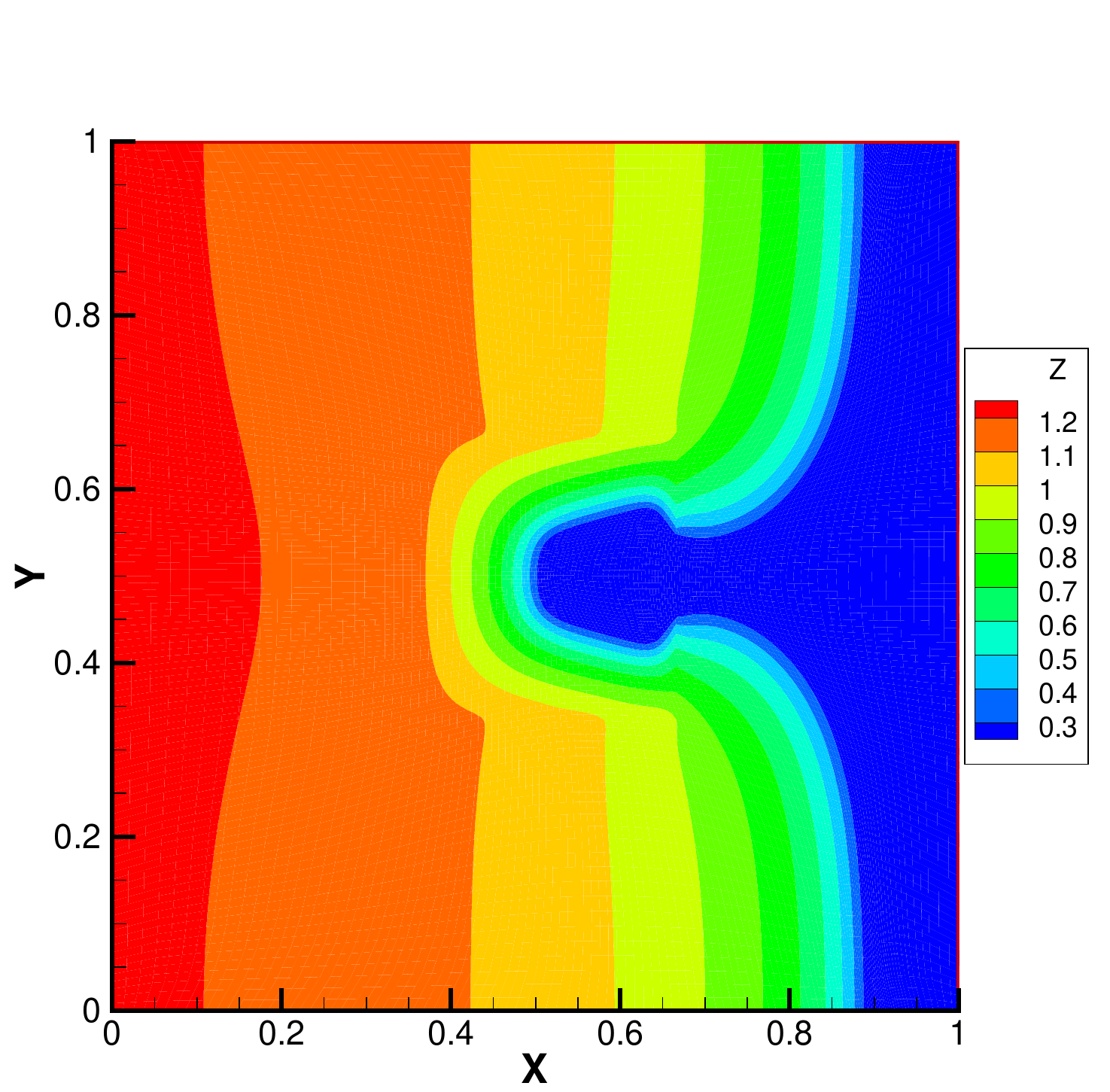}
\end{minipage}
}
\hbox{
\hspace{1in}
\begin{minipage}[t]{2.0in}
\centerline{\scriptsize (g):  with MM1 at $t=3.0$}
\includegraphics[width=2.0in]{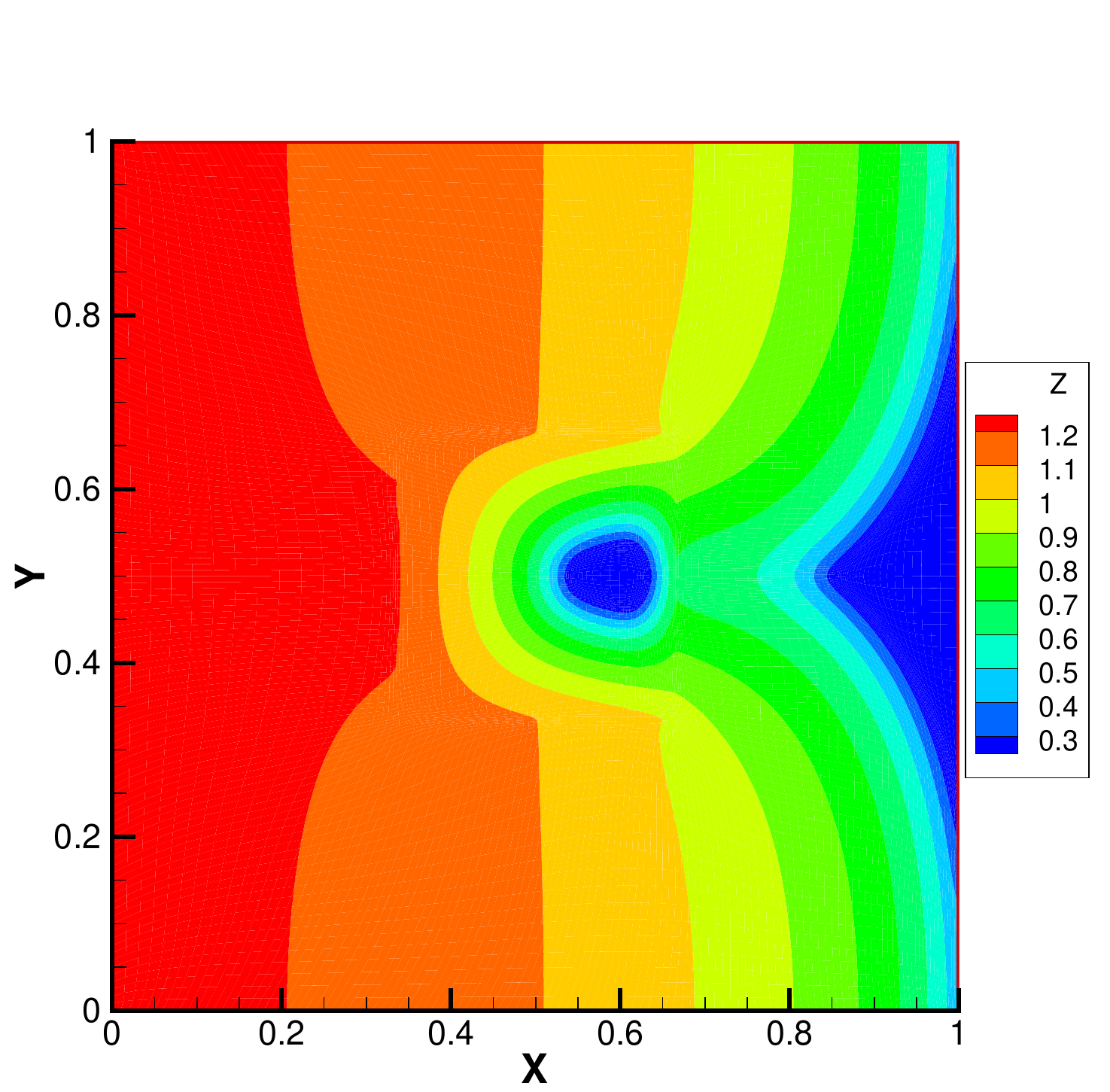}
\end{minipage}
\hspace{0.5in}
\begin{minipage}[t]{2.0in}
\centerline{\scriptsize (h):  with MM2 at $t=3.0$}
\includegraphics[width=2.0in]{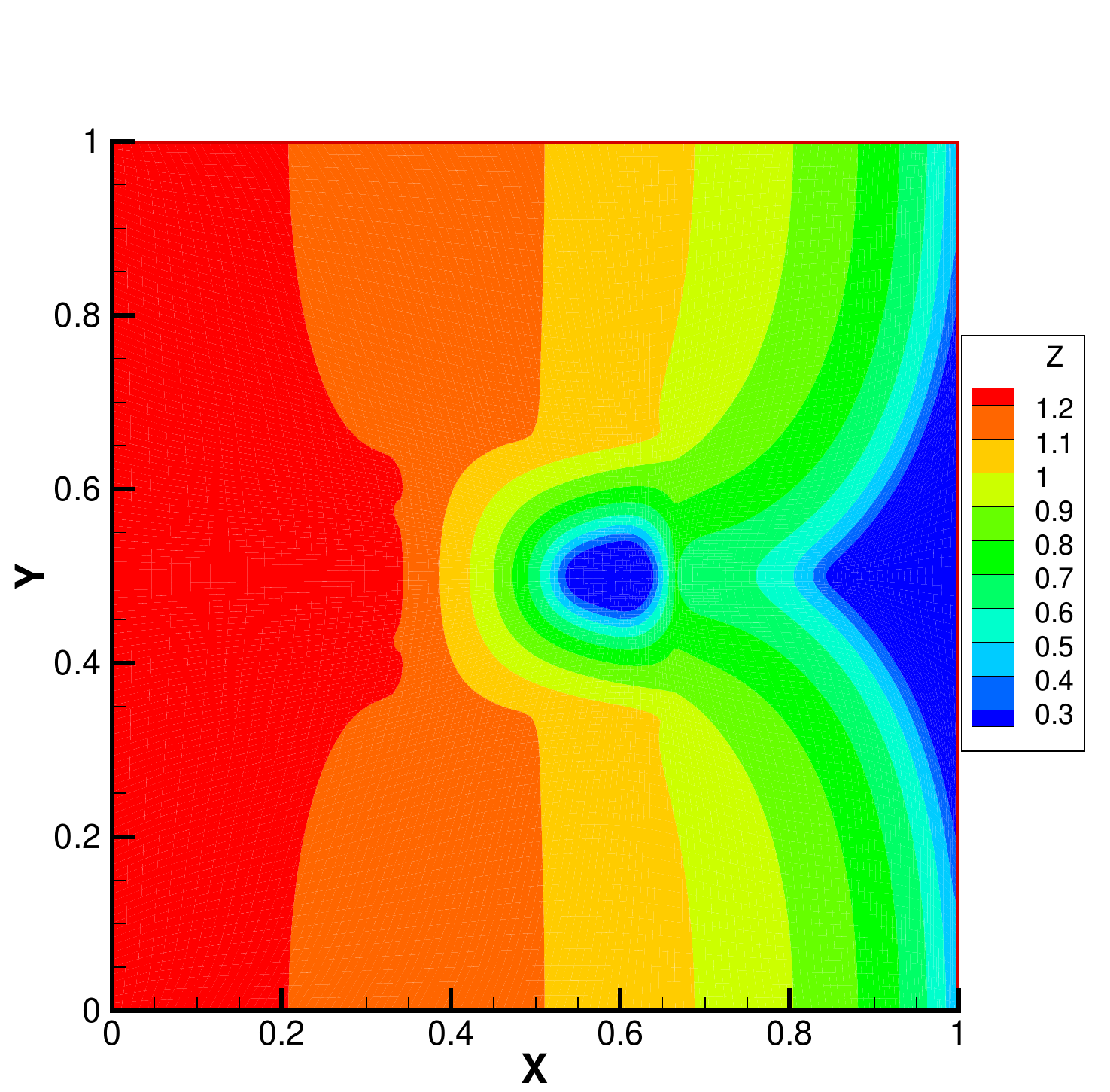}
\end{minipage}
}
\end{center}
\caption{Example~\ref{Example4.1}. The contours of the temperature obtained
with an MM1 of size  $121\times 121$ are compared with those
obtained with an MM2 of size  $41\times 41$ (with the physical PDE
being solved on a mesh of size $121\times 121$.}
\label{T7a}
\end{figure}

\begin{figure}
\begin{center}
\hbox{
\hspace{1in}
\begin{minipage}[t]{2.0in}
\centerline{\scriptsize (a):  with MM1 at $t=1.0$}
\includegraphics[width=2.0in]{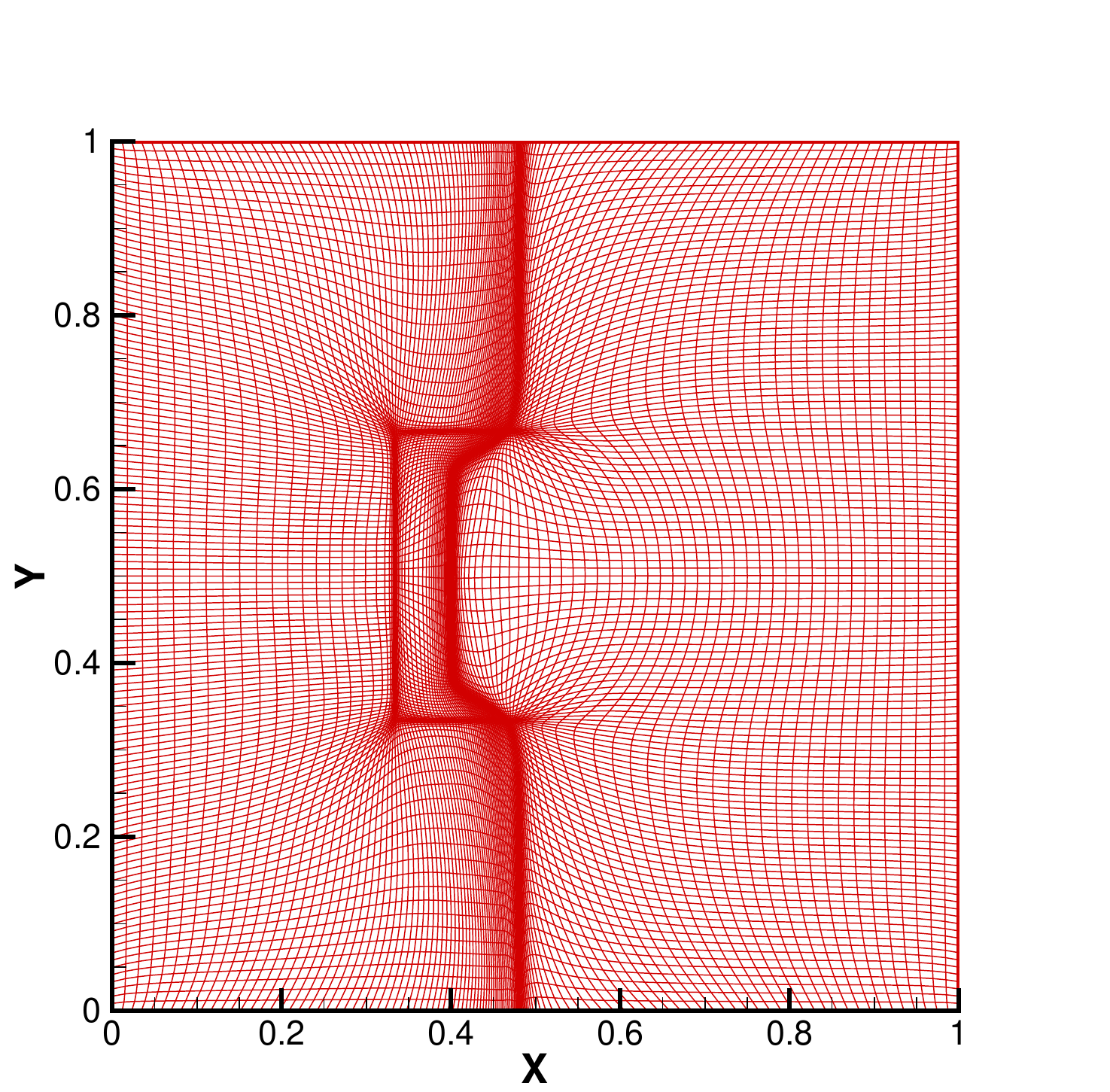}
\end{minipage}
\hspace{0.5in}
\begin{minipage}[t]{2.0in}
\centerline{\scriptsize (b):  with MM2 at $t=1.0$}
\includegraphics[width=2.0in]{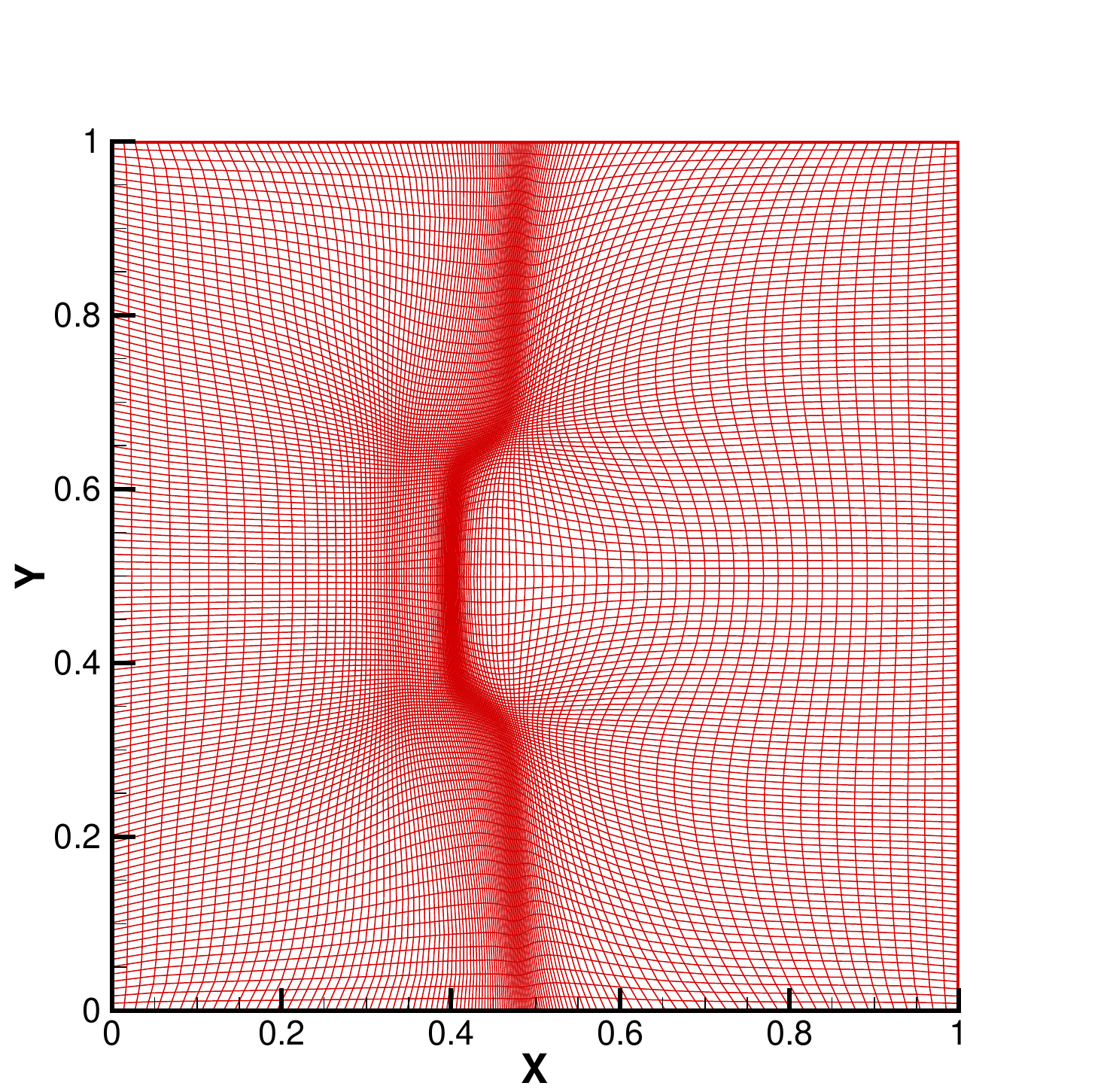}
\end{minipage}
}
\hbox{
\hspace{1in}
\begin{minipage}[t]{2.0in}
\centerline{\scriptsize (c):  with MM1 at $t=2.0$}
\includegraphics[width=2.0in]{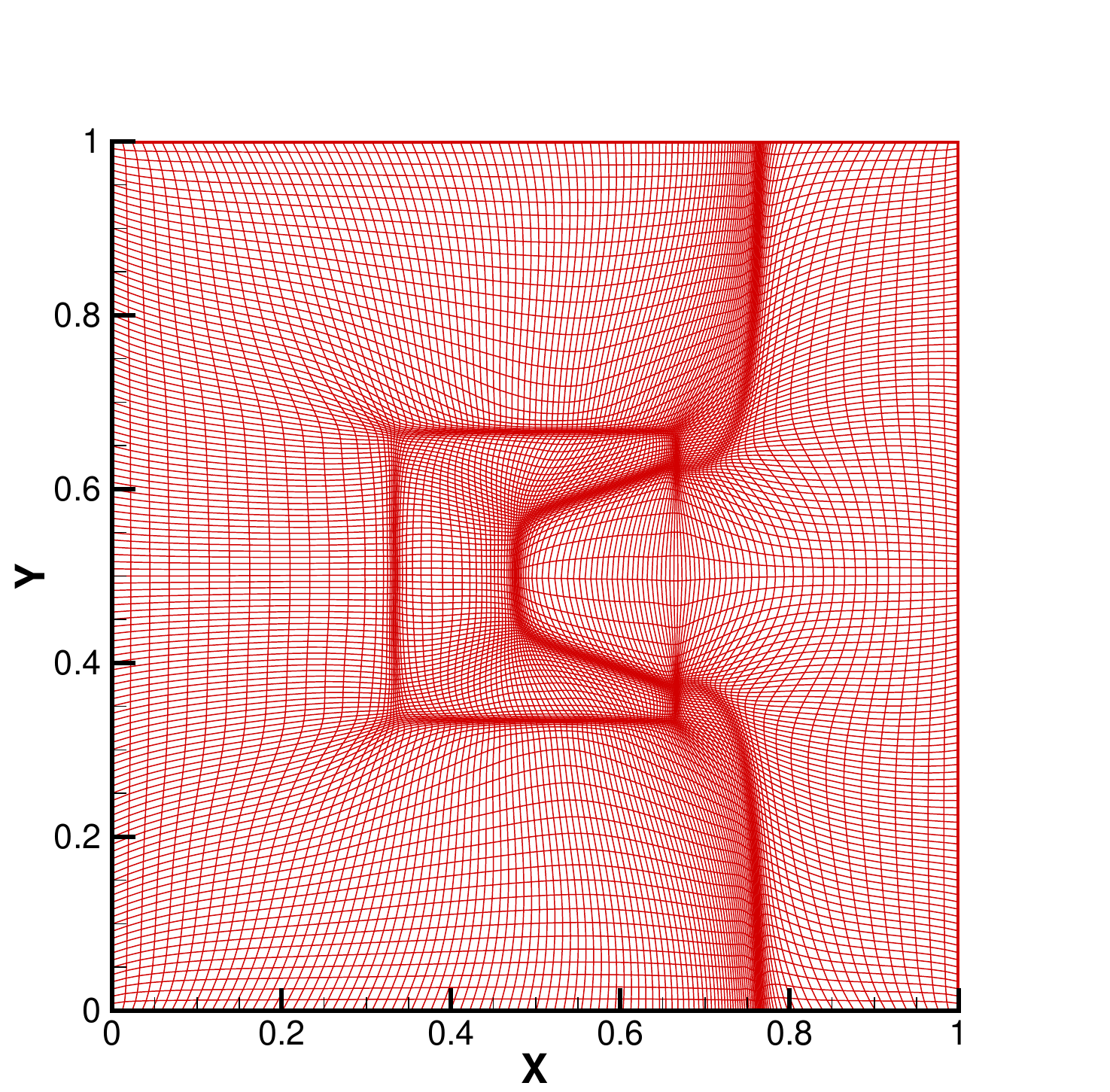}
\end{minipage}
\hspace{0.5in}
\begin{minipage}[t]{2.0in}
\centerline{\scriptsize (d):  with MM2 at $t=2.0$}
\includegraphics[width=2.0in]{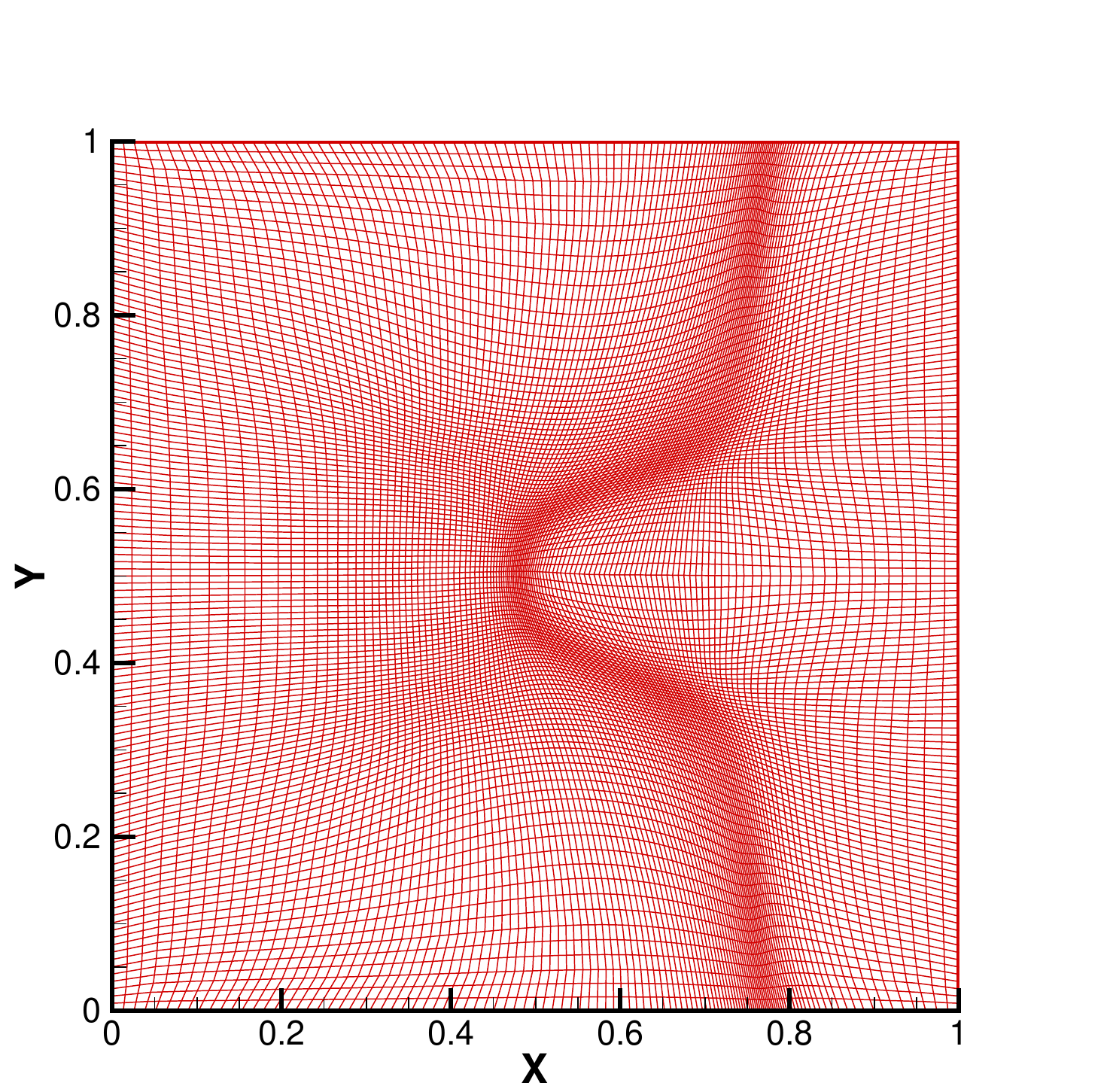}
\end{minipage}
}
\hbox{
\hspace{1in}
\begin{minipage}[t]{2.0in}
\centerline{\scriptsize (e):  with MM1 at $t=2.5$}
\includegraphics[width=2.0in]{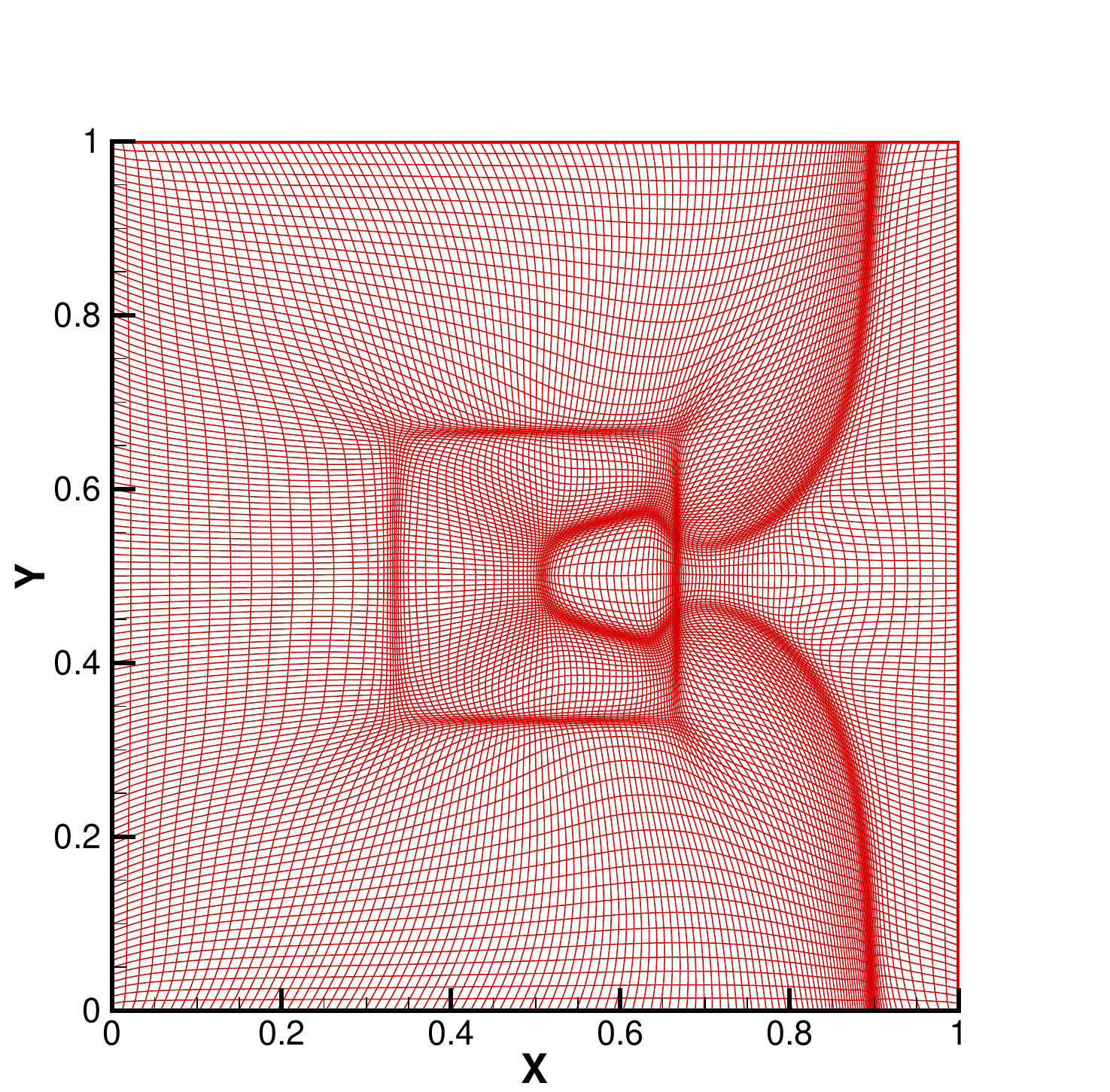}
\end{minipage}
\hspace{0.5in}
\begin{minipage}[t]{2.0in}
\centerline{\scriptsize (f):  with MM2 at $t=2.5$}
\includegraphics[width=2.0in]{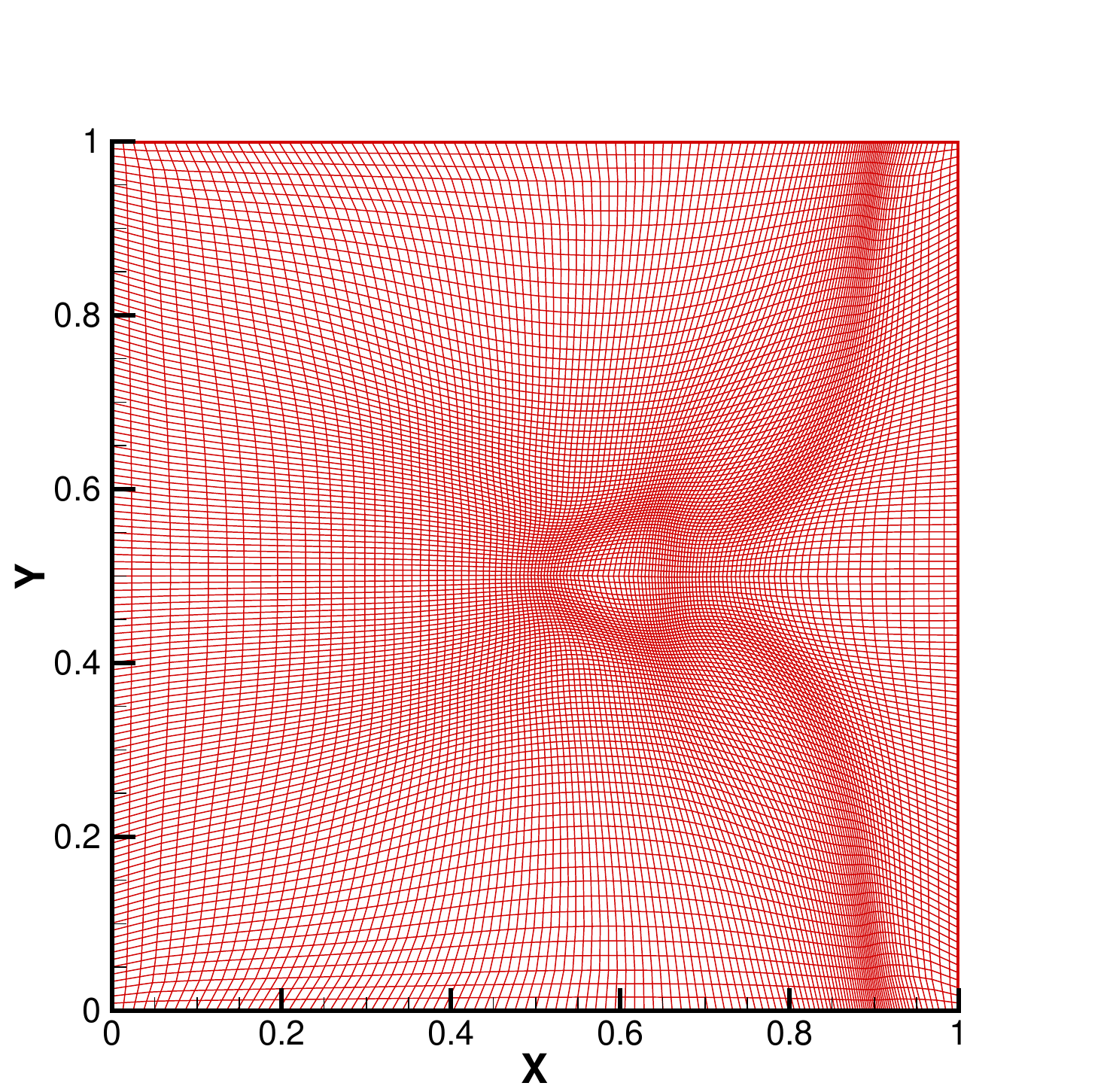}
\end{minipage}
}
\hbox{
\hspace{1in}
\begin{minipage}[t]{2.0in}
\centerline{\scriptsize (g):  with MM1 at $t=3.0$}
\includegraphics[width=2.0in]{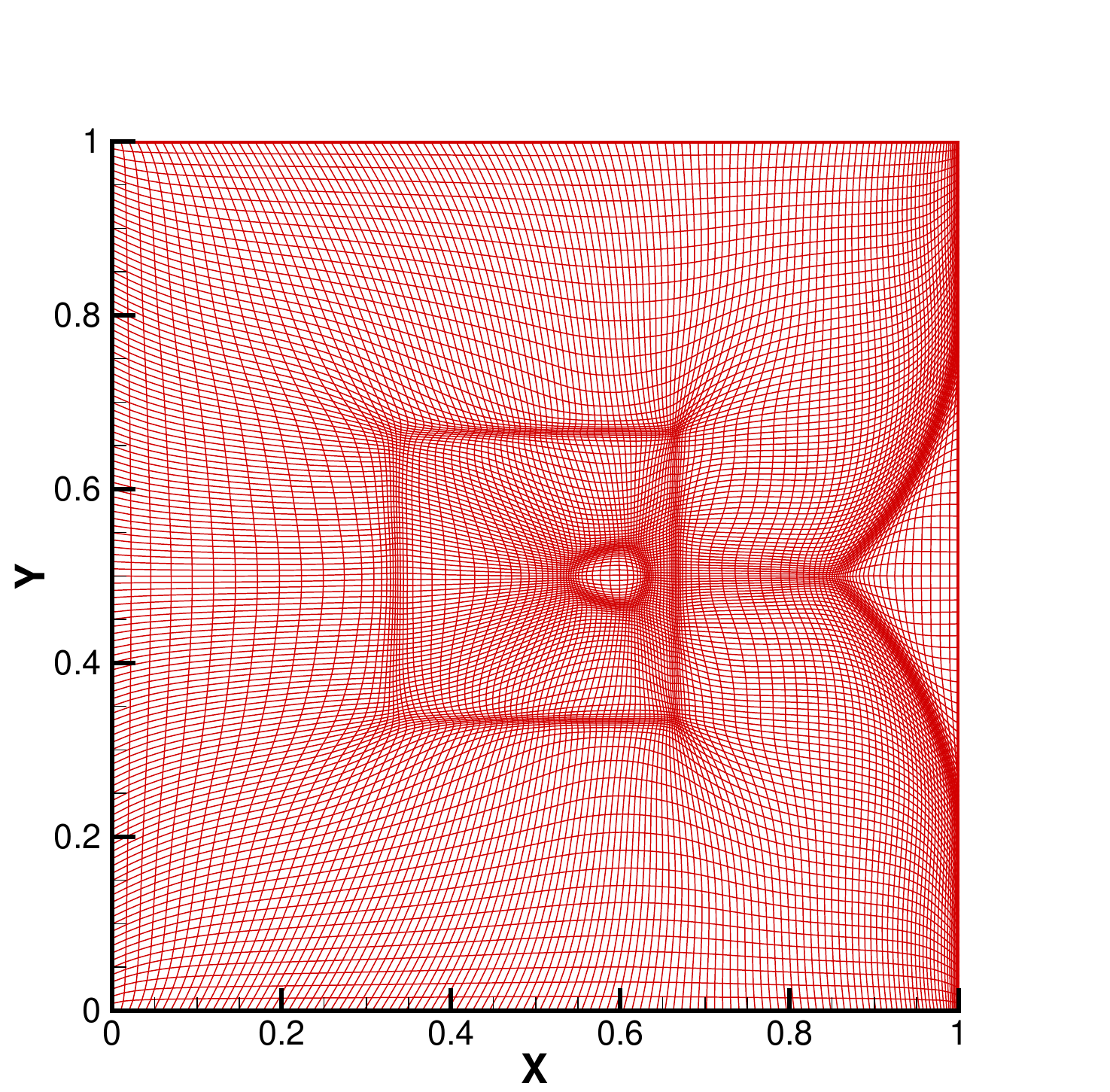}
\end{minipage}
\hspace{0.5in}
\begin{minipage}[t]{2.0in}
\centerline{\scriptsize (h):  with MM2 at $t=3.0$}
\includegraphics[width=2.0in]{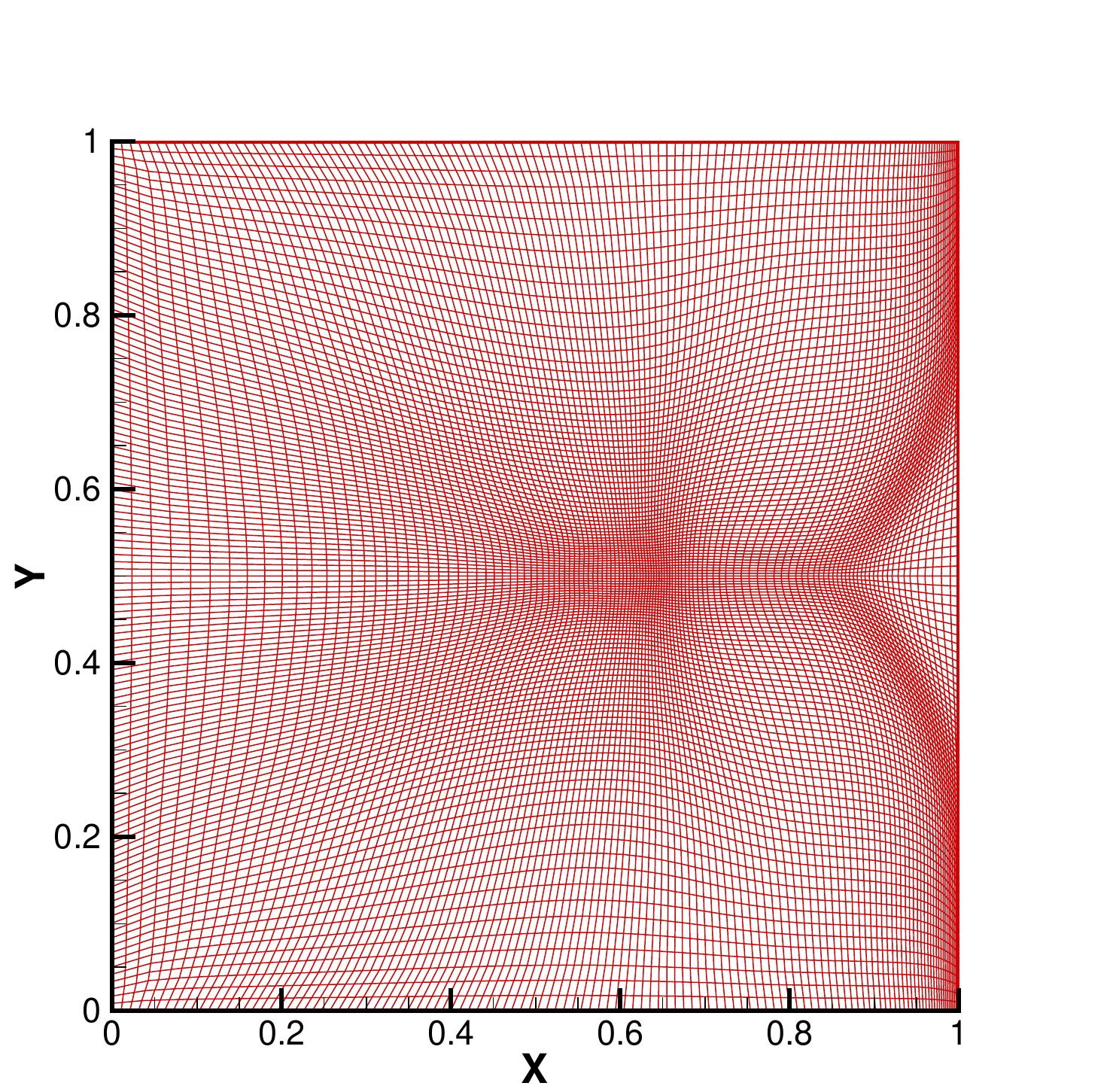}
\end{minipage}
}
\end{center}
\caption{Example~\ref{Example4.1}.  Moving meshes of size $121\times 121$ obtained with
MM1 and MM2 moving mesh strategies. MM2 is obtained by uniformly interpolating
a $41\times 41$ moving mesh.}
\label{T7b}
\end{figure}

\begin{exam}
\label{Example4.1}
For this example, the distribution of the atomic mass number is given by
\begin{equation}
\label{Z2}
z(x,y)=\begin{cases}
5,  & \text{for } (x,y) \in (\frac13,\frac23)\times  (\frac13,\frac23)\\
1,    & \text{otherwise} .
\end{cases}
\end{equation}
The initial and boundary conditions are given in (\ref{BC-1}) and (\ref{IC-2}), respectively.

A typical moving mesh of $81\times81$ and the computed solution thereon are shown in
Figs.~\ref{T5} and \ref{T6}. From the figures, we can see that the hot wave front
propagates from left to right and meets the central obstacle and then a Marshak wave is formed.
The profile of the Marshak wave has been captured accurately by the moving mesh
and the nodes concentrate around the front of the wave. This demonstrates the mesh
concentration ability of the moving mesh method.
Fig.~\ref{T7} shows the solutions obtained with a moving mesh of $61\times 61$ and a uniform mesh of
$121\times 121$, which are comparable.

The results obtained with a moving mesh of $121\times121$ are
compared in Figs.~\ref{T7a} and \ref{T7b} to those obtain with a two-level moving mesh strategy (MM2) \cite{Hua01}
where a mesh of size $41\times 41$ is moved using the moving mesh method but the physical PDEs are solved
on a mesh of $121\times121$ that is generated by uniformly refining the moving mesh.
Interestingly, MM2 leads to results with comparable accuracy but saves significant CPU time.
The CPU times for one-level and two-level moving meshes and uniform meshes are listed in Table~\ref{mytab1}.
From the table, one can see that the moving mesh method is more costly than the method with a uniform mesh
of the same size. This is not surprising since the moving mesh method solves more equations.
The efficiency of the moving mesh method can be improved significantly using the two-level moving mesh
strategy. For example, for the case with mesh $81\times81$, the CPU time
of MM2 (with the coarse mesh $41\times 41$) is about 25.3\% of that with the one-level moving mesh (MM1).
For the case $121\times 121$, the CPU time for MM2 is only about 5.7\% of that of MM1.
Moreover, when the mesh size increases from $41\times 41$ to $81\times 81$ the CPU time
increases about 13.3 times for MM1. This number is about 10.6 times when the mesh size
increases from $81\times 81$ to $121\times 121$. For MM2, the corresponding number
is only 3.36 and 2.38, respectively. Finally, we compare MM2 with the uniform mesh method.
From Table.~\ref{mytab1}, we can see that the difference between the two is getting smaller
as the mesh becomes finer.
\end{exam}

\begin{table}[htbp]
 \centering\small
\begin{threeparttable}
\caption{CPU time comparison among one-level and two-level moving mesh methods and
the uniform mesh method for Example \ref{Example4.1}. The CPU time is measured in seconds.
The last column is the ratio of the used CPU time to that used with a uniform mesh of the same size.}
  \label{mytab1}
 \medskip
 \begin{tabular}{|l|c|c|c|c|}
 \hline
   & Fine Mesh  & Coarse Mesh & Total CPU time & ratio \\
   \hline
  One-level MM & 41$\times$41    &   41$\times$41     &      2544       &  5.32  \\
  \hline
          & 81$\times$81    &    81$\times$81     &      33724        &  14.81      \\
  \hline
          & 121$\times$121  &     121$\times$121     &     356720    &   63.91   \\
  \hline

Tow-level MM   & 41$\times$41  &      41$\times$41     &    2544  &    5.32        \\
\hline
          & 81$\times$81  &    41$\times$41     &      8549      &    3.76        \\
\hline
          & 121$\times$121  &   41$\times$41     &      20325      &   3.64           \\
\hline

Fixed mesh     & 41$\times$41  &      n/a          &      478       &   1                 \\
\hline
          & 81$\times$81  &       n/a         &      2276      &    1                   \\
\hline
          & 121$\times$121  &     n/a          &     5581      &    1                \\
\hline
\end{tabular}
\end{threeparttable}
\end{table}

\begin{figure}
\begin{center}
\hbox{
\begin{minipage}[t]{2.3in}
\centerline{\scriptsize (a): t=1.0}
\includegraphics[width=2.3in]{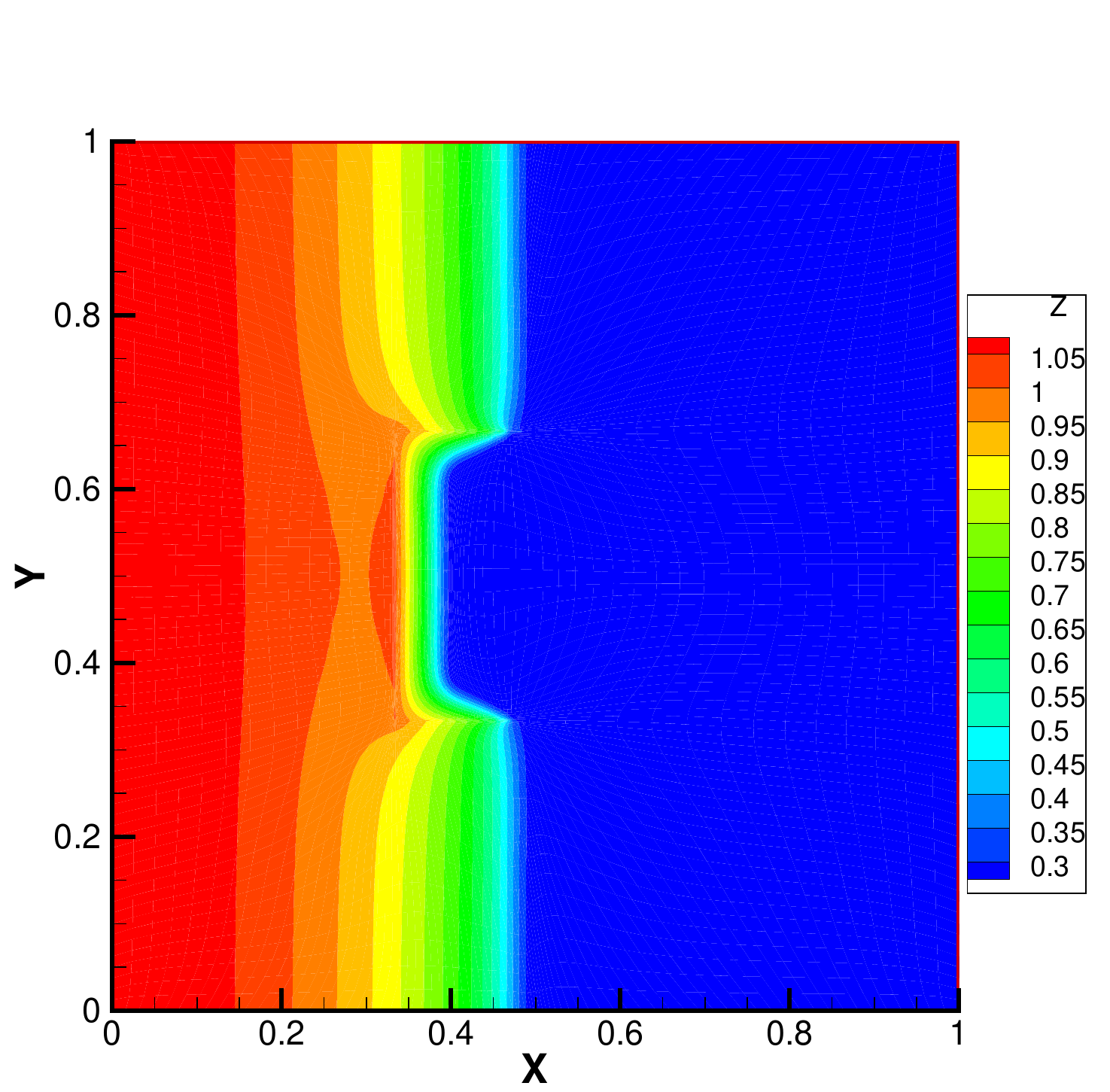}
\end{minipage}
\begin{minipage}[t]{2.3in}
\centerline{\scriptsize (b): t=1.5}
\includegraphics[width=2.3in]{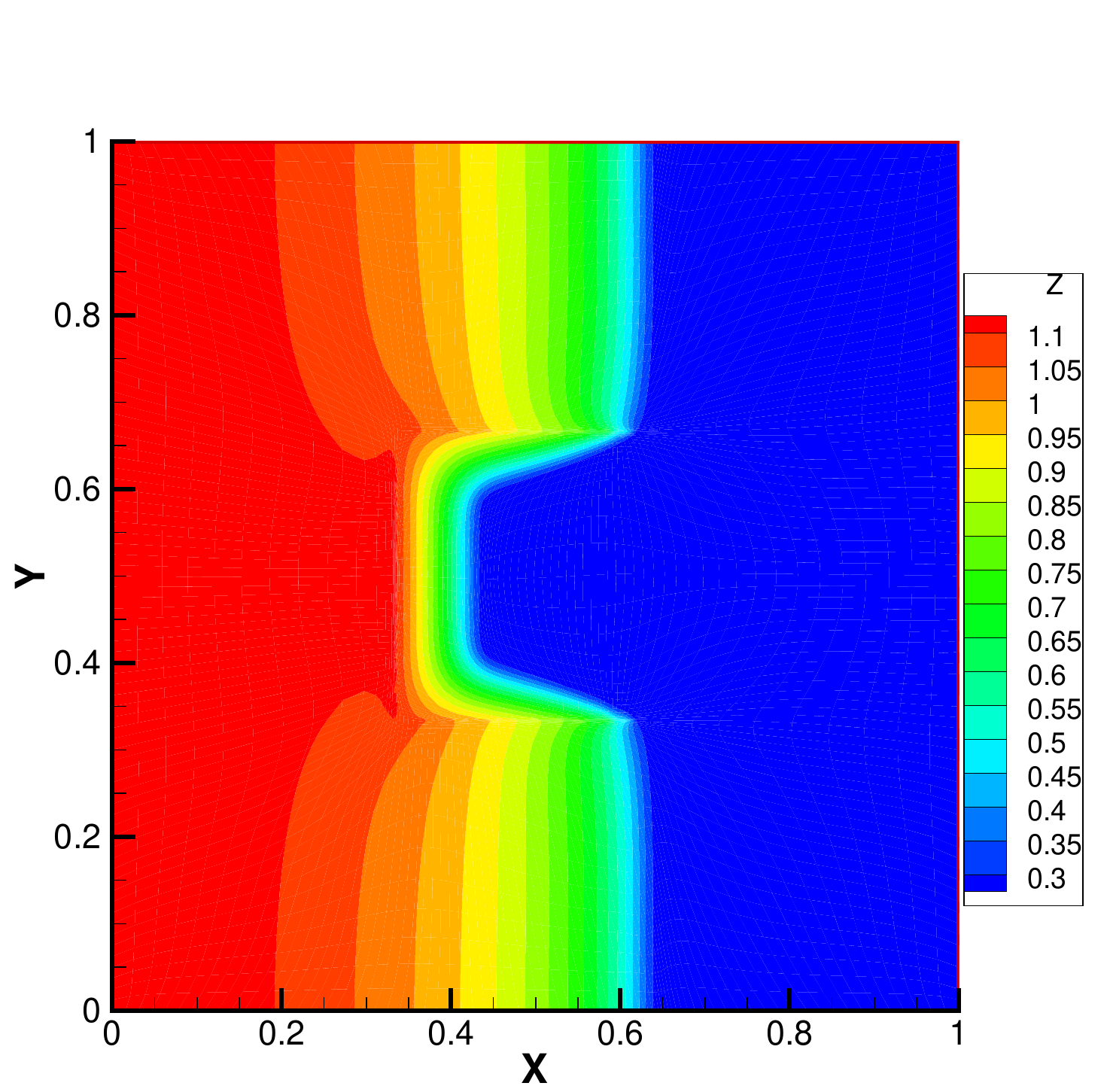}
\end{minipage}
\begin{minipage}[t]{2.3in}
\centerline{\scriptsize (c): t=2.0}
\includegraphics[width=2.3in]{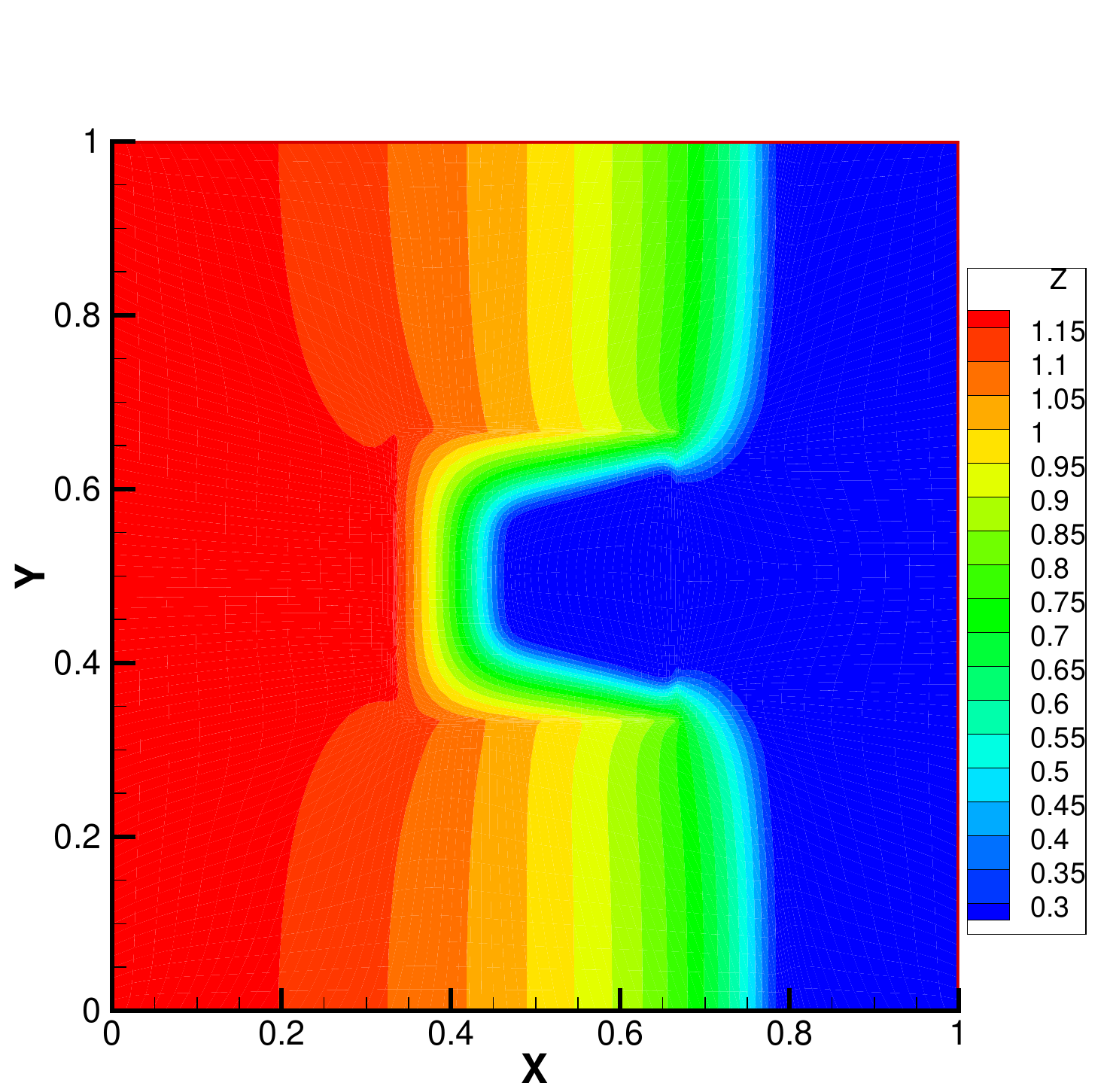}
\end{minipage}
}
\vspace{5mm}
\hbox{
\begin{minipage}[t]{2.3in}
\centerline{\scriptsize (d): t=2.4}
\includegraphics[width=2.3in]{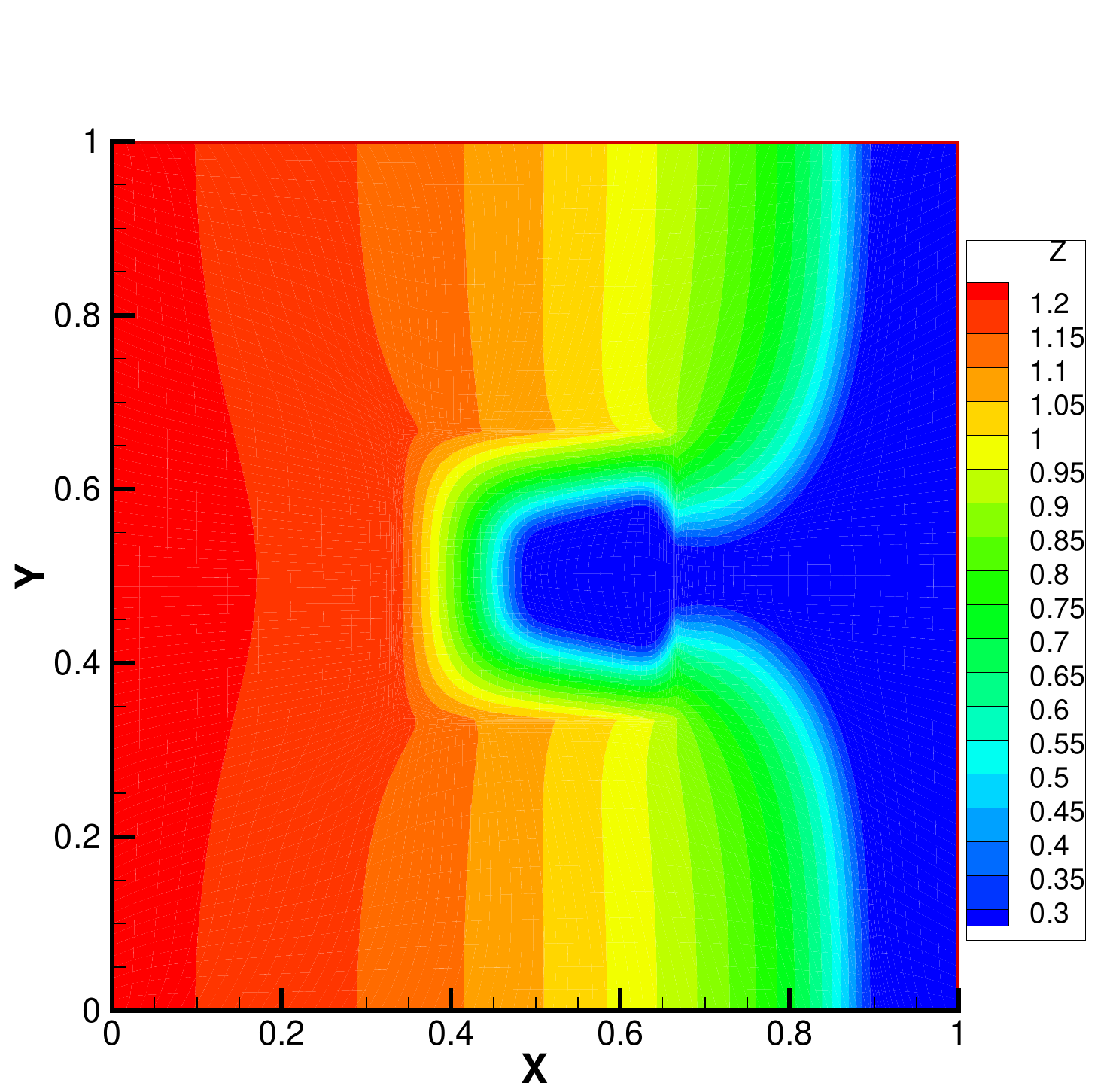}
\end{minipage}
\begin{minipage}[t]{2.3in}
\centerline{\scriptsize (e): t=2.8}
\includegraphics[width=2.3in]{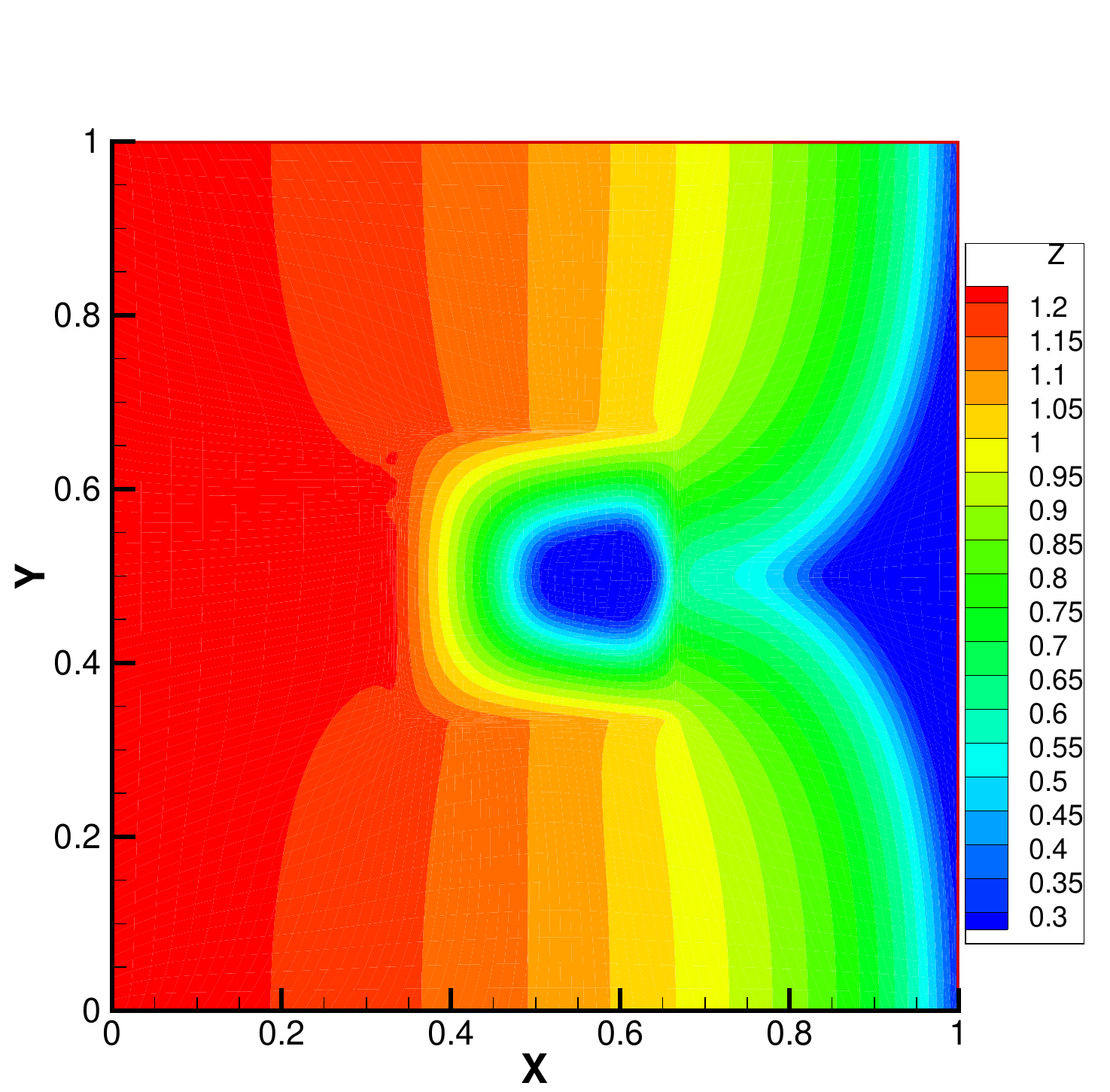}
\end{minipage}
\begin{minipage}[t]{2.3in}
\centerline{\scriptsize (f): t=3.0}
\includegraphics[width=2.3in]{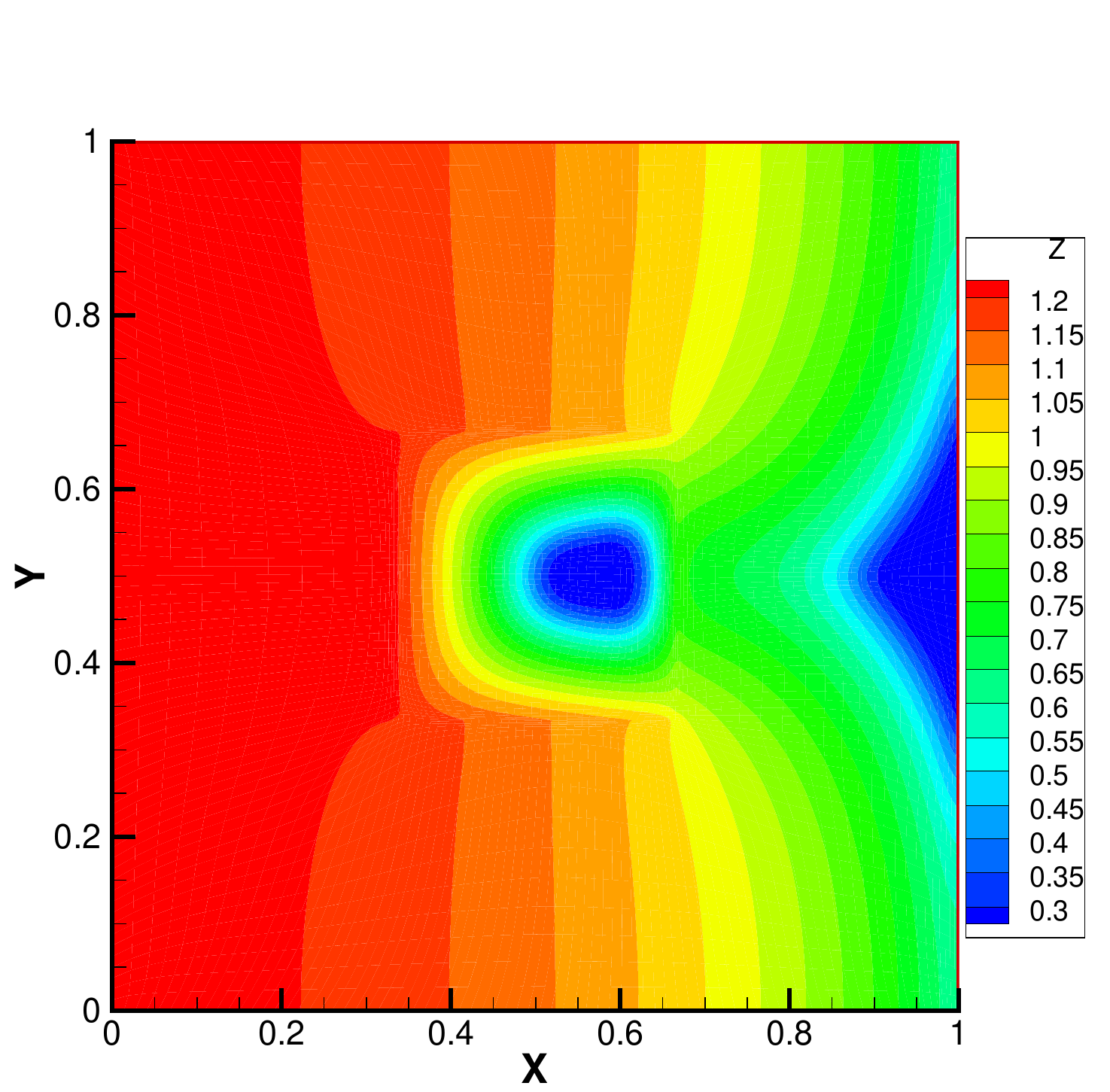}
\end{minipage}
}
\vspace{5mm}
\hbox{
\begin{minipage}[t]{2.3in}
\centerline{\scriptsize (g): t=3.5}
\includegraphics[width=2.3in]{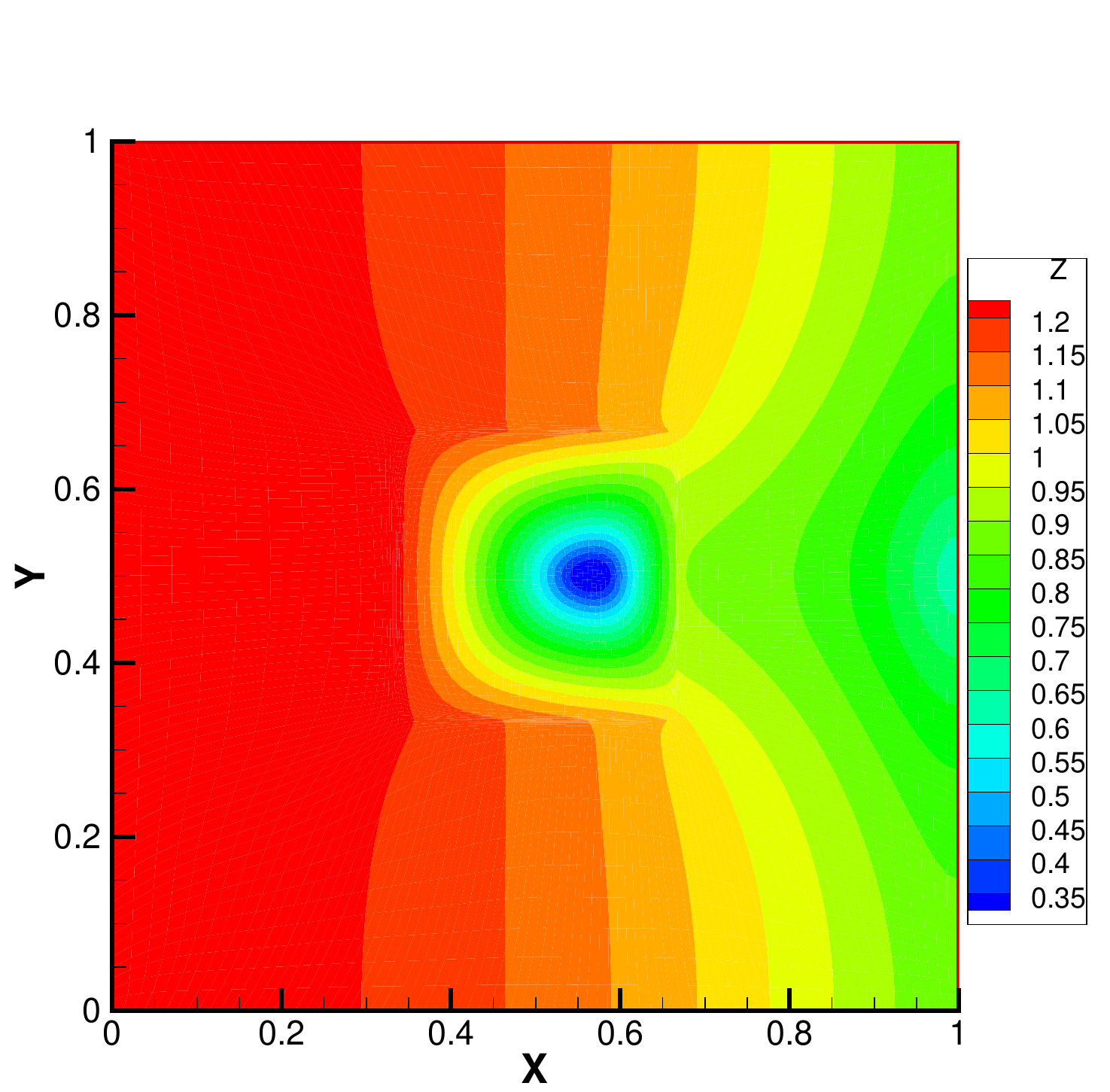}
\end{minipage}
\begin{minipage}[t]{2.3in}
\centerline{\scriptsize (h): t=4.0}
\includegraphics[width=2.3in]{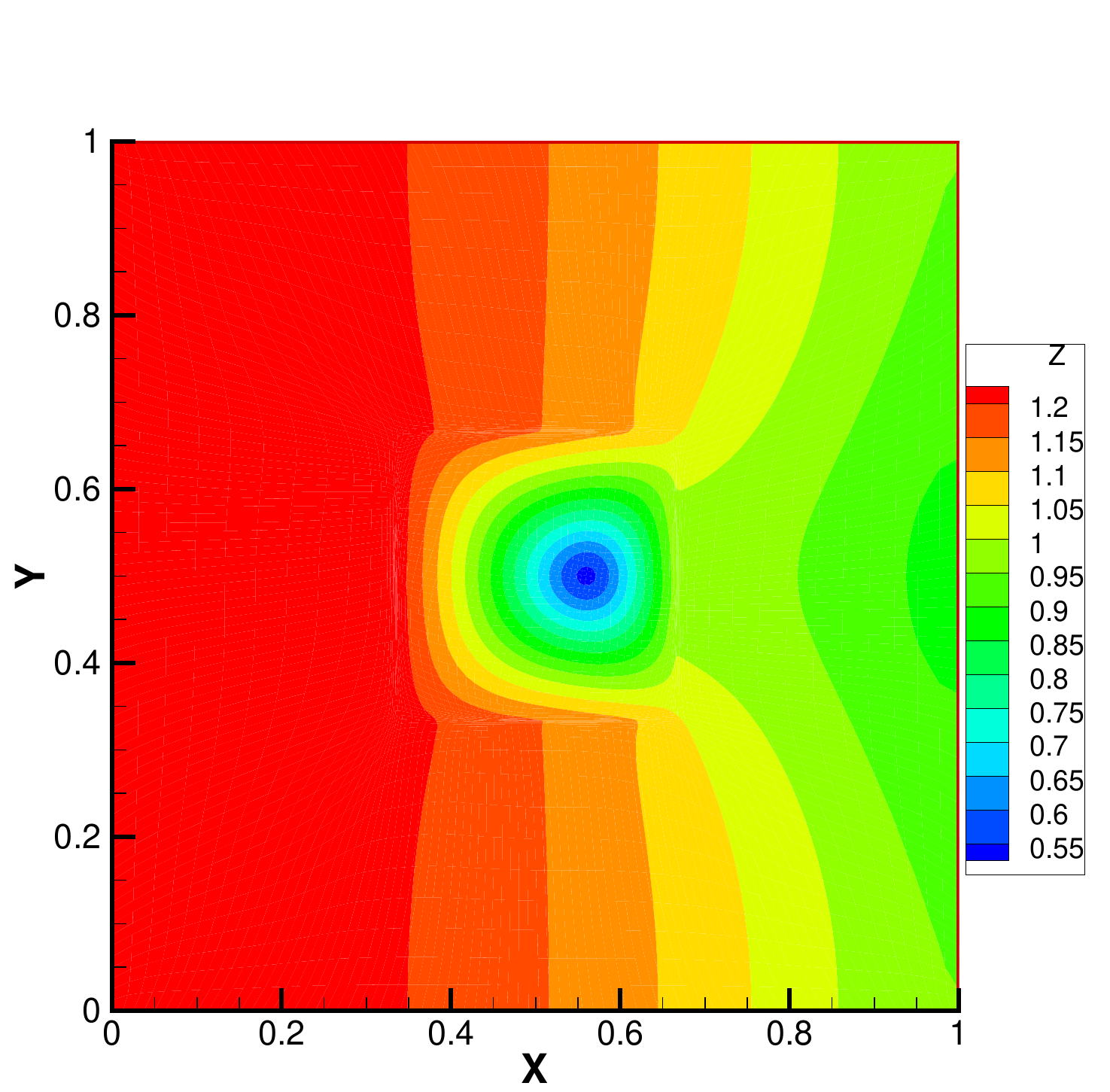}
\end{minipage}
\begin{minipage}[t]{2.3in}
\centerline{\scriptsize (i): t=5.0}
\includegraphics[width=2.3in]{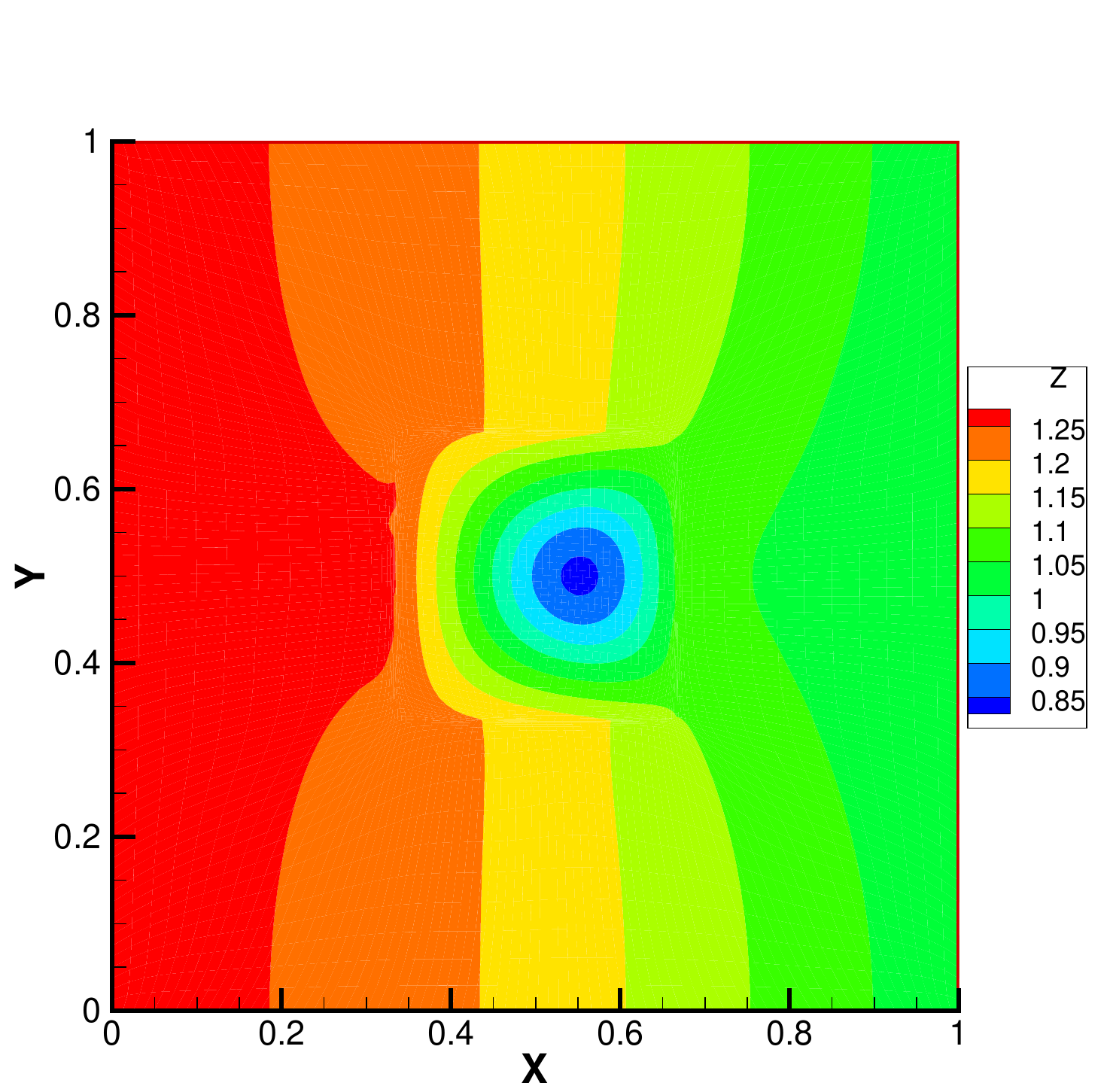}
\end{minipage}
}
\end{center}
\caption{Example~\ref{Example4.2}. The computed solution with a moving mesh of $81\times 81$
is shown at $t=1.0, 1.5, 2.0, 2.4, 2.8, 3.0, 3.5, 4.0, 5.0$.}
\label{T9}
\end{figure}

\begin{figure}
\begin{center}
\hbox{
\begin{minipage}[t]{2.3in}
\centerline{\scriptsize (a): t=1.0}
\includegraphics[width=2.3in]{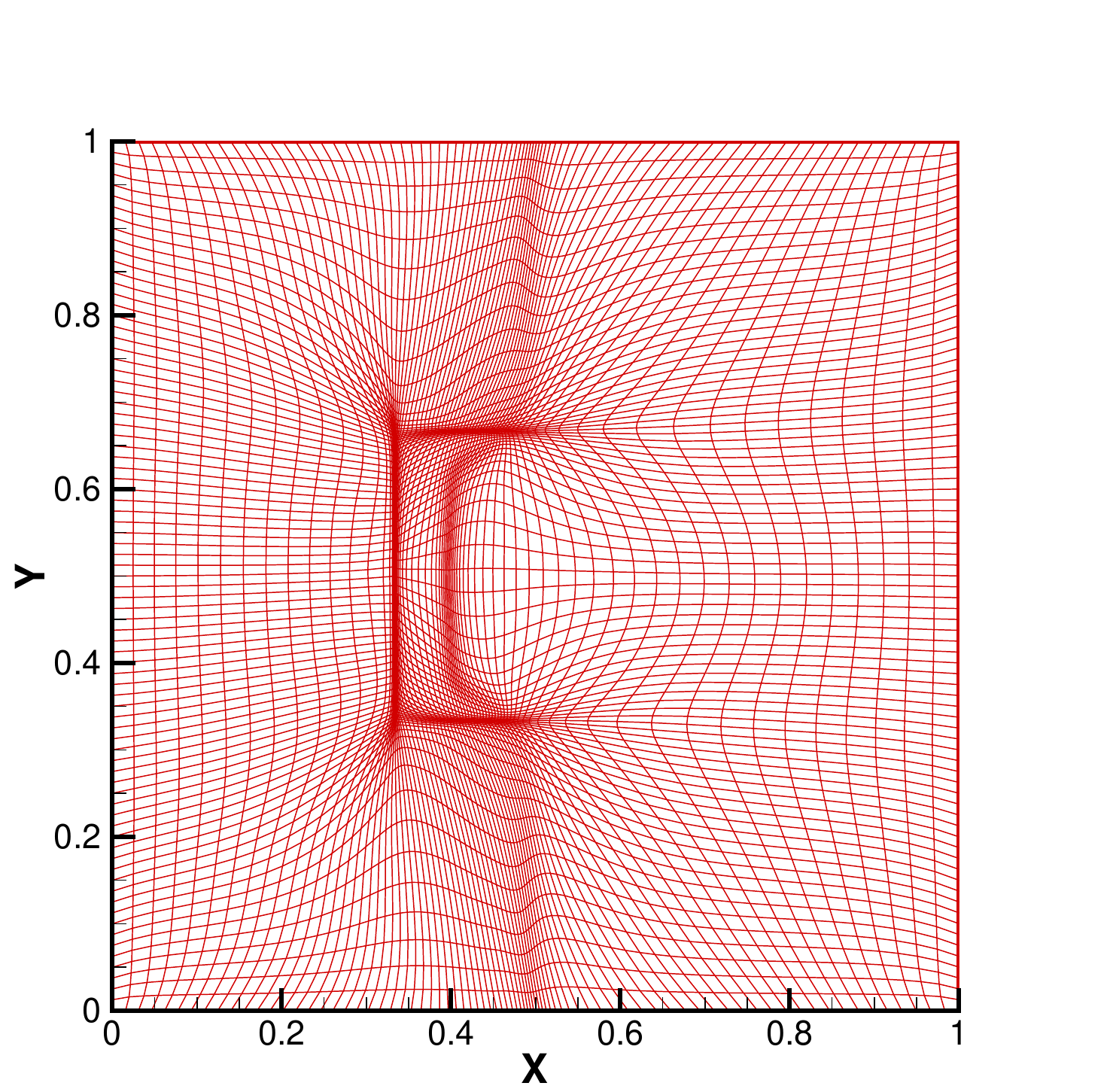}
\end{minipage}
\begin{minipage}[t]{2.3in}
\centerline{\scriptsize (b): t=1.5}
\includegraphics[width=2.3in]{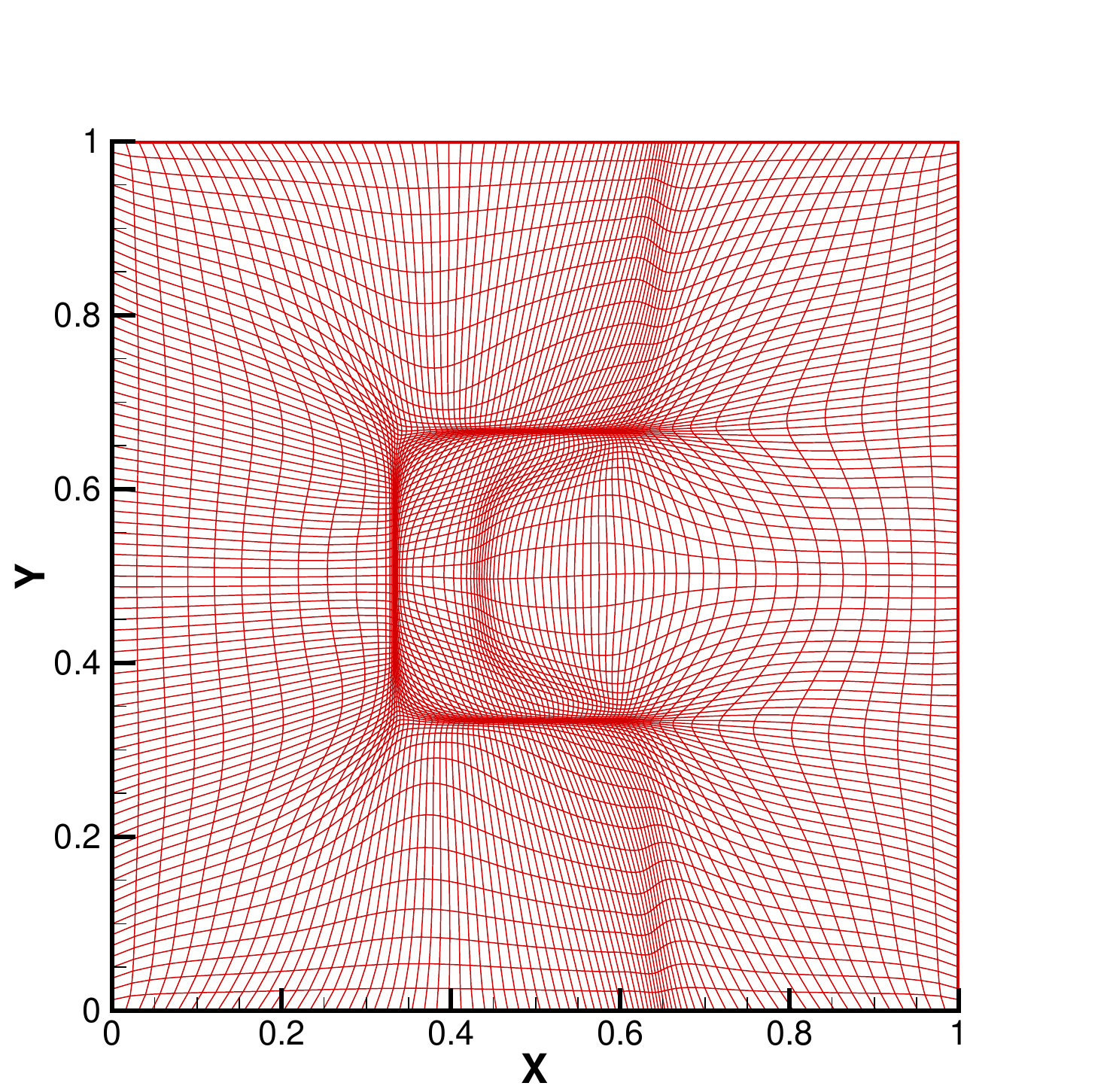}
\end{minipage}
\begin{minipage}[t]{2.3in}
\centerline{\scriptsize (c): t=2.0}
\includegraphics[width=2.3in]{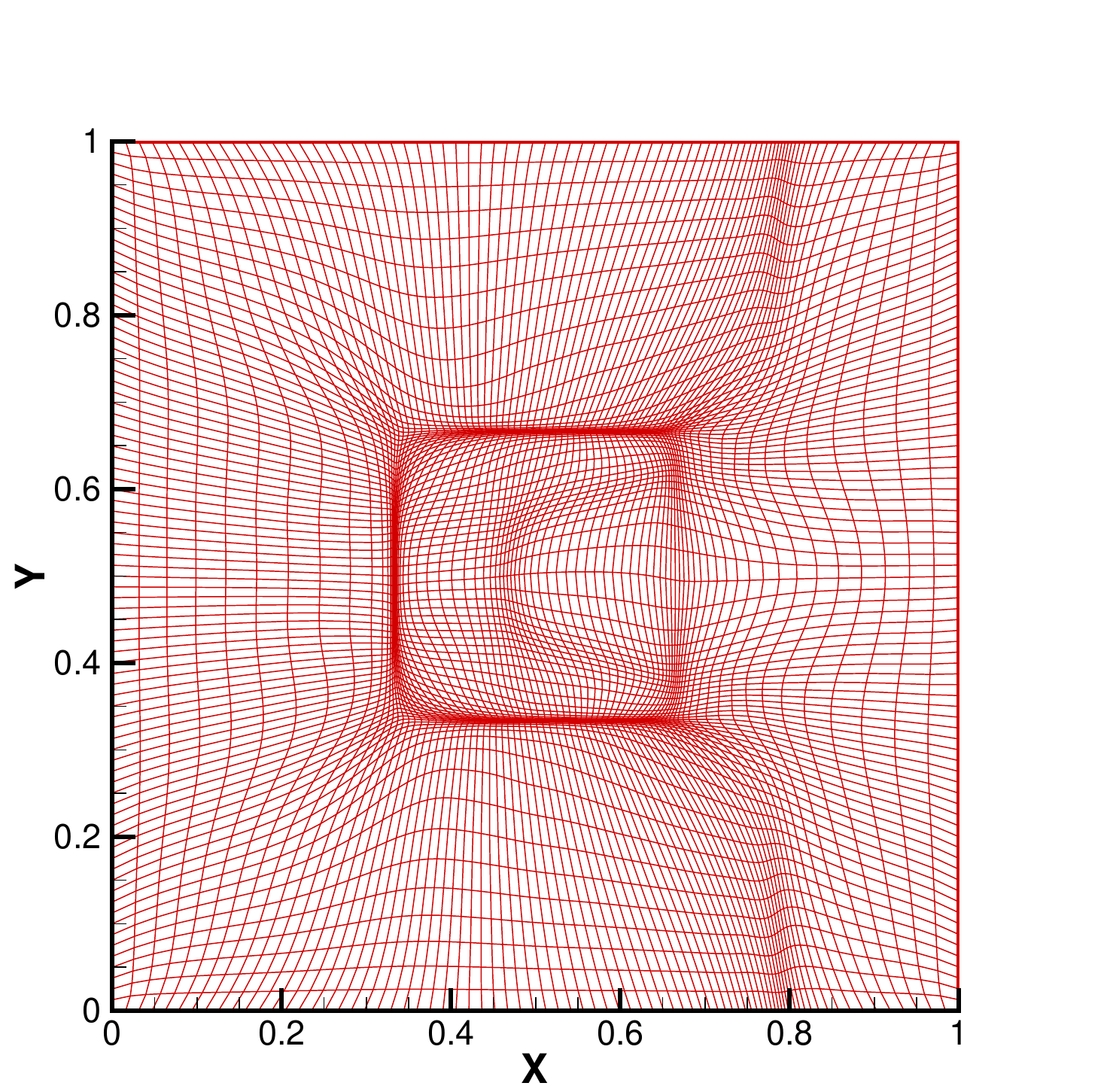}
\end{minipage}
}
\vspace{5mm}
\hbox{
\begin{minipage}[t]{2.3in}
\centerline{\scriptsize (d): t=2.4}
\includegraphics[width=2.3in]{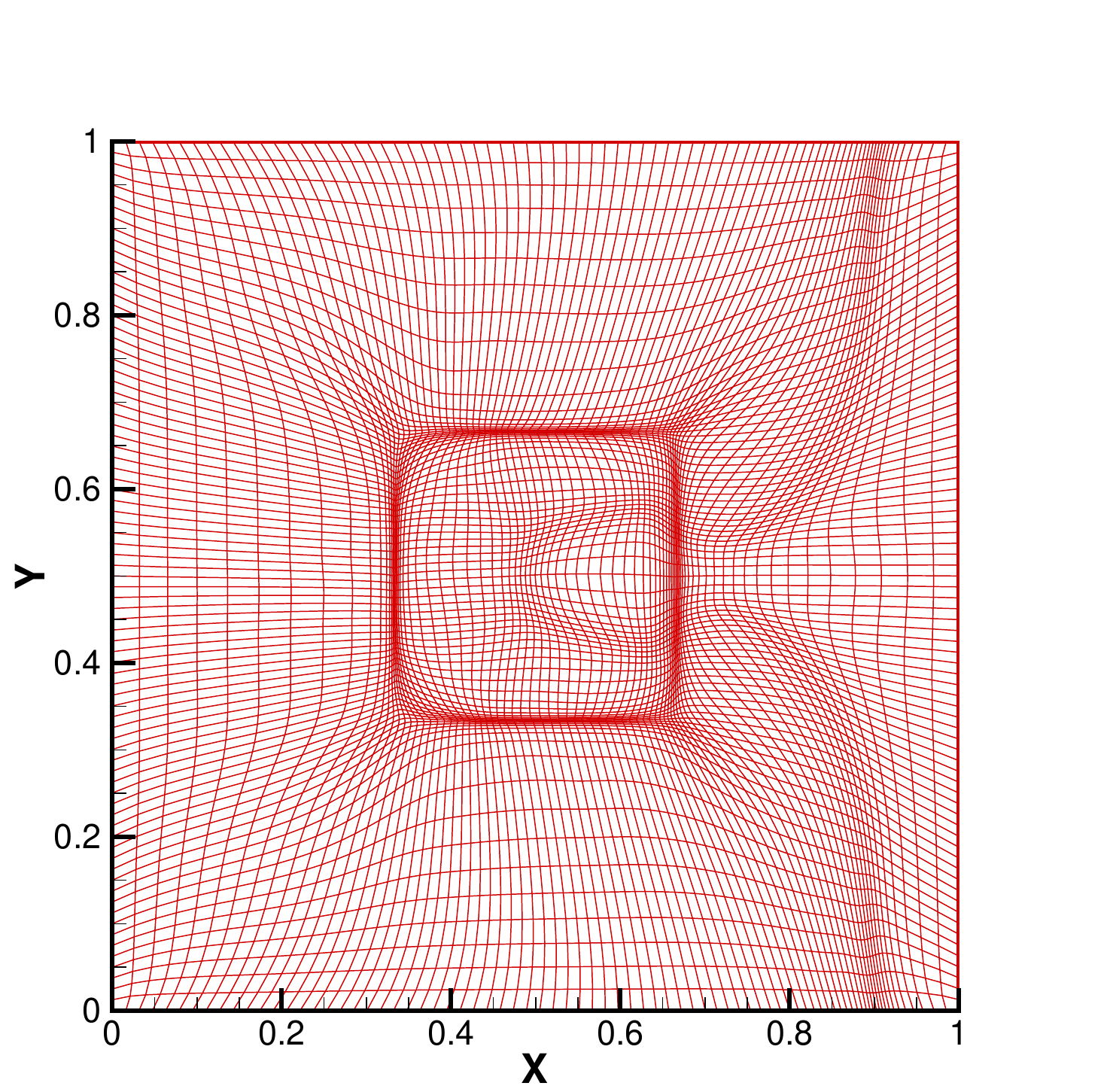}
\end{minipage}
\begin{minipage}[t]{2.3in}
\centerline{\scriptsize (e): t=2.8}
\includegraphics[width=2.3in]{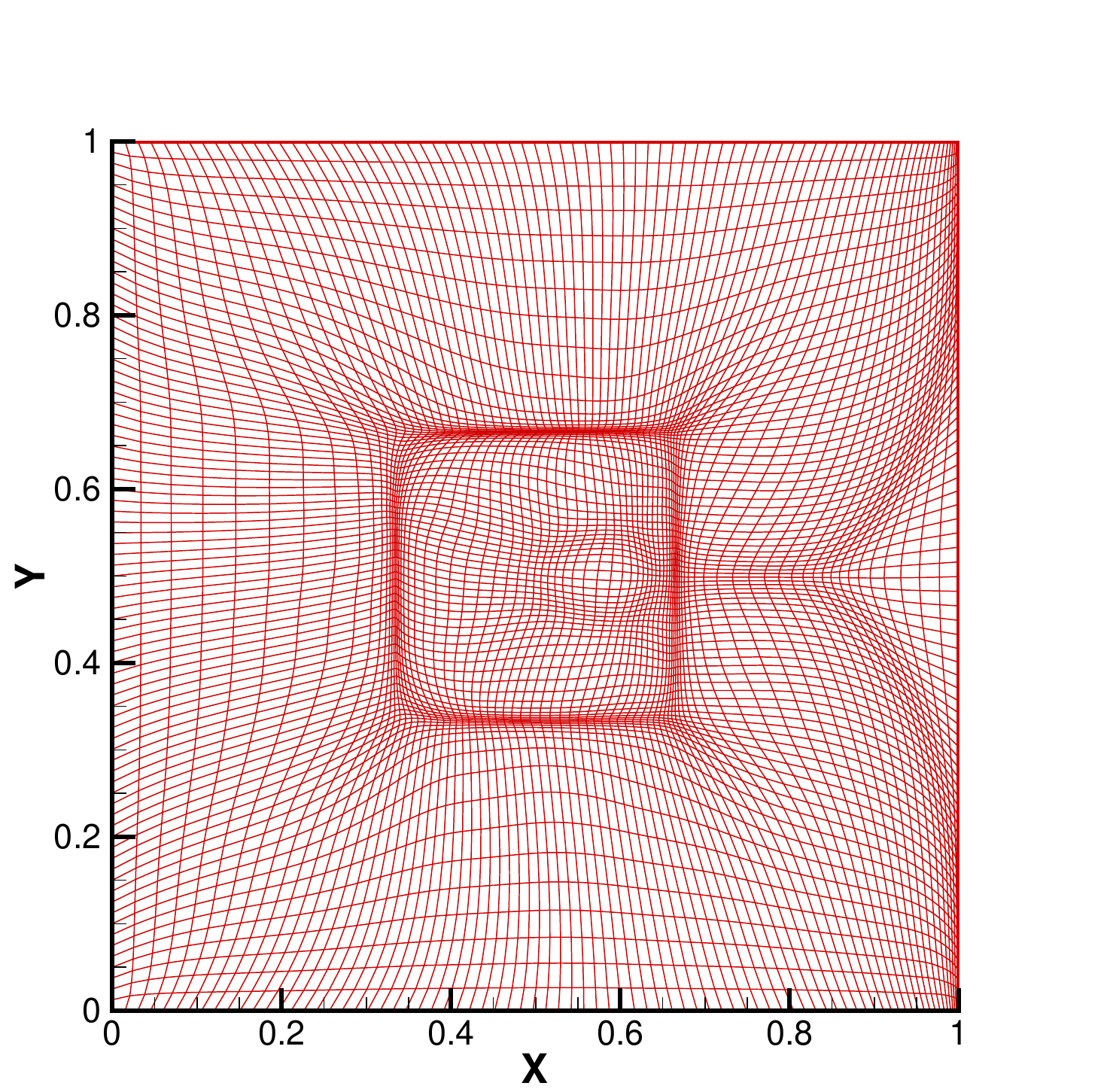}
\end{minipage}
\begin{minipage}[t]{2.3in}
\centerline{\scriptsize (f): t=3.0}
\includegraphics[width=2.3in]{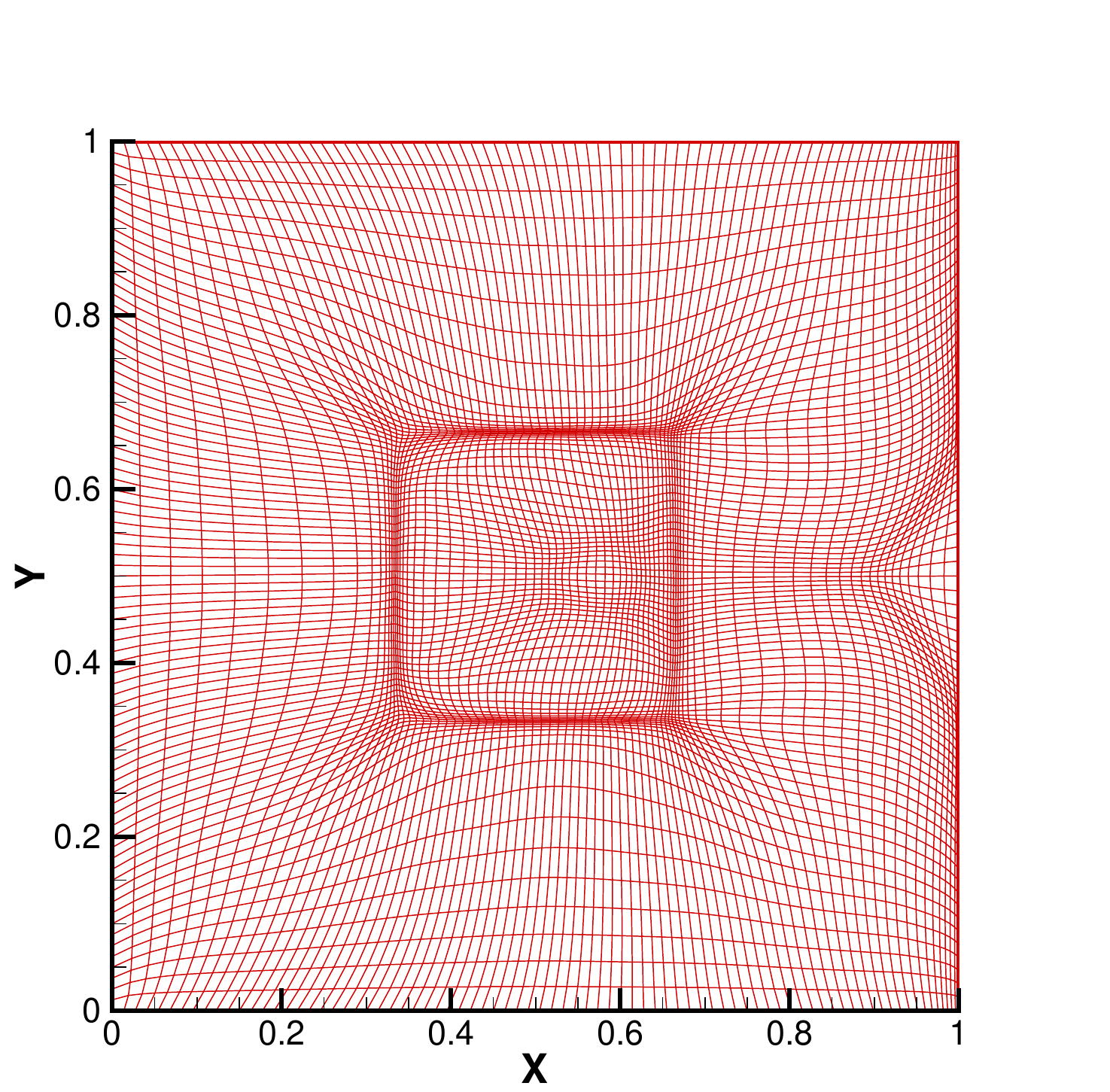}
\end{minipage}
}
\vspace{5mm}
\hbox{
\begin{minipage}[t]{2.3in}
\centerline{\scriptsize (g): t=3.5}
\includegraphics[width=2.3in]{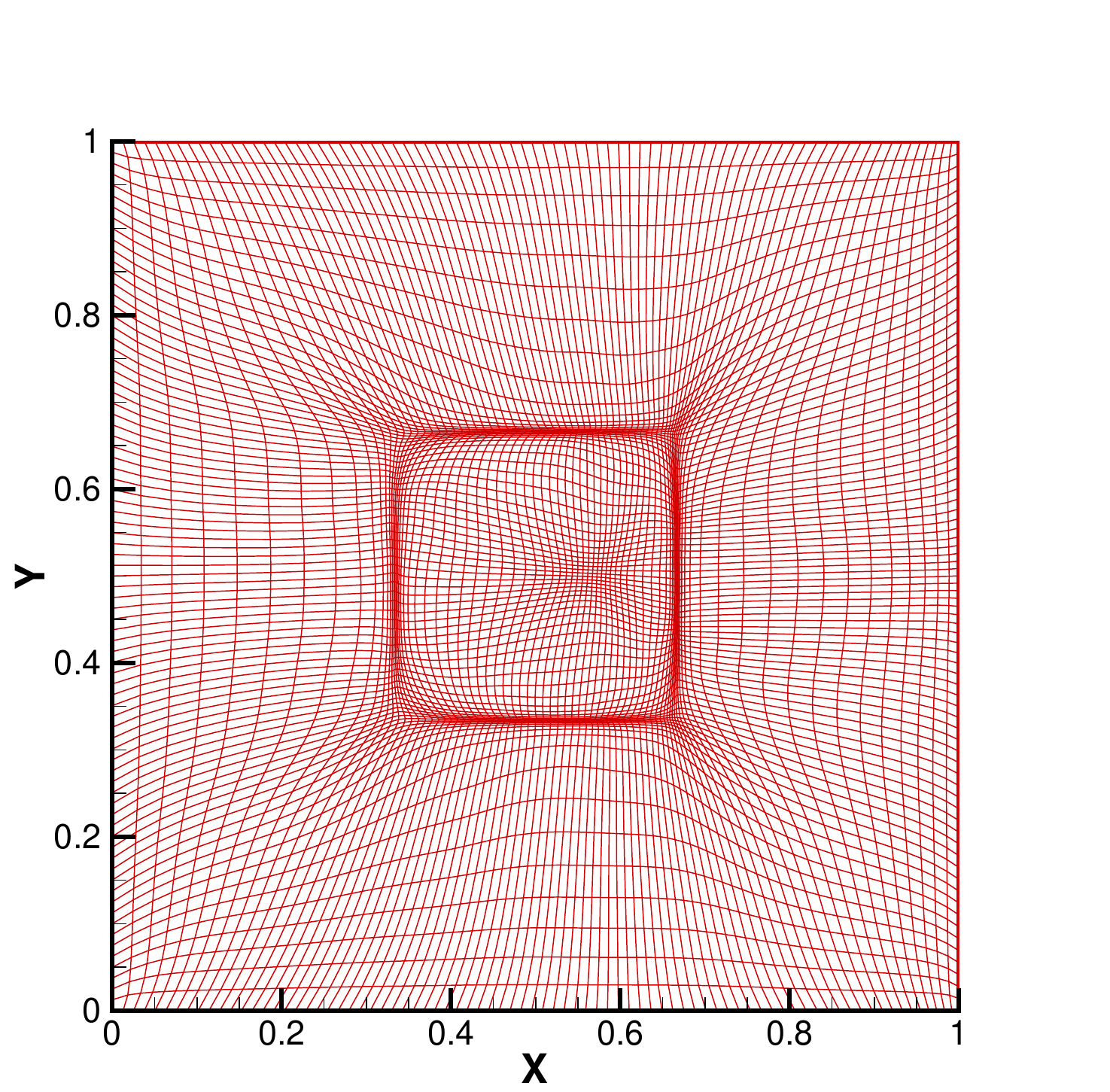}
\end{minipage}
\begin{minipage}[t]{2.3in}
\centerline{\scriptsize (h): t=4.0}
\includegraphics[width=2.3in]{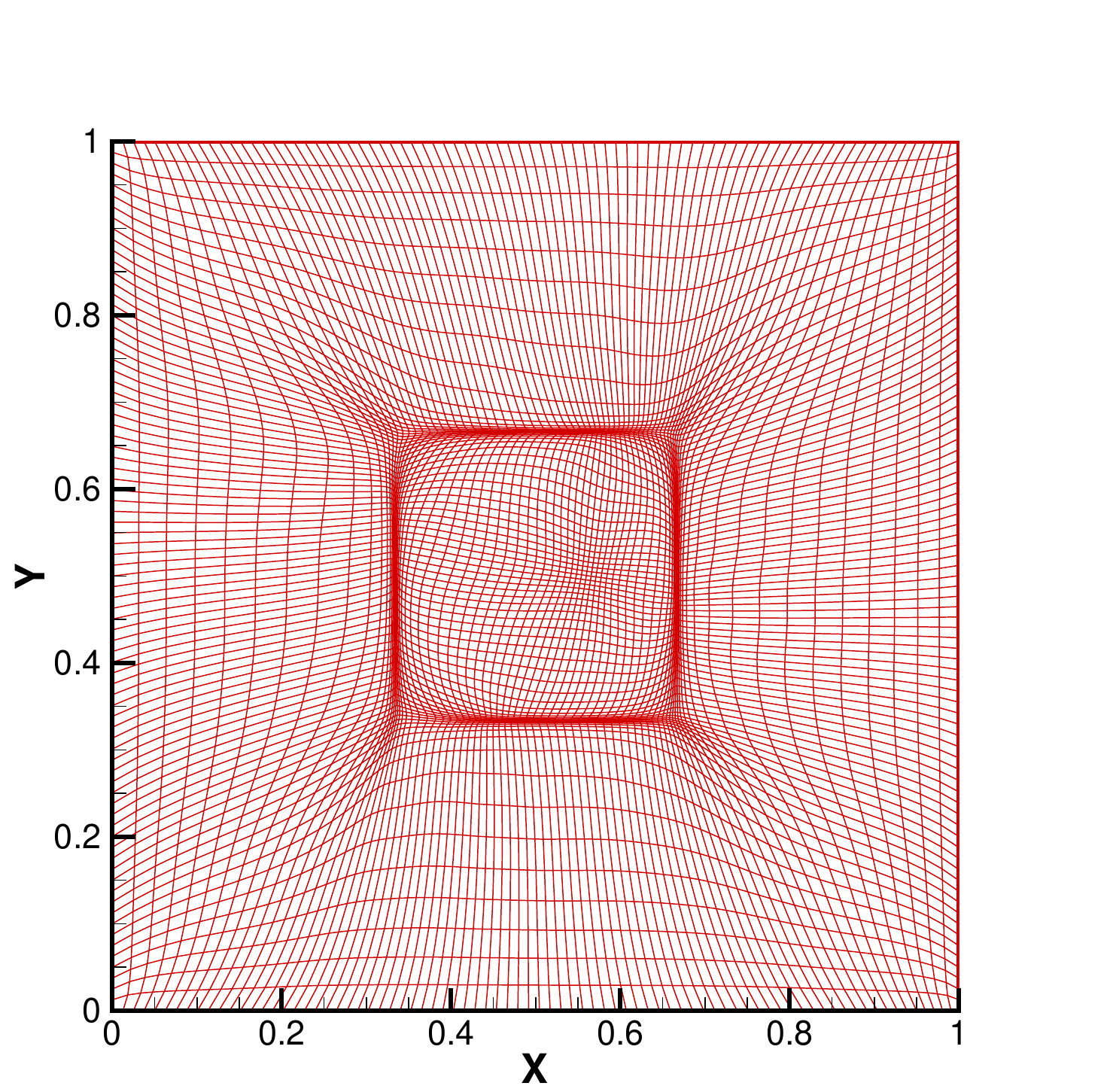}
\end{minipage}
\begin{minipage}[t]{2.3in}
\centerline{\scriptsize (i): t=5.0}
\includegraphics[width=2.3in]{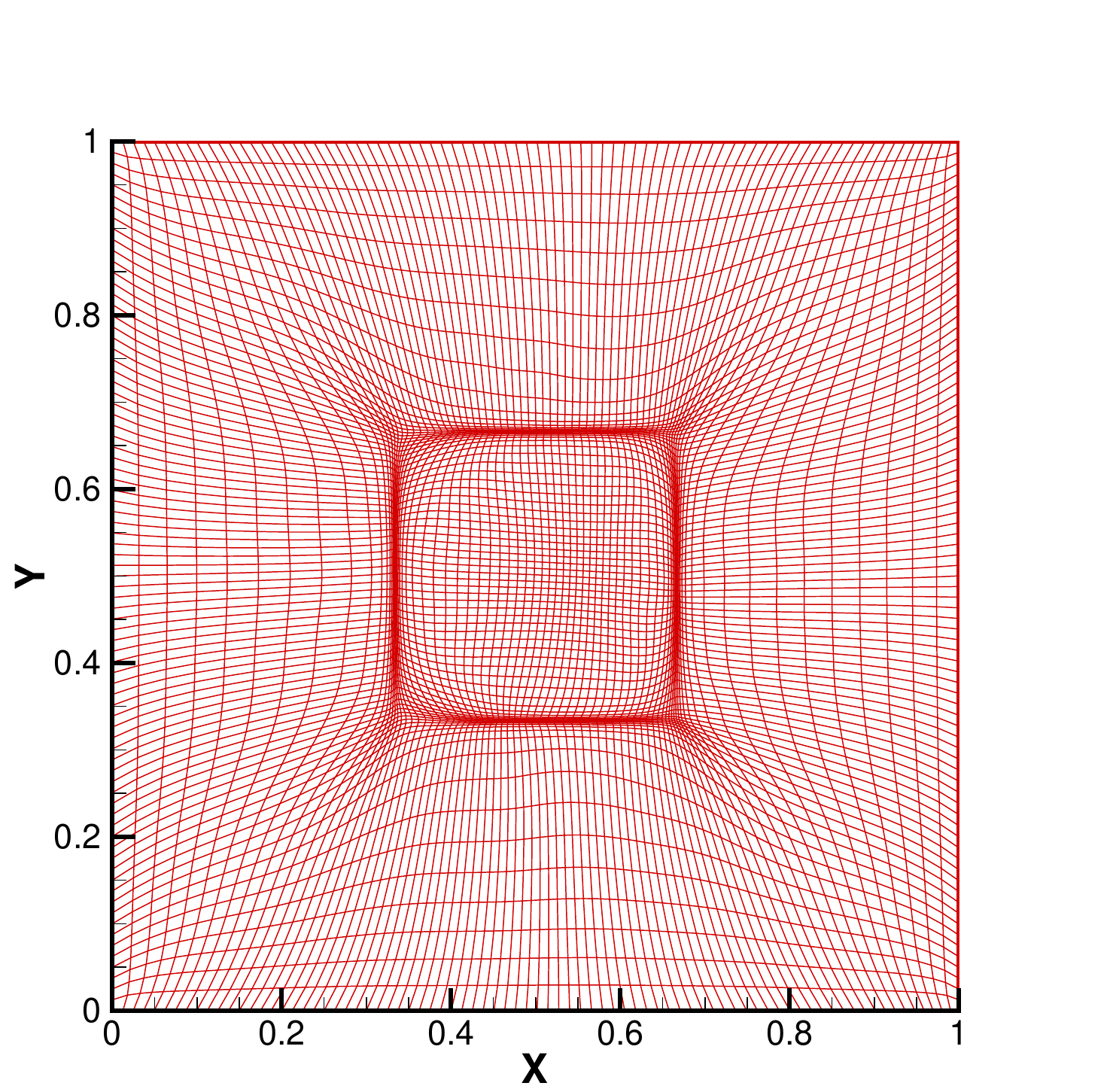}
\end{minipage}
}
\end{center}
\caption{Example~\ref{Example4.2}. The moving mesh of $81\times 81$
is shown at $t=1.0, 1.5, 2.0, 2.4, 2.8, 3.0, 3.5, 4.0, 5.0$.}
\label{T10}
\end{figure}

\begin{figure}
\begin{center}
\hbox{
\hspace{1in}
\begin{minipage}[t]{2.0in}
\centerline{\scriptsize (a):  with MM at $t=1.0$}
\includegraphics[width=2.0in]{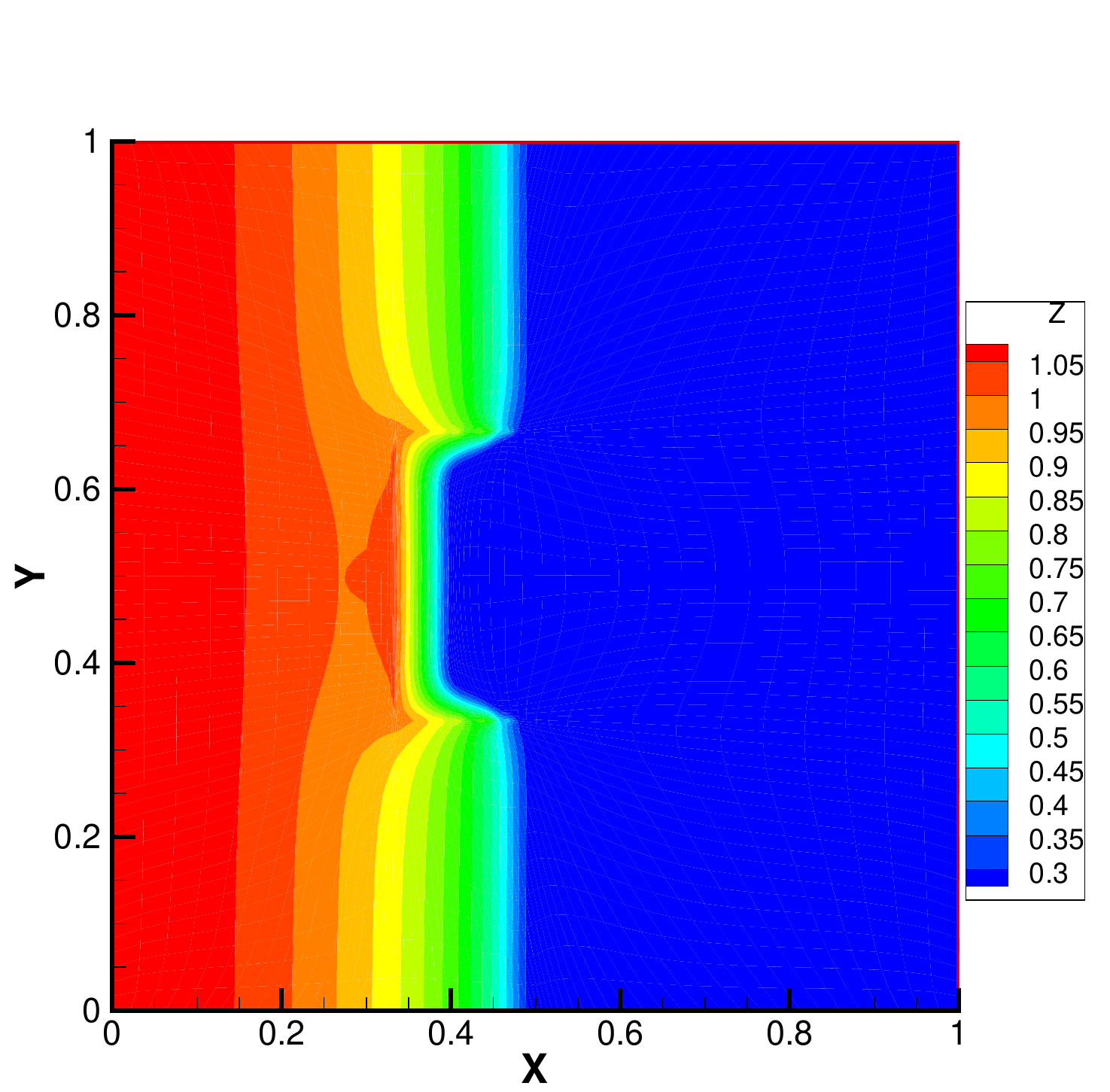}
\end{minipage}
\hspace{0.5in}
\begin{minipage}[t]{2.0in}
\centerline{\scriptsize (b):  with UM at $t=1.0$}
\includegraphics[width=2.0in]{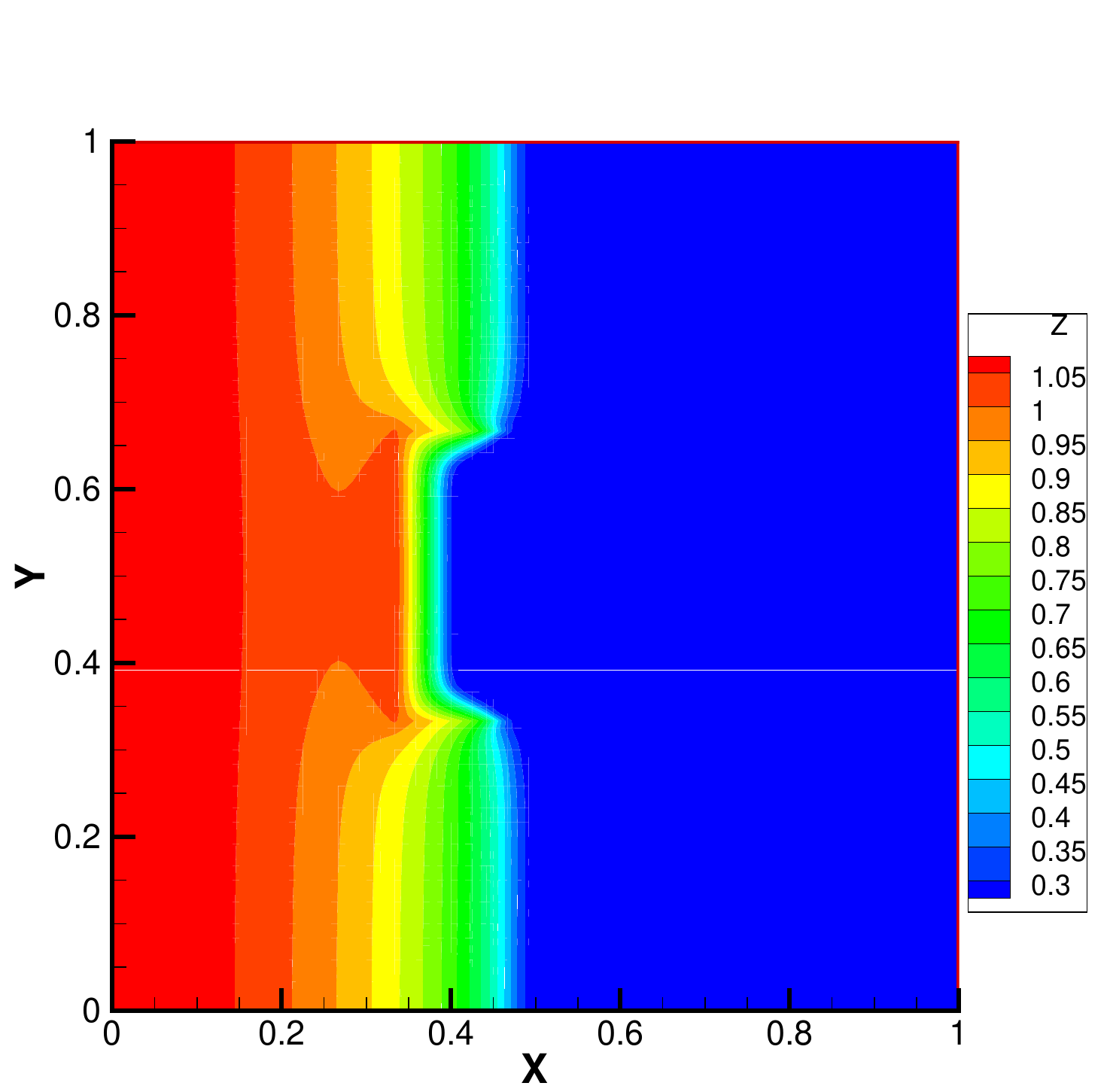}
\end{minipage}
}
\hbox{
\hspace{1in}
\begin{minipage}[t]{2.0in}
\centerline{\scriptsize (c):  with MM at $t=2.0$}
\includegraphics[width=2.0in]{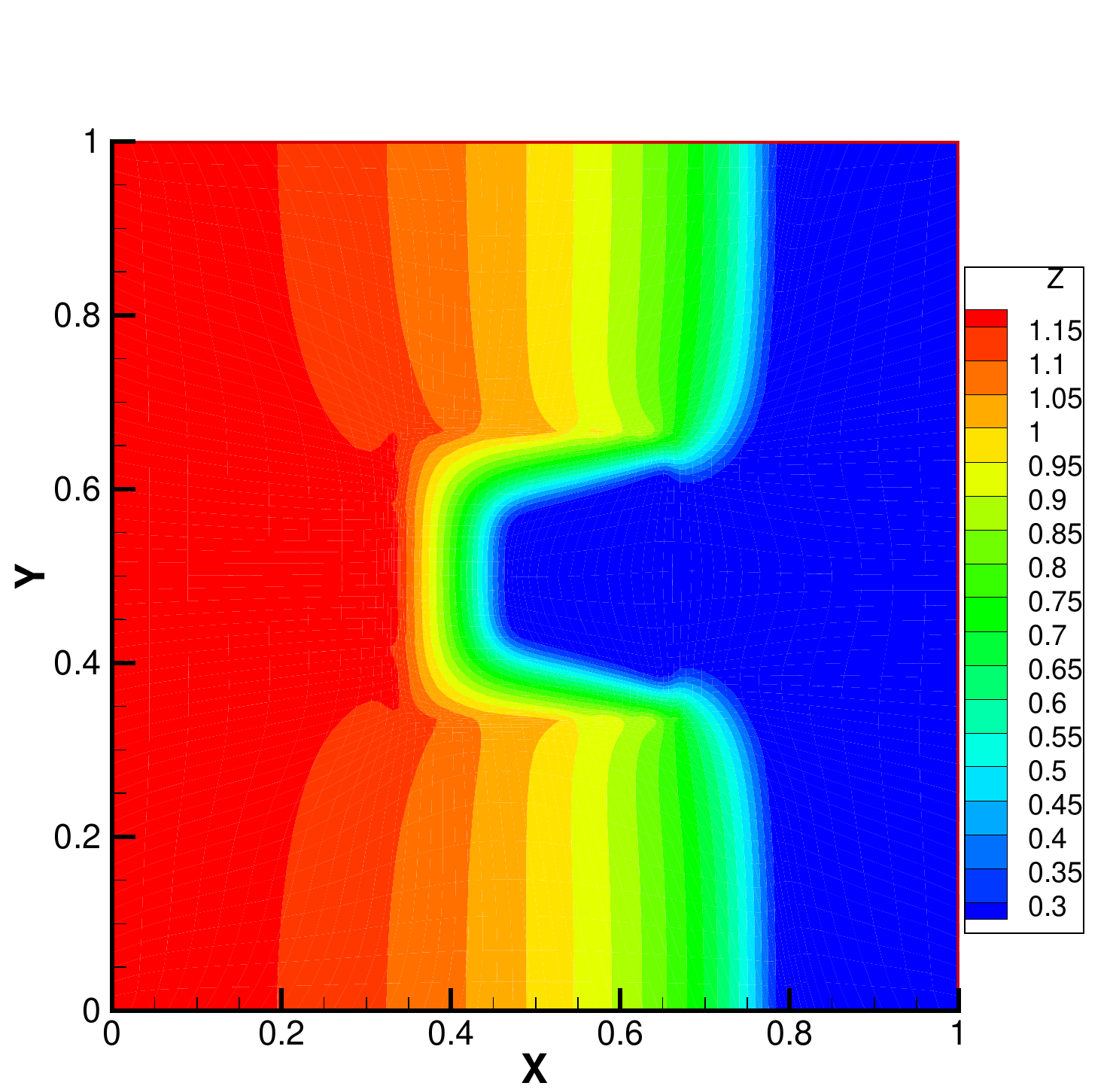}
\end{minipage}
\hspace{0.5in}
\begin{minipage}[t]{2.0in}
\centerline{\scriptsize (d):  with UM at $t=2.0$}
\includegraphics[width=2.0in]{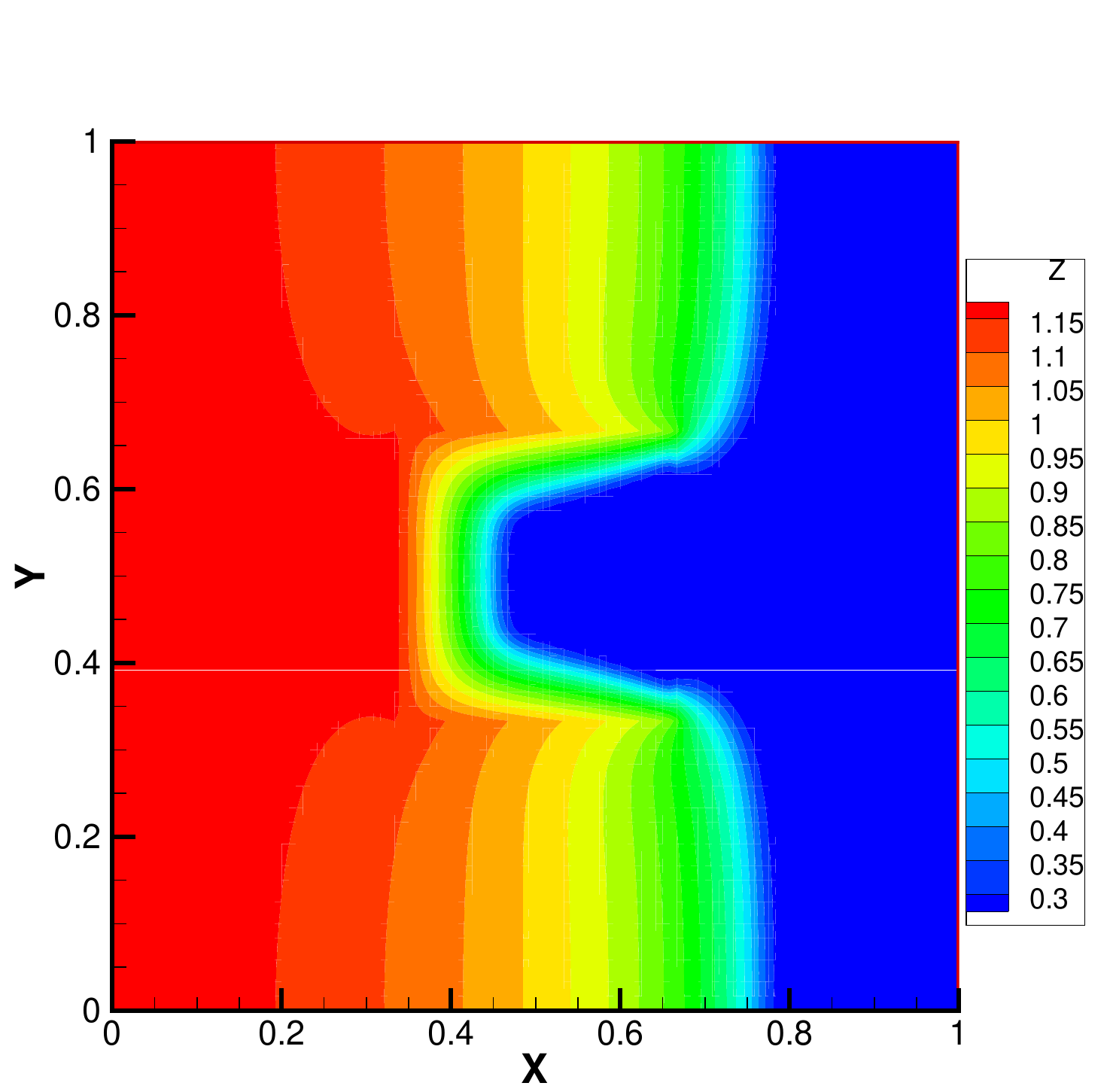}
\end{minipage}
}
\hbox{
\hspace{1in}
\begin{minipage}[t]{2.0in}
\centerline{\scriptsize (e):  with MM at $t=2.5$}
\includegraphics[width=2.0in]{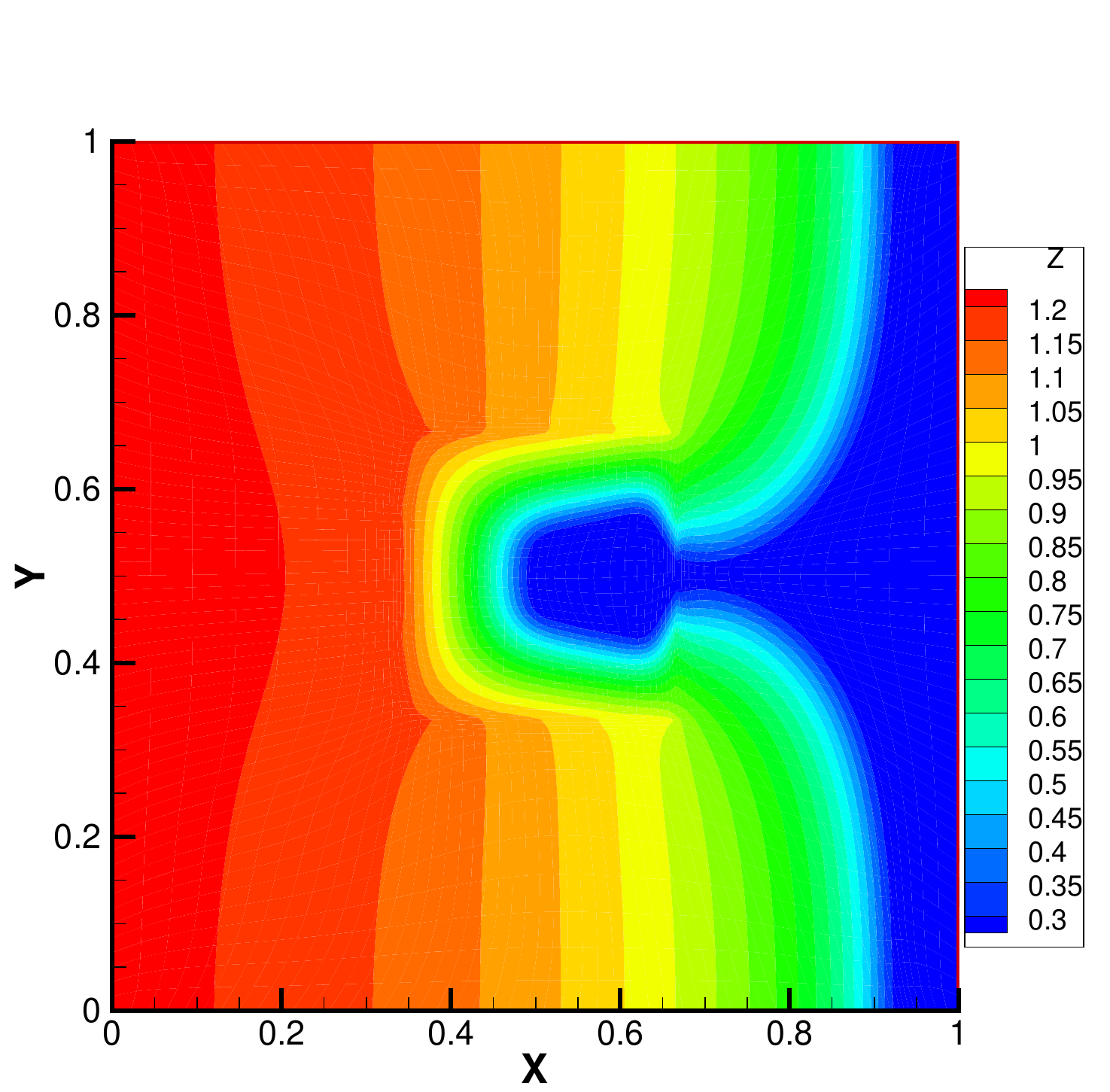}
\end{minipage}
\hspace{0.5in}
\begin{minipage}[t]{2.0in}
\centerline{\scriptsize (f):  with UM at $t=2.5$}
\includegraphics[width=2.0in]{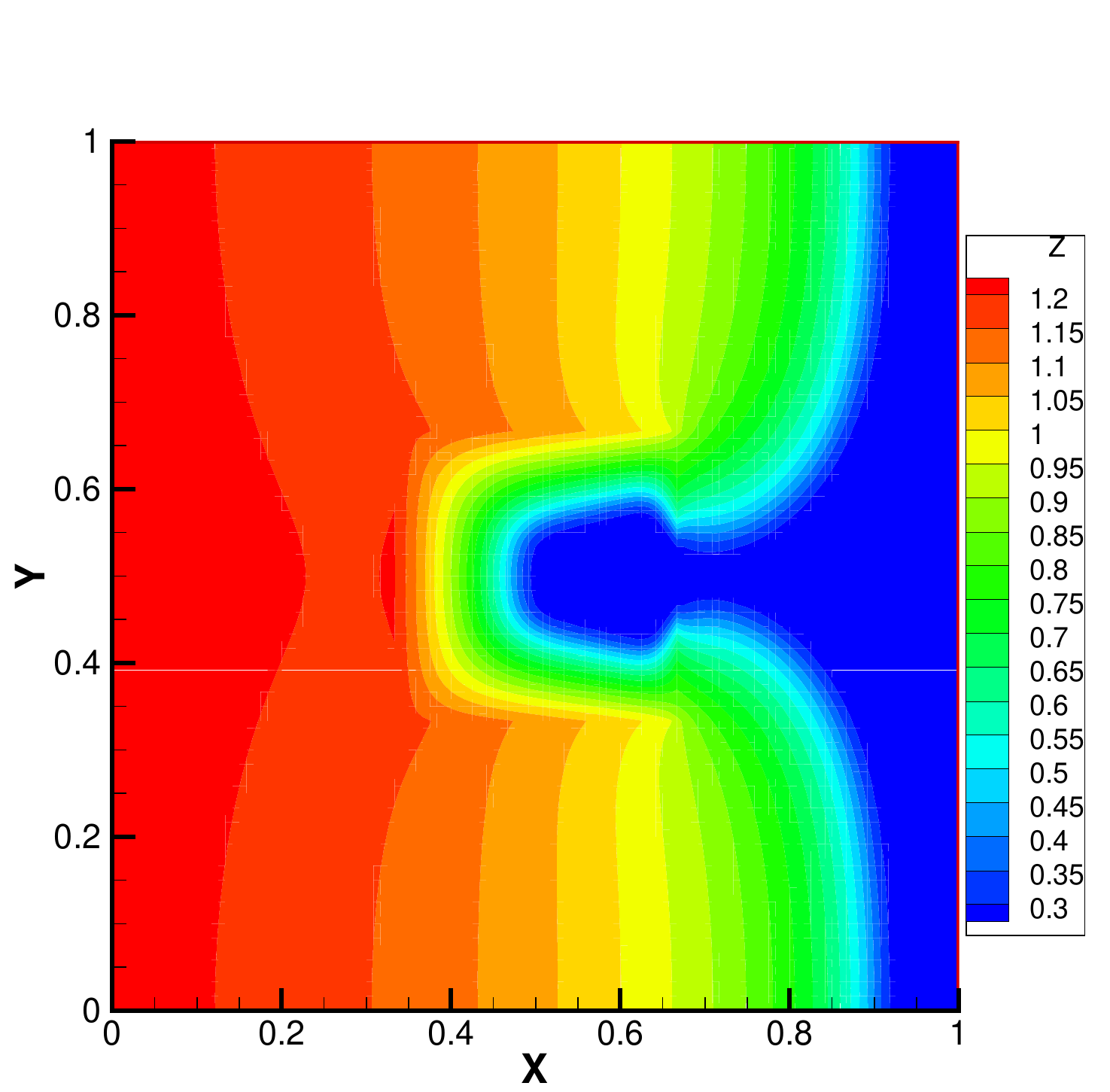}
\end{minipage}
}
\hbox{
\hspace{1in}
\begin{minipage}[t]{2.0in}
\centerline{\scriptsize (g):  with MM at $t=3.0$}
\includegraphics[width=2.0in]{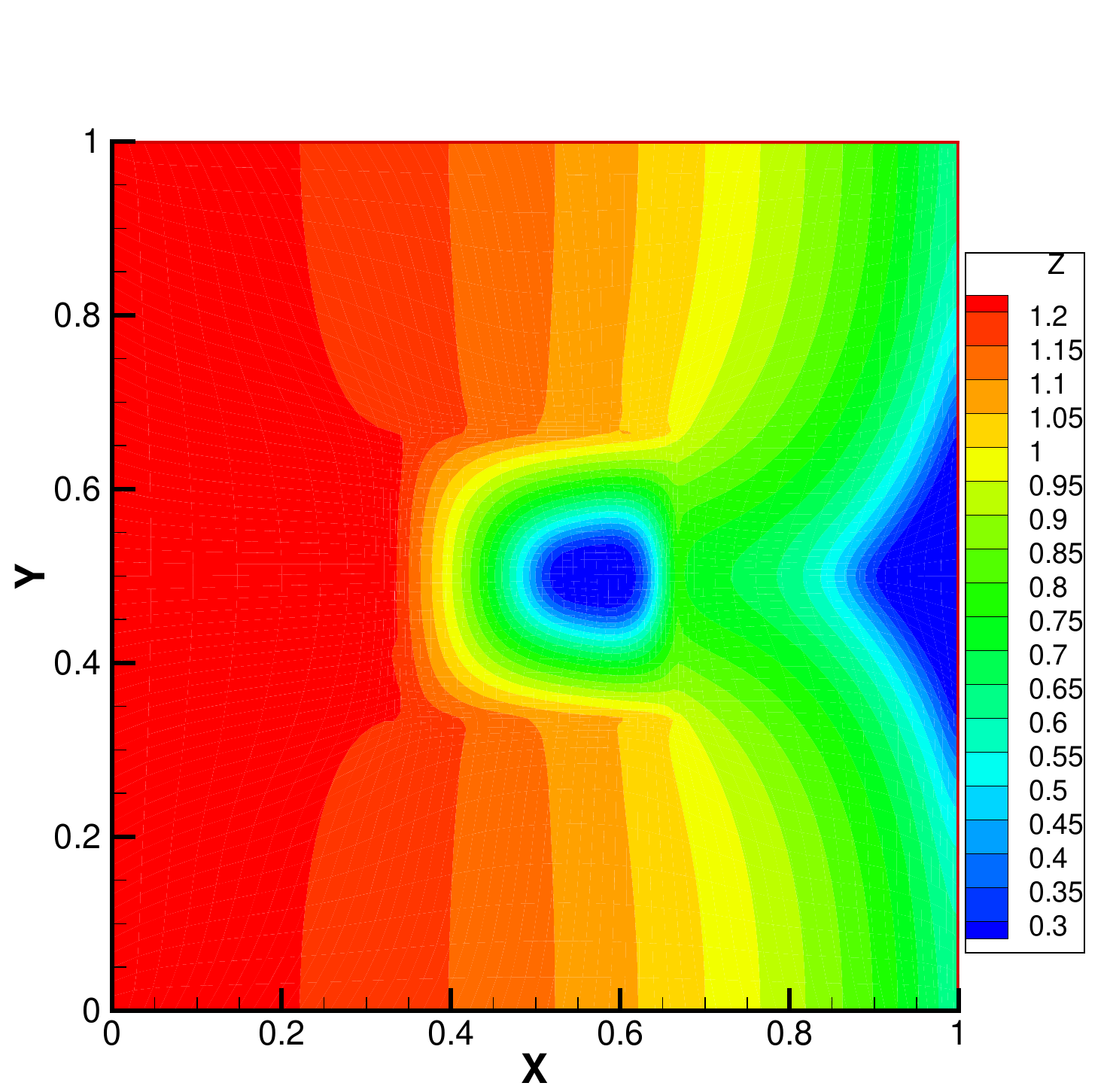}
\end{minipage}
\hspace{0.5in}
\begin{minipage}[t]{2.0in}
\centerline{\scriptsize (h):  with UM at $t=3.0$}
\includegraphics[width=2.0in]{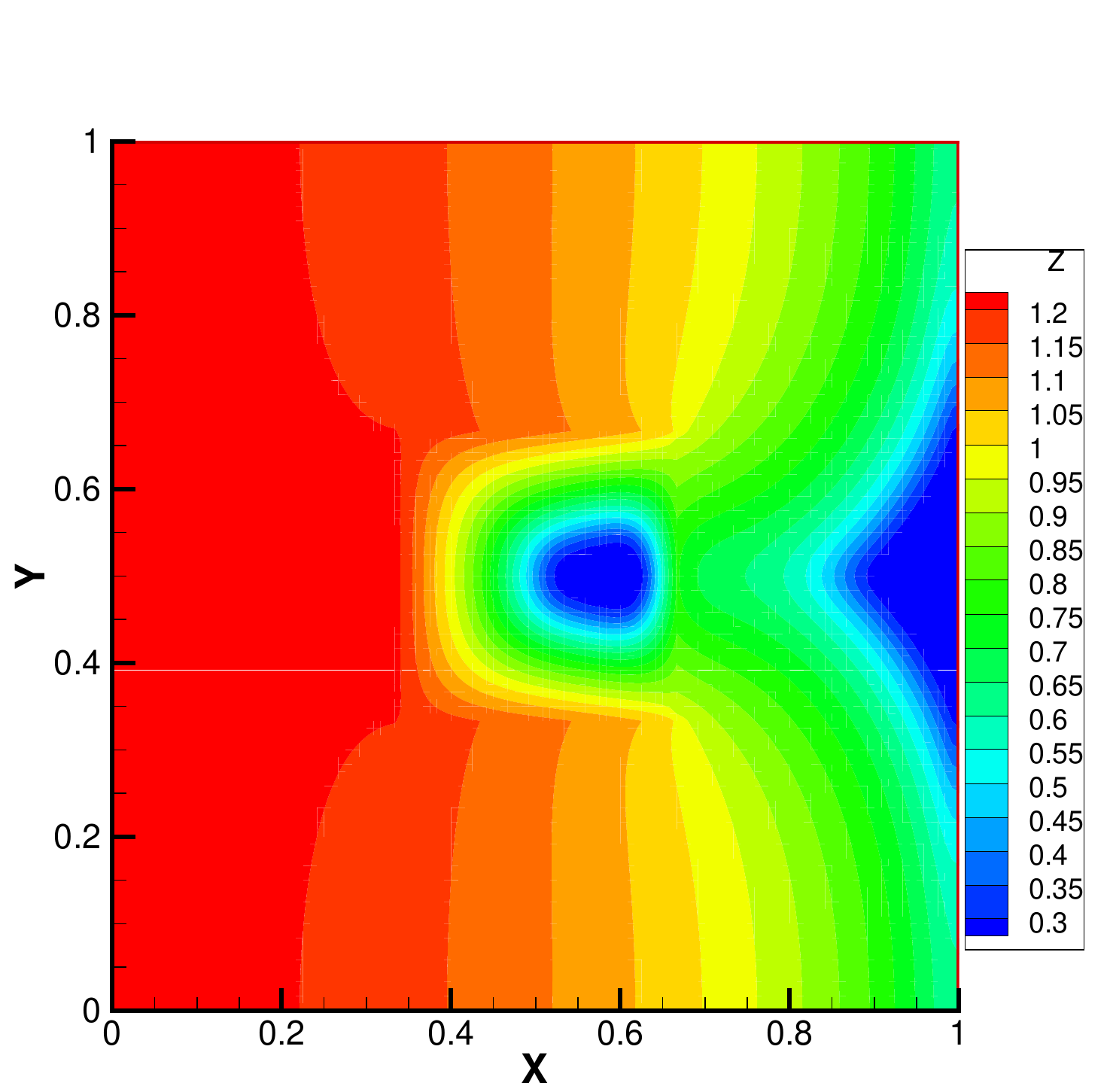}
\end{minipage}
}
\end{center}
\caption{Example~\ref{Example4.2}. The contours of the temperature obtained
with a moving mesh (MM) of size  $61\times 61$ are compared with those
obtained with a uniform mesh (UM) of size $121 \times121$.}
\label{T11}
\end{figure}

\begin{figure}
\begin{center}
\hbox{
\hspace{1in}
\begin{minipage}[t]{2.0in}
\centerline{\scriptsize (a):  with MM1 at $t=1.0$}
\includegraphics[width=2.0in]{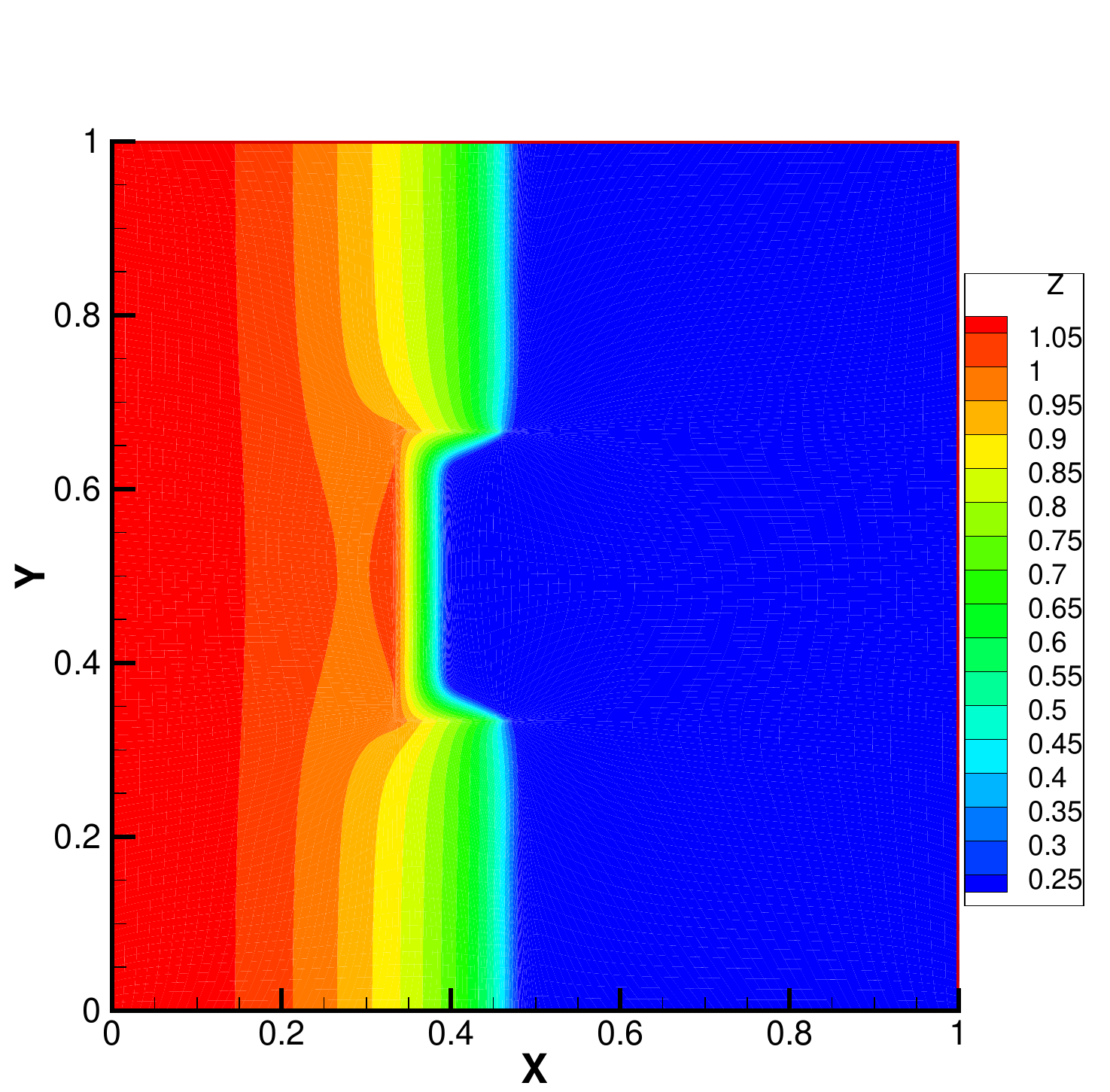}
\end{minipage}
\hspace{0.5in}
\begin{minipage}[t]{2.0in}
\centerline{\scriptsize (b):  with MM2 at $t=1.0$}
\includegraphics[width=2.0in]{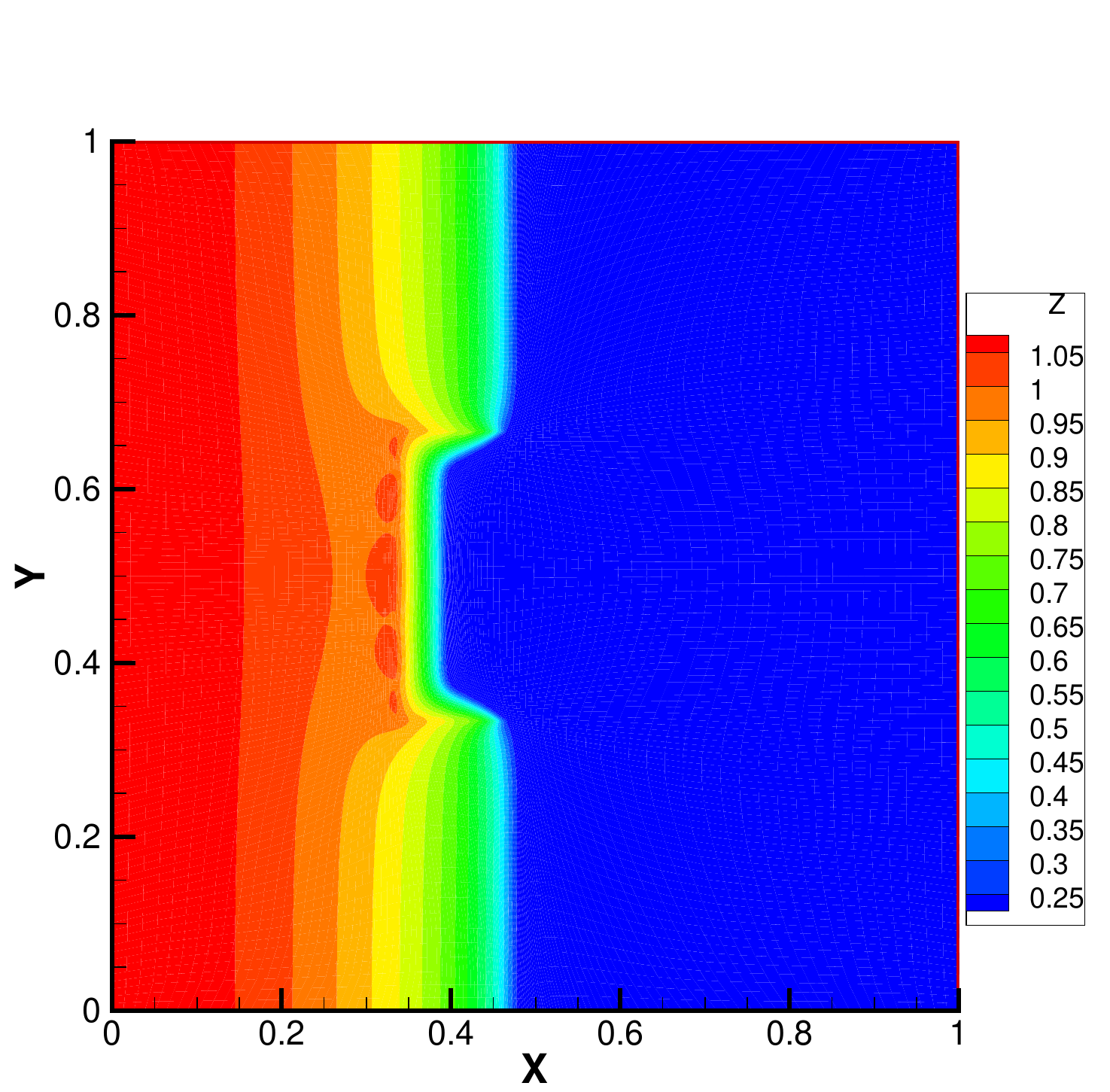}
\end{minipage}
}
\hbox{
\hspace{1in}
\begin{minipage}[t]{2.0in}
\centerline{\scriptsize (c):  with MM1 at $t=2.0$}
\includegraphics[width=2.0in]{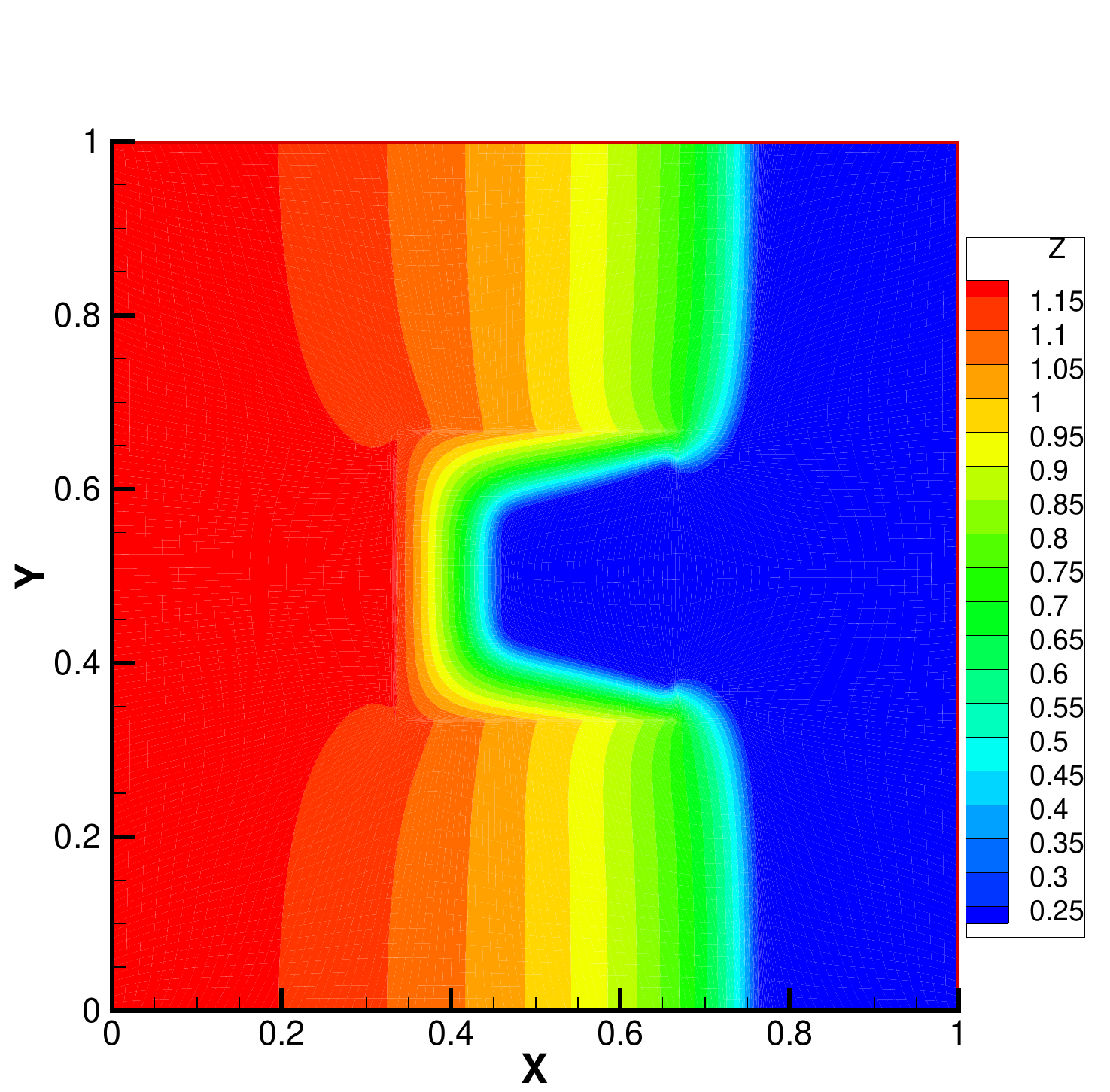}
\end{minipage}
\hspace{0.5in}
\begin{minipage}[t]{2.0in}
\centerline{\scriptsize (d):  with MM2 at $t=2.0$}
\includegraphics[width=2.0in]{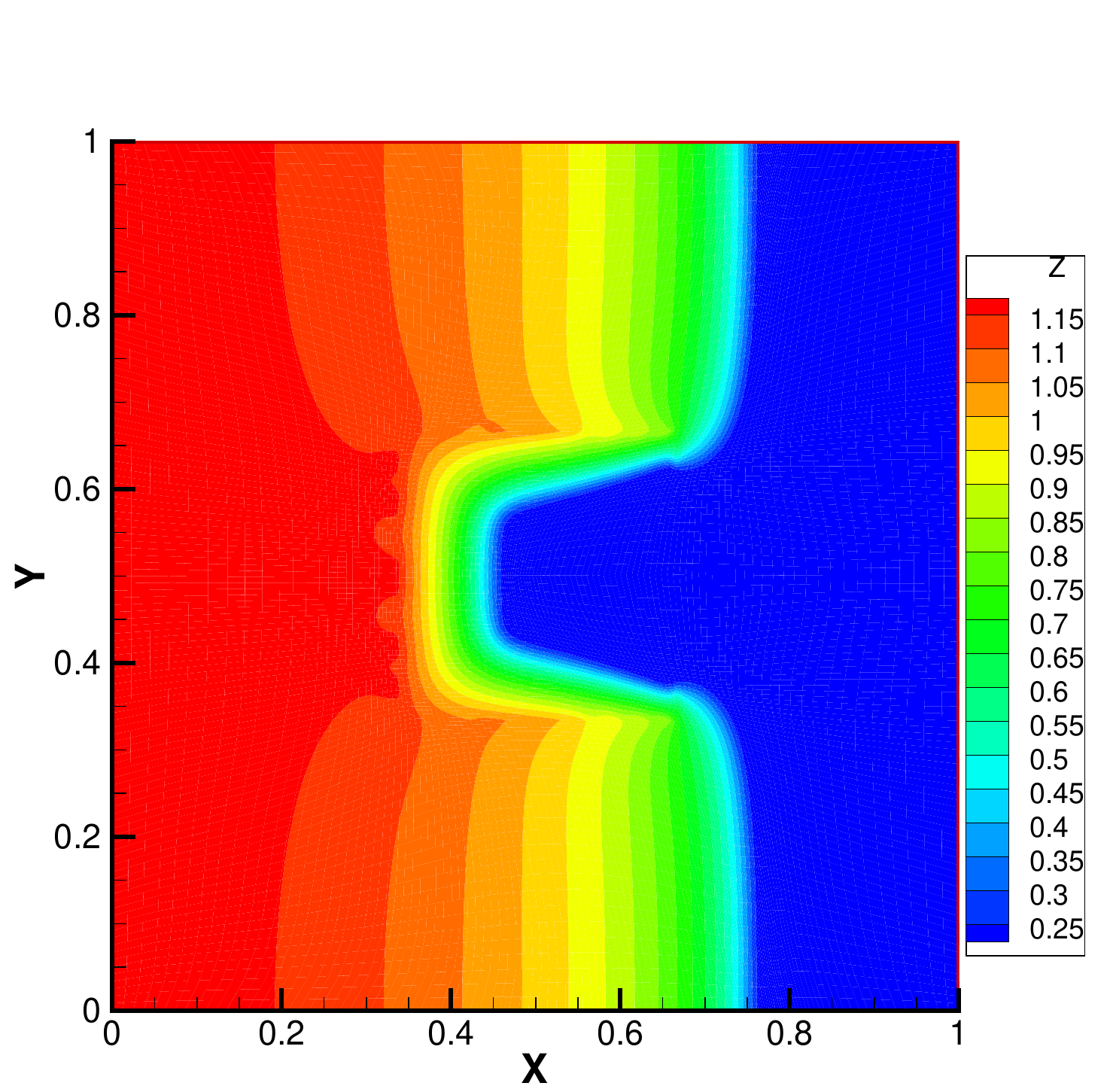}
\end{minipage}
}
\hbox{
\hspace{1in}
\begin{minipage}[t]{2.0in}
\centerline{\scriptsize (e):  with MM1 at $t=2.5$}
\includegraphics[width=2.0in]{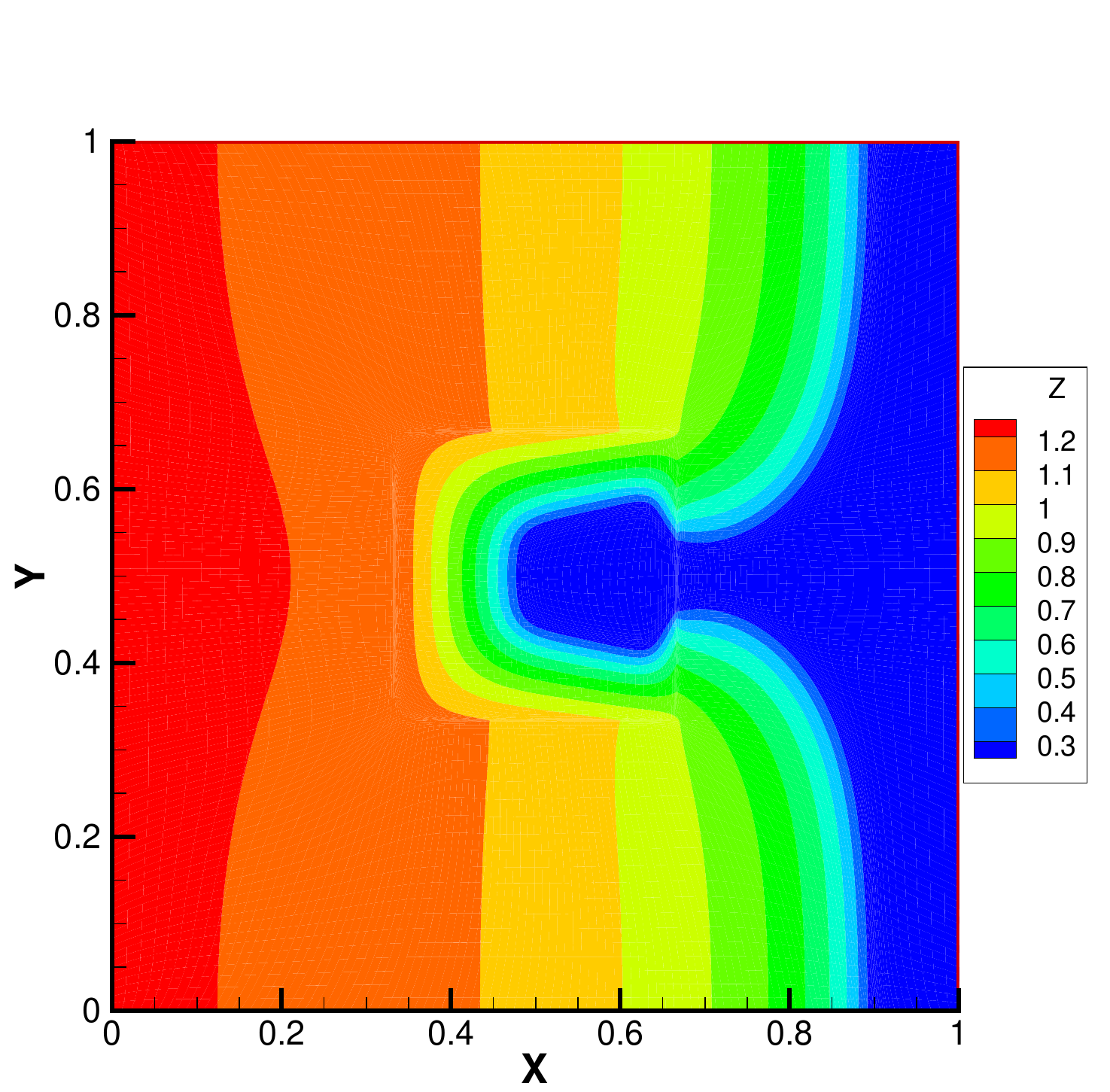}
\end{minipage}
\hspace{0.5in}
\begin{minipage}[t]{2.0in}
\centerline{\scriptsize (f):  with MM2 at $t=2.5$}
\includegraphics[width=2.0in]{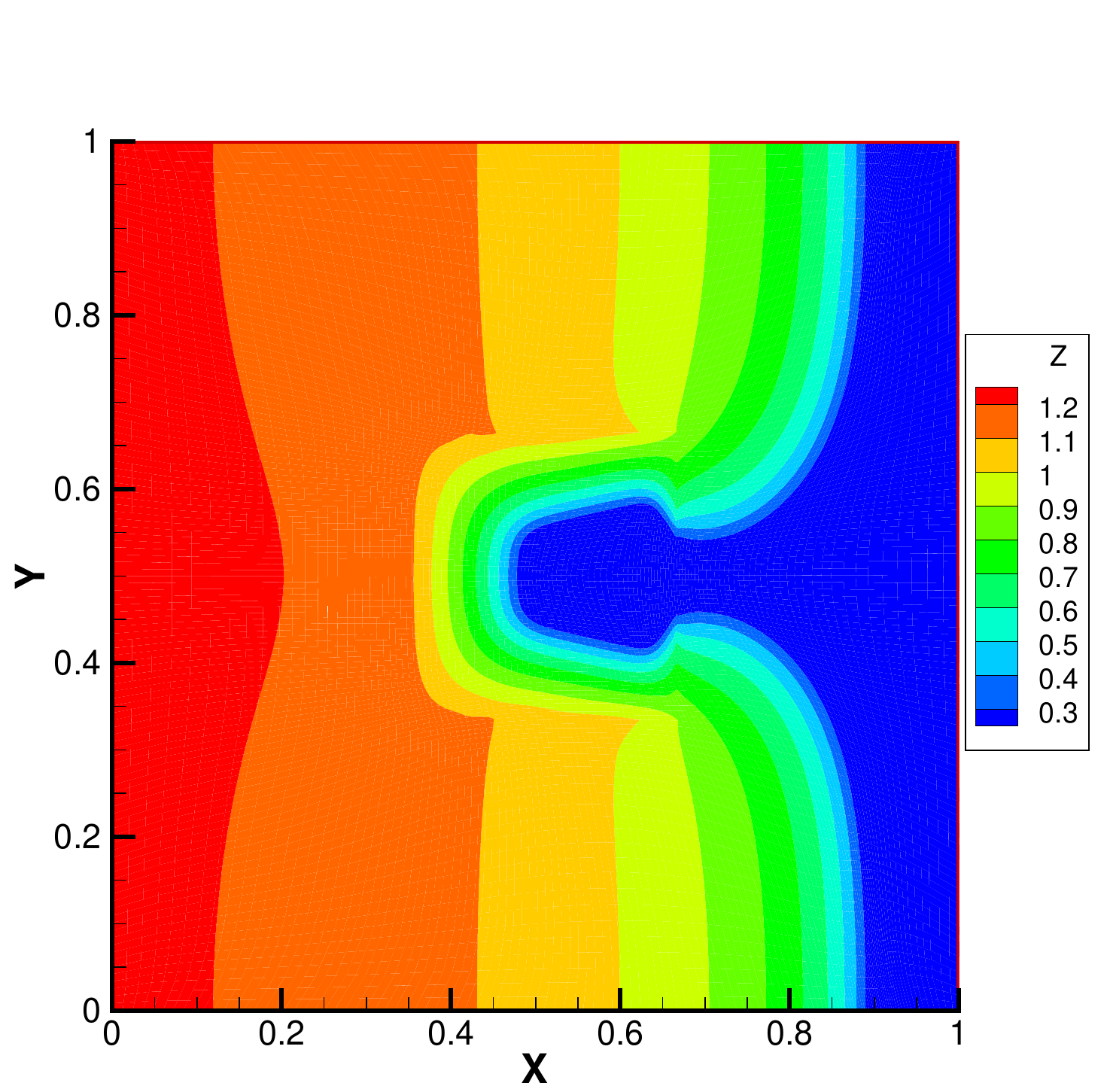}
\end{minipage}
}
\hbox{
\hspace{1in}
\begin{minipage}[t]{2.0in}
\centerline{\scriptsize (g):  with MM1 at $t=3.0$}
\includegraphics[width=2.0in]{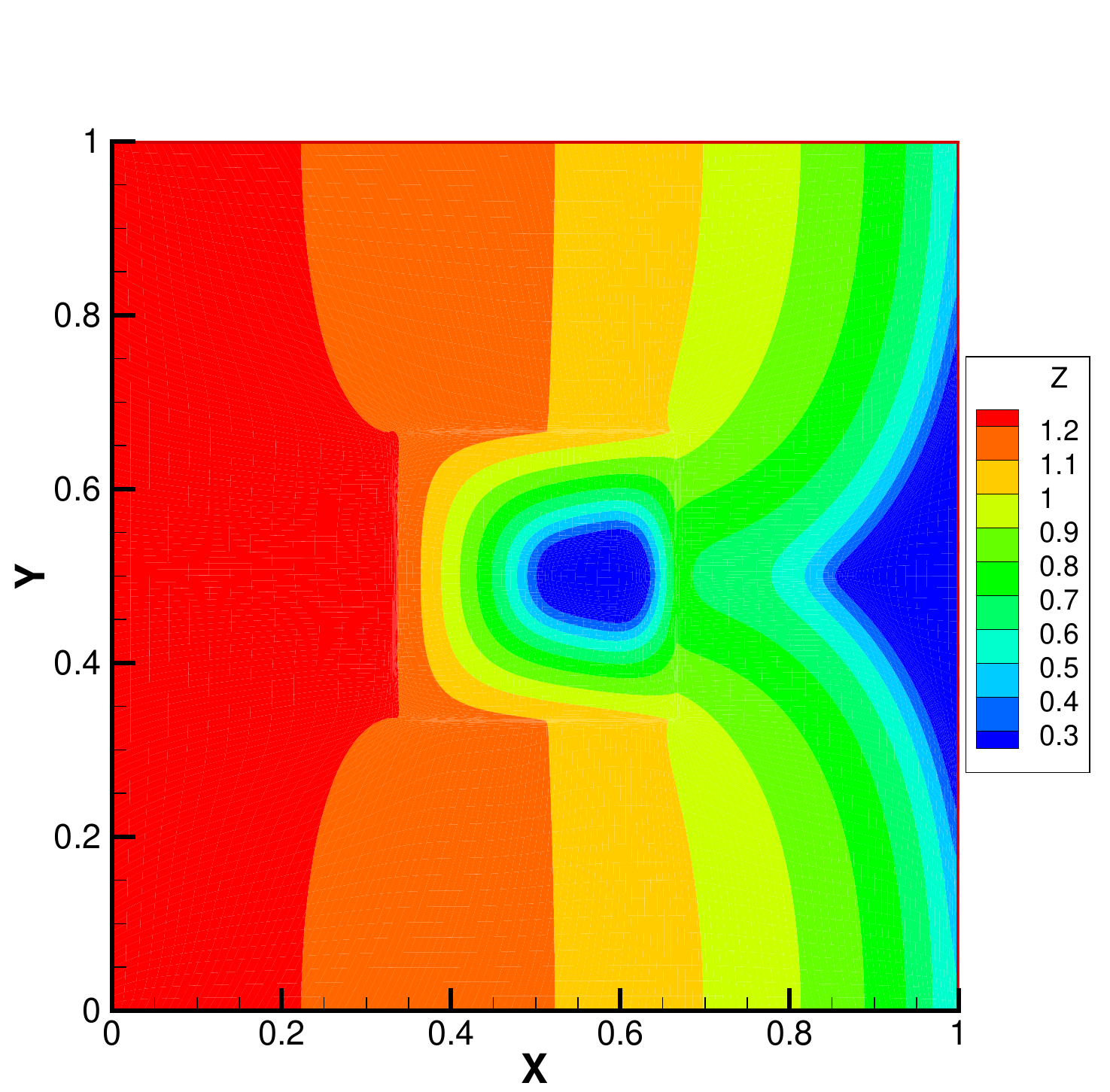}
\end{minipage}
\hspace{0.5in}
\begin{minipage}[t]{2.0in}
\centerline{\scriptsize (h):  with MM2 at $t=3.0$}
\includegraphics[width=2.0in]{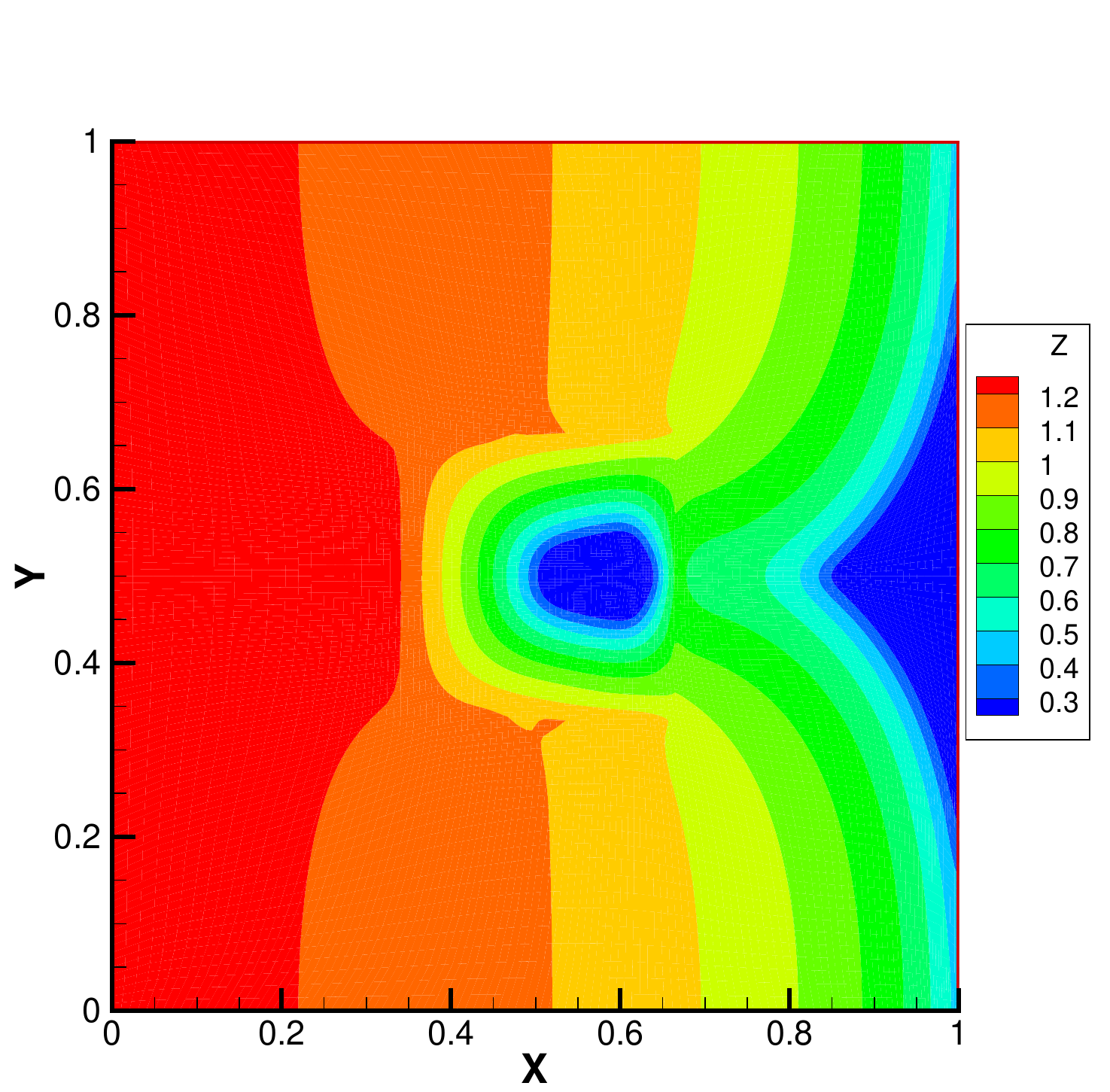}
\end{minipage}
}
\end{center}
\caption{Example~\ref{Example4.2}. The contours of the temperature obtained
with an MM1 of size  $121\times 121$ are compared with those
obtained with an MM2 of size  $41\times 41$ (with the physical PDE
being solved on a uniformly refined mesh of size $121\times 121$).}
\label{T12}
\end{figure}

\begin{figure}
\begin{center}
\hbox{
\hspace{1in}
\begin{minipage}[t]{2.0in}
\centerline{\scriptsize (a):  with MM1 at $t=1.0$}
\includegraphics[width=2.0in]{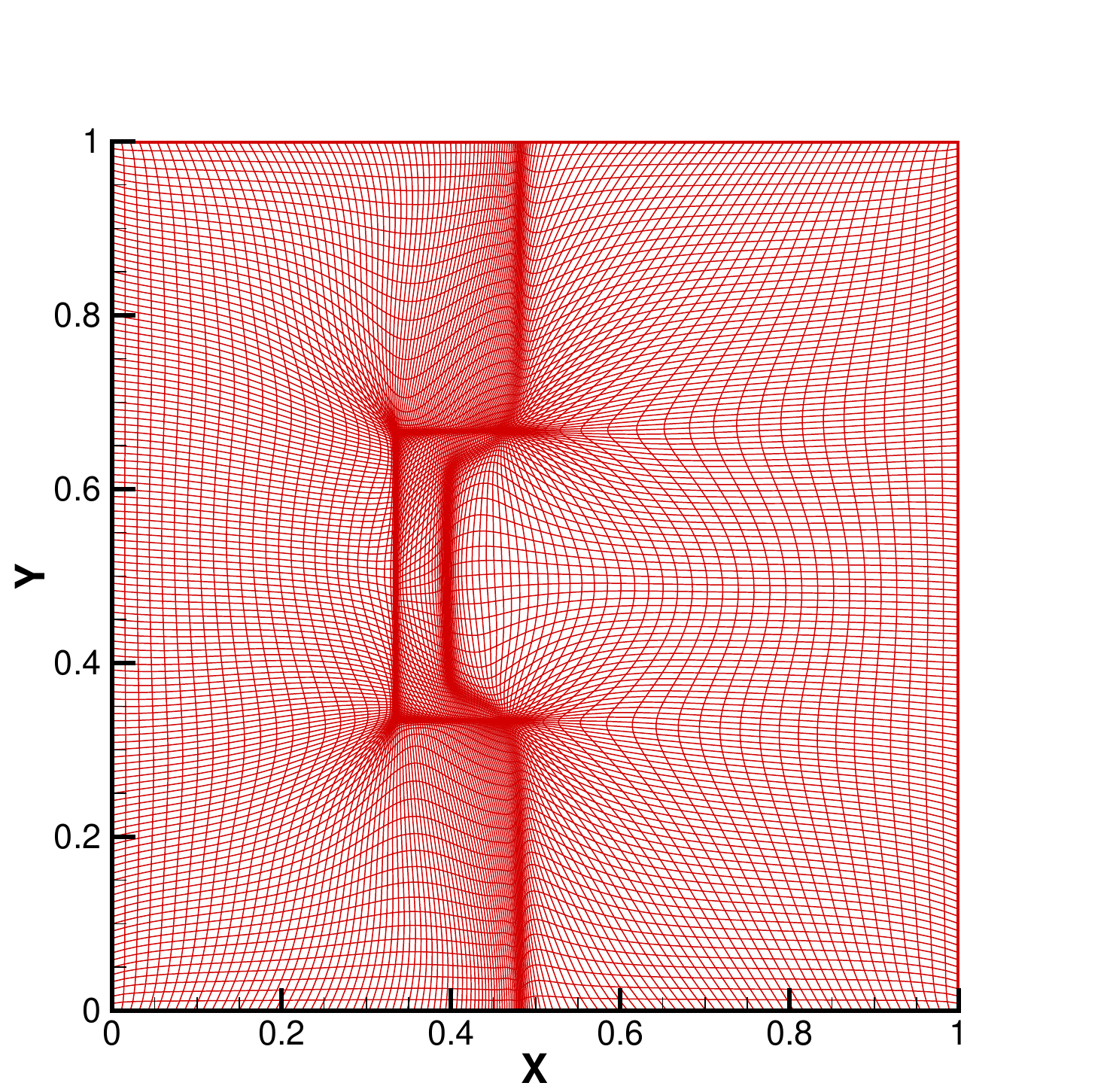}
\end{minipage}
\hspace{0.5in}
\begin{minipage}[t]{2.0in}
\centerline{\scriptsize (b):  with MM2 at $t=1.0$}
\includegraphics[width=2.0in]{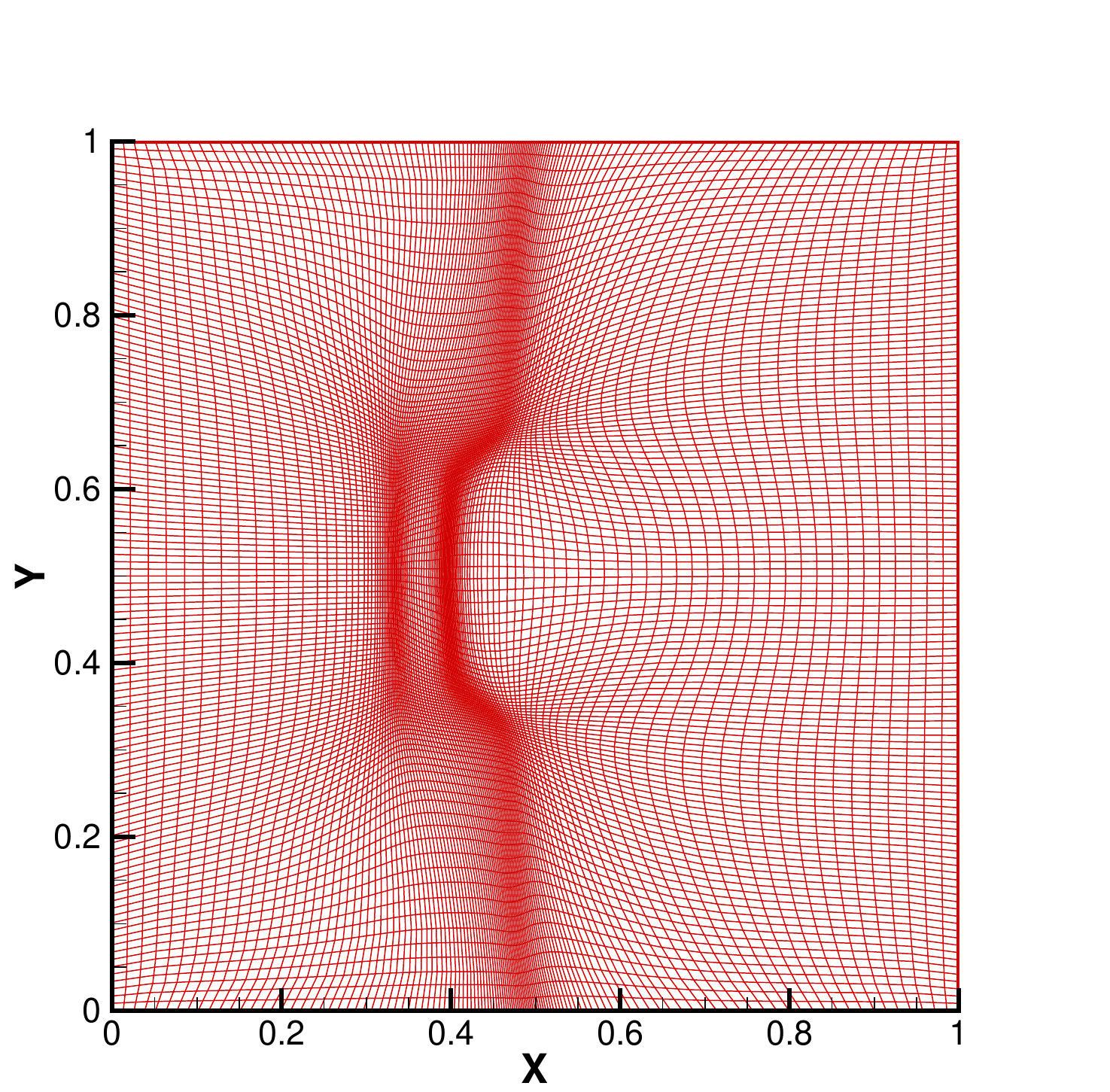}
\end{minipage}
}
\hbox{
\hspace{1in}
\begin{minipage}[t]{2.0in}
\centerline{\scriptsize (c):  with MM1 at $t=2.0$}
\includegraphics[width=2.0in]{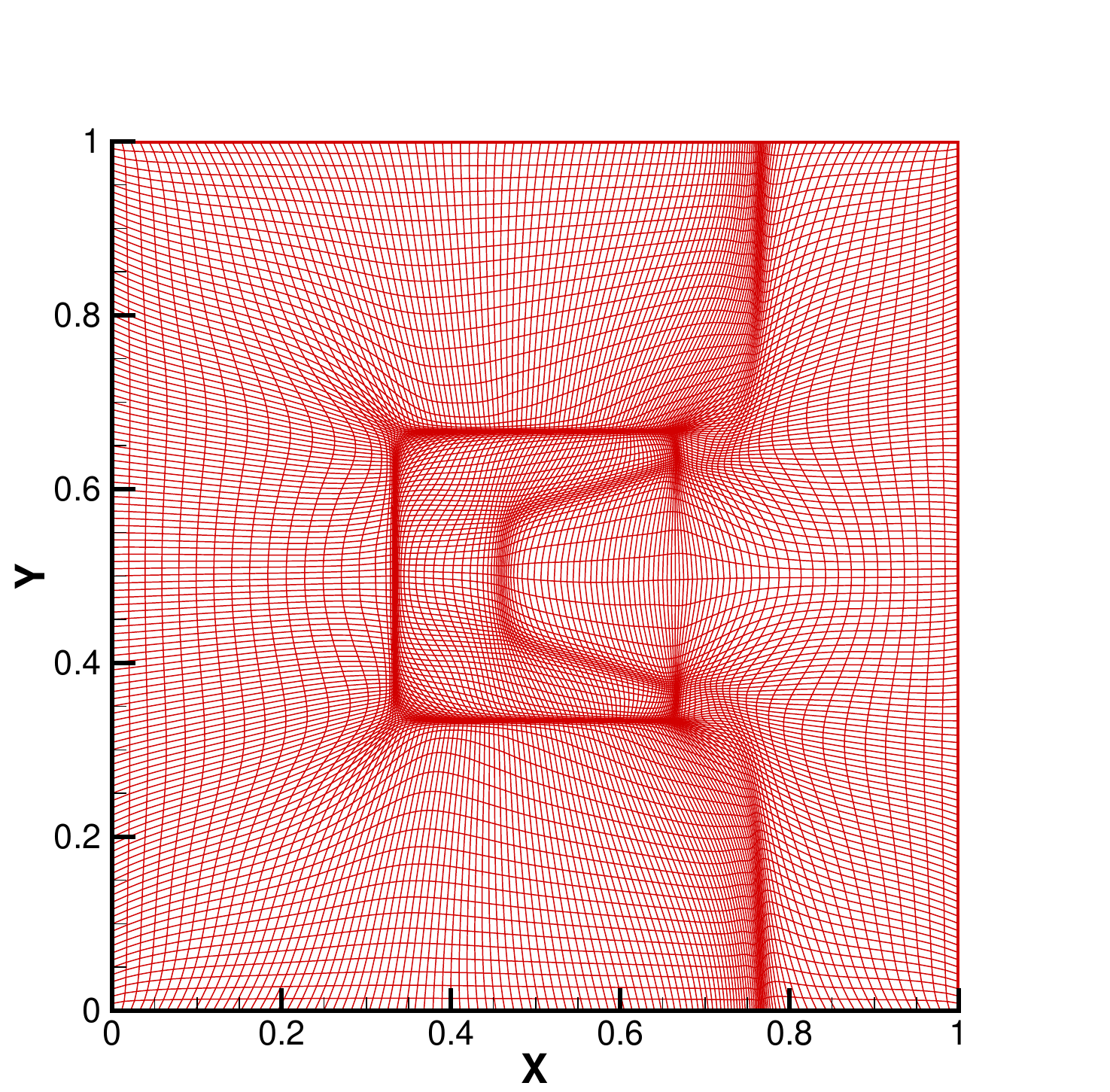}
\end{minipage}
\hspace{0.5in}
\begin{minipage}[t]{2.0in}
\centerline{\scriptsize (d):  with MM2 at $t=2.0$}
\includegraphics[width=2.0in]{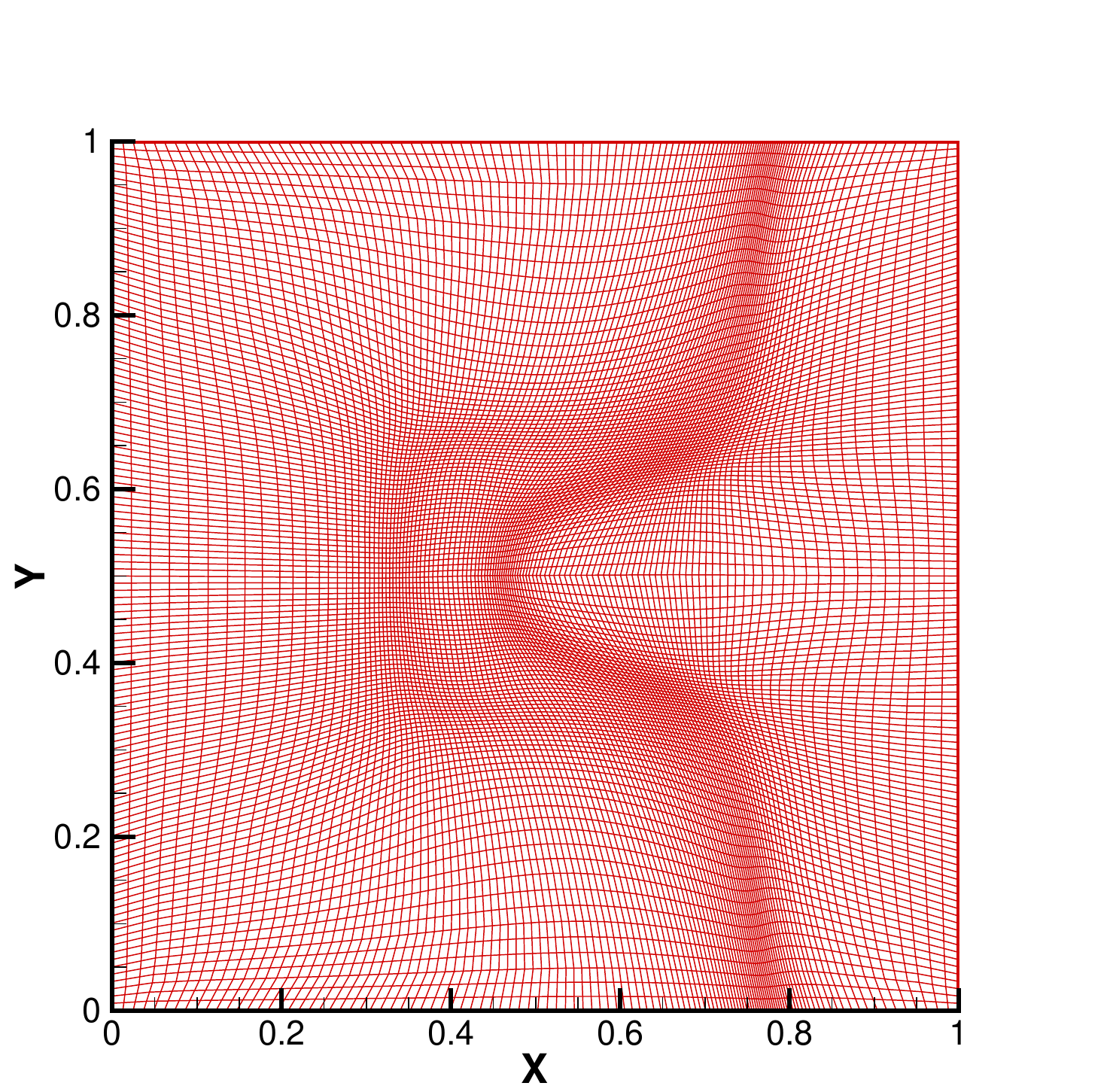}
\end{minipage}
}
\hbox{
\hspace{1in}
\begin{minipage}[t]{2.0in}
\centerline{\scriptsize (e):  with MM1 at $t=2.5$}
\includegraphics[width=2.0in]{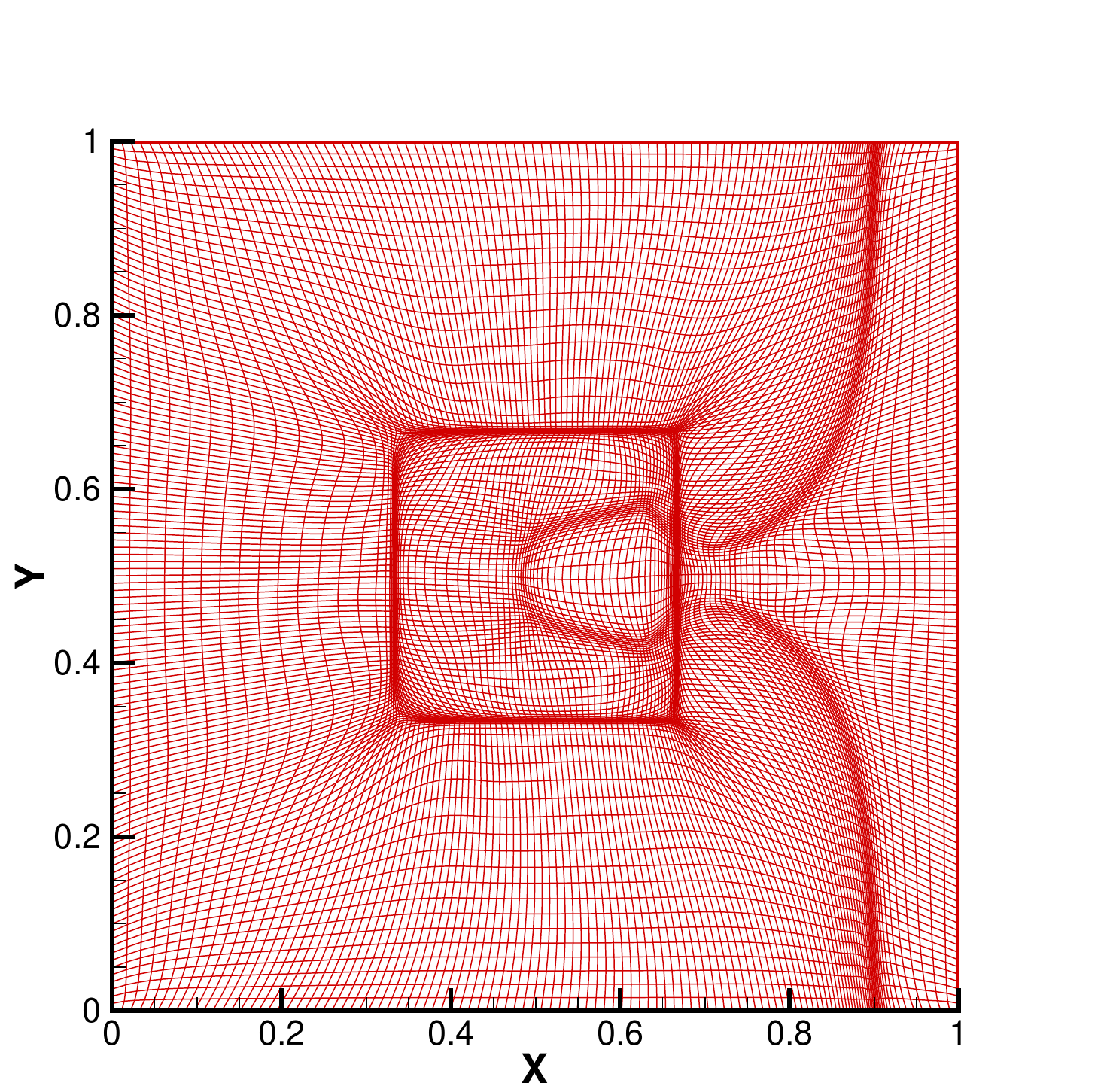}
\end{minipage}
\hspace{0.5in}
\begin{minipage}[t]{2.0in}
\centerline{\scriptsize (f):  with MM2 at $t=2.5$}
\includegraphics[width=2.0in]{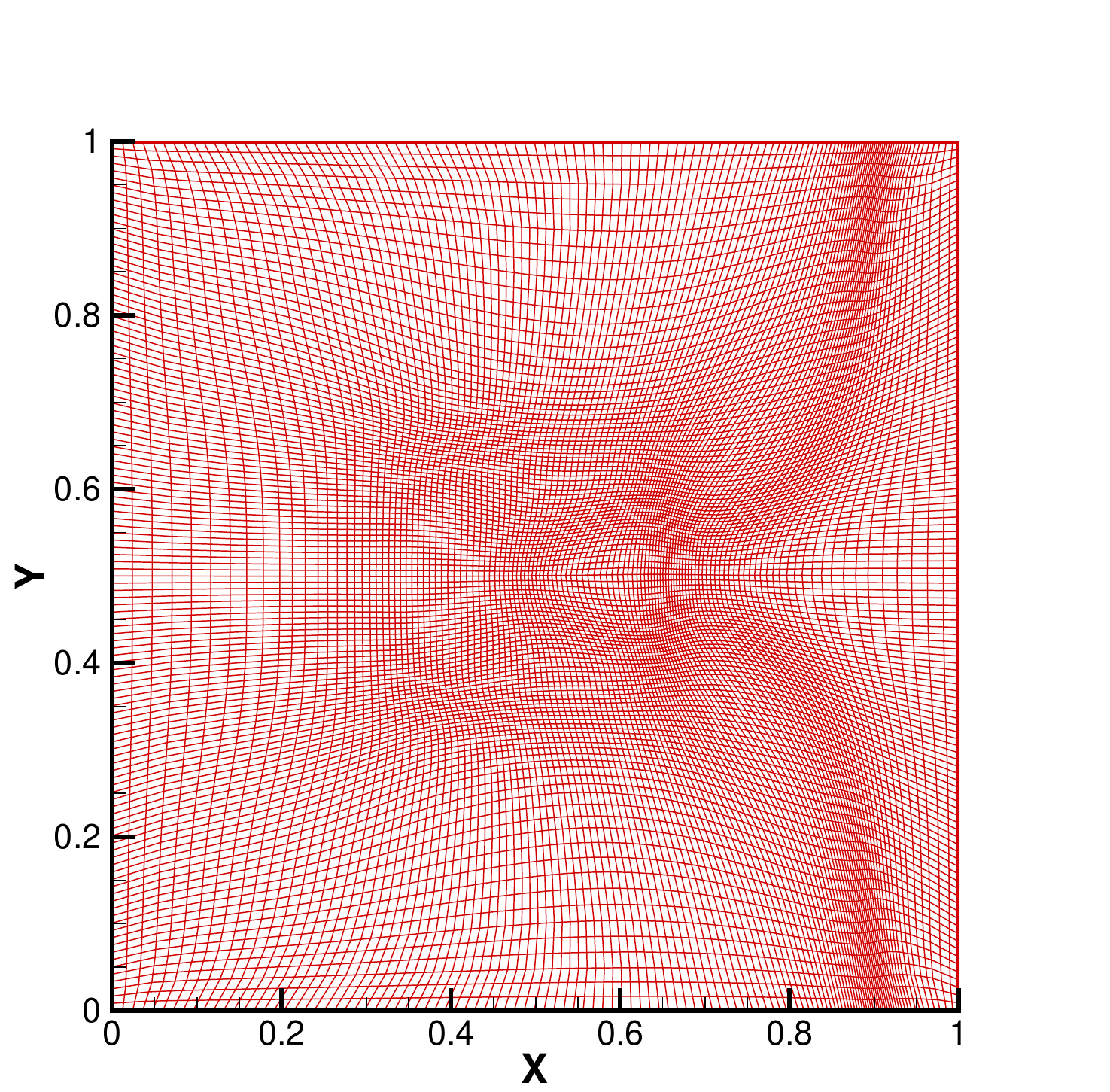}
\end{minipage}
}
\hbox{
\hspace{1in}
\begin{minipage}[t]{2.0in}
\centerline{\scriptsize (g):  with MM1 at $t=3.0$}
\includegraphics[width=2.0in]{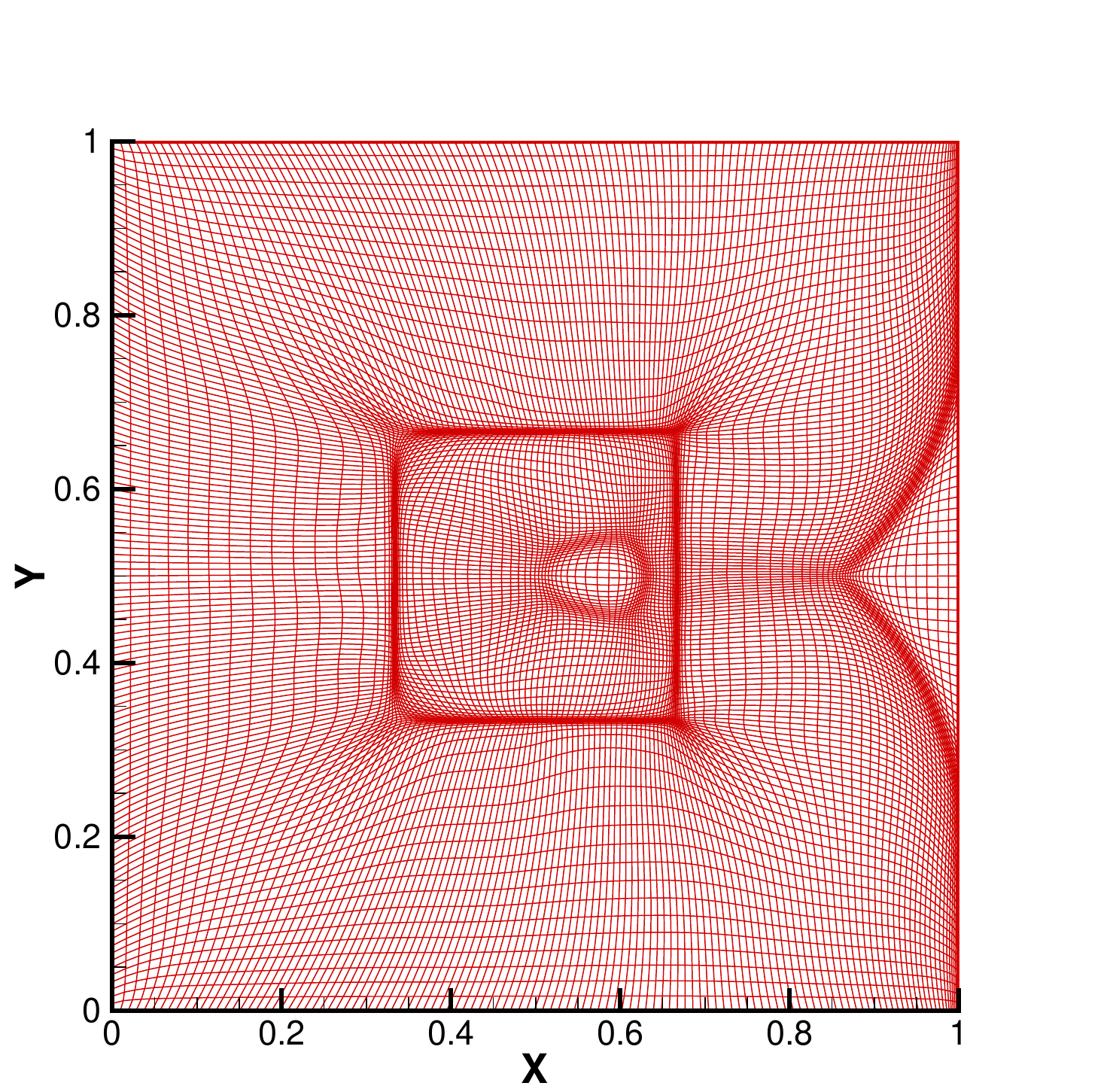}
\end{minipage}
\hspace{0.5in}
\begin{minipage}[t]{2.0in}
\centerline{\scriptsize (h):  with MM2 at $t=3.0$}
\includegraphics[width=2.0in]{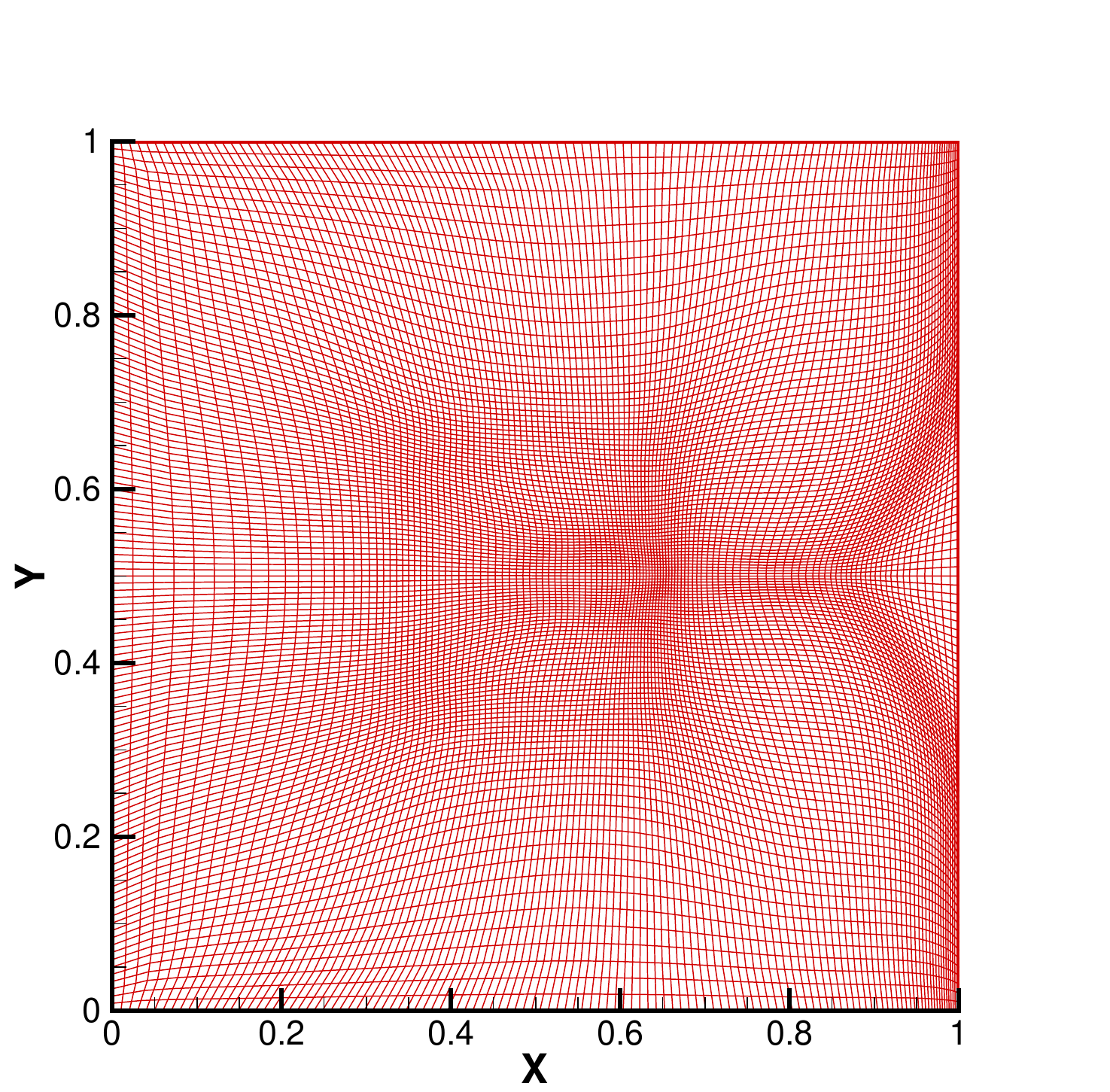}
\end{minipage}
}
\end{center}
\caption{Example~\ref{Example4.2}. Moving meshes of size $121\times 121$ obtained with
MM1 and MM2 moving mesh strategies. MM2 is obtained by uniformly interpolating
a $41\times 41$ moving mesh. }
\label{T13}
\end{figure}

\begin{figure}[htb]
\begin{center}
\includegraphics[width=2.5in]{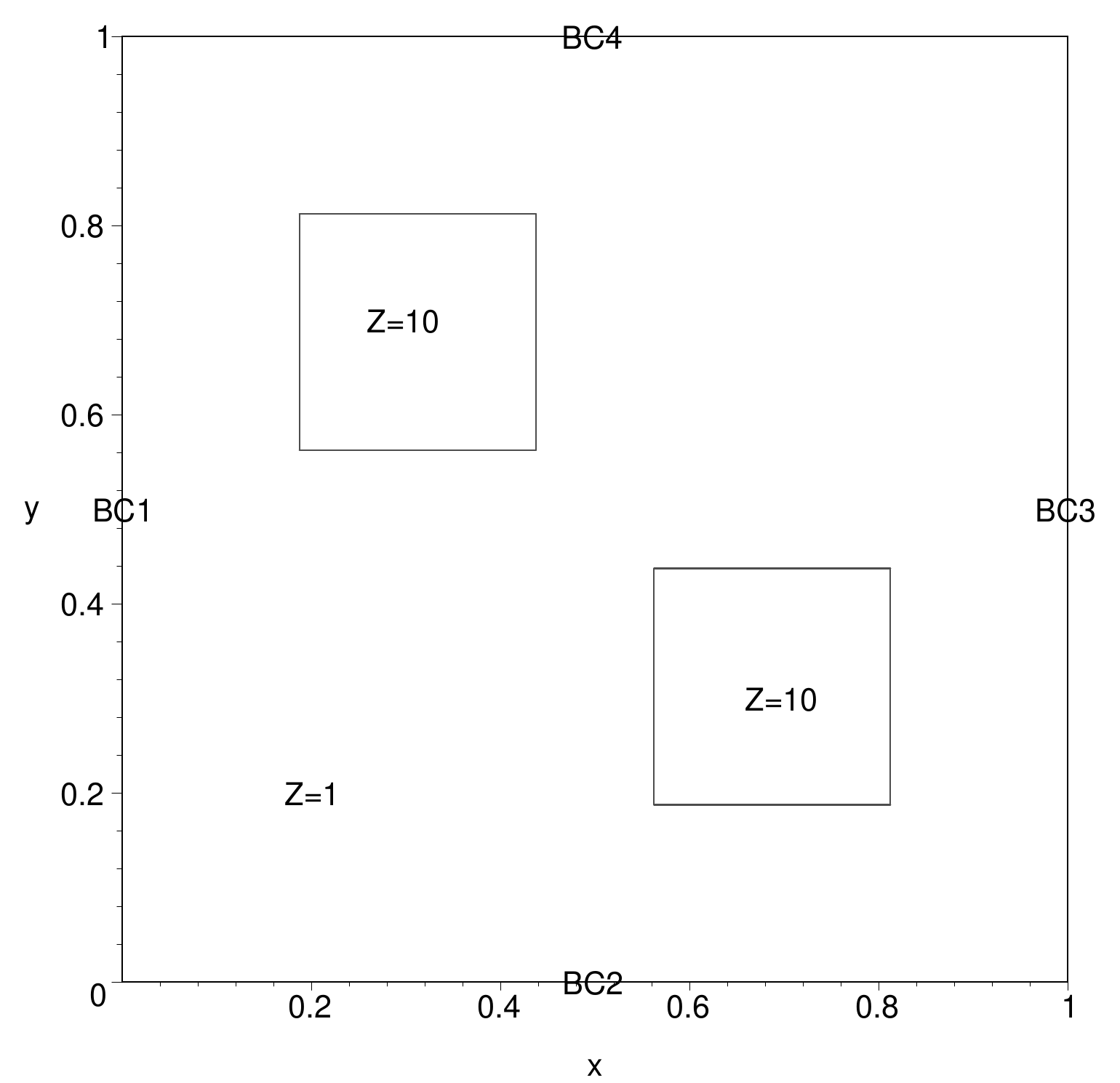}
\end{center}
\caption{The material configuration for Example~\ref{Example4.3}.
The walls BC1, BC2, BC3, BC4 are insulated with respect to radiation diffusion and material
conduction.}
\label{BCD3}
\end{figure}

\begin{exam}
\label{Example4.2}
The setting of this example is the same as the previous example except that the distribution of
the atomic mass number is given by
\begin{equation}
\label{Z3}
Z(x,y)=\begin{cases}
10,  & \text{for }(x,y) \in (\frac13,\frac23)\times (\frac13,\frac23)\\
1,    & \text{otherwise}.
\end{cases}
\end{equation}
Note that the jump in the values of $z$ is more significant than that in the previous example.

The moving mesh of $81\times 81$ and the solution are shown
in Figs.~\ref{T9} and \ref{T10}. From the figures, we can see that
the shape of the central obstacle and the profile of the Marshak wave have been captured
and reflected accurately by the mesh concentration.
It is also worth mentioning that the solutions obtained here are comparable to those
obtained by Kang \cite[Fig. 6 on Page 15]{Kang2002} but with more mesh points.
Comparison results are shown in Fig. \ref{T11} for a moving mesh of $61\times 61$ versus
a uniform mesh of $121\times 121$ and in Figs.~\ref{T12} and \ref{T13} for a one-level moving mesh
of $121\times121$ versus a two-level moving mesh of $41\times 41$ (with the physical PDE
being solved on a uniformly refined mesh of $121\times121$). The results are all comparable.
Moreover, the CPU time is listed in Table~\ref{mytab2}. It can be seen that the two-level
moving mesh strategy can significantly improve the efficiency of the moving mesh method
without compromising the accuracy.

\end{exam}

\begin{table}[htbp]
 \centering\small
\begin{threeparttable}
  \caption{CPU time comparison among one-level, two-level moving mesh
  and uniform mesh methods for Example \ref{Example4.2}. The CPU time is measured in seconds.
  The last column is the ratio of the used CPU time to that used with a uniform mesh of the same size.}
  \label{mytab2}
 \medskip
 \begin{tabular}{|l|c|c|c|c|}
 \hline
   & Fine Mesh  & Coarse Mesh & Total CPU time & ratio \\
   \hline
  One-level MM & 41$\times$41    &   41$\times$41     &      2951       &  5.85 \\
  \hline
          & 81$\times$81    &    81$\times$81     &      139374        &  58.76      \\
  \hline
          & 121$\times$121  &     121$\times$121     &     $1786732$   &  $333.03$   \\
  \hline

Tow-level MM   & 41$\times$41  &      41$\times$41     &    2951  &   5.85       \\
\hline
          & 81$\times$81  &    41$\times$41     &      9581      &   4.03         \\
\hline
          & 121$\times$121  &   41$\times$41     &     21888      &   4.08           \\
\hline

Fixed mesh     & 41$\times$41  &      n/a          &     498       &   1                 \\
\hline
          & 81$\times$81  &       n/a         &      2372      &    1                   \\
\hline
          & 121$\times$121  &     n/a          &     5365      &    1                \\
\hline
\end{tabular}
\end{threeparttable}

\end{table}

\begin{figure}
\begin{center}
\hbox{
\begin{minipage}[t]{2.3in}
\centerline{\scriptsize (a): t=0.5}
\includegraphics[width=2.3in]{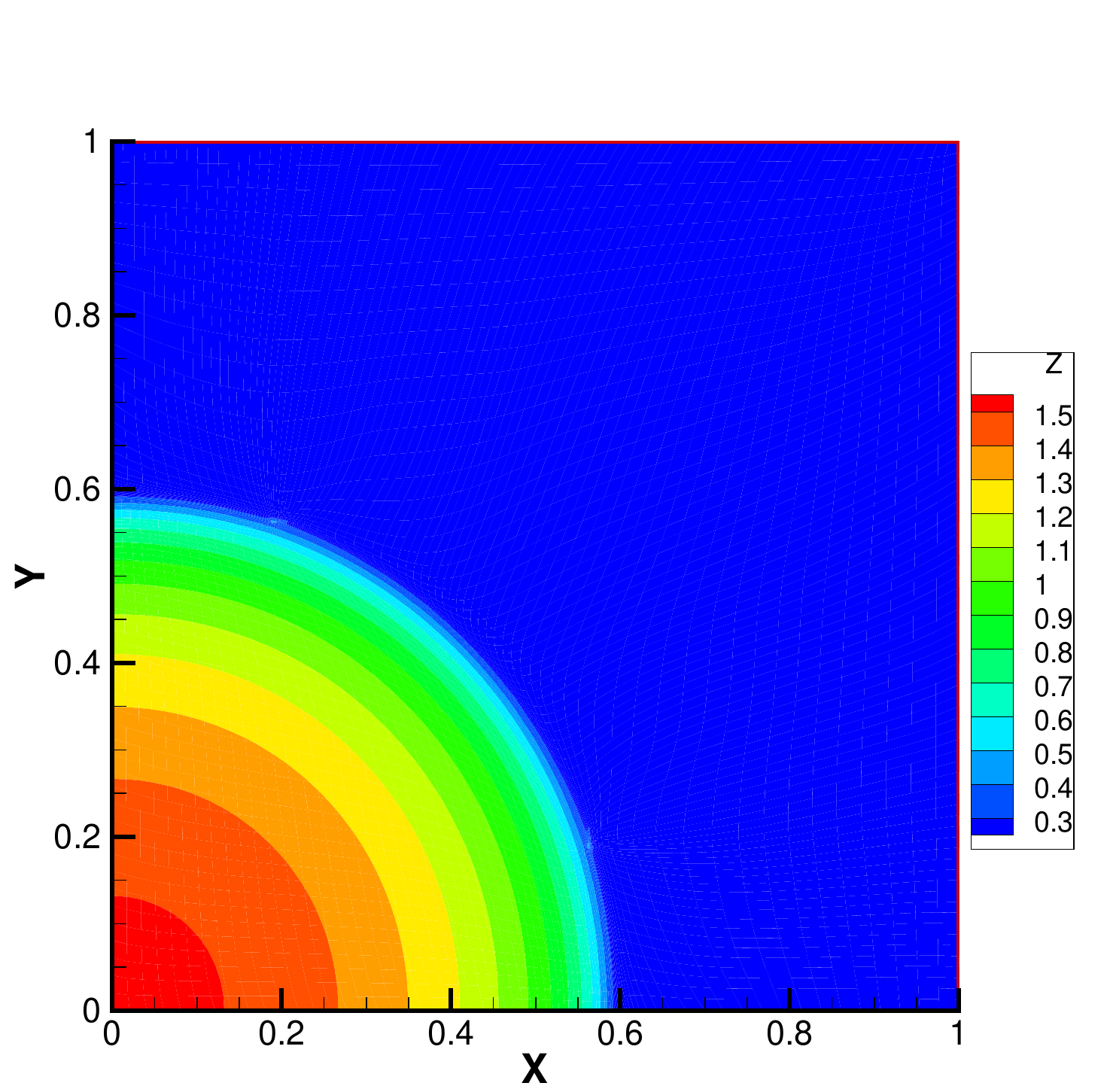}
\end{minipage}
\begin{minipage}[t]{2.3in}
\centerline{\scriptsize (b): t=0.7}
\includegraphics[width=2.3in]{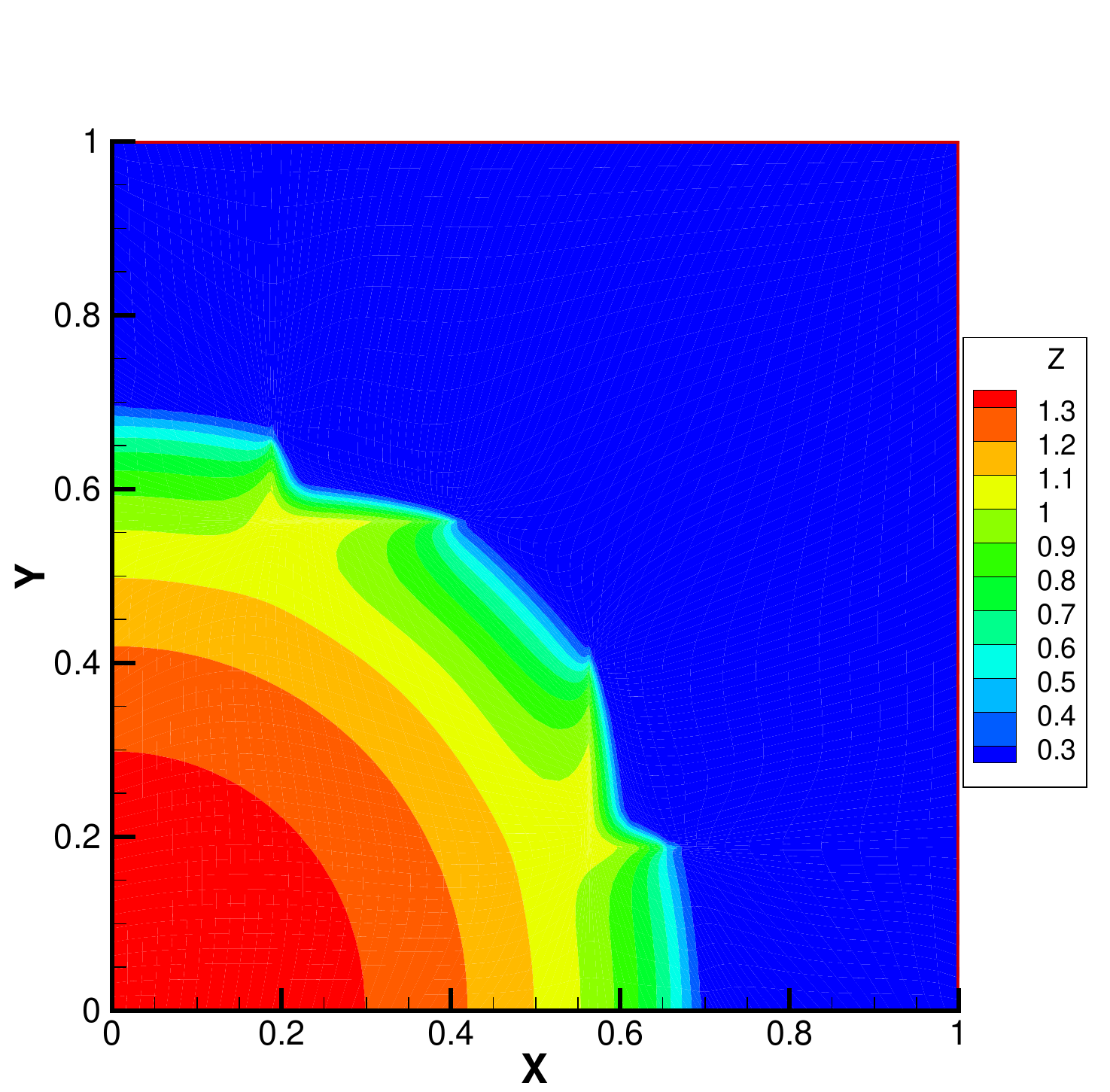}
\end{minipage}
\begin{minipage}[t]{2.3in}
\centerline{\scriptsize (c): t=0.8}
\includegraphics[width=2.3in]{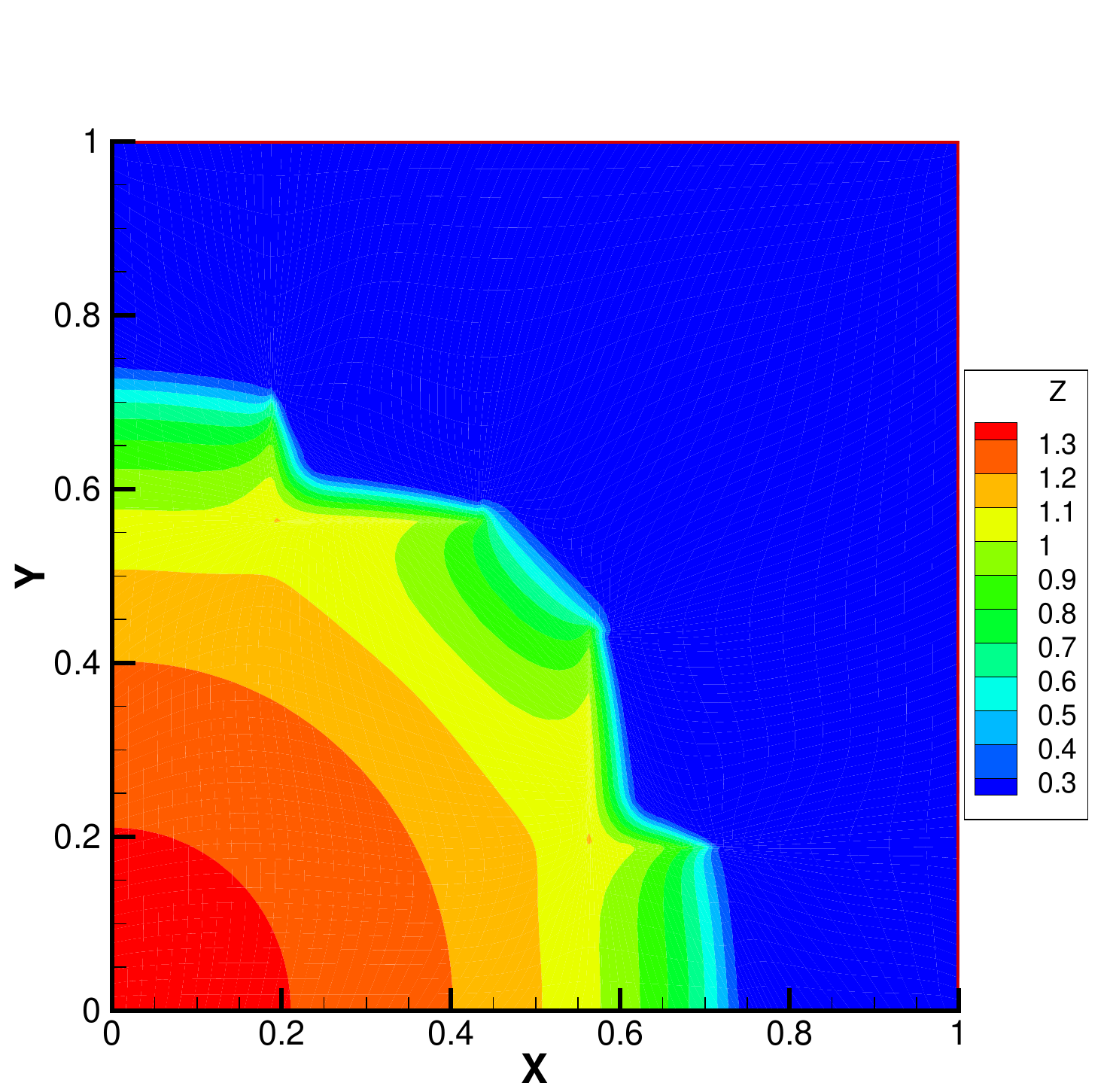}
\end{minipage}
}
\vspace{5mm}
\hbox{
\begin{minipage}[t]{2.3in}
\centerline{\scriptsize (d): t=0.9}
\includegraphics[width=2.3in]{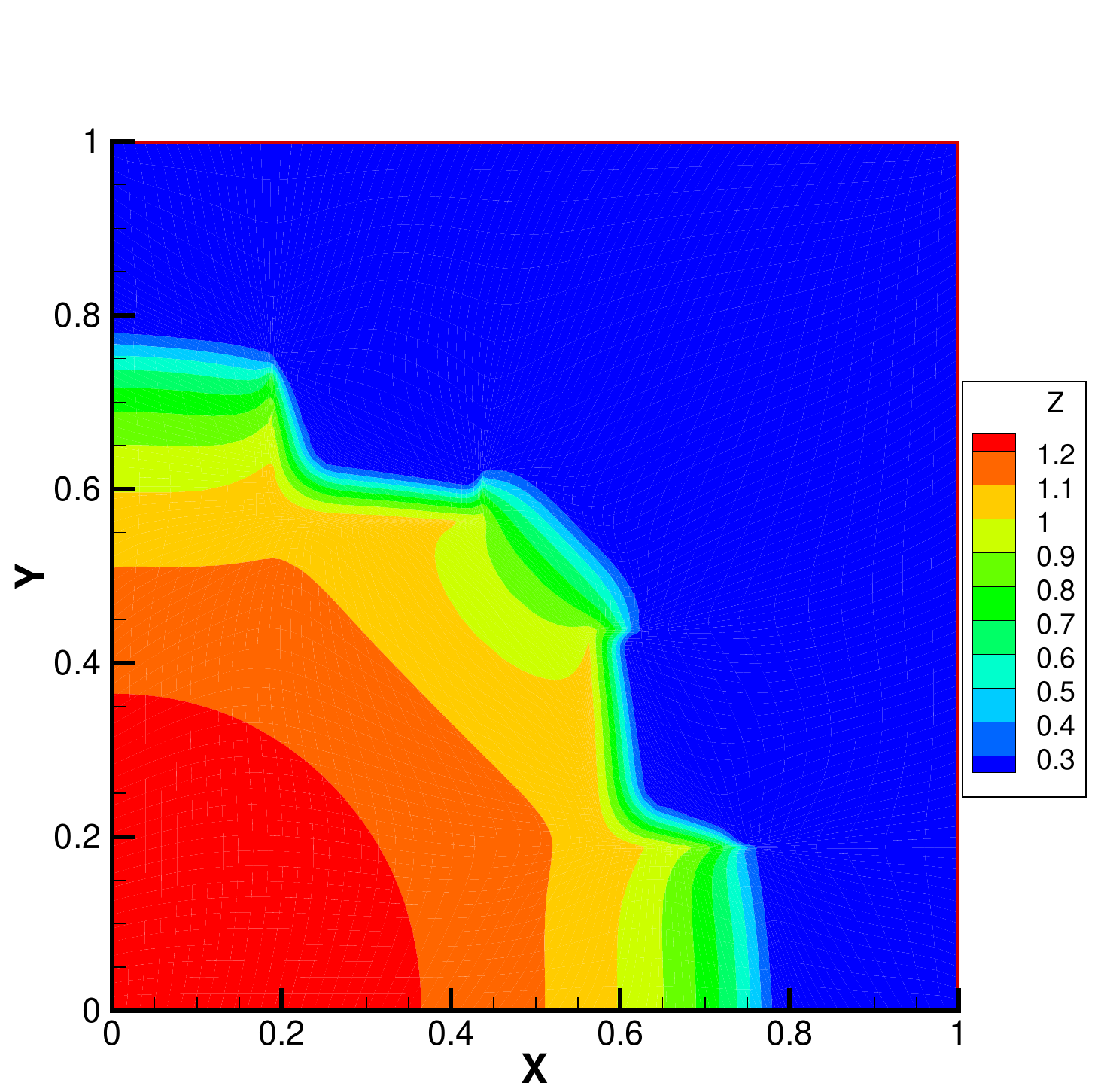}
\end{minipage}
\begin{minipage}[t]{2.3in}
\centerline{\scriptsize (e): t=1.0}
\includegraphics[width=2.3in]{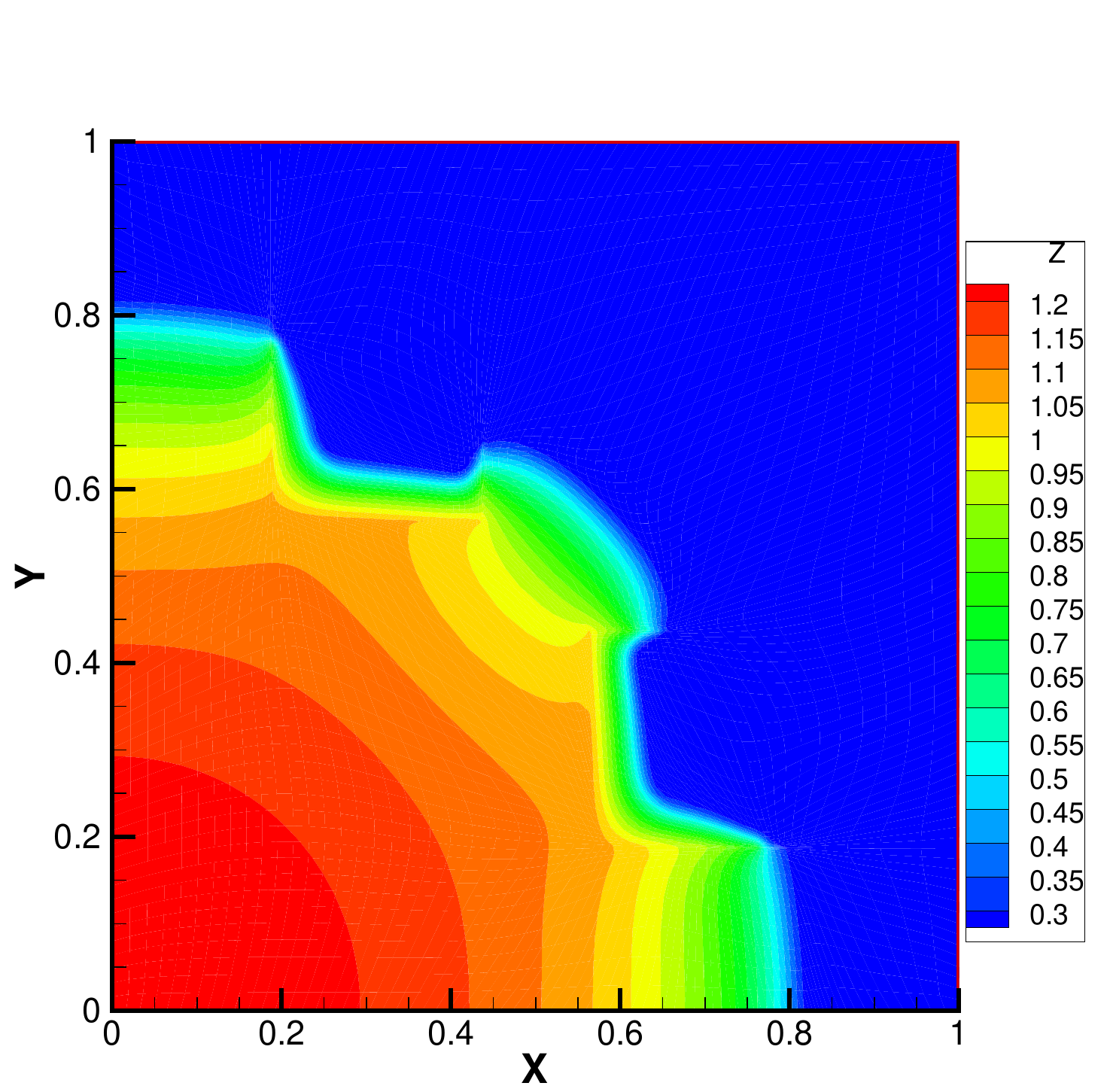}
\end{minipage}
\begin{minipage}[t]{2.3in}
\centerline{\scriptsize (f): t=1.5}
\includegraphics[width=2.3in]{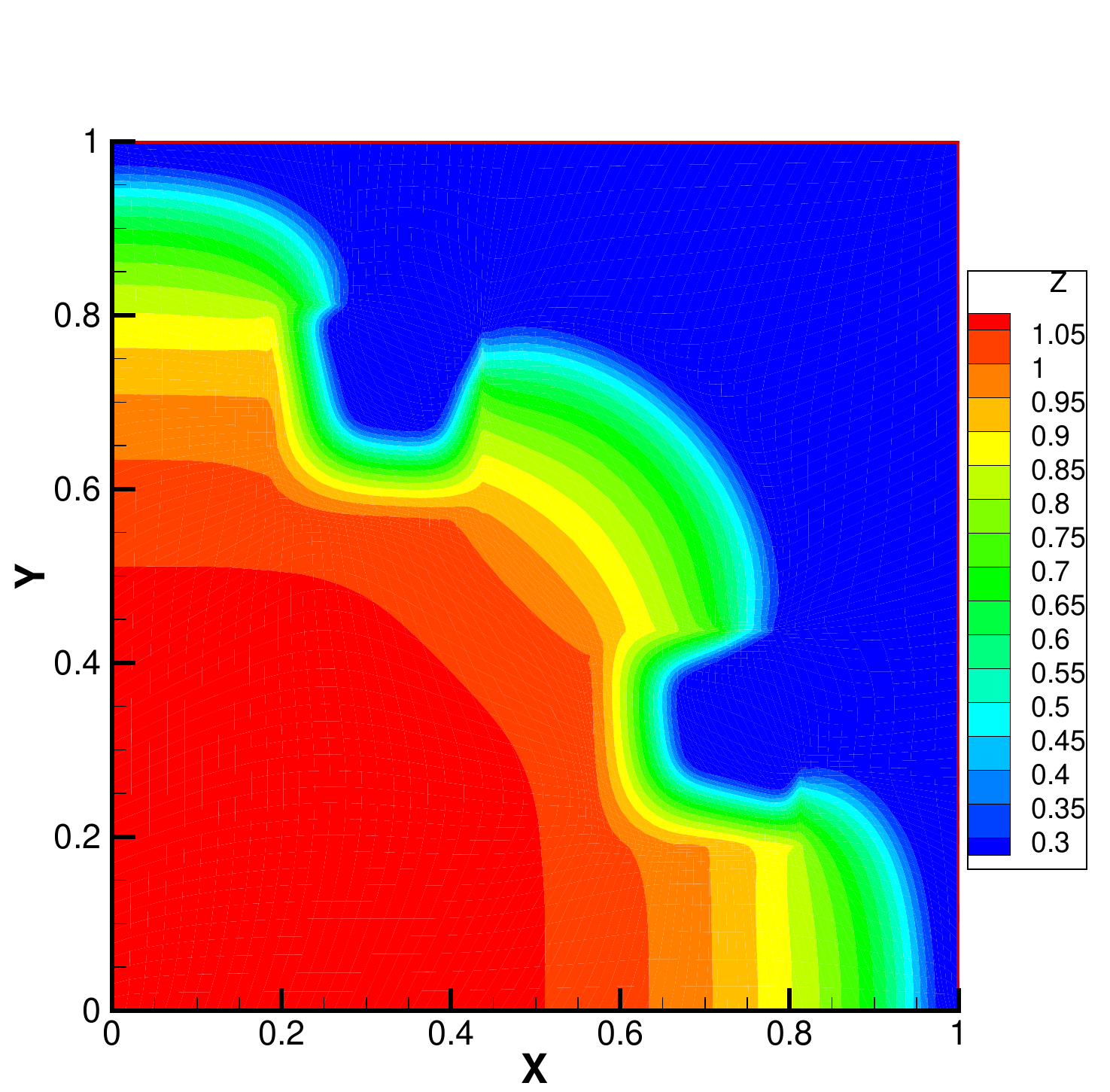}
\end{minipage}
}
\vspace{5mm}
\hbox{
\begin{minipage}[t]{2.3in}
\centerline{\scriptsize (g): t=2.0}
\includegraphics[width=2.3in]{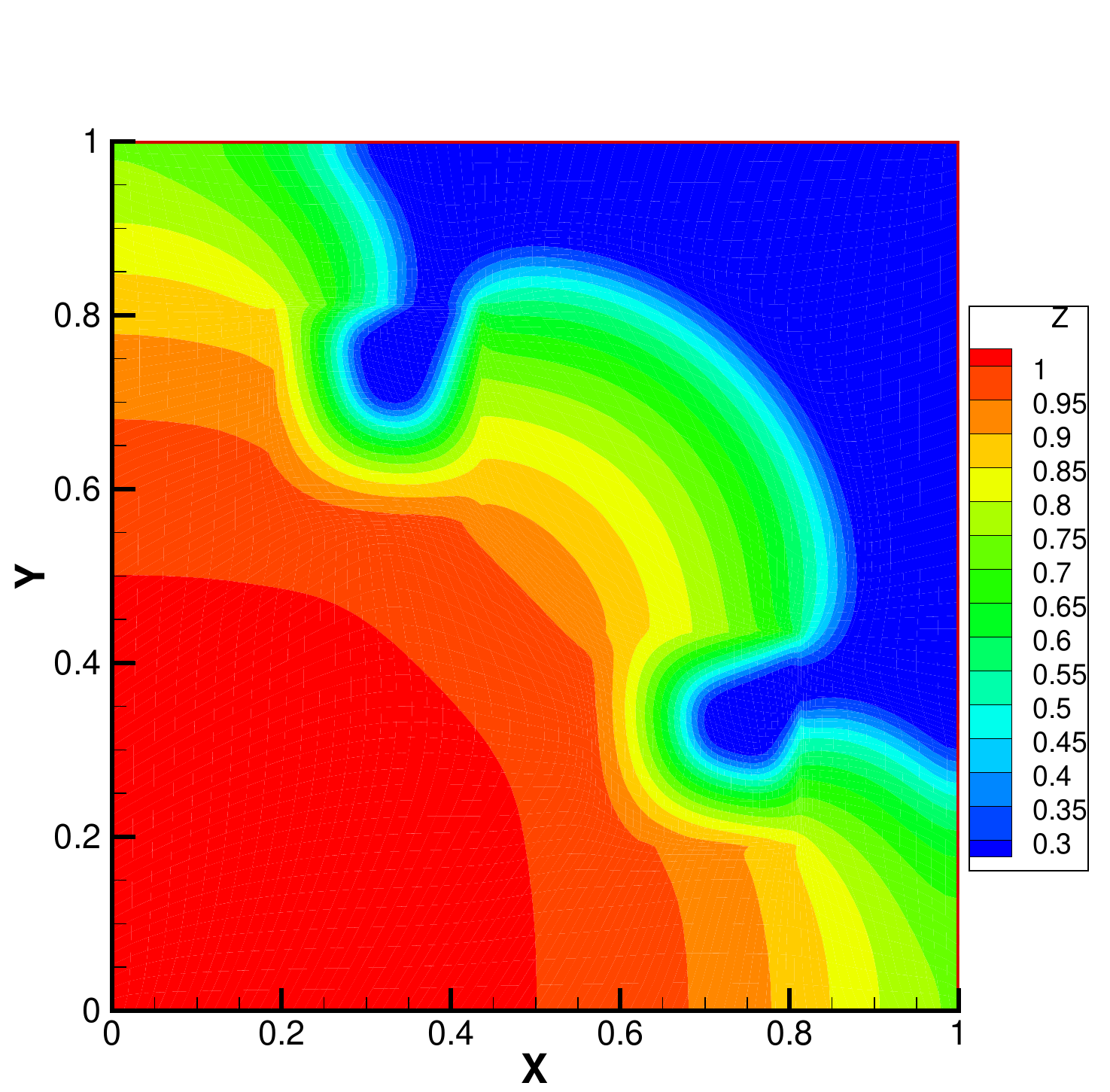}
\end{minipage}
\begin{minipage}[t]{2.3in}
\centerline{\scriptsize (h): t=2.5}
\includegraphics[width=2.3in]{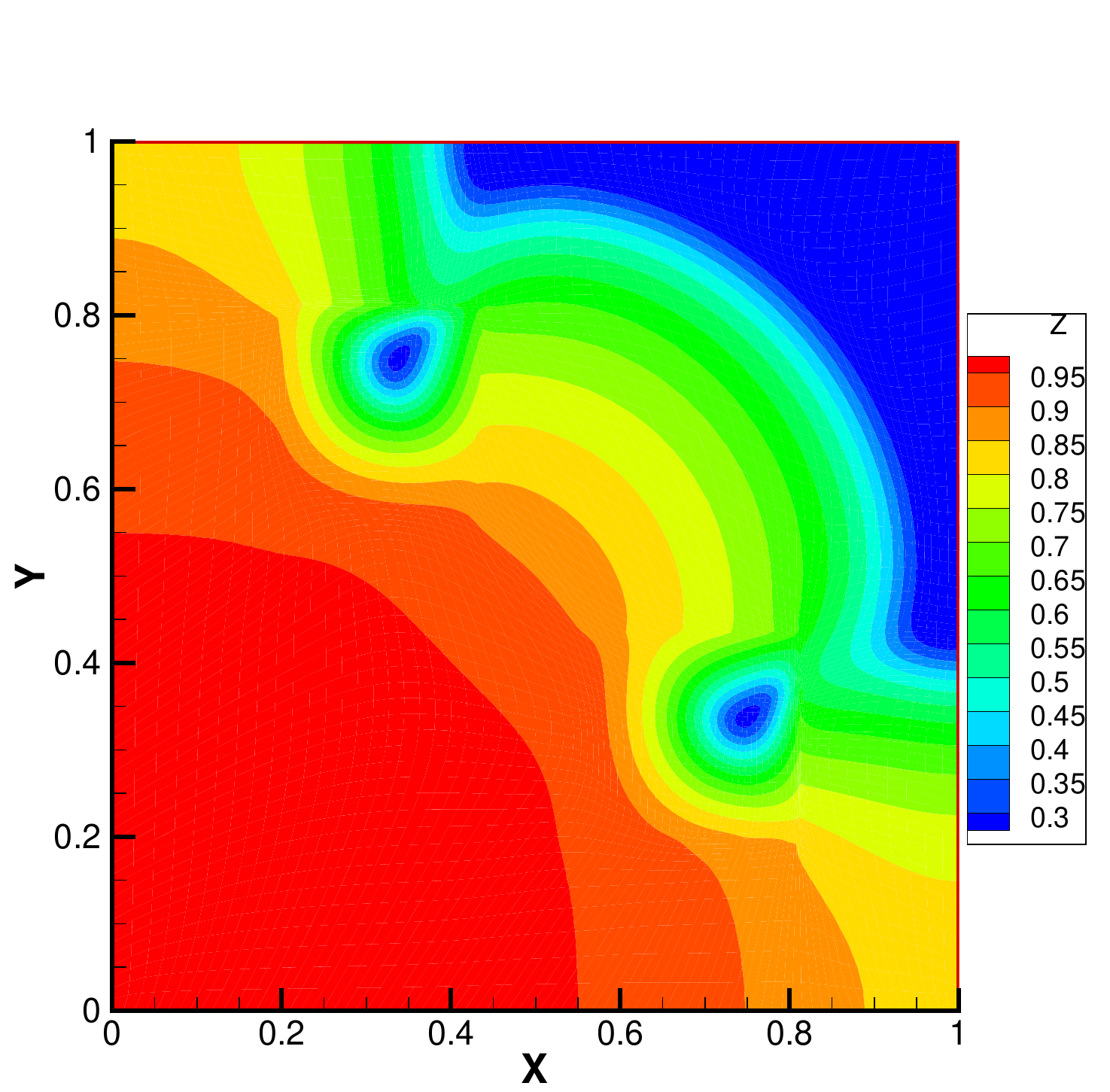}
\end{minipage}
\begin{minipage}[t]{2.3in}
\centerline{\scriptsize (i): t=3.0}
\includegraphics[width=2.3in]{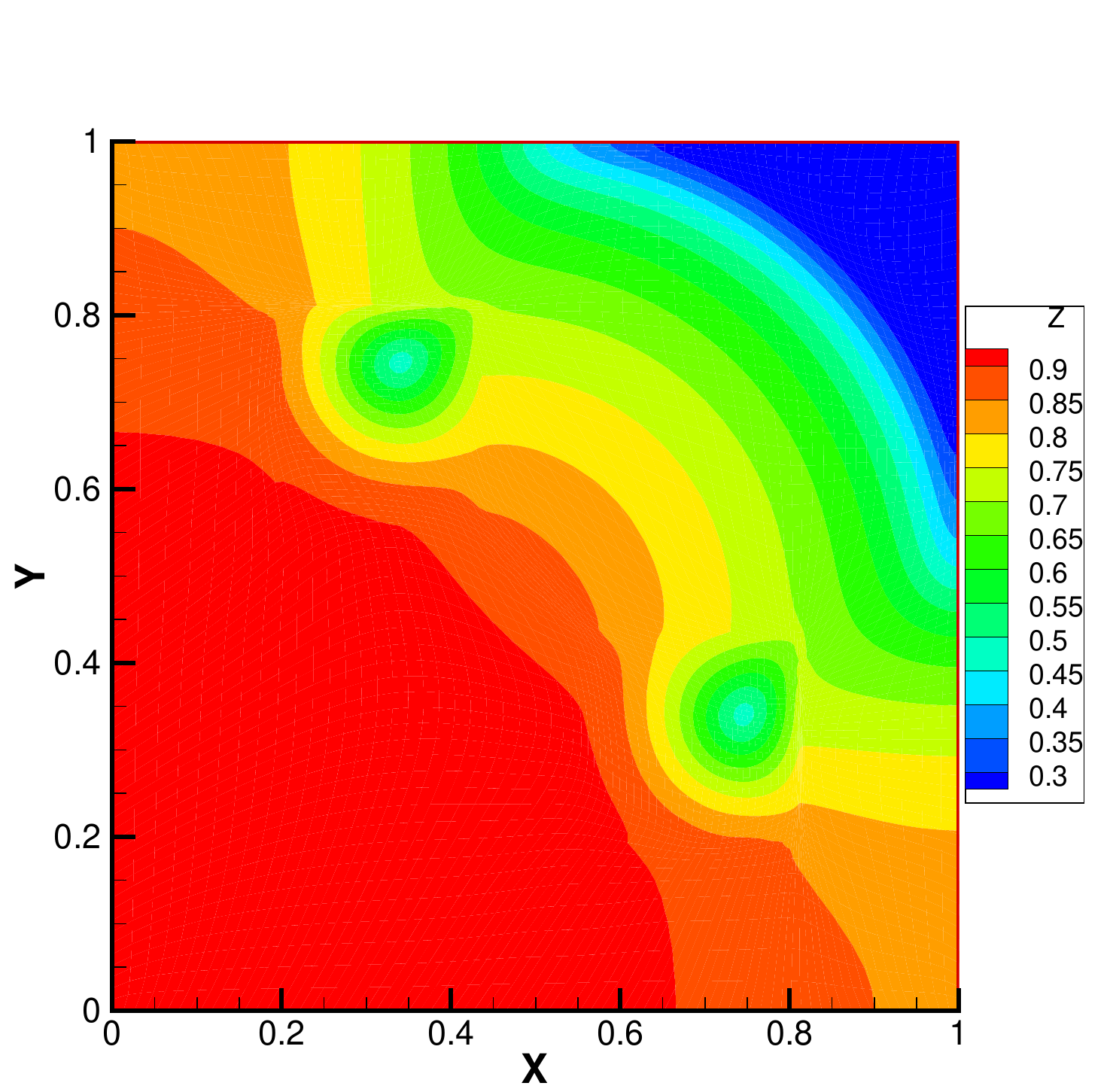}
\end{minipage}
}
\end{center}
\caption{Example~\ref{Example4.3}. The computed solution on a moving mesh of $81\times 81$
is shown at $t=0.5, 0.7, 0.8, 0.9, 1.0, 1.5, 2.0, 2.5, 3.0$.}
\label{T28}
\end{figure}

\begin{figure}
\begin{center}
\hbox{
\begin{minipage}[t]{2.3in}
\centerline{\scriptsize (a): t=0.5}
\includegraphics[width=2.3in]{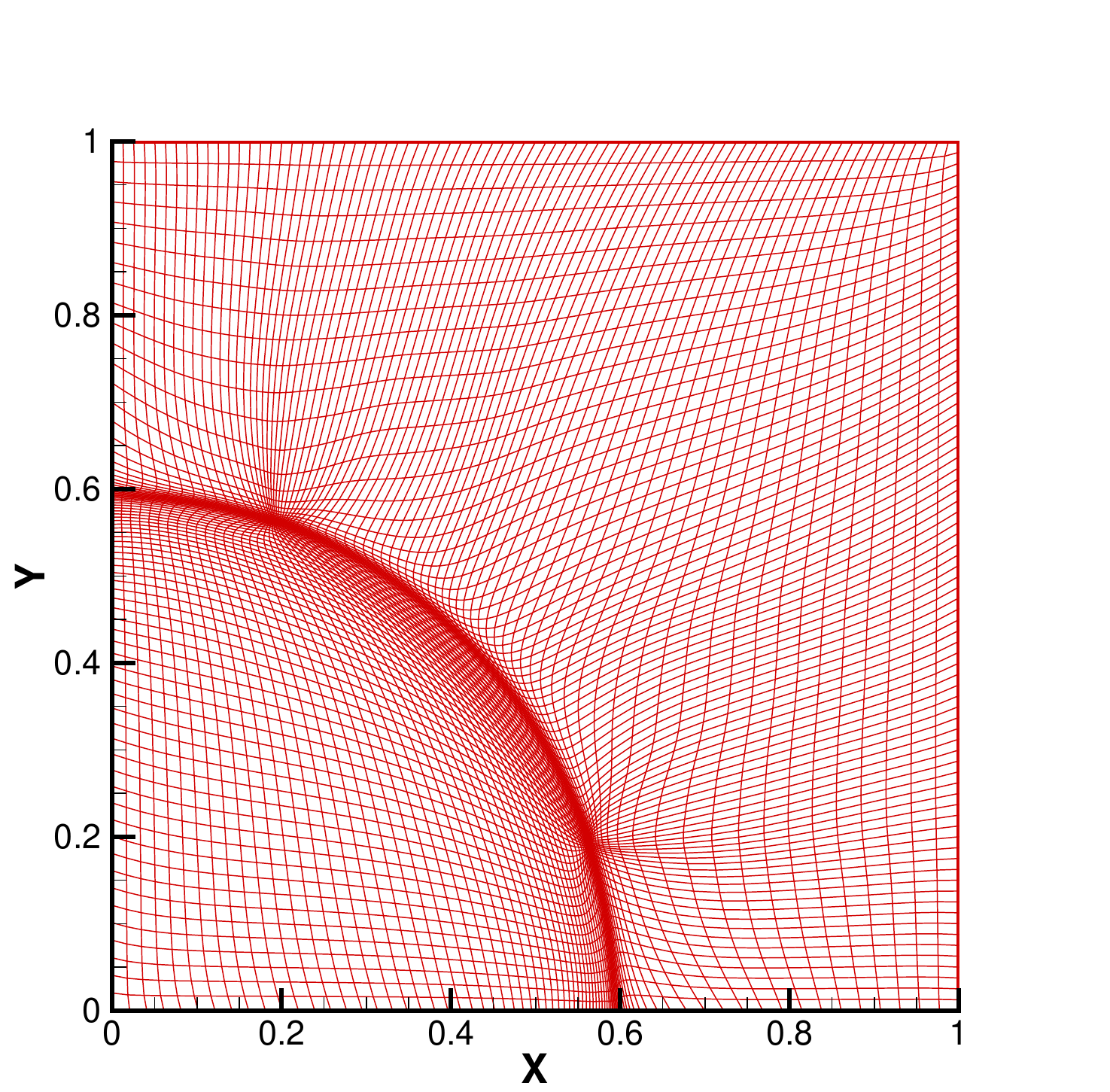}
\end{minipage}
\begin{minipage}[t]{2.3in}
\centerline{\scriptsize (b): t=0.7}
\includegraphics[width=2.3in]{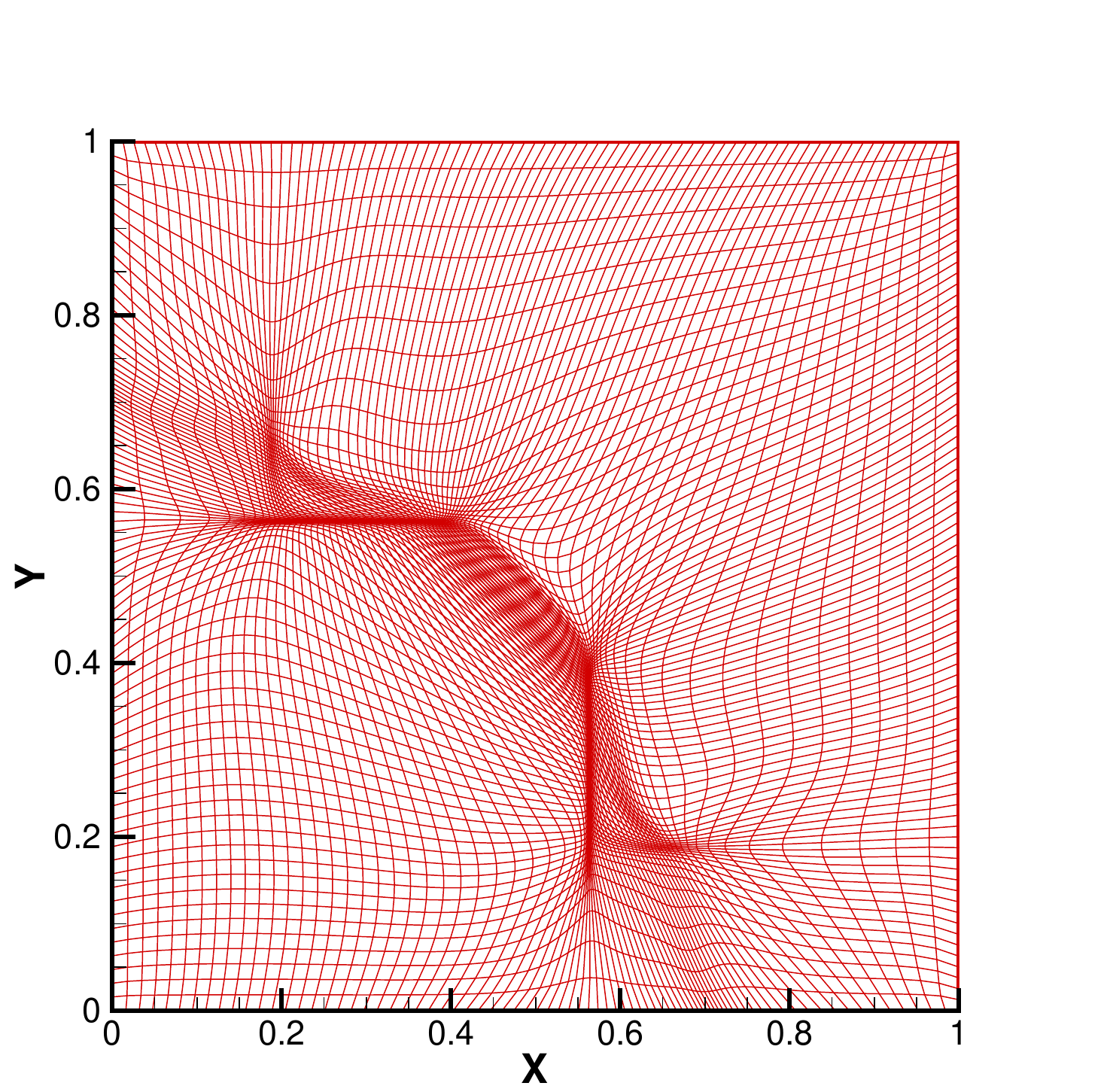}
\end{minipage}
\begin{minipage}[t]{2.3in}
\centerline{\scriptsize (c): t=0.8}
\includegraphics[width=2.3in]{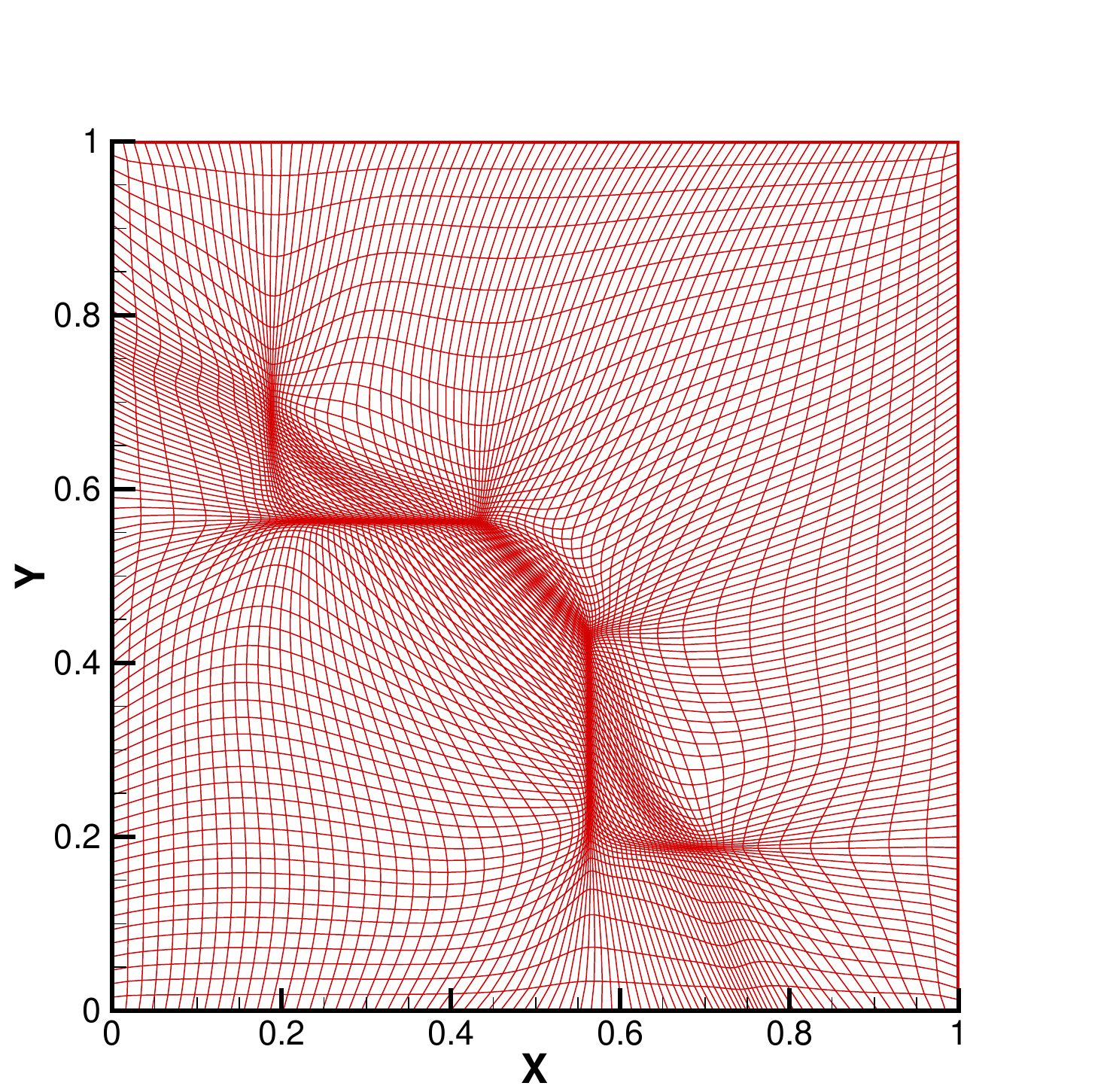}
\end{minipage}
}
\vspace{5mm}
\hbox{
\begin{minipage}[t]{2.3in}
\centerline{\scriptsize (d): t=0.9}
\includegraphics[width=2.3in]{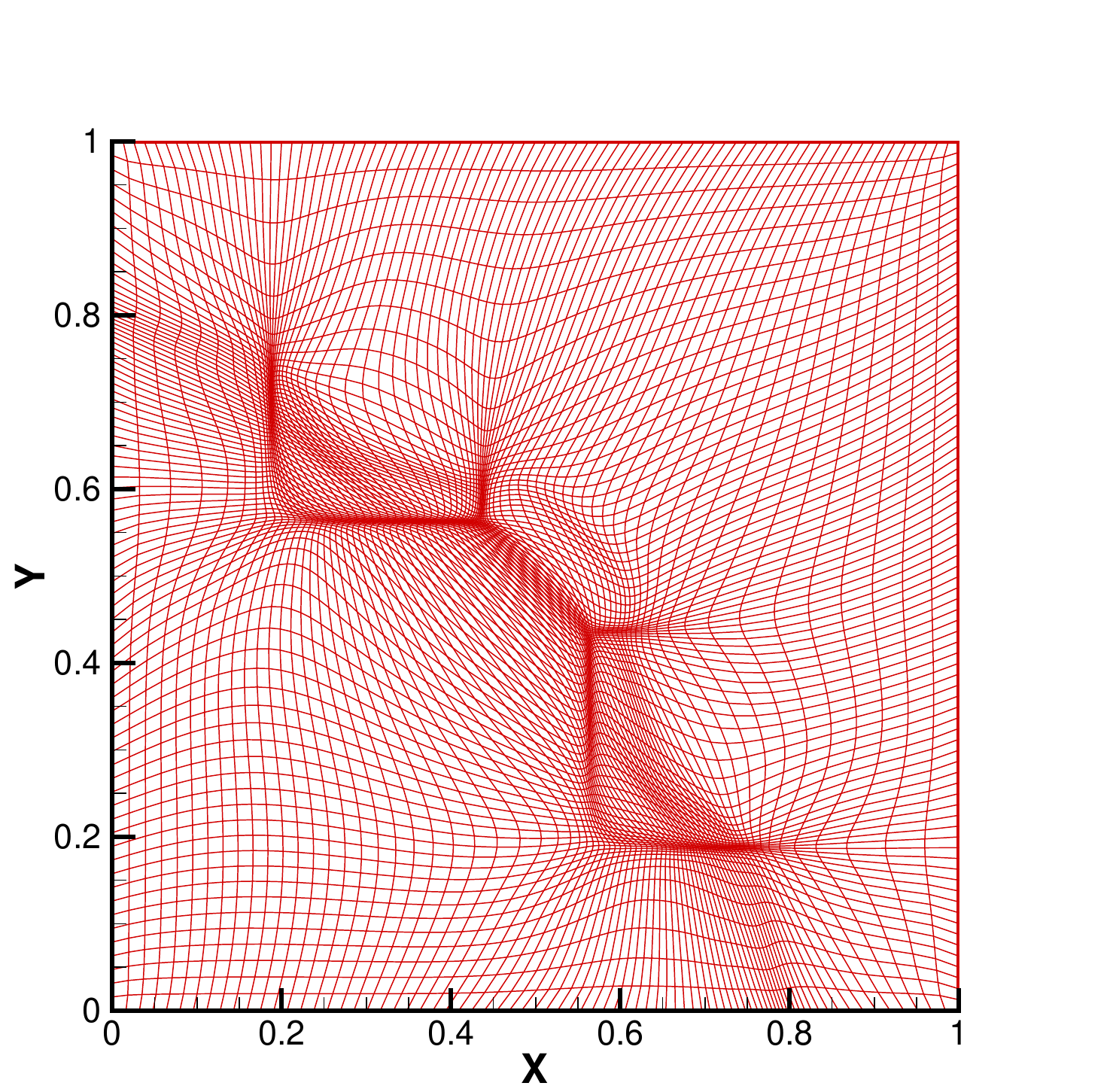}
\end{minipage}
\begin{minipage}[t]{2.3in}
\centerline{\scriptsize (e): t=1.0}
\includegraphics[width=2.3in]{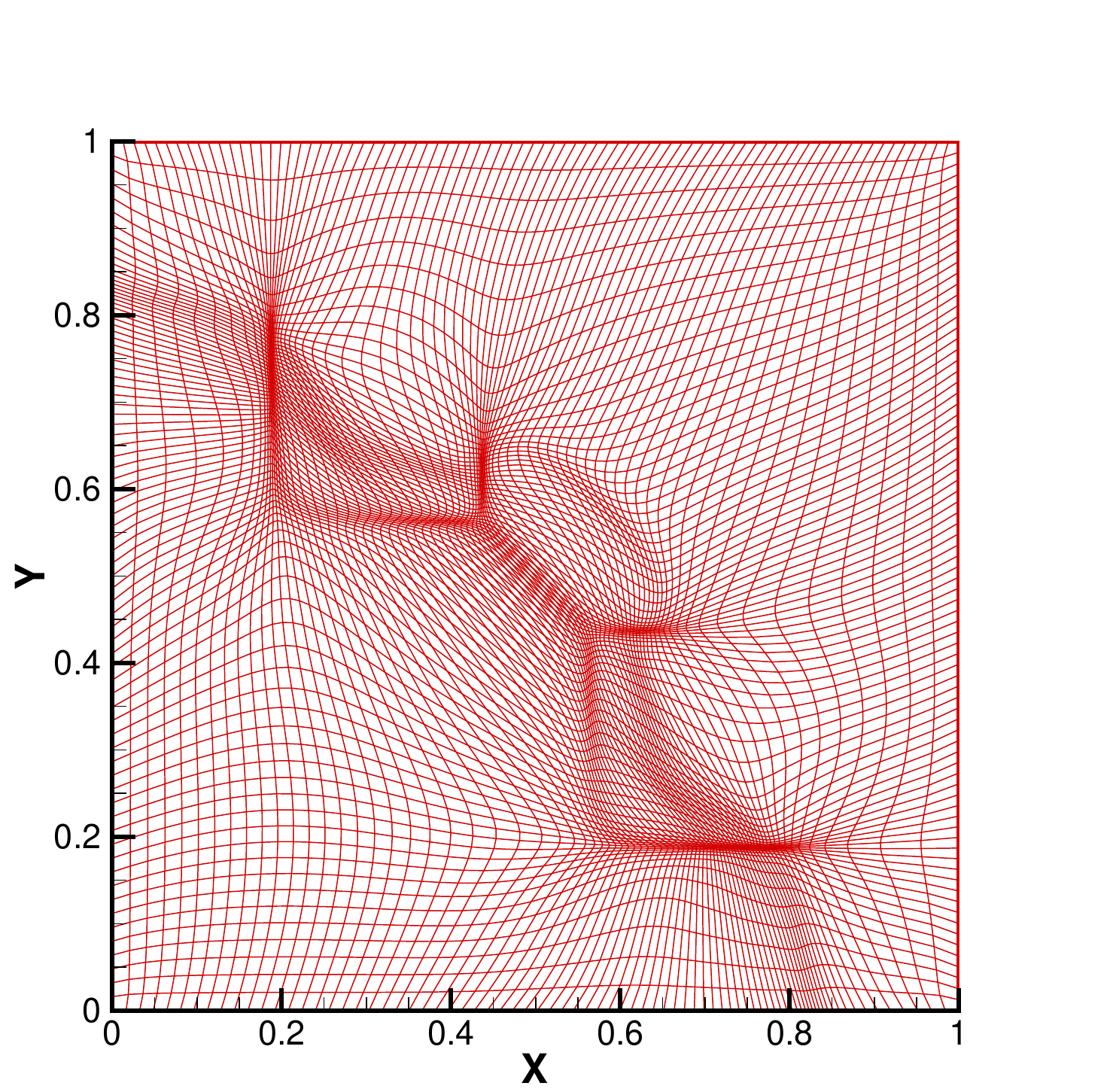}
\end{minipage}
\begin{minipage}[t]{2.3in}
\centerline{\scriptsize (f): t=1.5}
\includegraphics[width=2.3in]{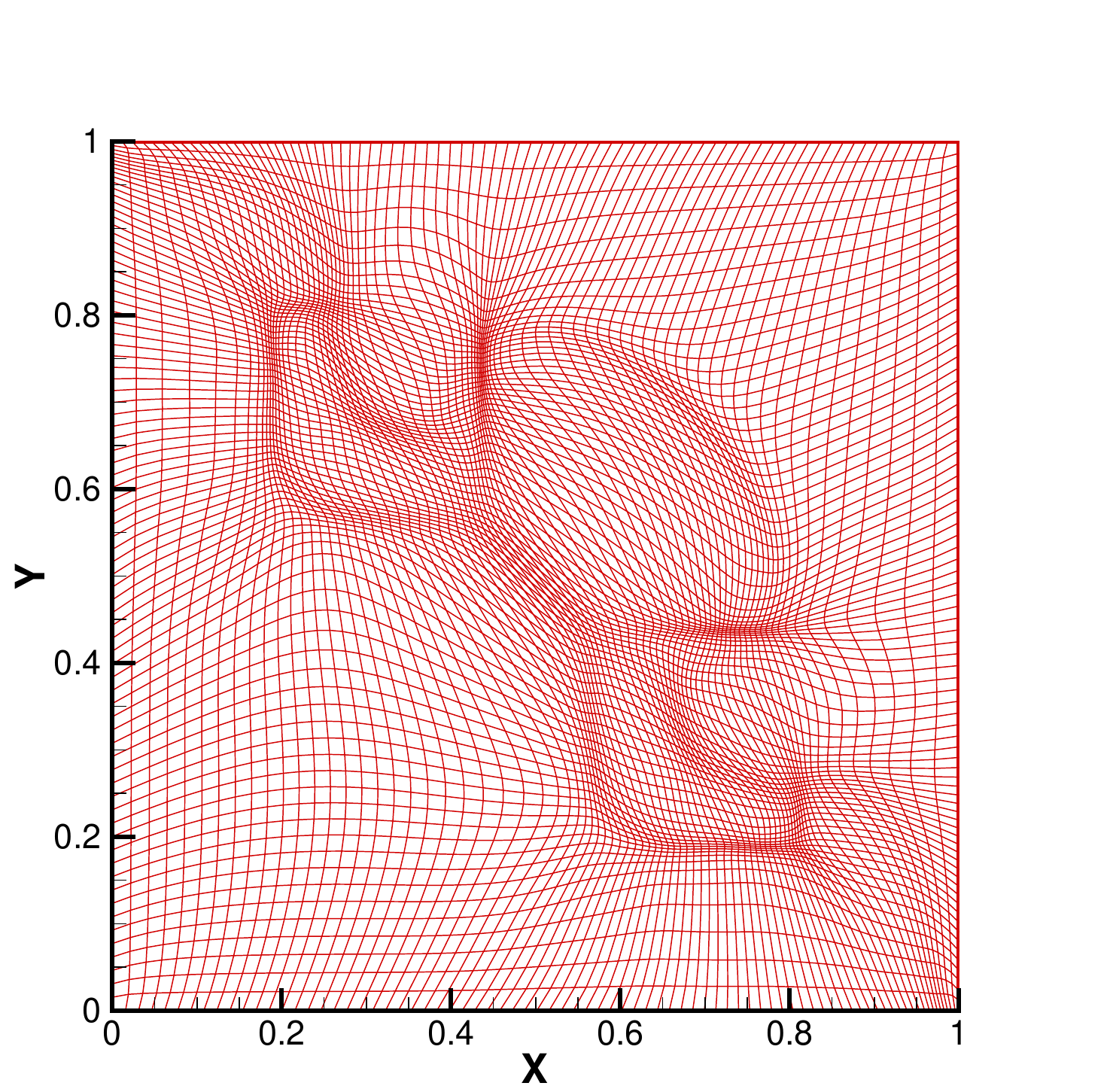}
\end{minipage}
}
\vspace{5mm}
\hbox{
\begin{minipage}[t]{2.3in}
\centerline{\scriptsize (g): t=2.0}
\includegraphics[width=2.3in]{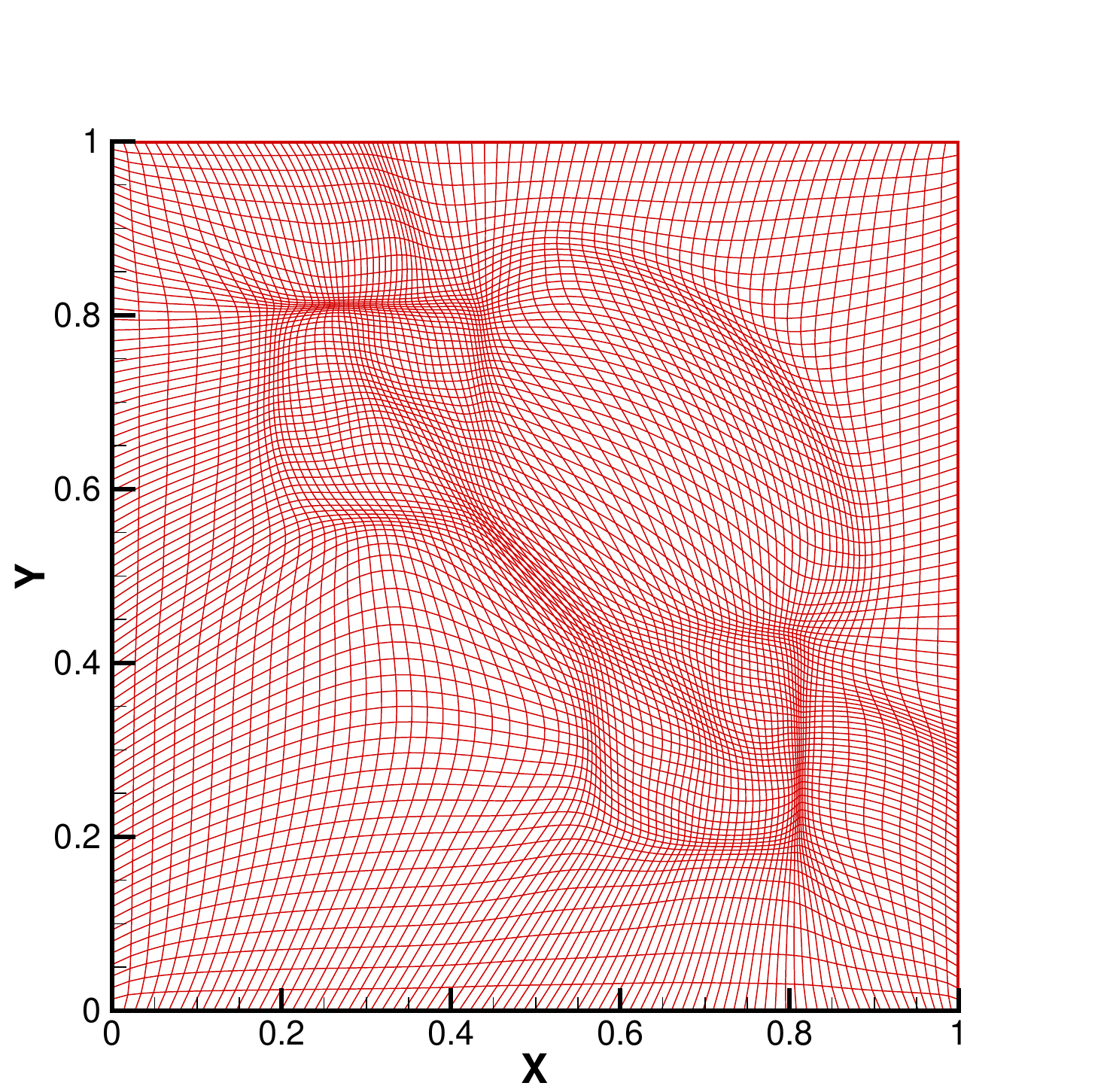}
\end{minipage}
\begin{minipage}[t]{2.3in}
\centerline{\scriptsize (h): t=2.5}
\includegraphics[width=2.3in]{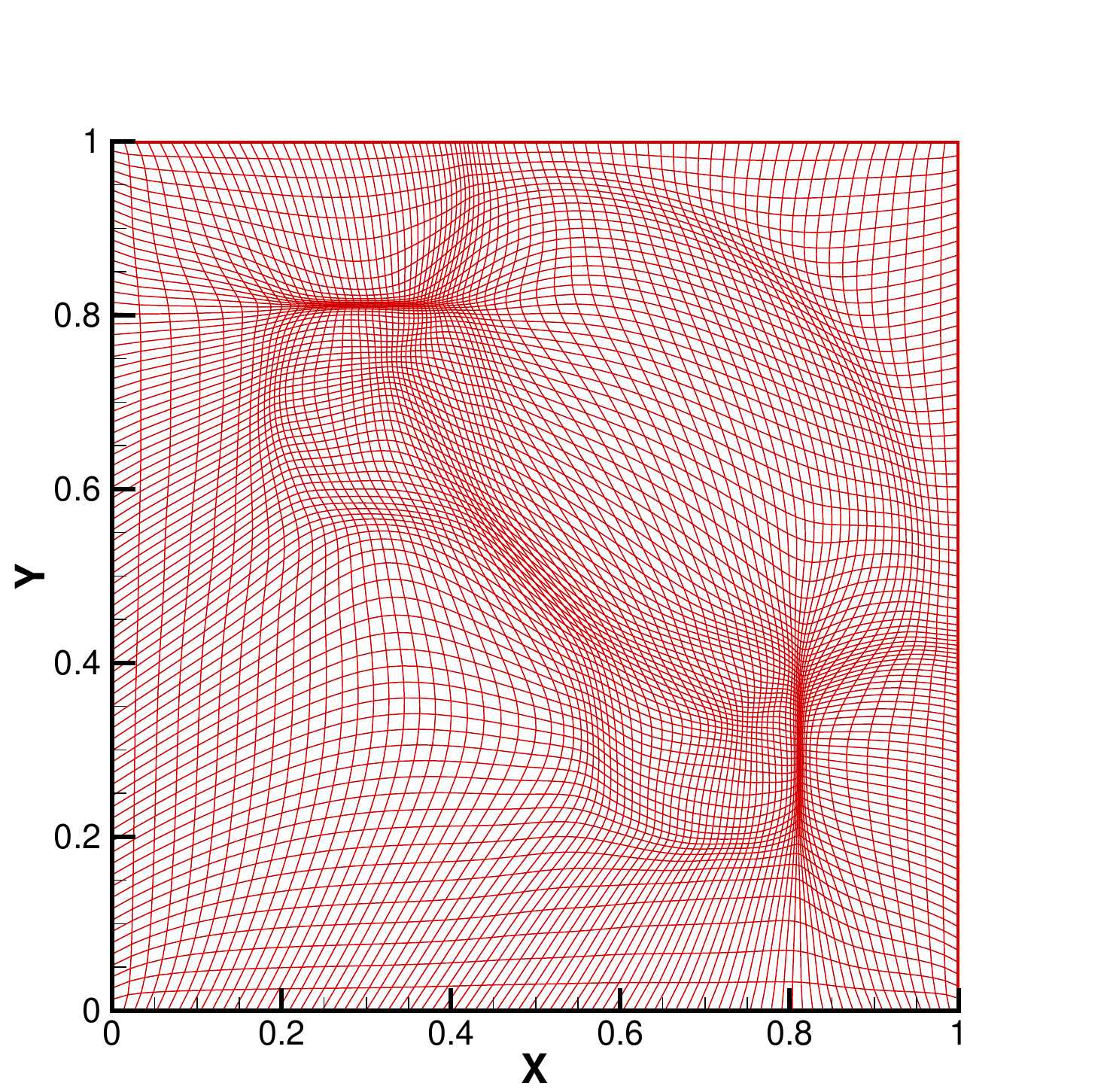}
\end{minipage}
\begin{minipage}[t]{2.3in}
\centerline{\scriptsize (i): t=3.0}
\includegraphics[width=2.3in]{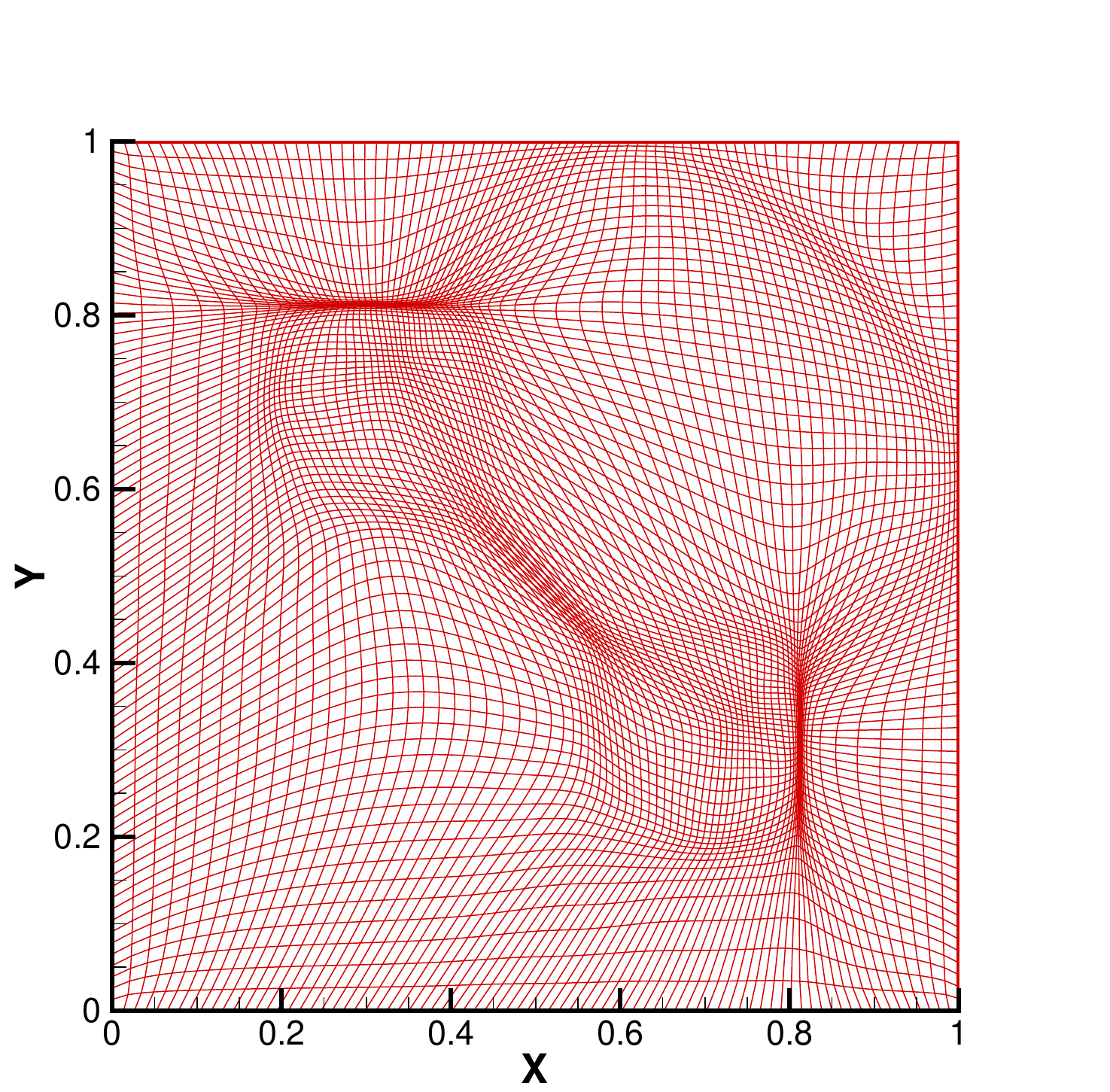}
\end{minipage}
}
\end{center}
\caption{Example~\ref{Example4.3}. The moving mesh of $81\times 81$ is shown at
$t=0.5, 0.7, 0.8, 0.9, 1.0, 1.5, 2.0, 2.5, 3.0$. }
\label{T29}
\end{figure}

\begin{figure}
\begin{center}
\hbox{
\hspace{1in}
\begin{minipage}[t]{2.0in}
\centerline{\scriptsize (a):  with MM at $t=1.0$}
\includegraphics[width=2.0in]{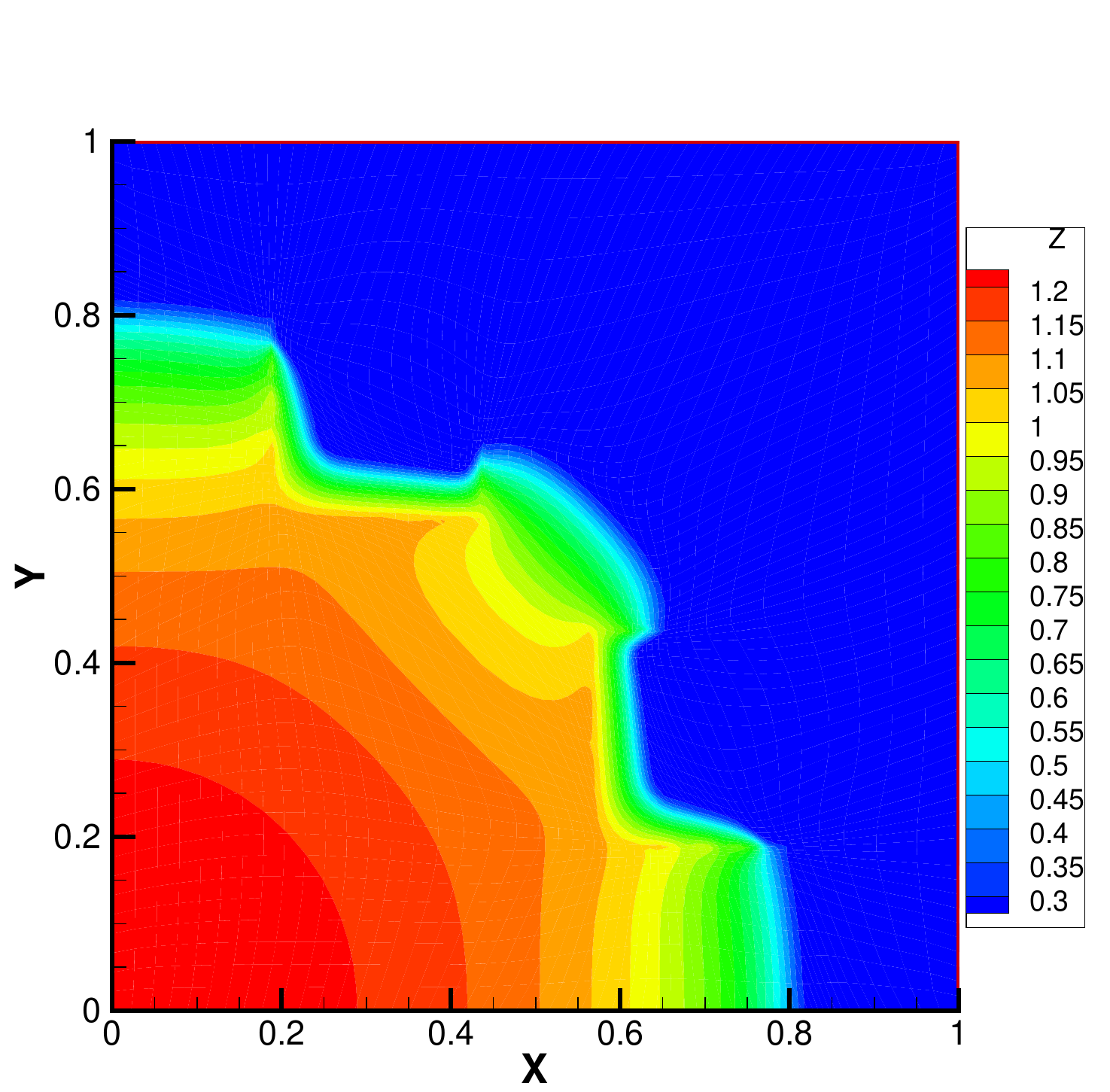}
\end{minipage}
\hspace{0.5in}
\begin{minipage}[t]{2.0in}
\centerline{\scriptsize (b):  with UM at $t=1.0$}
\includegraphics[width=2.0in]{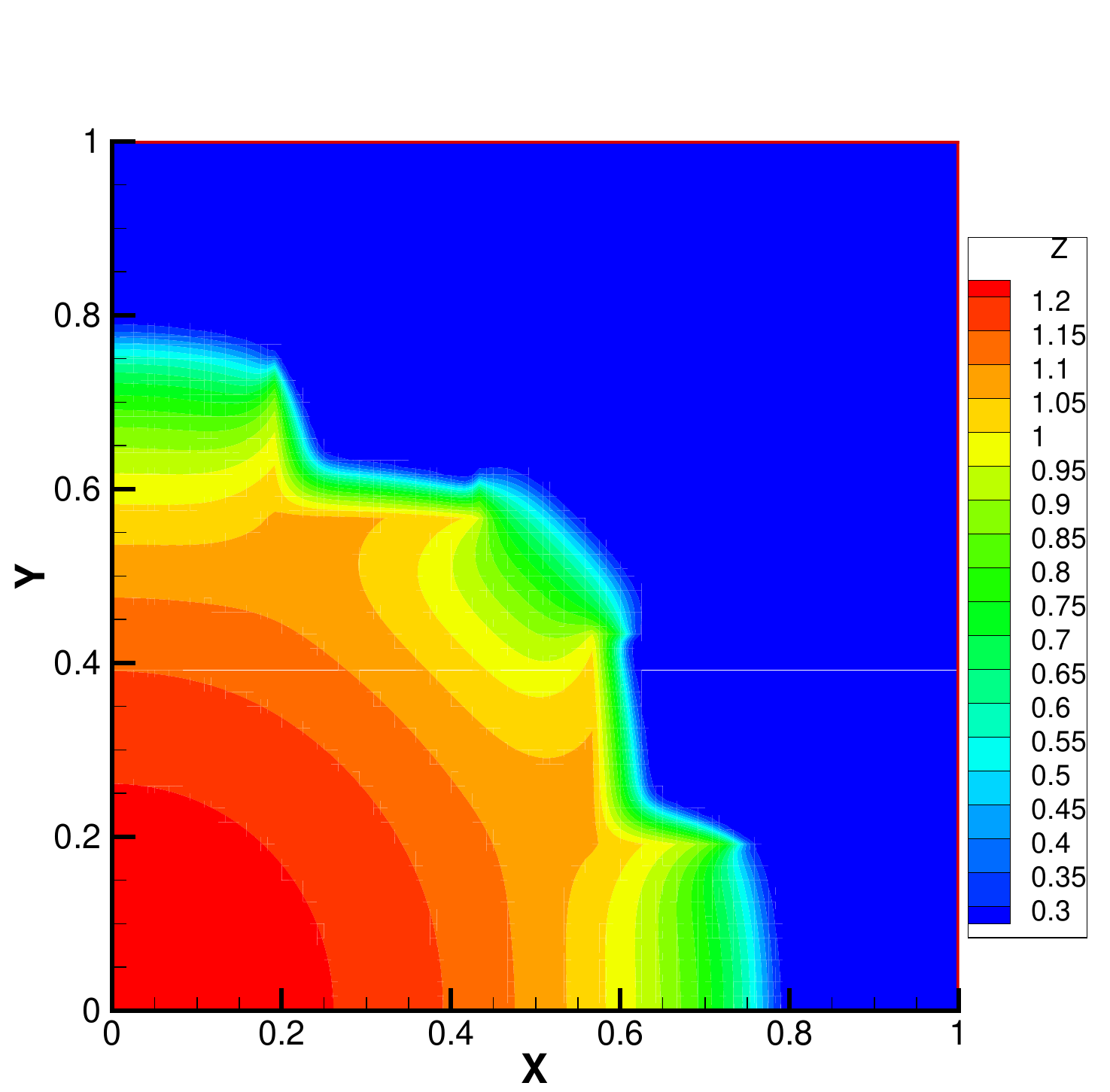}
\end{minipage}
}
\hbox{
\hspace{1in}
\begin{minipage}[t]{2.0in}
\centerline{\scriptsize (c):  with MM at $t=2.0$}
\includegraphics[width=2.0in]{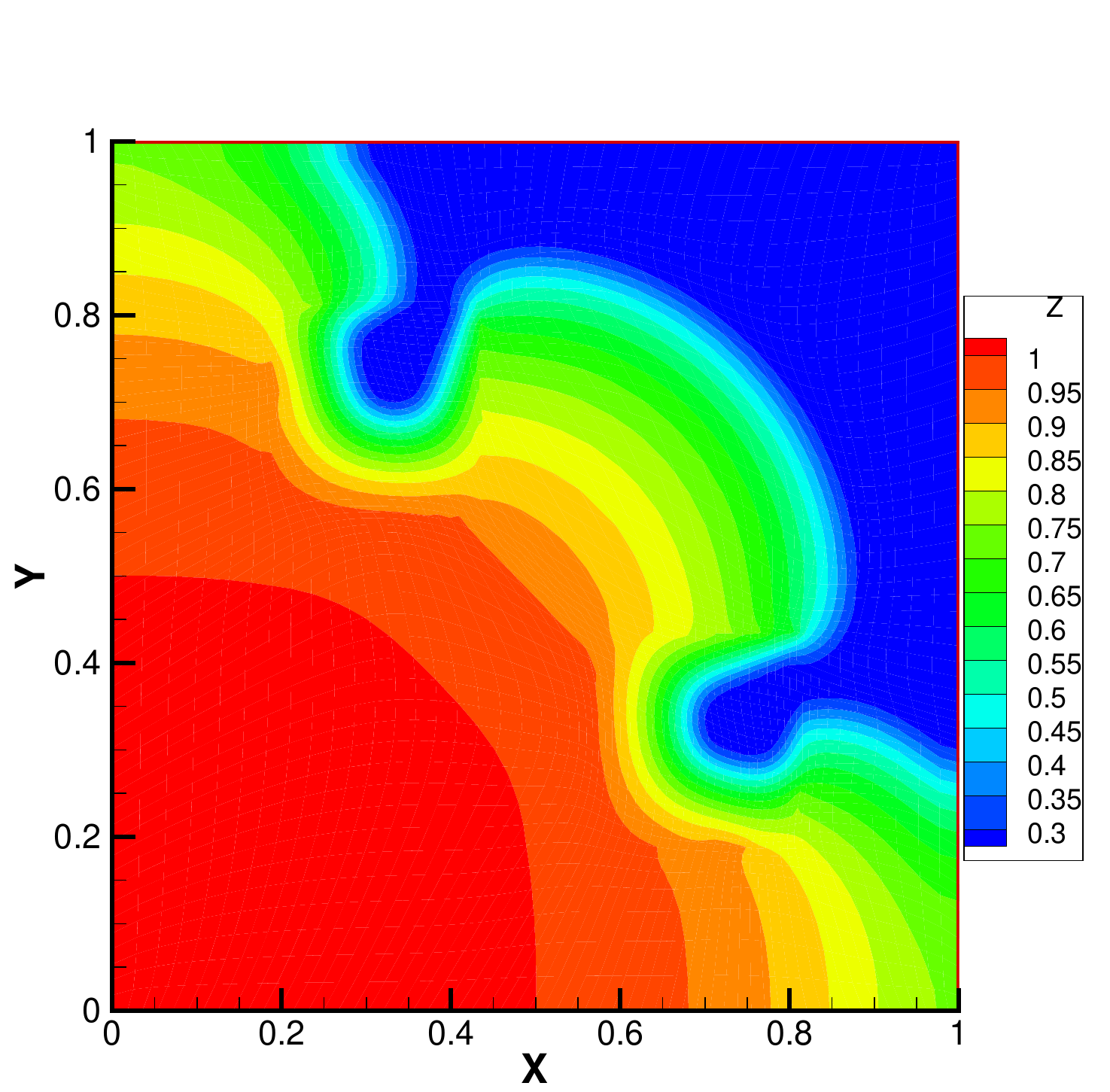}
\end{minipage}
\hspace{0.5in}
\begin{minipage}[t]{2.0in}
\centerline{\scriptsize (d):  with UM at $t=2.0$}
\includegraphics[width=2.0in]{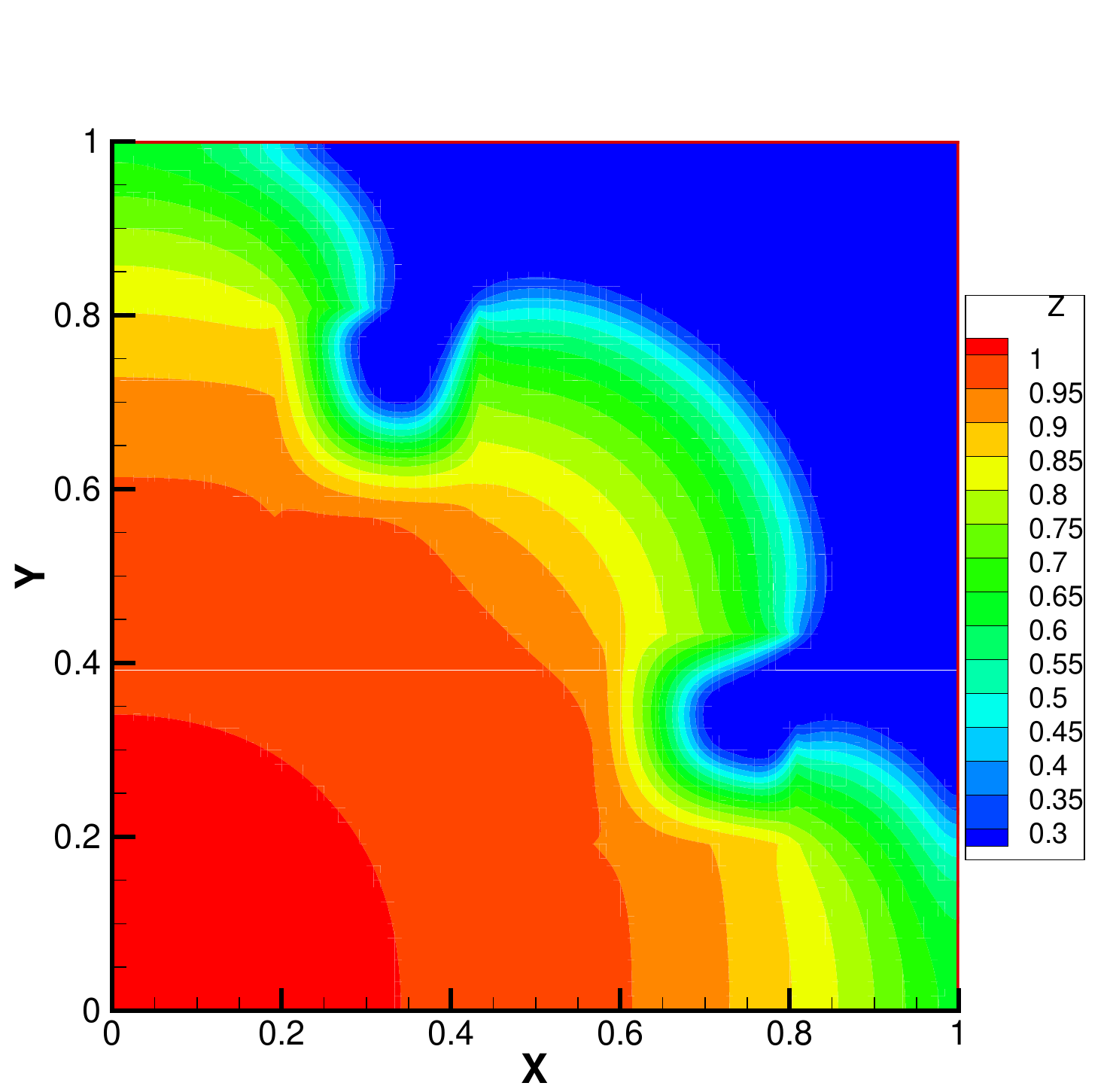}
\end{minipage}
}
\hbox{
\hspace{1in}
\begin{minipage}[t]{2.0in}
\centerline{\scriptsize (e):  with MM at $t=2.5$}
\includegraphics[width=2.0in]{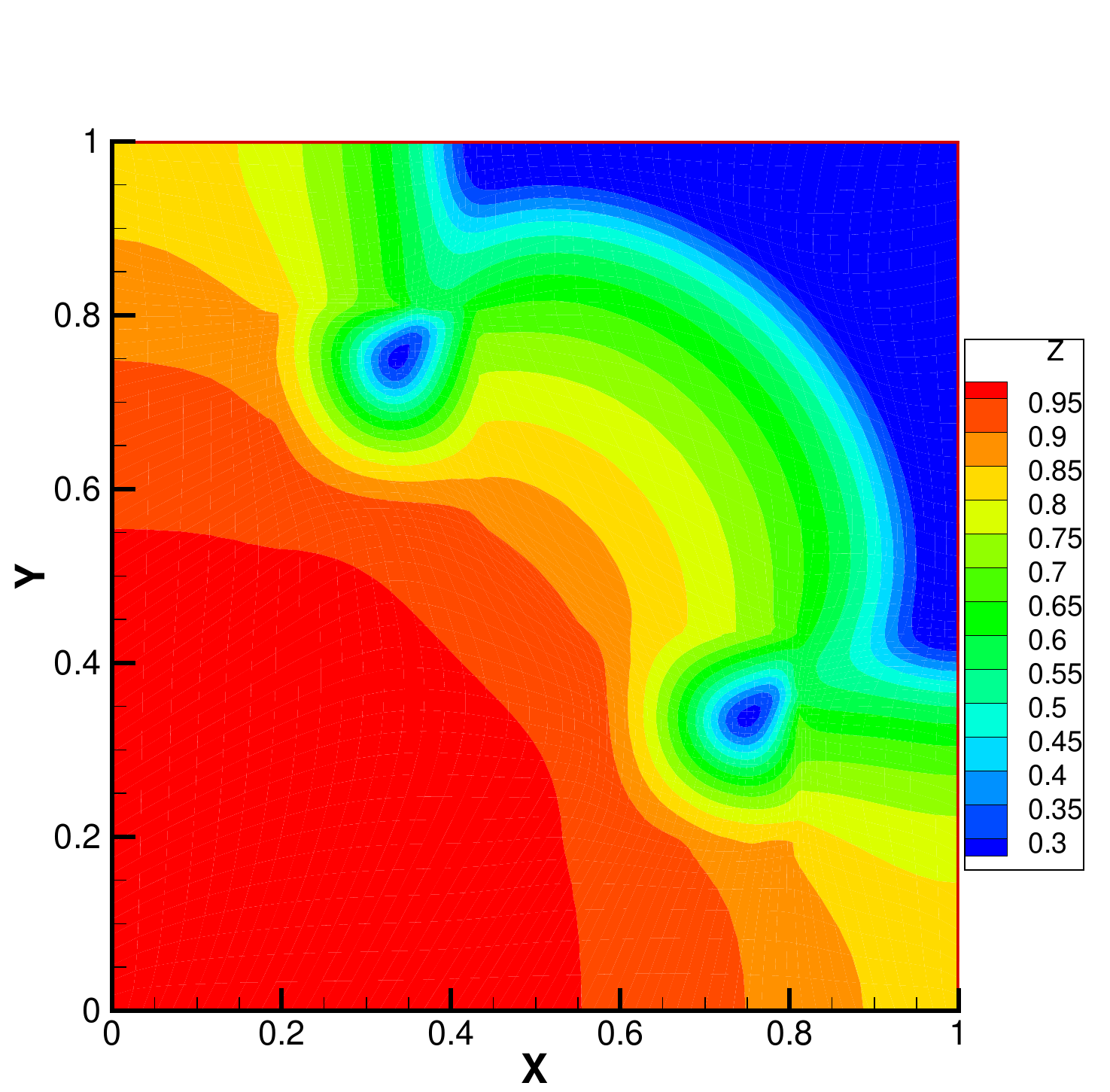}
\end{minipage}
\hspace{0.5in}
\begin{minipage}[t]{2.0in}
\centerline{\scriptsize (f):  with UM at $t=2.5$}
\includegraphics[width=2.0in]{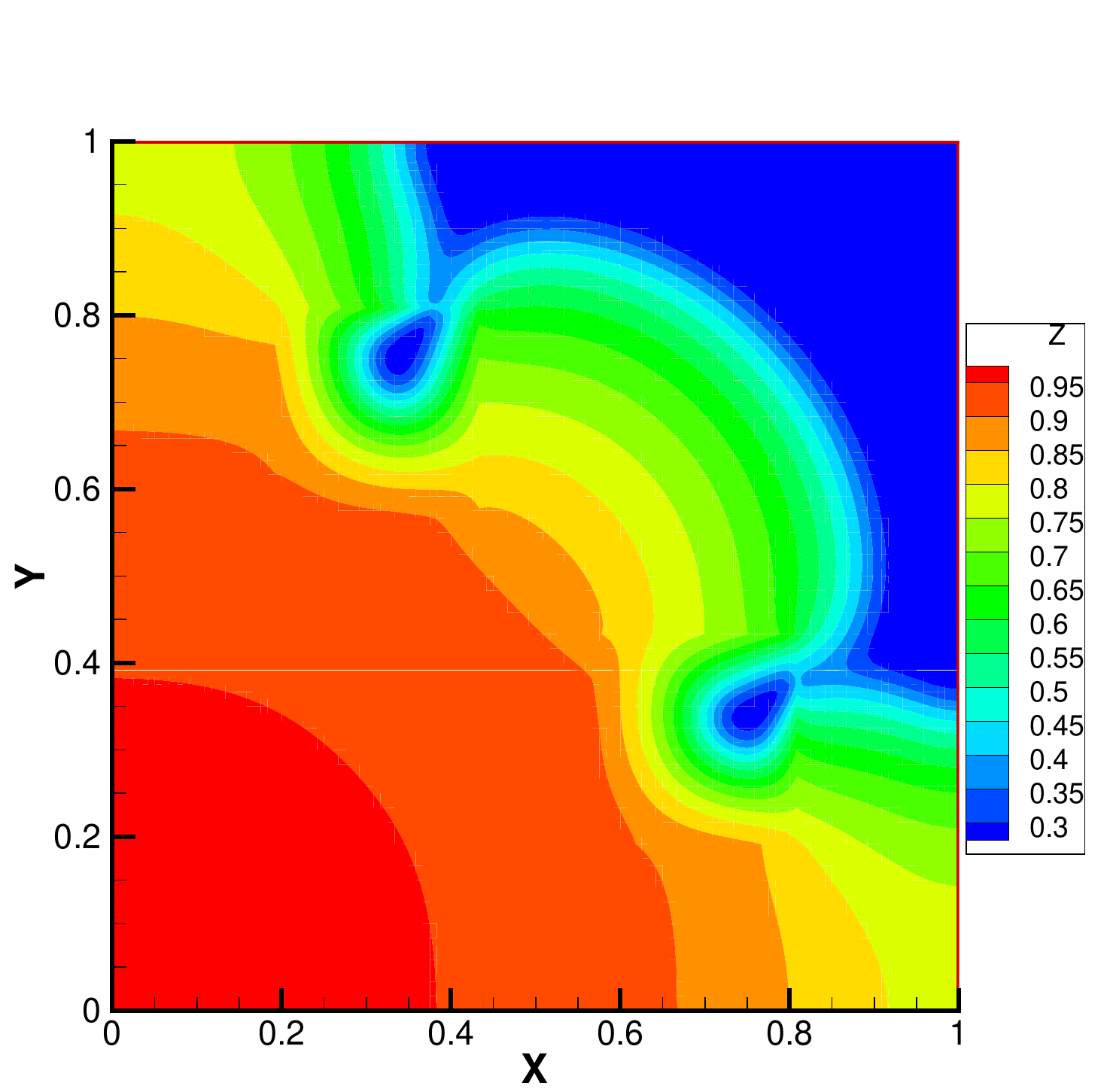}
\end{minipage}
}
\hbox{
\hspace{1in}
\begin{minipage}[t]{2.0in}
\centerline{\scriptsize (g):  with MM at $t=3.0$}
\includegraphics[width=2.0in]{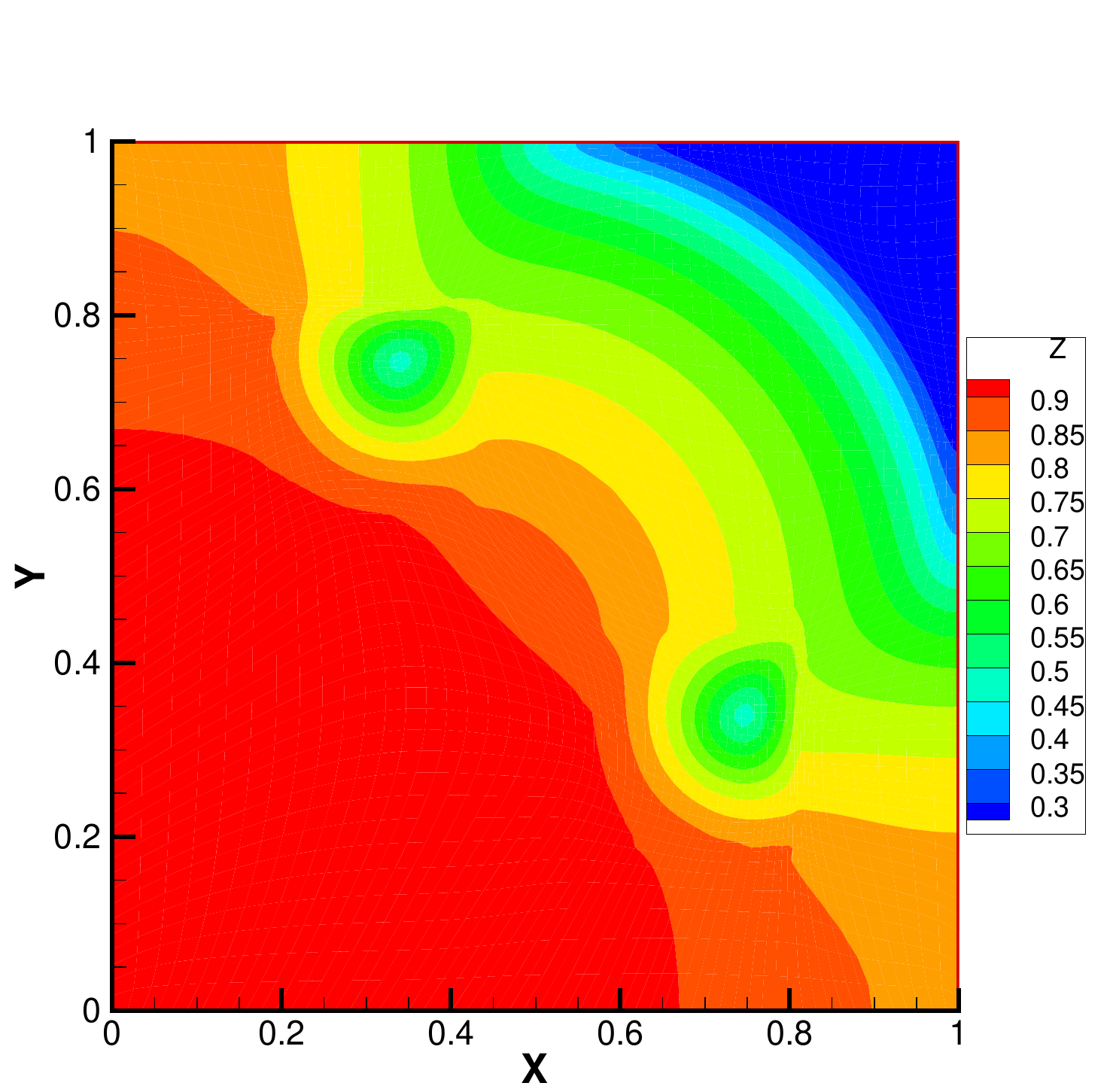}
\end{minipage}
\hspace{0.5in}
\begin{minipage}[t]{2.0in}
\centerline{\scriptsize (h):  with UM at $t=3.0$}
\includegraphics[width=2.0in]{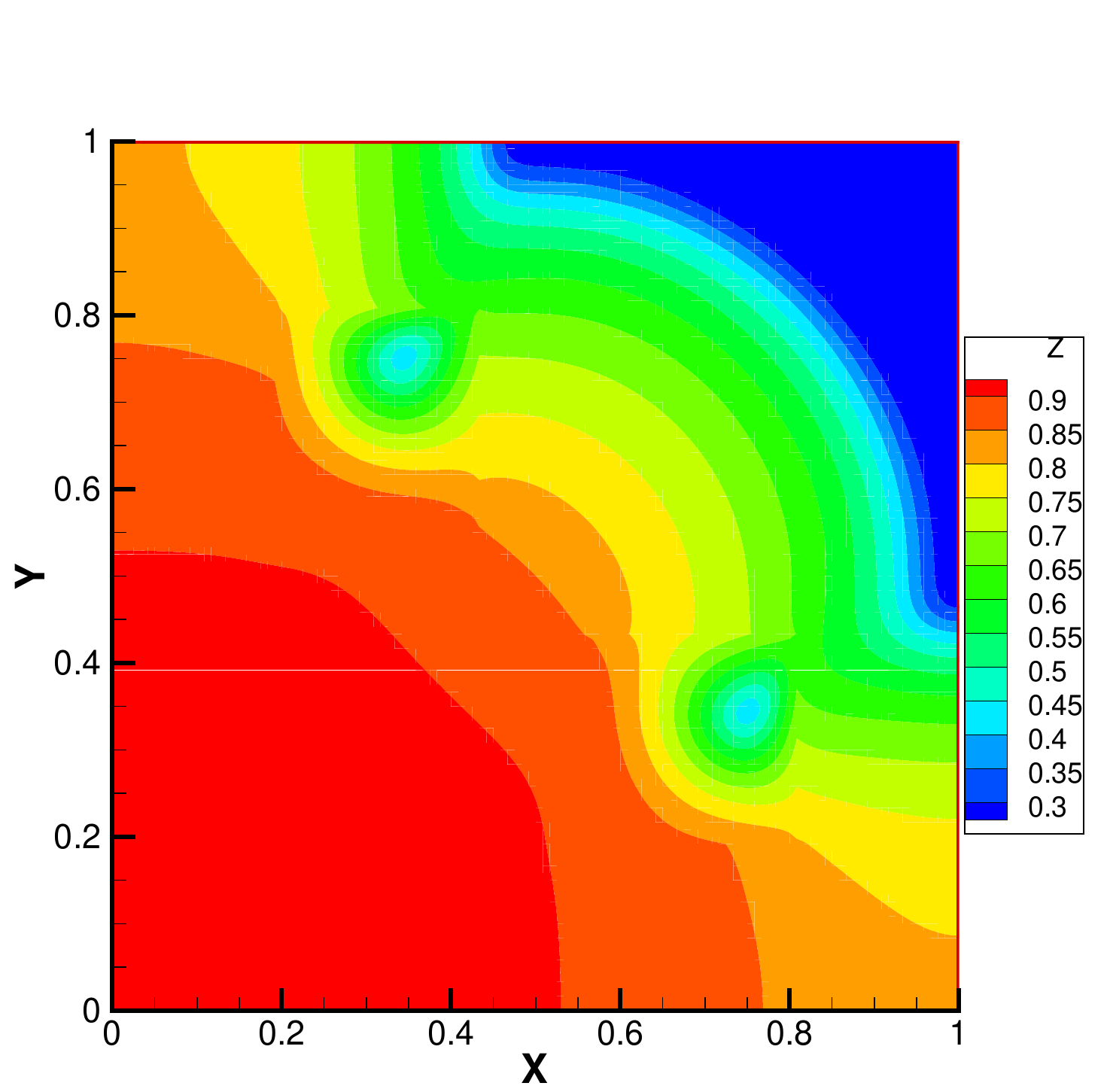}
\end{minipage}
}
\end{center}
\caption{Example~\ref{Example4.3}. The contours of the temperature obtained
with a moving mesh (MM) of size  $61\times 61$ are compared with those
obtained with a uniform mesh (UM) of size $121 \times121$.}
\label{T30}
\end{figure}

\begin{figure}
\begin{center}
\hbox{
\hspace{1in}
\begin{minipage}[t]{2.0in}
\centerline{\scriptsize (a):  with MM1 at $t=1.0$}
\includegraphics[width=2.0in]{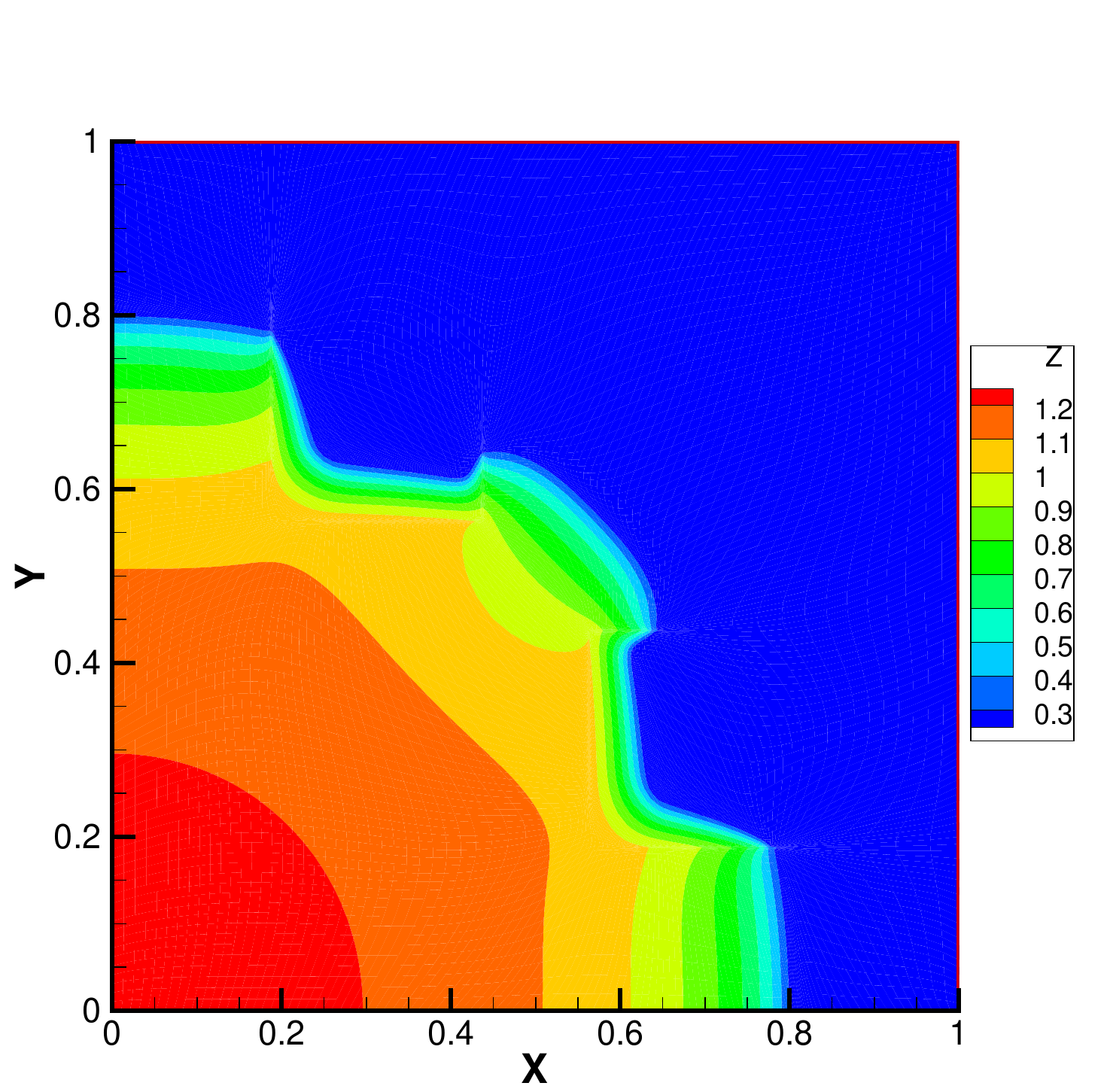}
\end{minipage}
\hspace{0.5in}
\begin{minipage}[t]{2.0in}
\centerline{\scriptsize (b):  with MM2 at $t=1.0$}
\includegraphics[width=2.0in]{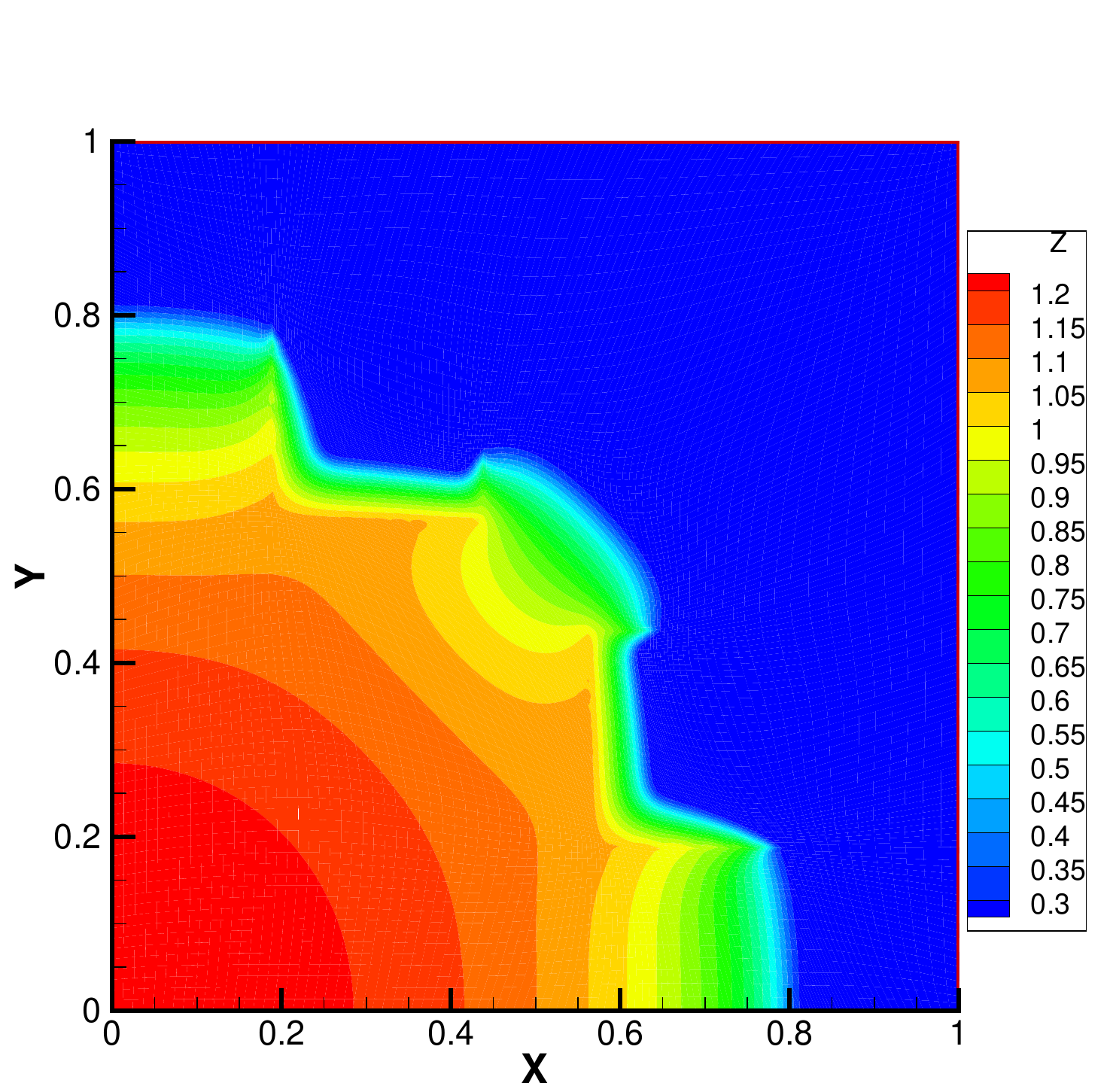}
\end{minipage}
}
\hbox{
\hspace{1in}
\begin{minipage}[t]{2.0in}
\centerline{\scriptsize (c):  with MM1 at $t=2.0$}
\includegraphics[width=2.0in]{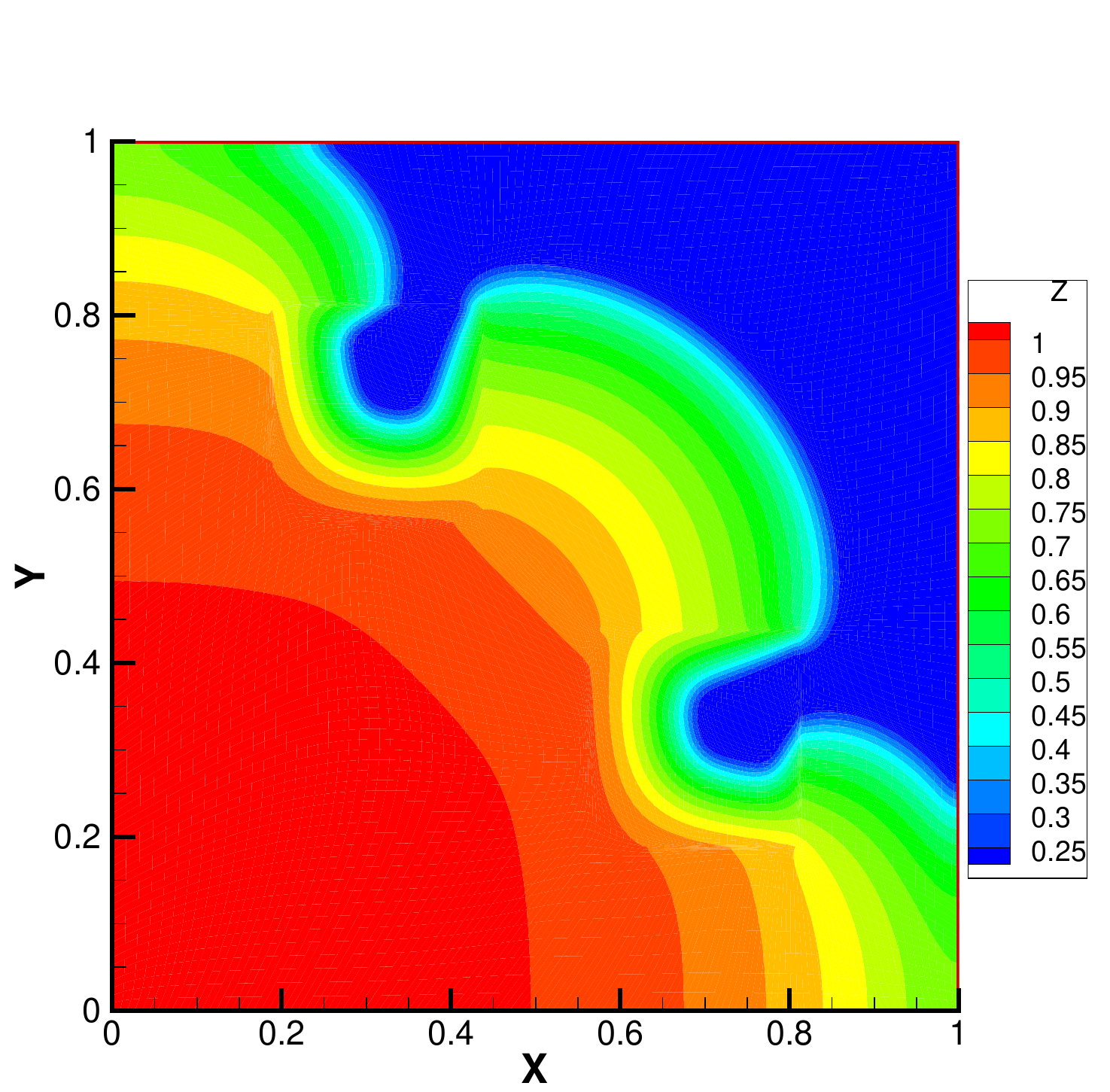}
\end{minipage}
\hspace{0.5in}
\begin{minipage}[t]{2.0in}
\centerline{\scriptsize (d):  with MM2 at $t=2.0$}
\includegraphics[width=2.0in]{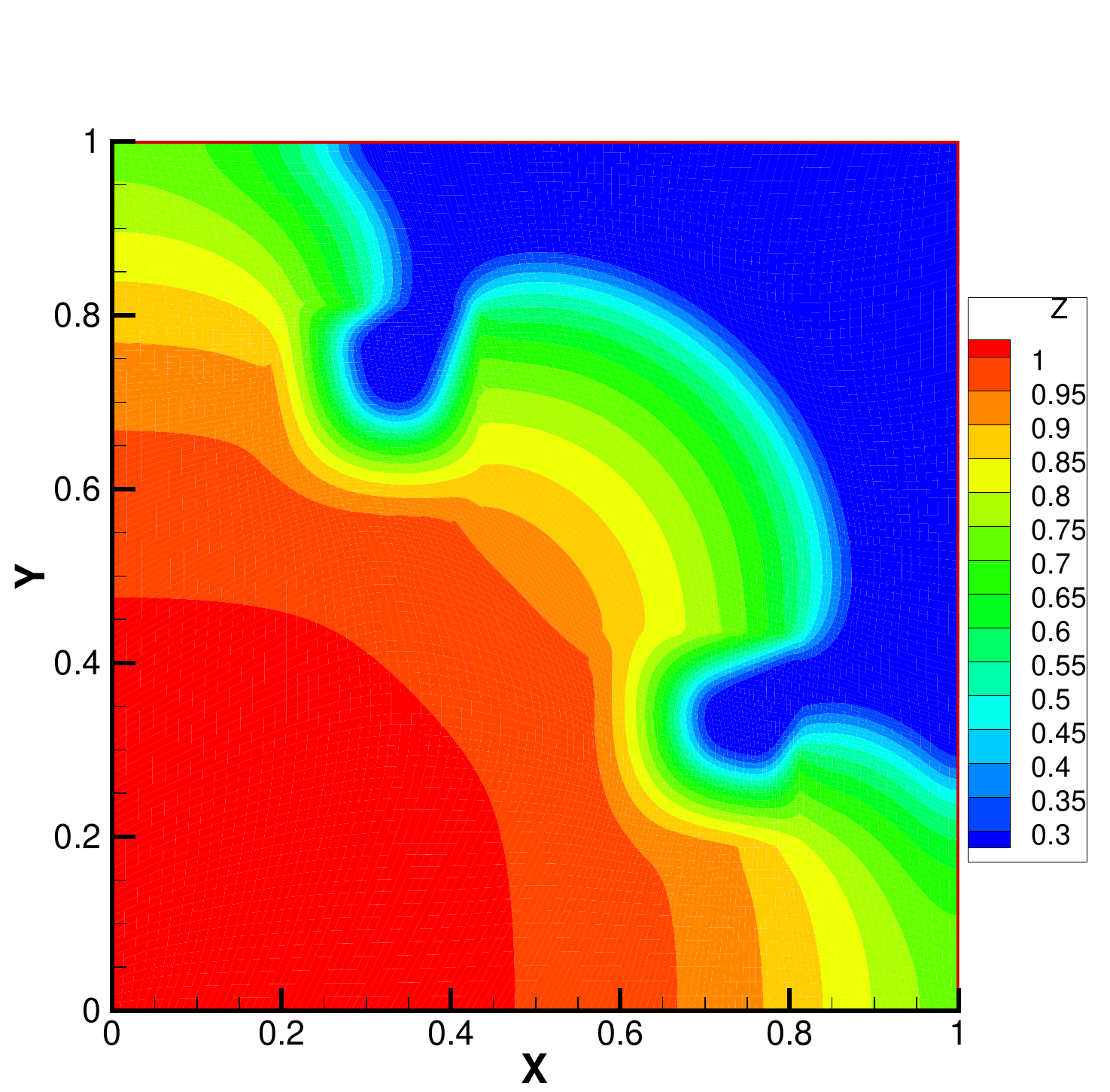}
\end{minipage}
}
\hbox{
\hspace{1in}
\begin{minipage}[t]{2.0in}
\centerline{\scriptsize (e):  with MM1 at $t=2.5$}
\includegraphics[width=2.0in]{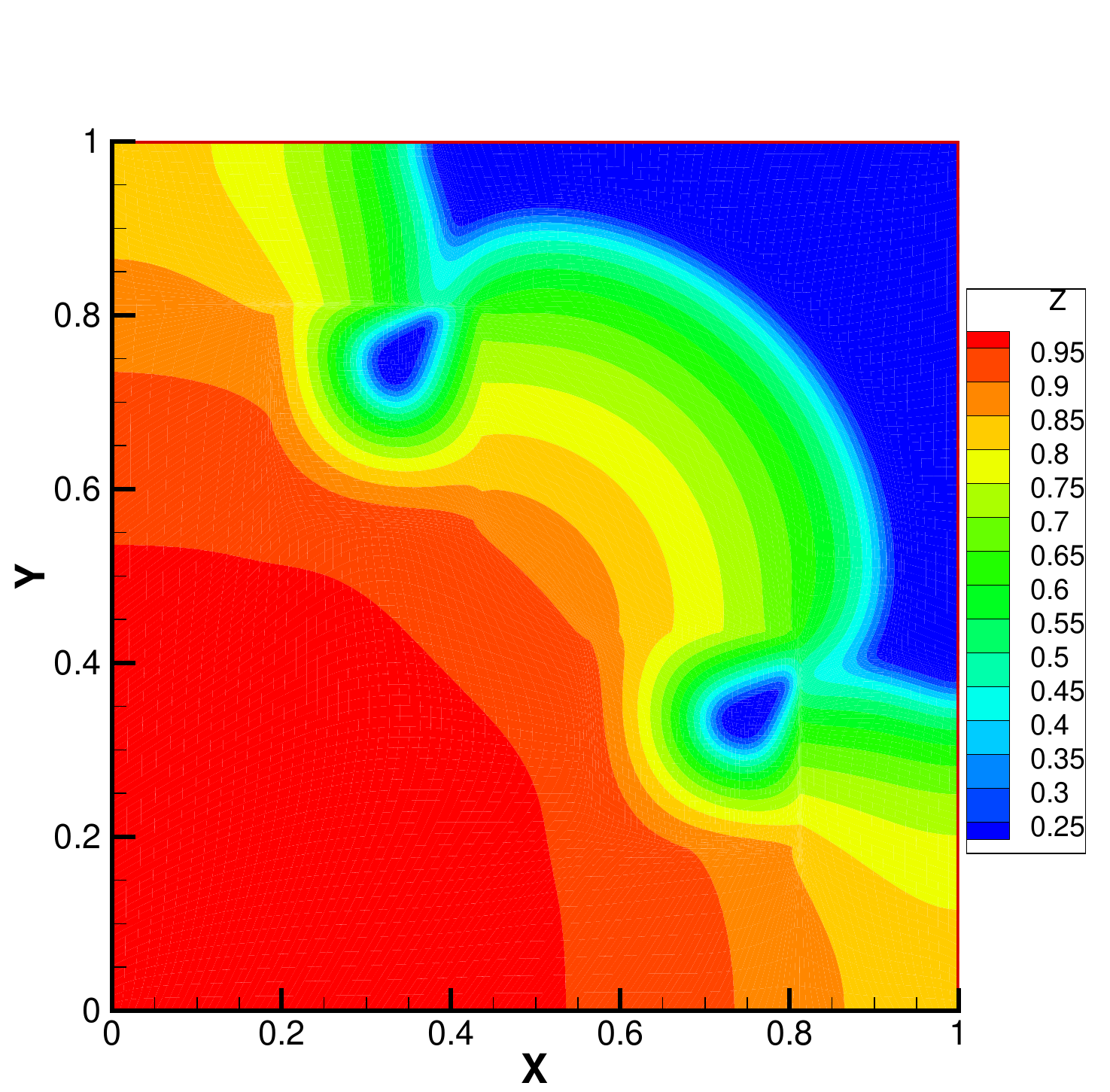}
\end{minipage}
\hspace{0.5in}
\begin{minipage}[t]{2.0in}
\centerline{\scriptsize (f):  with MM2 at $t=2.5$}
\includegraphics[width=2.0in]{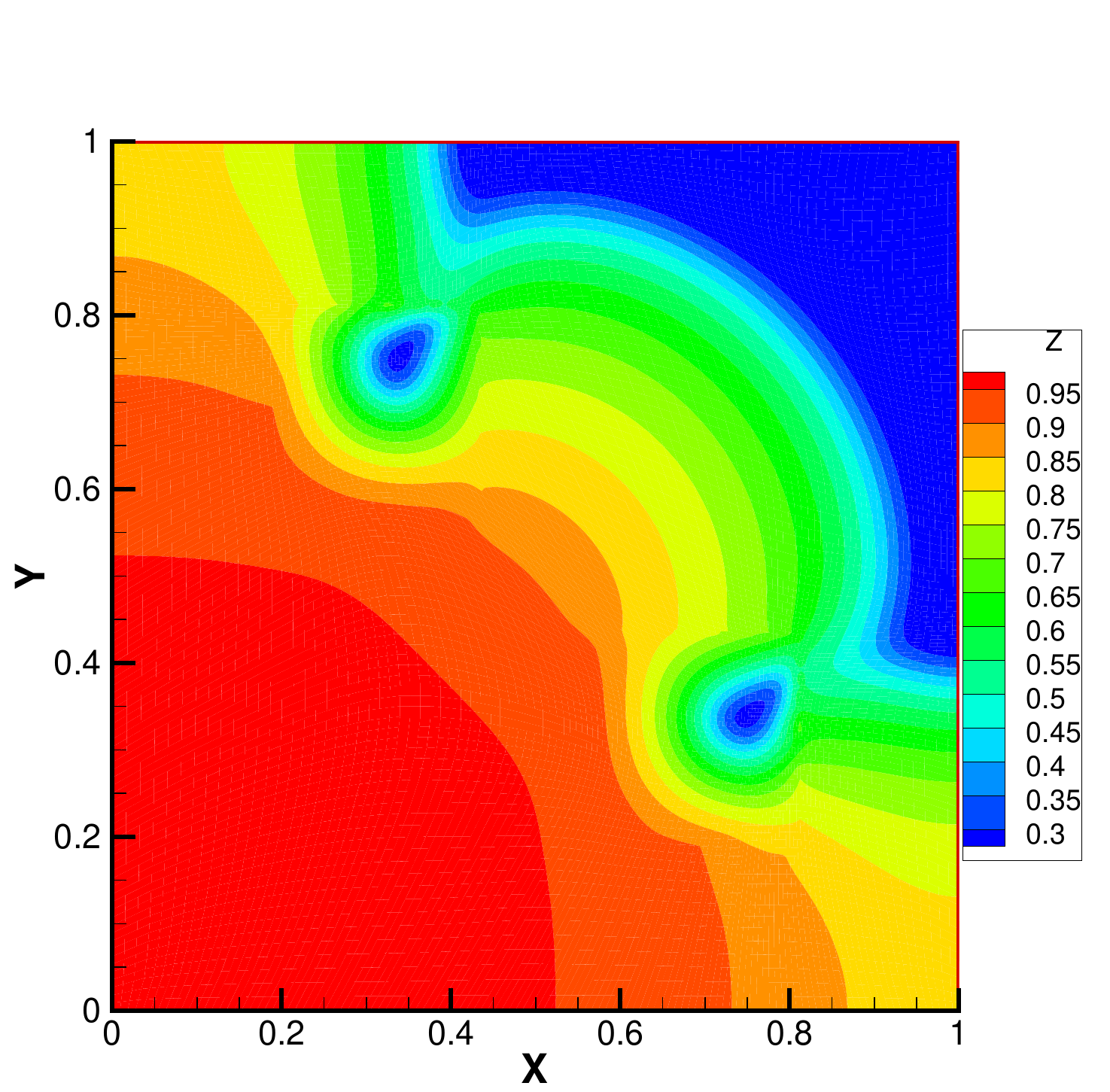}
\end{minipage}
}
\hbox{
\hspace{1in}
\begin{minipage}[t]{2.0in}
\centerline{\scriptsize (g):  with MM1 at $t=3.0$}
\includegraphics[width=2.0in]{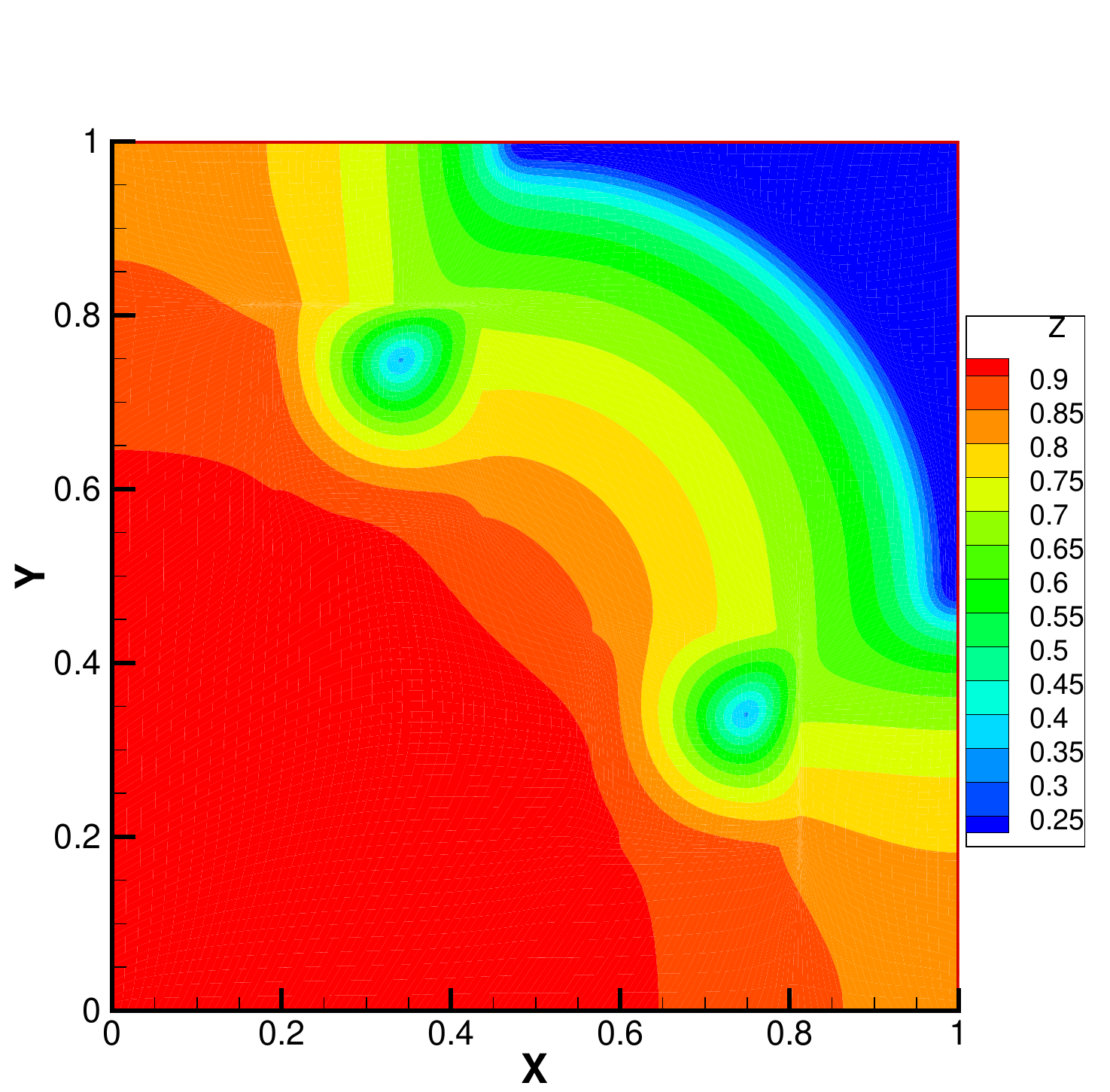}
\end{minipage}
\hspace{0.5in}
\begin{minipage}[t]{2.0in}
\centerline{\scriptsize (h):  with MM2 at $t=3.0$}
\includegraphics[width=2.0in]{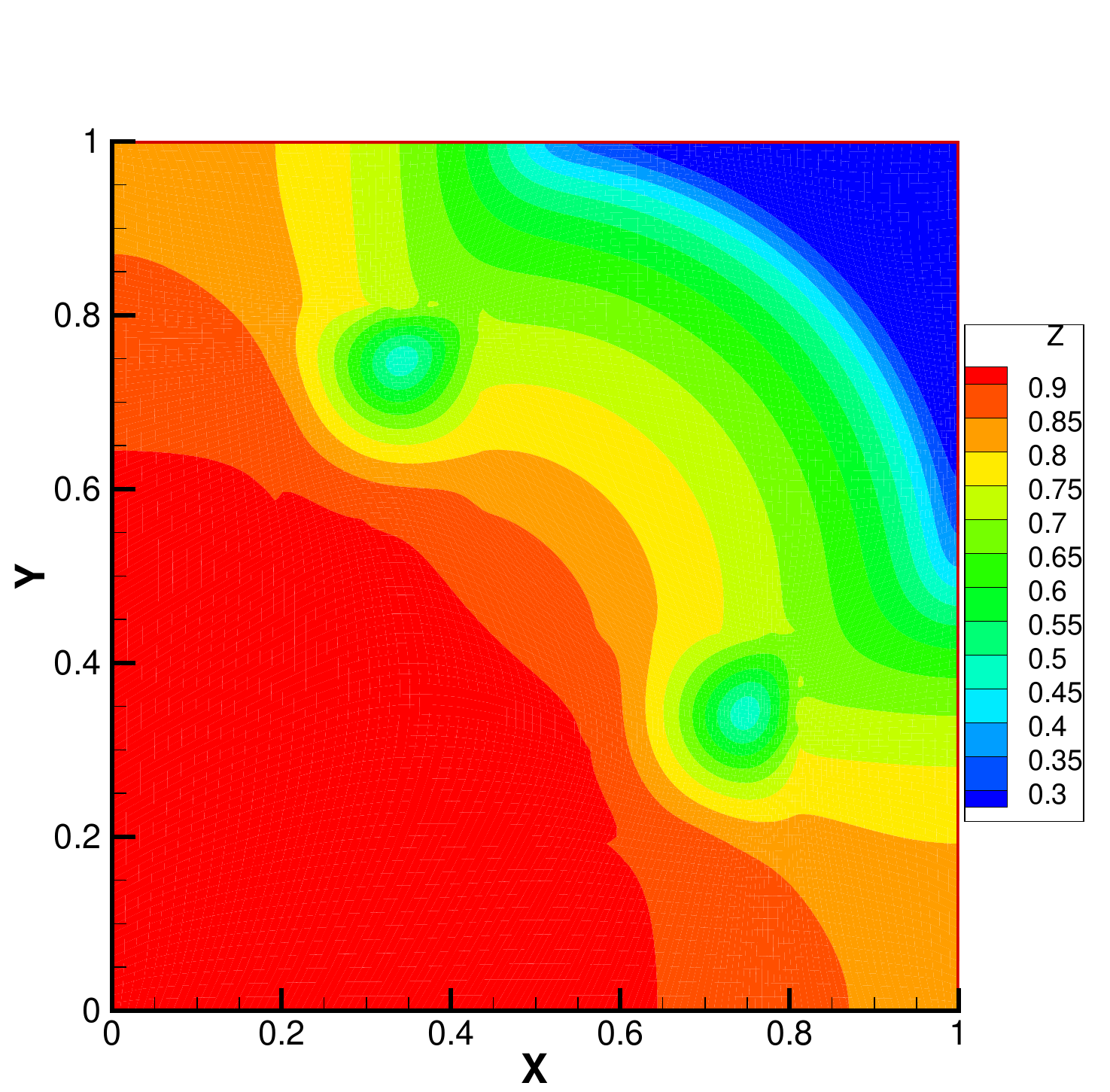}
\end{minipage}
}
\end{center}
\caption{Example~\ref{Example4.3}.The contours of the temperature obtained
with an MM1 of size  $121\times 121$ are compared with those
obtained with an MM2 of size  $41\times 41$ (with the physical PDE
being solved on a uniformly refined mesh of size $121\times 121$).}
\label{T31}
\end{figure}

\begin{figure}
\begin{center}
\hbox{
\hspace{1in}
\begin{minipage}[t]{2.0in}
\centerline{\scriptsize (a):  with MM1 at $t=1.0$}
\includegraphics[width=2.0in]{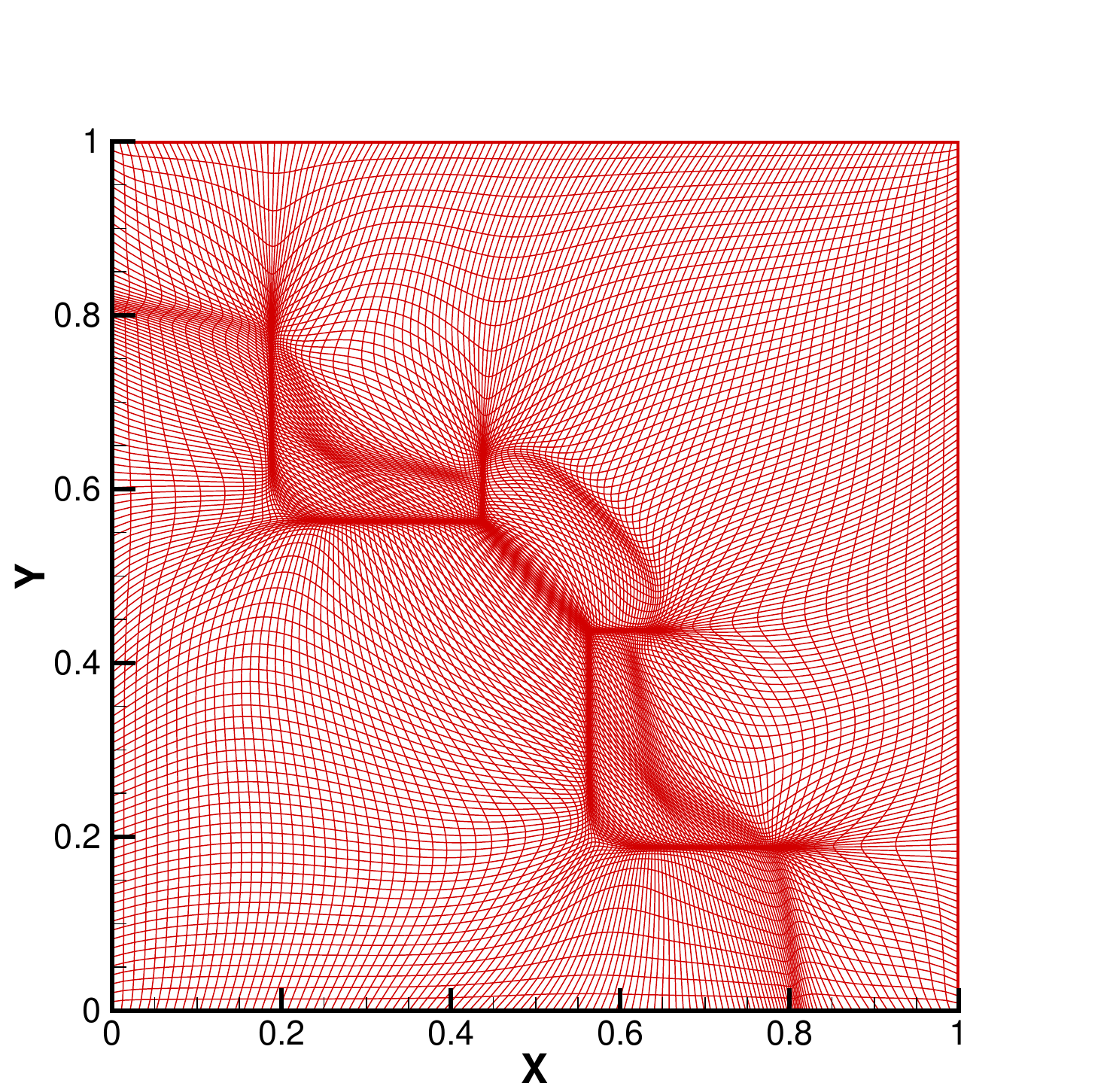}
\end{minipage}
\hspace{0.5in}
\begin{minipage}[t]{2.0in}
\centerline{\scriptsize (b):  with MM2 at $t=1.0$}
\includegraphics[width=2.0in]{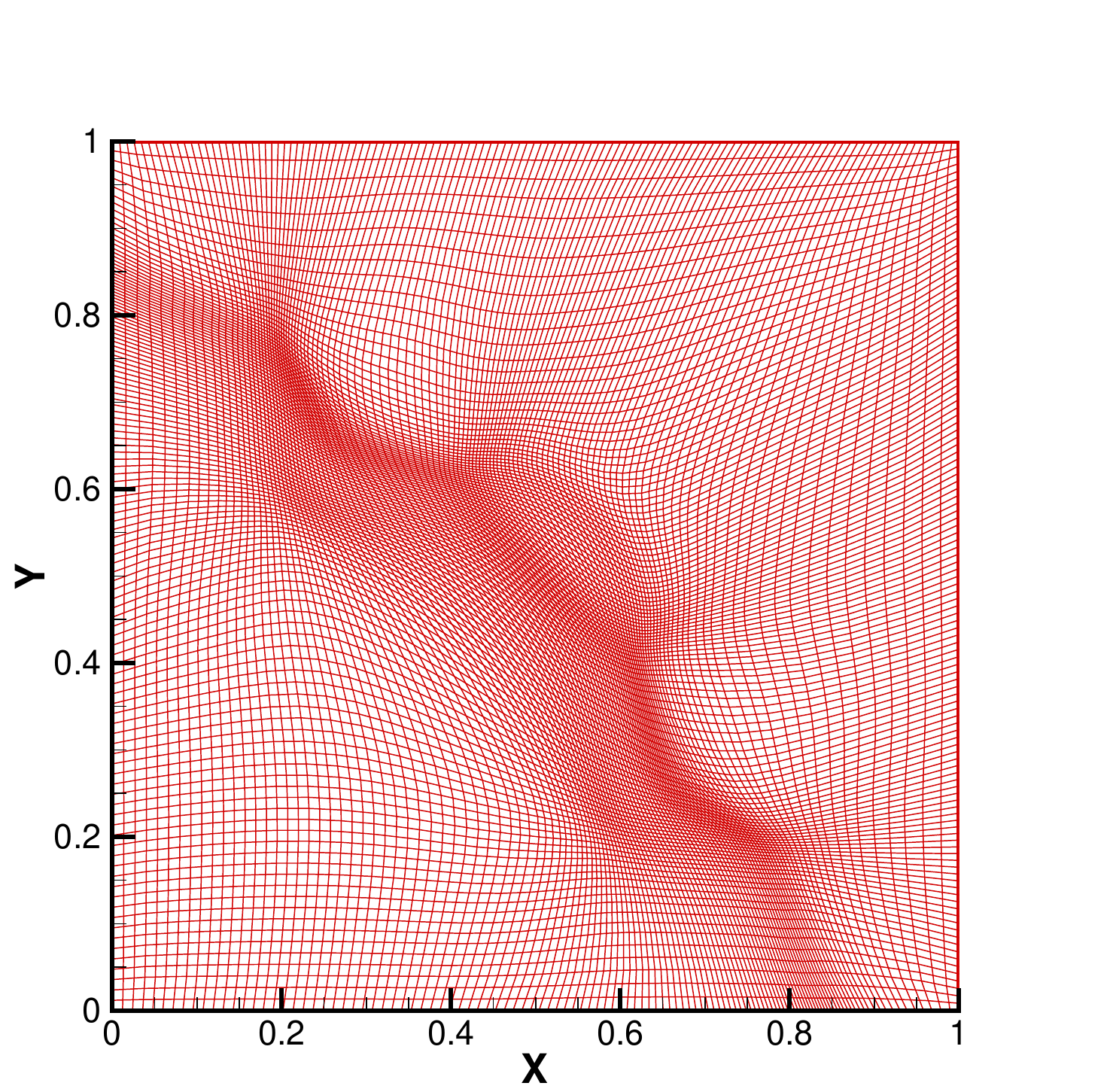}
\end{minipage}
}
\hbox{
\hspace{1in}
\begin{minipage}[t]{2.0in}
\centerline{\scriptsize (c):  with MM1 at $t=2.0$}
\includegraphics[width=2.0in]{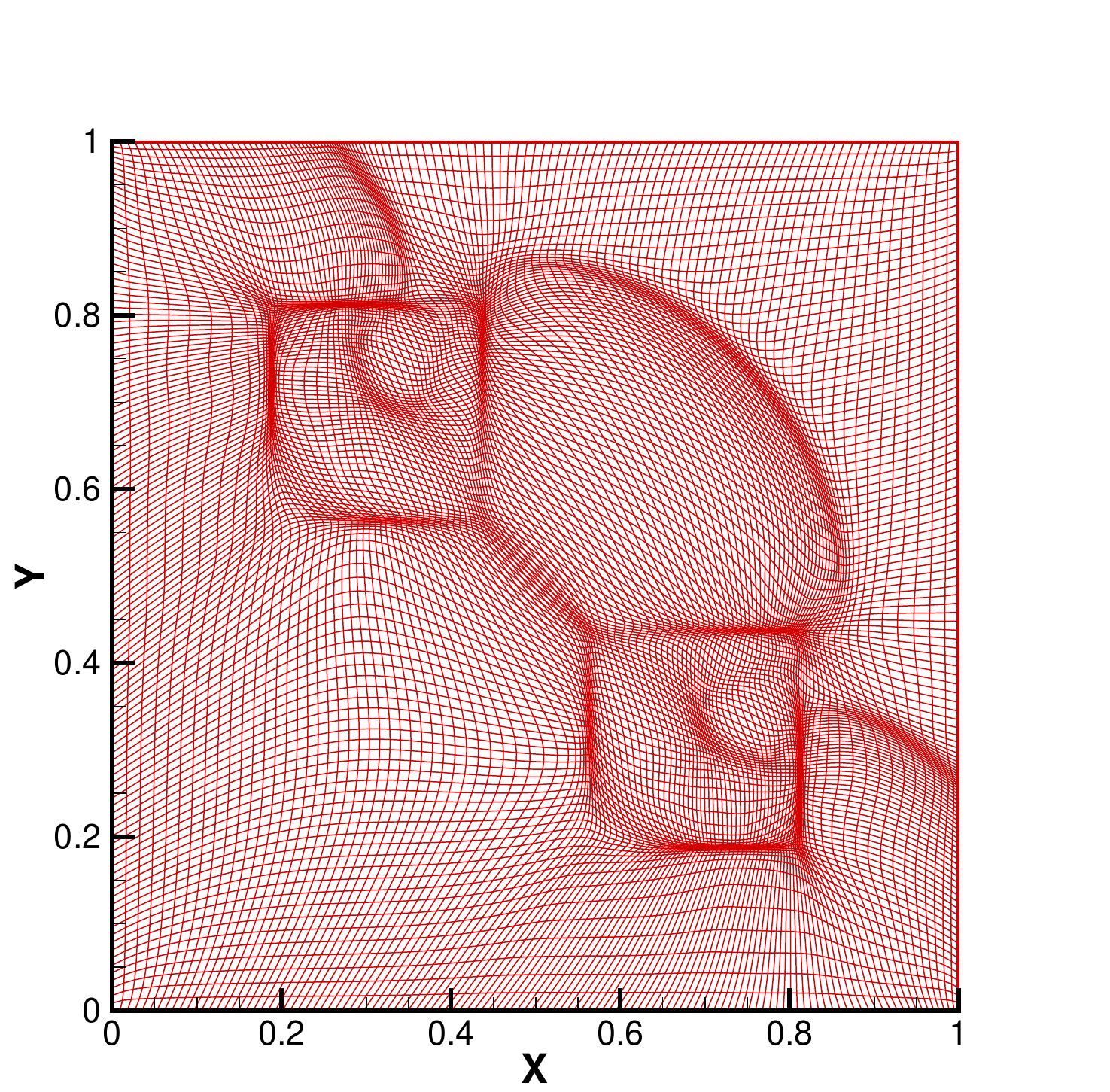}
\end{minipage}
\hspace{0.5in}
\begin{minipage}[t]{2.0in}
\centerline{\scriptsize (d):  with MM2 at $t=2.0$}
\includegraphics[width=2.0in]{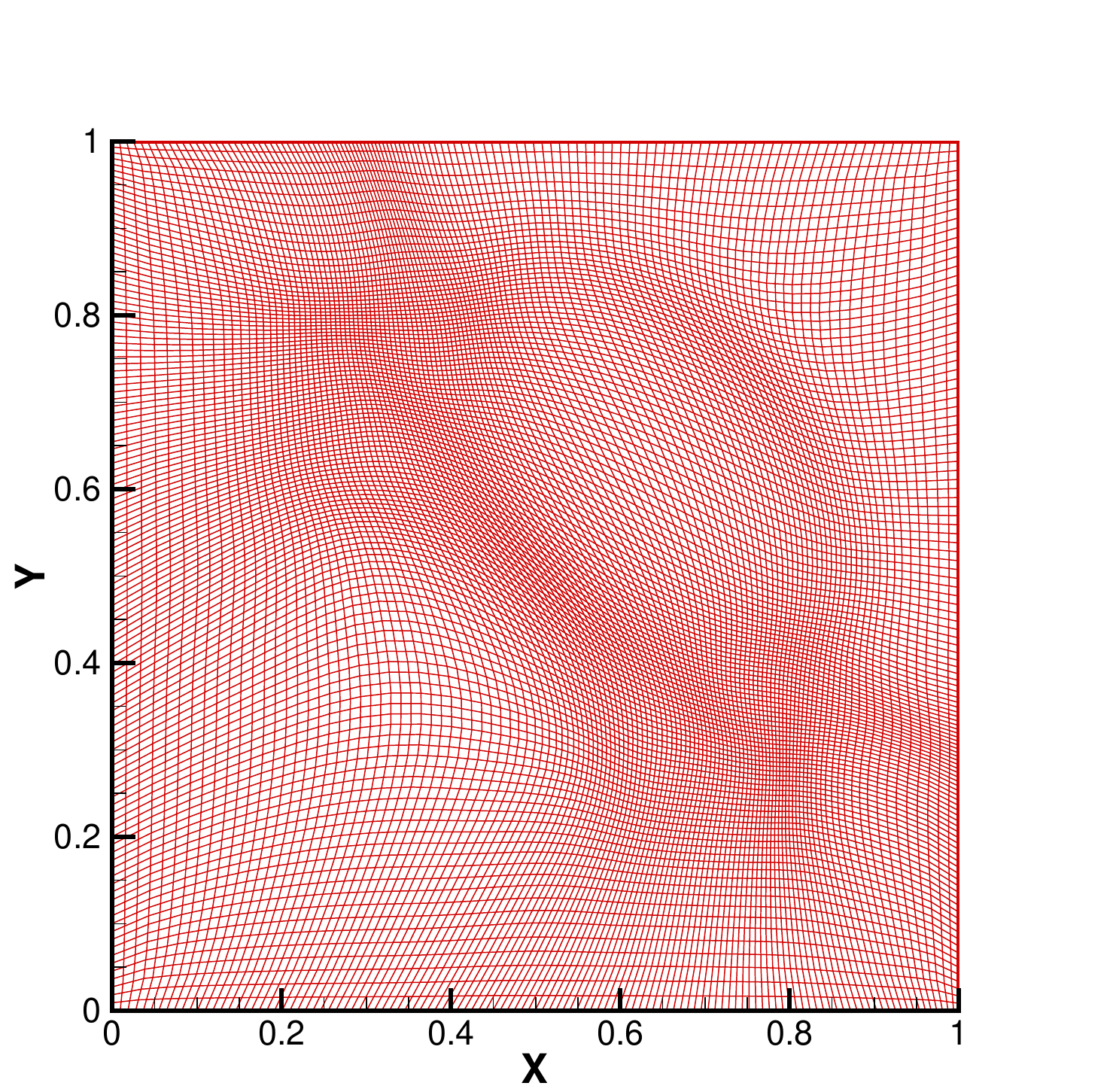}
\end{minipage}
}
\hbox{
\hspace{1in}
\begin{minipage}[t]{2.0in}
\centerline{\scriptsize (e):  with MM1 at $t=2.5$}
\includegraphics[width=2.0in]{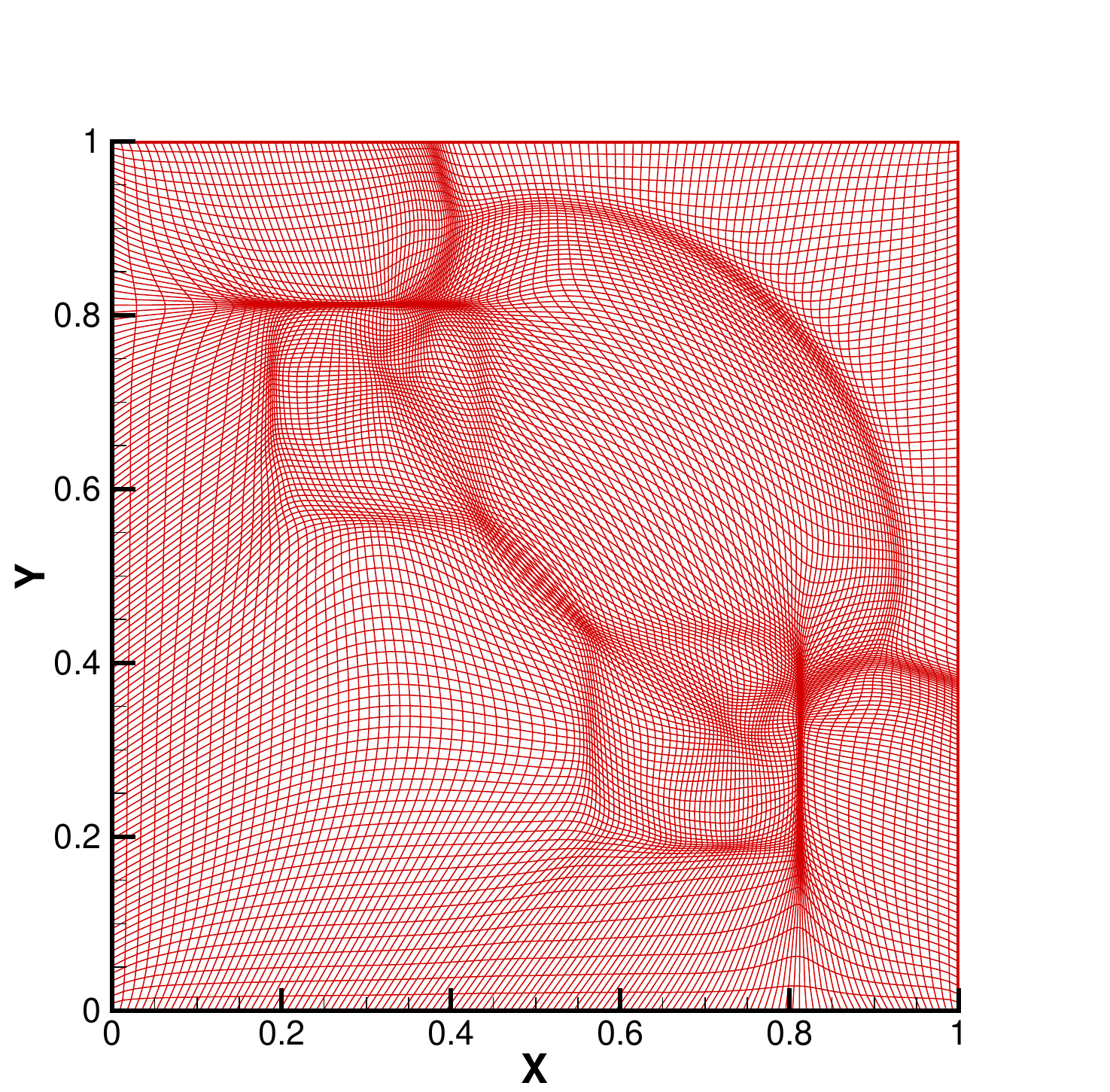}
\end{minipage}
\hspace{0.5in}
\begin{minipage}[t]{2.0in}
\centerline{\scriptsize (f):  with MM2 at $t=2.5$}
\includegraphics[width=2.0in]{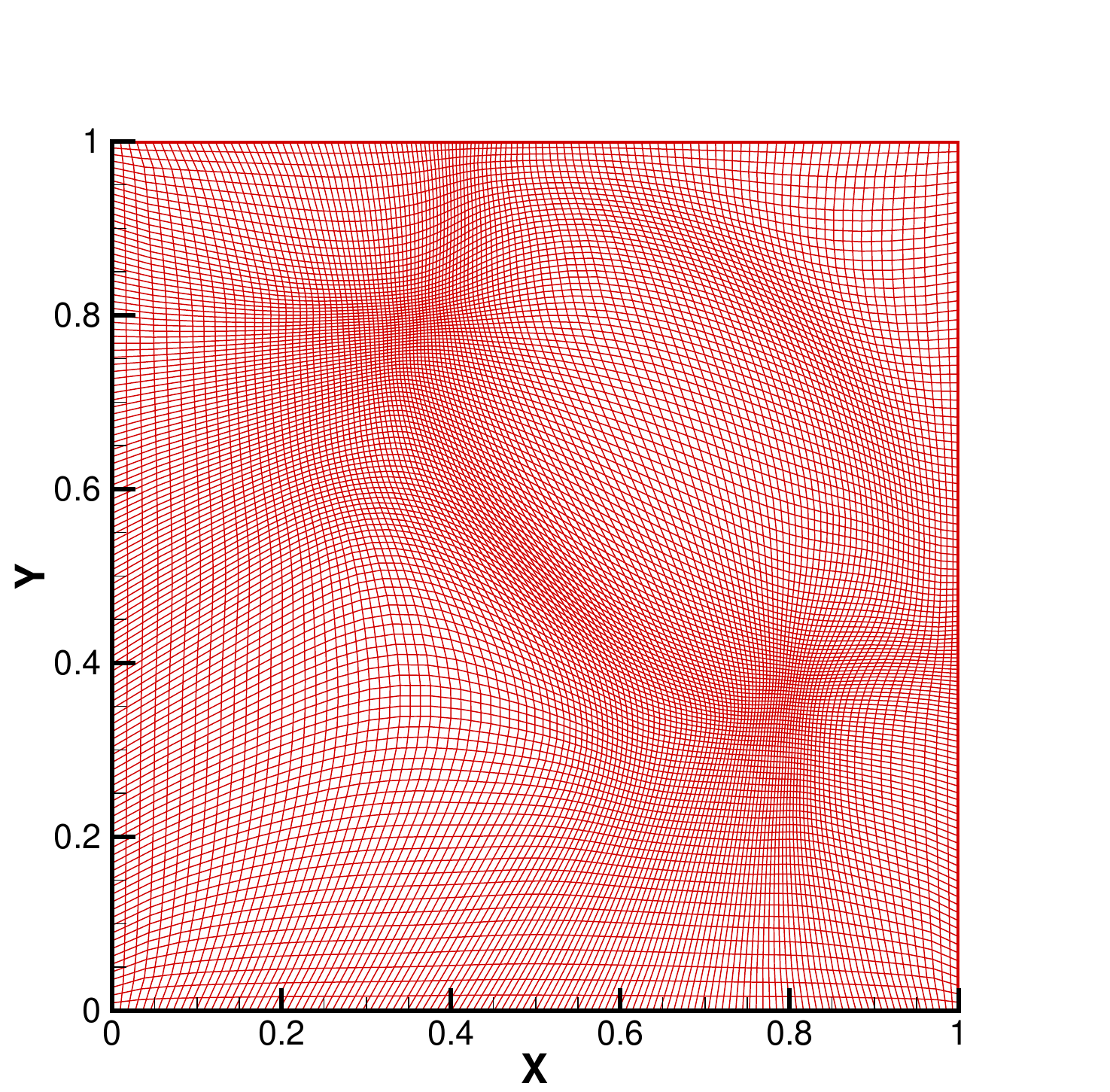}
\end{minipage}
}
\hbox{
\hspace{1in}
\begin{minipage}[t]{2.0in}
\centerline{\scriptsize (g):  with MM1 at $t=3.0$}
\includegraphics[width=2.0in]{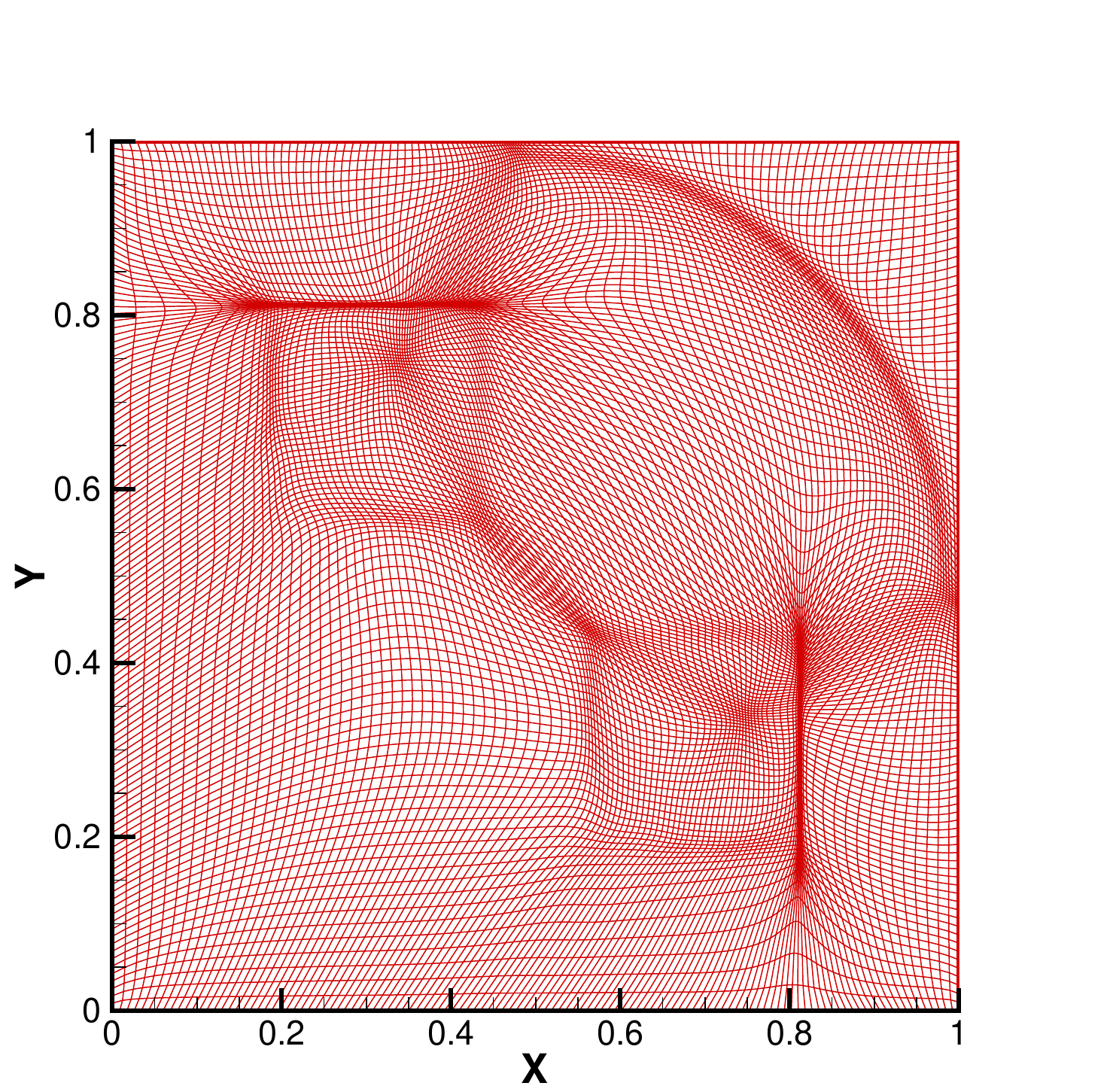}
\end{minipage}
\hspace{0.5in}
\begin{minipage}[t]{2.0in}
\centerline{\scriptsize (h):  with MM2 at $t=3.0$}
\includegraphics[width=2.0in]{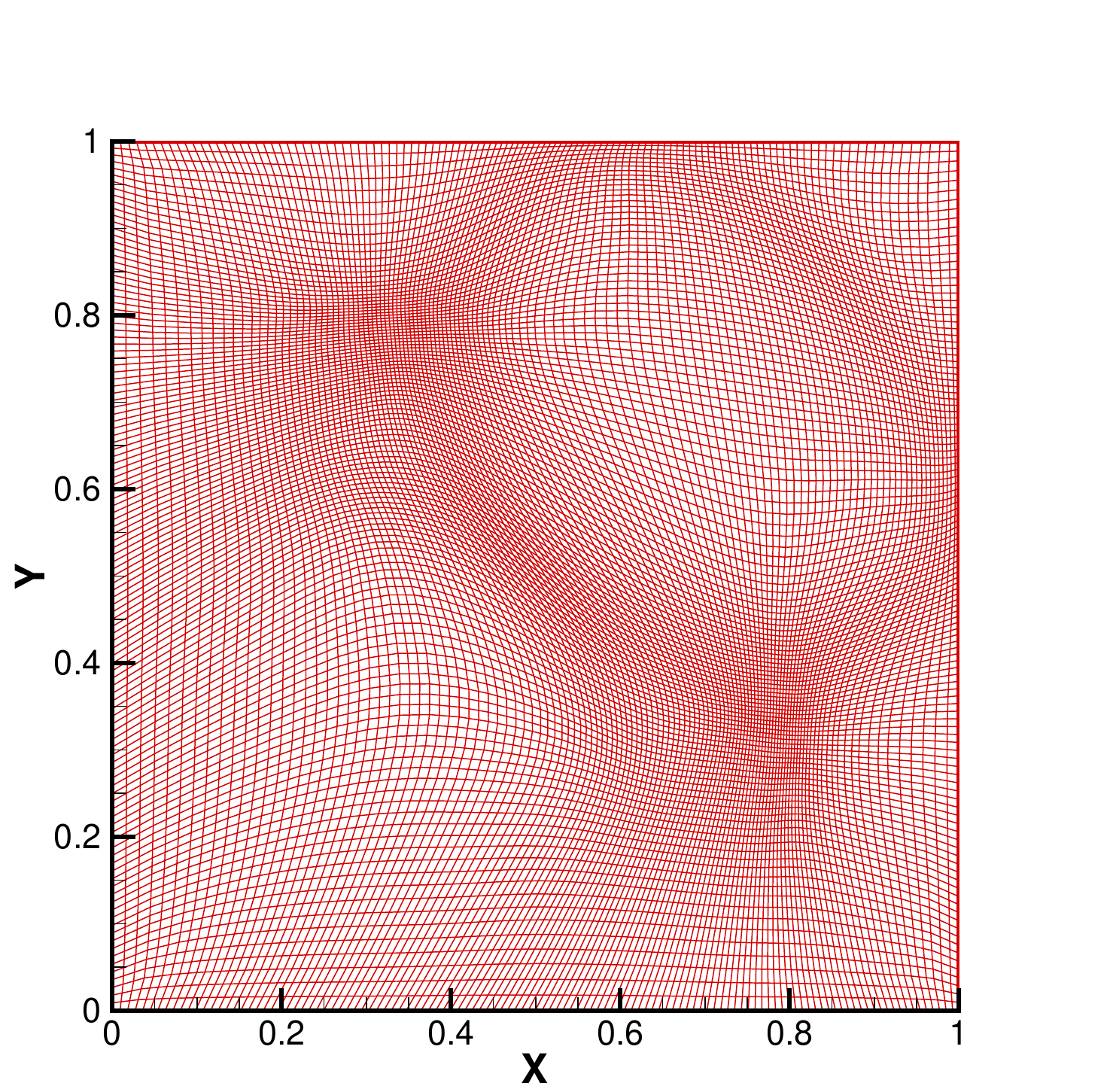}
\end{minipage}
}
\end{center}
\caption{Example~\ref{Example4.3}. Moving meshes of size $121\times 121$ obtained with
MM1 and MM2 moving mesh strategies. MM2 is obtained by uniformly interpolating
a $41\times 41$ moving mesh.}
\label{T32}
\end{figure}

\begin{exam}
\label{Example4.3}

The material configuration for this example is shown in Figs. \ref{BCD3}.
The insets are $(3/16,7/16)\times (9/16,13/16)$ and $(9/16,13/16)\times (3/16,7/16)$
and the distribution of the atomic mass number is given as
\begin{equation}
\label{Z4}
z(x,y)=\begin{cases}
10,   & \text{for }(x,y) \in (\frac{3}{16},\frac{7}{16})\times (\frac{9}{16},\frac{13}{16}) \\
10,   & \text{for } (x,y) \in (\frac{9}{16},\frac{13}{16})\times (\frac{3}{16},\frac{7}{16}) \\
1,    & \text{otherwise}.
\end{cases}
\end{equation}
The boundary of $\Omega = (0,1)\times (0,1)$ is considered as insulated with respect to both radiation
and material conduction, i.e.,
\begin{equation}
\label{BC1}
\frac{\partial E}{\partial n} = \frac{\partial T}{\partial n} = 0, \quad \text{ on } \partial \Omega .
\end{equation}
The initial condition is taken as (cf.  \cite{yuan20092})
\begin{equation}
E(x, y, 0)=0.001+100 \; \mathrm{exp}\left (-100 (x^2+y^2)\right ),
\quad T(x,y, 0) = E(x, y, 0)^{\frac 1 4}.
\label{IC1}
\end{equation}

A typical moving mesh of size $81\times 81$ and the computed solution thereon are shown
in Figs.~\ref{T28} and \ref{T29}. Once again, we can see that our moving mesh method
is able to capture the Marshak wave accurately. The results are in good agreement with
those by Sheng et al. \cite{yuan20092}.
Comparison results are shown in Fig.~\ref{T30} for a moving mesh of $61\times61$
and a uniform mesh of $121\times121$ and in Figs.~\ref{T31} and \ref{T32}
for one-level and two-level moving meshes of $121\times121$. The CPU time is recorded
in Table~\ref{mytab3}.
\end{exam}

\begin{table}[htbp]
 \centering\small
\begin{threeparttable}
  \caption{CPU time comparison among one-level, two-level moving mesh
  and uniform mesh methods for Example \ref{Example4.3}. The CPU time is measured in seconds.
  The last column is the ratio of the used CPU time to that used with a uniform mesh of the same size.}
  \label{mytab3}
 \medskip
 \begin{tabular}{|l|c|c|c|c|}
 \hline
   & Fine Mesh  & Coarse Mesh & Total CPU time & ratio \\
   \hline
  One-level MM & 41$\times$41    &   41$\times$41     &     2359       &  4.72 \\
  \hline
          & 81$\times$81    &    81$\times$81     &      38828        &  17.25      \\
  \hline
          & 121$\times$121  &     121$\times$121     &     $983286$   & $174.16$      \\
  \hline

Tow-level MM   & 41$\times$41  &      41$\times$41     &      2359  &   4.72       \\
\hline
          & 81$\times$81  &    41$\times$41     &     8912      &    3.96        \\
\hline
          & 121$\times$121  &   41$\times$41     &      22836      &  4.04          \\
\hline

Fixed mesh     & 41$\times$41  &      n/a          &      500       &   1                 \\
\hline
          & 81$\times$81  &       n/a         &      2251     &    1                   \\
\hline
          & 121$\times$121  &     n/a          &     5646     &    1                \\
\hline
\end{tabular}
\end{threeparttable}
\end{table}

\section{Conclusions}
\label{SEC:conclusions}

In the previous sections we have studied the moving mesh finite difference solution of the 2T model
for multi-material, non-equilibrium radiation diffusion equations based on the MMPDE moving mesh strategy. 
The model involves nonlinear diffusion coefficients and its solutions stay positive for all time when
they are positive initially. Nonlinear diffusion and preservation of solution positivity pose challenges
in the numerical solution of the model. A coefficient-freezing predictor-corrector
method has been used for treating nonlinear diffusion while a cutoff strategy with a positive threshold
\cite{LuHuVV2012} has been employed to keep the solutions positive. A two-level moving mesh
strategy and the sparse matrix solver UMFPACK with the MAC OSX acceleration have been used to
improve the efficiency of the computation.

The method has been applied to three examples of multi-material non-equilibrium radiation diffusion.
The numerical results show that the method is able to capture the profiles and local structures
of Marshak waves with adequate mesh concentration. The numerical solutions are in good agreement
with those in the existing literature. Comparison studies have also been made between uniform and
adaptive moving meshes and between one-level and two-level moving meshes. It is shown that
the two-level moving mesh strategy can significantly improve the computational efficiency with
only a mild accuracy compromise. Extending the current method to three-dimensional radiation diffusion
models \cite{LSY2017} and more realistic three-temperature models \cite{AJW2017} will be an interesting
research topic for near future.

\vspace{20pt}

\noindent
{\bf Acknowledgmnents.}\\
The work was supported in part by NSFC (China) (Grant No. 11701555),
NSAF (China) (Grant No. U1630247), and Science Challenge Project (China)
(Grant~No.~JCKY2016212A502).


\end{document}